\documentclass[11pt,letter,reqno]{amsart}
\usepackage{amsthm,amsfonts,amsmath,amscd,amssymb}
\usepackage[english]{babel}
\usepackage[T1]{fontenc}
\usepackage[utf8]{inputenc}
\usepackage[
  margin=2cm,includefoot,footskip=30pt]{geometry}
\usepackage{lmodern}
\usepackage{bm}
\usepackage{booktabs}
\usepackage{array}
\usepackage{graphics,graphicx}
\usepackage{subcaption}
\usepackage{float}
\usepackage[font=small]{caption}
\usepackage{cite}
\usepackage[dvipsnames]{xcolor}
\usepackage{tikz-cd}
\usetikzlibrary{decorations.markings,angles}
\usepackage{stmaryrd} \SetSymbolFont{stmry}{bold}{U}{stmry}{m}{n}
\usepackage{mathrsfs}
\usepackage{extarrows}
\usepackage{mathtools}
\usepackage{tensor}
\usepackage{tabularx}
\usepackage{enumerate}
\usepackage{enumitem}
\usepackage{adjustbox}
\usepackage[mathcal]{eucal}
\usepackage{bbold}
\usepackage{bm}
\usepackage{tabularx}
\usepackage{framed}
\usepackage{CJKutf8}
\usepackage{lscape} 

\usepackage{etoolbox}
\patchcmd{\part}{\normalfont\bfseries}{\normalfont\bfseries\centering}{}{}

\newcolumntype{P}[1]{>{\centering\arraybackslash}p{#1}}
\newcolumntype{M}[1]{>{\centering\arraybackslash}m{#1}}

\usepackage[strict]{changepage}
\usepackage{tensor}
\usepackage[all]{xy}
\usepackage{ulem}

\usepackage{fancybox}  
\usepackage[breakable]{tcolorbox}  
\usepackage{graphicx}
\usepackage{tcolorbox}  

\usepackage{verbatim}%%%%%

\usepackage{hyperref}
\hypersetup{
	colorlinks,
	urlcolor={red!55!black},
	citecolor={green!60!black},
	linkcolor={red!55!black}
}

\mathchardef\mhyphen="2D
\def\on{\operatorname}

\makeatletter
\providecommand{\leftsquigarrow}{%
  \mathrel{\mathpalette\reflect@squig\relax}%
}
\newcommand{\reflect@squig}[2]{%
  \reflectbox{$\m@th#1\rightsquigarrow$}%
}
\makeatother

\definecolor{ao}{rgb}{0.0, 0.5, 0.0}

\usepackage{mdframed}
\newtheorem{ftheo}{Theorem}

\newenvironment{fthm}
  {\begin{mdframed}[innertopmargin = 3pt, innerbottommargin=3pt,skipabove=5pt,skipbelow=5pt,linewidth=0.25pt,nobreak=true,align=center]\begin{ftheo}}
  {\end{ftheo}\end{mdframed}}

\usepackage{lipsum}
\newtheorem{conjecture}{Conjecture}

\newenvironment{fconj}
  {\begin{mdframed}[innertopmargin = 3pt, innerbottommargin=3pt,skipabove=5pt,skipbelow=5pt,linewidth=0.25pt,nobreak=true,align=center]\begin{conjecture}}
  {\end{conjecture}\end{mdframed}}
  
\usepackage[capitalize, noabbrev, nosort]{cleveref}
\newtheorem{theorem}{Theorem}[section]
\newtheorem{lemma}[theorem]{Lemma}
\newtheorem{proposition}[theorem]{Proposition}
\newtheorem{corollary}[theorem]{Corollary}

\newtheorem{corollaryalpha}[theorem]{Corollary}

\theoremstyle{definition}

\newtheorem{definition}[theorem]{Definition}

\newtheorem{remark}[theorem]{Remark}
\newtheorem{example}[theorem]{Example}

\usepackage{mdframed}
\usepackage{lipsum}

\newenvironment{fdef}
  {\begin{mdframed}[innertopmargin = 3pt, innerbottommargin=3pt,skipabove=5pt,skipbelow=5pt,linewidth=0.25pt,nobreak=true,align=center]\begin{definition}}
  {\end{definition}\end{mdframed}}

\mdfdefinestyle{style1}{leftmargin=5cm,rightmargin=5cm,skipbelow=1cm,skipabove=1cm,
   innertopmargin=0pt}

   \newenvironment{fcor}
  {\begin{mdframed}[innertopmargin = 3pt, innerbottommargin=3pt,skipabove=5pt,skipbelow=5pt,linewidth=0.25pt,nobreak=true,align=center]\begin{corollary}}
  {\end{corollary}\end{mdframed}}

\usepackage{varwidth}
\newcommand{\B}{\mathcal{B}}
\newcommand{\C}{\mathbb{C}}
\newcommand{\bG}{\mathbb{G}}
\newcommand{\D}{\mathbb{D}}

\newcommand{\M}{\mathcal{M}}
\newcommand{\N}{\mathbb{N}}
\newcommand{\T}{\mathcal{T}}
\newcommand{\LL}{\mathcal{L}}

\newcommand{\kk}{\mathbb{k}}
\newcommand{\R}{\mathbb{R}}
\newcommand{\cC}{\mathbb{C}}
\newcommand{\Z}{\mathbb{Z}}
\newcommand{\lL}{\mathbb{L}}

\newcommand{\g}{\gamma}
\renewcommand{\r}{\rho}
\newcommand{\w}{\mathfrak{w}}
\renewcommand{\tt}{\mathfrak{t}}

\renewcommand{\b}{\beta}
\newcommand{\lr}{\longrightarrow}
\newcommand{\s}{\sigma}

\newcommand{\SA}{\mathcal{A}}
\newcommand{\SX}{\mathcal{X}}
\newcommand{\SW}{\mathcal{W}}
\newcommand{\RW}{\mathcal{R}\SW}
\newcommand{\sse}{\subseteq}

\newcommand{\La}{\Lambda}
\newcommand{\la}{\lambda}

\newcommand{\dd}{\partial}
\newcommand{\std}{\text{st}}
\newcommand{\Br}{\text{Br}}
\newcommand{\gr}{\operatorname{gr}}
\newcommand{\SL}{\operatorname{SL}}

\newcommand{\sbs}{\hfill\rule{1.2ex}{1.2ex}}

    \newcommand{\borel}[0]{\mathsf{B}}
    
    \renewcommand{\T}[0]{\mathsf{T}}
    
    \newcommand{\FF}[0]{\mathsf{F}}
\newcommand{\uni}[0]{\mathsf{U}}

\newcommand{\dynkin}[0]{\mathsf{D}}
\newcommand{\G}[0]{\mathsf{G}}

\newcommand{\fW}{\mathfrak{w}}

\newcommand{\ff}{\mathfrak{f}}

\DeclareMathOperator{\Gr}{Gr} 
\def\Le{\reflectbox{$\mathrm{L}$}}

\newcommand{\ord}{\operatorname{ord}}
\newcommand{\tht}{\vartheta}

\DeclareRobustCommand{\ateb}{\text{\reflectbox{$\beta$}}}
\newcommand{\SR}{\mathscr{R}}
\newcommand{\PP}{\mathbb{P}}
\renewcommand{\L}{\mathbb{L}}
\renewcommand{\top}{\operatorname{top}}

\newcommand{\Spec}{\operatorname{Spec}}
\newcommand{\BS}{\operatorname{BS}}
\newcommand{\hH}{\operatorname{H}}
\newcommand{\hol}{\operatorname{Hol}}
\newcommand{\GL}{\operatorname{GL}}
\newcommand{\brn}{\mbox{Br}_n^+}
\newcommand{\wb}{\w:\beta\lr\delta(\beta)}
\newcommand{\mL}{\operatorname{\mu_{\mathscr{D}}(L)}}

\title{Algebraic and geometric aspects of Demazure weaves}
\author{Roger Casals}
\date{\today}

\begin{document}
\maketitle

\begin{abstract}
This is an expository account of recent work in the study of braid varieties and their relation to low-dimensional contact and symplectic topology. A driving theme is the development and use of Demazure weaves to construct cluster algebra structures on the coordinate rings of such braid varieties.
\end{abstract}

\setcounter{tocdepth}{1}
\tableofcontents

\section*{Introduction}\label{sec:intro}

The object of these lecture notes is to provide an expository account of recent work in the study of braid varieties and their relation to low-dimensional contact and symplectic topology. A driving theme will be the development and use of Demazure weaves to construct cluster algebra structures on the coordinate rings of such braid varieties.\\

\noindent The first part of the manuscript addresses algebraic and combinatorial aspects of braid varieties building on the diagrammatic calculus of weaves. Its main goal is to explain how to use Demazure weaves to construct cluster algebra structures on the coordinate rings of braid varieties. The second part is focused on the study of Legendrian links and their Lagrangian fillings using weaves, which helps illustrate one fruitful connection between these braid varieties and their cluster structures, weaves, and symplectic topology.\\

%%%%%%%%%%%%%%%%%%%%%%%%%%%%%%%%%%%%%
%%%%%%%%%%%%%%%%%%%%%%%%%%%%%%%%%%%%%
%%%%%%%%%%%%%%%%%%%%%%%%%%%%%%%%%%%%%

\noindent {\bf Scientific context for Part I.} The study of Lie groups $\G$, their Lie algebras $\mathfrak{g}$, and their representation theory is useful in mathematics, see e.g.~\cite{MR1920389,MR2069175,MR2440737,MR3331229,MR1240056,MR1345386,MR3443860}. In many of its incarnations, it builds on the Weyl groups $W(\G)$ and the manner in which they help organize different structures around $\G$, e.g.~via the Bruhat decomposition, or by encoding symmetries of roots. Among other features, braid varieties $X(\beta)$ and their associated weave diagrams $\wb$, where $\beta\in\Br^+(W(G))$ is a positive braid word, provide new ground to explore a version of {\bf non-reduced Lie theory}, where the role of elements of the Weyl group $W(\G)$ is played by the positive braid monoid $\Br^+(W(\G))$, and that of reduced expressions for $w\in W(\G)$ is played by positive braid words $\beta=\s_{i_1}\s_{i_2}\cdot\ldots\cdot\s_{i_l}$. Some instances of this phenomenon are as follows:\\

\begin{enumerate}
    \item The intersection theory of Schubert cells in flag varieties $\G/\borel$ and Grassmannians. Historically rooted in enumerative algebraic geometry, these also naturally arise when studying the Bruhat decomposition
        \begin{equation}\label{eq:Bruhat_decomposition}\G=\bigsqcup_{w\in W}\borel w\borel
        \end{equation}
    of $\G$ and projecting to $\G/\borel$. For instance, it follows from \eqref{eq:Bruhat_decomposition} that elements in $W(\G)$ dictate the possible relative positions between two flags. Overall, if $\borel_-$ is the opposite Borel, intersections of the form
        \begin{equation}\label{eq:intersection_shape}
            \fbox{$(\borel u\borel)\cap(\borel_- v\borel_-)$}\qquad u,v\in W(\G)
        \end{equation}
and their analogues in flag varieties and Grassmannians have a central role in Lie theory and algebraic combinatorics. Double Bruhat cells, Richardson varieties, and intersections of Schubert varieties can be understood as instances of \eqref{eq:intersection_shape}, and the variants $(\borel u\borel)\cap\uni_-$ provide unipotent cells. Braid varieties $X(\beta)$ generalize these and, coarsely put, can be understood as being of the form \eqref{eq:intersection_shape} if we wrote $(\borel\beta\borel)\cap(\borel_-\delta(\beta)\borel_-)$, where we need to make rigorous sense of $\borel\beta\borel$ for a positive braid word $\beta$, and $\delta(\beta)\in W(\G)$ is the Demazure product of $\beta$, which is still an element of the Weyl group.\\

\item The coordinate rings of varieties arising as \eqref{eq:intersection_shape} tend to have remarkable combinatorial structures, central to them being the fact that they admit the structure of a cluster algebra, cf.~\cite{FWZ16,FZ2002I,FZ2003II,BFZ2005,FZ2007IV,BZ2005}. This is also the case for many of their cousins, such as the positroid subvarieties of Grassmannians, cf.~\cite{GL2023,Scott06,SSBW2019,FSB2022,JKS2016,BKM2016}. An important fact about braid varieties $X(\beta)$, proven in \cite{CGGLSS25}, and see also \cite{GalashinLamShermanBennettSpeyer2022TypeA,GalashinLamShermanBennett2026GeneralType,casals2025comparingclusteralgebrasbraid}, is that their coordinate rings $\Z[X(\beta)]$ admit the structure of a cluster algebra, and each Demazure weave $\wb$ provides a cluster seed for such a cluster algebra. It does so in an explicit manner, so that both the quiver and the cluster variables are readily computable via weaves.\\

\noindent Here the diagrammatic calculus of weaves, cf.~\cite{CasalsZaslow,CasalsWeng22,CGGJ24,CGGLSS25}, could be seen as generalizing the role played by the appropriate algebraic combinatorial objects in the reduced cases, e.g.~reduced plabic graphs in the case of positroids, or (tagged) triangulations in the case of cluster algebras associated to surfaces \cite{FG2006,FST08,FT12}. Thus an important feature of many Lie theoretic varieties, namely the existence of cluster algebra structures and diagrammatic combinatorics used to manipulate them, seems to still persist when we enter the non-reduced land through the lens of braid varieties.\\

\item Building on the introduction of global bases in representation theory, the coordinate rings of varieties of the form \eqref{eq:intersection_shape} provide fruitful soil for algebraic combinatorial models describing global bases. In particular, the classical theory of crystal bases, which captures part of the structure of such global bases, has $W(\G)$ and the roots of $\G$ as building ingredients. To wit, the Lusztig data that parametrize different crystals of interest depend on a choice of a reduced expression for a certain $w\in W(\G)$, depending on the coordinate ring being studied. In line with this, the monoidal categorification of cluster coordinate rings of varieties of the form \eqref{eq:intersection_shape}, where such a ring is realized as a Grothendieck ring $K_0$ of a monoidal category, was pioneered by D.~Hernandez and B.~Leclerc, cf.~\cite{MR2682185,MR3077685,MR3331727,MR3500832,MR3576119,MR4404353}.\\

In the non-reduced case of braid varieties $X(\beta)$, monoidal categorifications have recently been constructed and their corresponding theory of global bases and crystals is now rather well-established thanks to the work of M.~Kashiwara and his co-authors, cf.~\cite{KangKashiwaraKimOh2018Monoidal,KKOP2020Monoidal,KKOP2024MonoidalII,KKOP2025MonoidalIII,KKOP2024PBWAffine,KashiwaraPark2022CategoricalCrystals,Park2023ExtendedCrystals,OhPark2025BosonicPBW,KKOP2025GlobalBases,KKOP2024BraidSymmetries,KKOP2025Faithful} and see also \cite{Bi2026Quantization,Bi2026TwistedProducts,HuhJungKimPark2026Comparison}. In particular, there is a precise sense in which the classical $w$-Lusztig data, with $w\in W(\G)$, can now be generalized to $\beta$-Lusztig data, where $\beta$ is a positive braid word. The two key ingredients to make sense of this non-reduced version of Lusztig data are the extended crystal $\widehat{B}(\infty)$ and the braid group symmetries that they construct on the monoidal category, the extended crystal $\widehat{B}(\infty)$, and their associated objects, including the bosonic extension algebra, which in a sense contains $\Z[X(\beta)]$ for all positive braid words $\beta$.\\

\item The additive categorification of coordinate rings of varieties of the form \eqref{eq:intersection_shape} leads to a rich interplay between the classical tilting theory for Abelian categories of modules, e.g.~over preprojective algebras, and the combinatorics of quiver mutation and cluster variables, see e.g.~\cite{MR2242628,GLS08,GeissLeclercSchroer2013QuantumCoordinate,MR3534832} and \cite{JKS2016,BKM2016,MR3694688,MR4662273,MR4456104,MR2844757}. In addition to providing enlightening results on the representation theory of quivers (and their potentials, see \cite{DWZ08,DWZ10_2}), these additive categorifications nicely interact with classical structures from Lie theory, such as highest weight structures, cf.~\cite{CPS88,Kra17,BS24}. While highest weight structures were in part motivated by the study of Verma modules and the representation theory of the Lie algebras $\mathfrak{g}$, they naturally appear when additively categorifying the cluster algebras for the coordinate rings of varieties of the form \eqref{eq:intersection_shape}, see e.g.~\cite{IR11} and \cite{MR2242628,GLS08,BIRS09}.\\

We have now recently established in \cite{CasalsChrist} the non-reduced generalization of such additive categorifications for braid varieties $X(\beta)$ and their Demazure weaves. A key fact is that non-reducedness naturally leads to considering derived categories, rather than just their Abelian hearts.\footnote{Or more generally, pretriangulated dg categories, or stable $\infty$-categories, which is the viewpoint in \cite{CasalsChrist}.} An important feature of these additive categorifications of braid varieties $X(\beta)$ is the emergence of a derived analogue for highest weight theory, with full exceptional collections acquiring the role of standard objects. The role of universal extensions and quotients is now played by a truncated version of spherical twist functors, which allow for the construction of both simple-minded and silting collections from the exceptional collection, much in line with how simple and indecomposable projective objects are classically built from standard objects. In short, we have a glimpse of the non-reduced world and see how important categorical results and theories in the Abelian (reduced) setting seem to be part of, and sometimes better understood via, a derived (non-reduced) theory based on the positive braid monoid.\\

    \item Seeded in G.~Lusztig's study of positivity for reductive Lie groups \cite{Lusztig,Lusztig1997PositivityBases,Lusztig1998,Lusztig2019PositivityII,Lusztig2025Survey}, positive structures are starting to have a more unified theory around them. Examples include the positivity phenomena for cluster algebras \cite{FZ1999,Fomin2010,FominLusztig2022,FZ2000Tests,BFZ1996,BZ1997,LS2015,MSW2011,Davison2018,DM2021}, the positive representations in higher Teichm\"uller theory \cite{FG2006,MR2221137}, positive geometries \cite{ABL2017,ABHY2018,ALS2021,Lam2024Invitation,Lam2025Moduli}, and positivity in matroids \cite{ARW2017}, see also \cite{KW2014}. Varieties appearing as \eqref{eq:intersection_shape} do exhibit such positivity, e.g.~existence of a positive atlas, and positive properties for some of their global bases, cf.~G.~Lusztig's works above. Such positivity features are also shared by braid varieties $X(\beta)$, see e.g.~\cite{BaoHe2022TwistedPositivity}. These positivity properties have an important role when discussing tropicalization and parametrizing global bases. In the context of Demazure weaves, even the most basic definitions, such as those of Lusztig cycles and their associated quiver, are grounded on tropical identities arising from this underlying positivity.\\
\end{enumerate}

\noindent The above are a few instances of this reduced to non-reduced phenomenon, selected just from personal knowledge. See also \cite{trinh2021heckecategoryunipotentlocus,trinh2022unipotentelementstwistinglink} for another interesting example, where identities of R.~Steinberg \cite[Theorem 15.1]{Steinberg1968Endomorphisms} and N.~Kawanaka \cite[Theorem 4.1]{Kawanaka1975} on certain counts over finite fields in the reduced setting are enhanced to $\ell$-adic cohomological isomorphisms and established in the braid variety non-reduced setting, cf.~\cite[Corollary 4]{trinh2021heckecategoryunipotentlocus} and see also \cite[Conjectures 1 \& 17]{trinh2022unipotentelementstwistinglink}. It should be emphasized that there are many previous works that point toward such non-reduced extensions in and around Lie theory. These include works of P.~Deligne, G.~Lusztig and others, see e.g.~\cite{Deligne1972Tresses,Deligne1997Action,DeligneLusztig1976Representations,BroueMichel1997Regular,DigneMichelRouquier2007Cohomologie,Lusztig2011Varieties} and cf.~\cite[Appendix B]{trinh2021heckecategoryunipotentlocus}. Note as well that varieties that arise by considering (not necessarily reduced) positive braid words have also a prominent standing in geometry and topology, see \cite{Bott1956Morse,BottSamelson1958Morse,Demazure1974Desingularisation,LakshmibaiLittelmannMagyar2002SMT} and the second part of these notes, especially its relation to \cite{FG2006}.

\begin{table}[h!]
  \centering

\begin{tabular}{ |P{4cm}|P{5cm}|P{5cm}|}
	\hline
     & {\bf Reduced expression $w\in W(\G)$} & {\bf Positive braid word $\beta\in\on{Br}^+(W(\G))$}\\
	\hline

    Algebraic variety & $(\borel u\borel)\cap(\borel_- v\borel_-)$ & $X(\beta)\approx(\borel \beta\borel)\cap(\borel_- \delta(\beta)\borel_-)$\\
	\hline

    \rule{0pt}{0.5cm}Coordinate ring & Cluster algebra & Cluster algebra\\
	\hline

    \rule{0pt}{0.5cm} Some algebraic combinatorial objects & double wiring diagrams, plabic graphs, and more & Demazure weaves $\wb$\\
	\hline

    \rule{0pt}{0.5cm} Positive global bases and crystals & canonical bases and $w$-Lusztig data & canonical bases and $\beta$-Lusztig data\\
	\hline

    \rule{0pt}{0.5cm}Constructing additive categorifications & Via simples and tilting in an Abelian category & Simple-minded and silting in derived category\\
	\hline

    \rule{0pt}{0.5cm}Constructing monoidal categorifications & Via subcategories of $\mathcal{U}_q'(\widehat{\mathfrak{g}})$-mod& By incorporating braid group action and symmetries\\
    \hline
\end{tabular}
  \caption{Summary of analogies between objects and results when using a reduced expression for an element $w\in W(\G)$ of the Weyl group, and their putative counterparts when using a (not necessarily reduced) positive braid word $\beta\in\Br^+(W(\G))$.}  \label{table:summary_nonreduced}
\end{table}
\vspace{0.1cm}

\noindent The comments above are summarized in Table \ref{table:summary_nonreduced}. There are many additional aspects to Table \ref{table:summary_nonreduced} that are unfortunately omitted, e.g.~the work of B.~Keller, C.~Amiot, Y.~Wu and others on cluster and Higgs categories, see e.g.~\cite{MR2640929,amiot2011generalizedclustercategories,MR2764906,MR3024264,MR4289034,MR4157113,MR2762547,Keller12_DT2,MR2795754,MR2385670,MR2316979,MR4581969,MR4576399}, for the row on additive categorifications, or the work of P.~Galashin, T.~Lam, M.~Sherman-Bennett and D.~Speyer on 3D plabic graphs \cite{GalashinLamShermanBennettSpeyer2022TypeA} for another instance of Type A combinatorics for non-reduced $\beta$.\footnote{I apologize for other involuntary sins of omission: should the reader think a citation to another work must be added, please contact the author.} The content of these notes will focus on the study and construction of cluster algebras in and around braid varieties, i.e.~$(1)$ and $(2)$ above. Items $(3)$ and $(4)$, and in a sense $(5)$, require an understanding of such cluster structures in $(2)$, and thus it seems appropriate to focus on $(2)$ for now.\\

%%%%%%%%%%%%%%%%%%%%%%%%%%%%%%%%%%%%%
%%%%%%%%%%%%%%%%%%%%%%%%%%%%%%%%%%%%%
%%%%%%%%%%%%%%%%%%%%%%%%%%%%%%%%%%%%%
\noindent {\bf Main Results related to Part I.} Let $\G$ be a simple algebraic Lie group and $\beta\in\Br^+(\G)$ a positive braid word. The first part of these notes focuses on the combinatorial and algebraic aspects of braid varieties $X(\beta)$, \cite{CGGJ24,CGGLSS25} and see \cite[Appendix B]{trinh2021heckecategoryunipotentlocus}. Intuitively, the points in the algebraic varieties $X(\beta)$ parametrize linear tuples of flags in $\G/\borel$ in certain relative positions. The relative positions are given by the crossings, i.e.~the letters, of the positive braid word $\beta$. Specifically, let $\beta = \sigma_{i_1}\cdots \sigma_{i_{\ell}}$ be a positive braid word and $\delta(\beta) \in W(\G)$ its Demazure product. By definition, the {\it braid variety} associated with $\beta$ is
\begin{equation}\label{eq:braid_variety}
X(\beta) := \{(\borel_1, \dots, \borel_{\ell+1}) \in (\G/\borel)^{\ell + 1} \mid \borel_1 = \borel, \borel_{k} \buildrel s_{i_{k}} \over \longrightarrow \borel_{k+1}, \borel_{\ell + 1} = \delta(\beta) \borel \},
\end{equation}
\noindent where $\delta(\beta) \in W(\G) \cong N_{\G}(\T)/\T$ has been lifted to $N_{\G}(\T)$. These are studied in detail in \Cref{sec:weaves_braid_varieties}. In particular, they are smooth irreducible varieties over $\Z$ and their isomorphism type only depends on the underlying braid element $[\beta]\in\B(\G)$.

\begin{center}
	\begin{figure}[h!]
		\centering
		\includegraphics[scale=1]{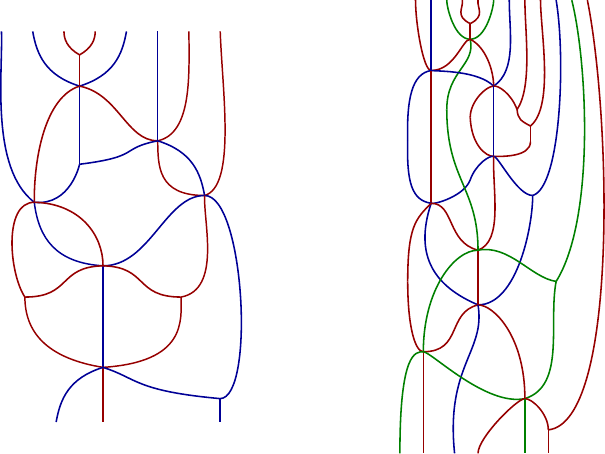}
		\caption{Two instances of Demazure weaves. These diagrams are read top-to-bottom, and can be intuitively understood as ways to simplify a positive braid word $\beta$, starting at the top, down to a permutation, at the bottom. (Left) A Demazure weave $\w:\beta\lr\delta(\beta)$ in $\SL_3$, from the 3-stranded braid word $\beta=\s_1^2\s_2^2\s_1^2\s_2^2$ to $\delta(\beta)=\s_1\s_2\s_1$, where blue weave edges are labeled by $\s_1\in S_3$ and red weave edges by $\s_2\in S_3$. (Right) A Demazure weave $\w:\beta\lr\delta(\beta)$ in $\SL_4$, from the 4-stranded braid word $\beta=\s_2\s_1\s_3\s_2^2\s_1\s_2^2\s_1\s_3\s_2$ to its Demazure product $\delta(\beta)$, where blue weave edges are labeled by $\s_1\in S_3$, and red and green weave edges by $\s_2$ and $s_3\in S_4$ respectively.}
		\label{fig:Example_TwoWeavesSL3SL4}
	\end{figure}
\end{center}

A central contribution of Part I is to develop the diagrammatic calculus of Demazure weaves, which provides an algebraic and combinatorial framework to study the geometry of $X(\beta)$ and the algebraic properties of its coordinate ring of regular functions $\Z[X(\beta)]$. \Cref{fig:Example_TwoWeavesSL3SL4} depicts two instances of Demazure weaves for $\G=\SL_3$ and $\SL_4$. In a nutshell, a Demazure weave $\wb$ is a specific sequence that simplifies a positive braid word $\beta$ to its Demazure product $\delta(\beta)$ by using either braid moves, which contribute to 4 and 6-valent vertices, or the operation $\s_i^2\to\s_i$, which contributes a trivalent vertex. We will introduce the notion of weave equivalences, which generalize relations between syzygies in the positive braid monoid, and the operation of weave mutation.

\begin{table}
\begin{center}
\caption{Some of the concepts being studied in Part I.}\label{table:PartI_introduction}
\small
\begin{tabular}{ |P{3.8cm}|P{5.5cm}|P{4.8cm}|  }
	\hline
	\multicolumn{3}{|c|}{ {\bf First algebraic and combinatorial aspects of Demazure weaves}} \\
	\hline

{\bf Weaves} & {\bf Braid varieties} & {\bf Cluster theory}\\
\hline

Positive braid word $\beta$ & Algebraic variety $X(\beta)$ & Cluster algebra $\Z[X(\beta)]$\\
\hline
Demazure weave $\w:\beta\to\delta(\beta)$ & Open embedding of algebraic torus $\iota_w:T_\w\hookrightarrow X(\beta)$ & Chart morphism $\Z[X(\beta)]\to\Z[u_1^{\pm1},\ldots,u_r^{\pm1}]$\\
\hline
Lusztig cycles $\{\g_v\}$ & Rational functions $X^\w_v\in\Z(X(\beta))$ & cluster $\SX$-variables $\{X_v^\w\}$\\
\hline
Dual Lusztig cycles $\{\r_v\}$ & Regular functions $A^\w_v\in\Z[X(\beta)]$ & cluster $\SA$-variables $\{A_v^\w\}$\\
\hline
Intersecting weave cycles & Algebraic 2-form $\omega_{cl}\in\Gamma(T_\w,\Omega^2 T_\w)$ & Quiver $Q_\w$\\
\hline
Weave mutation $\w\stackrel{\mu}{\sim}\w'$ & Change of torus $\iota_\w(T_\w)\neq \iota_{\w'}(T_{\w'})$ & Cluster mutation\\
\hline
\end{tabular}
\normalsize
\end{center}
\end{table}

\begin{comment}

\begin{enumerate}
    \item Any pair of Demazure weaves $\wb$ can be connected by a sequence of weave equivalences and weave mutations.\\
    \item The intersection theory of Lusztig cycles and dual Lusztig cycles of a Demazure weave $\wb$ is invariant under weave equivalences, and it changes according to quiver mutation for weave mutations.\\
    
    \item A Demazure weave $\wb$ provides a cluster $\SA$-seed $(T^\SA_\w,Q_\w,\{A_i^\w\})$ so that $A_i^\w\in\Z[X(\beta)]$, $Q_\w$ encodes the intersections of dual Lusztig cycles, and we have an isomorphism
    \begin{equation}\label{eq:braid_variety_is_cluster}
        \Z[X(\beta)]\cong\SA(Q_\w)
    \end{equation}
    where $\SA(Q_\w)$ denotes the cluster algebra associated to the quiver $Q_\w$.\\

    \item Similarly, a Demazure weave $\wb$ provides a cluster $\SX$-seed $(T^\SX_\w,Q_\w,\{X_i^\w\})$ so that $X_i^\w\in\Z(X(\beta))$, $Q_\w$ encodes the intersections of Lusztig cycles, and $X(\beta)$ is isomorphic to the cluster $\SX$-scheme associated to $Q_\w$.\\
    
    \item Several geometric features of $X(\beta)$, including a defining ideal, a torus action and a cyclic symmetry, an algebraic closed $2$-form compatible with the cluster structures, and a natural Landau-Ginzburg potential associated to a positive braid word $\beta$, can all be expressed explicitly in terms of Reeb coordinates of $X(\beta)$.\\
\end{enumerate}
\end{comment}

From a Demazure weave $\wb$, we will explain how to construct a number of combinatorial objects, including its Lusztig cycles, its dual Lusztig cycles, and an intersection theory for them. In particular, we will obtain a quiver $Q_\w$ for each Demazure weave $\wb$, which encodes the intersections between Lusztig cycles. The vertices of the quiver are indexed by the trivalent vertices of $\w$, i.e.~one for each time we applied the operation $\s_i^2\to\s_i$ when simplifying $\beta$ to $\delta(\beta)$. The first combinatorial properties of these Demazure weaves are as follows:

\begin{fthm}\label{thm:weave_intro1}
Let $\beta\in\brn$ be a positive braid word. Then the following hold:

\begin{enumerate}[label=$(\arabic*)$]
    \item Two Demazure weaves $\w,\w':\beta\lr\delta(\beta)$ are connected by a sequence of weave equivalences and weave mutations.\\

    \item If $\w_1\simeq\w_2$ is a weave equivalence, the two intersection quivers $Q_{\w_1}$ and $Q_{\w_2}$ encoding the intersections of Lusztig cycles for $\w_1$ and $\w_2$ coincide, i.e.~$Q_{\w_1}=Q_{\w_2}$ after relabeling vertices. The same holds for the dual Lusztig cycles.\\

    \item If $\w_1\stackrel{\mu}{\sim}\w_2$ is a weave mutation at the Lusztig cycle $\gamma_{v_k}$ of $\w_1$, then the intersection quiver $Q_{\w_2}$ is obtained from $Q_{\w_1}$ by quiver mutation at its $k$th vertex, i.e.~$Q_{\w_2}=\mu_k(Q_{\w_1})$. The same holds for the dual Lusztig cycles.\\
\end{enumerate}
\end{fthm}
\vspace{0.2cm}
\noindent From a Demazure weave $\wb$, we will also construct a number of objects related to the braid variety $X(\beta)$, including an open torus chart $T_\w\sse X(\beta)$, and several functions on $X(\beta)$. This is summarized in Table \ref{table:PartI_introduction}. An important instance of this is the construction of a rational function $A_v\in\Z(X(\beta))$ for each trivalent vertex of a Demazure weave $\w$. In a nutshell, the first key results regarding cluster algebras on braid varieties can be summarized as follows.

\begin{fthm}\label{thm:cluster_variables_weave_intro2}
Let $\beta\in\brn$ be a positive braid word, $X(\beta)$ its braid variety and $\w:\beta\lr\delta(\beta)$ a Demazure weave with $d$ trivalent vertices $v_1,\ldots,v_d$. Then the following hold:

\begin{enumerate}[label=$(\arabic*)$]
    \item For every trivalent vertex $v\in\w$, the rational function $A_v\in\Z(X(\beta))$ is a regular function. That is, $A_v\in\Z[X(\beta)]$, so $A_v:X(\beta)\lr\mathbb{A}^1_\Z$ has no poles.\\

    \item For every trivalent vertex $v\in\w$, the rational function $A_v\in\Z[X(\beta)]$ is non-vanishing if and only if there is a non-vanishing weight for the Lusztig cycle $\g_v$ at an edge of the south boundary of $\w$.\\

    \item The restrictions of $A_{v_1},\ldots,A_{v_d}$ to the weave torus $T_\w\sse X(\beta)$ are a coordinate system, i.e.~
    $$T_\w\cong \Spec\Z[A_{v_1}^{\pm1},\ldots,A_{v_d}^{\pm1}].$$
    That is $T_\w\sse X(\beta)$ is cut-out as an open subvariety as the non-vanishing of all the regular functions $A_{v_i}$.\\

    \item If $\w_1\simeq\w_2$ is a weave equivalence, the two sets of regular functions $\{A^{\w_1}_{v_1},\ldots,A^{\w_1}_{v_d}\}$ and $\{A^{\w_2}_{v_1},\ldots,A^{\w_2}_{v_d}\}$ coincide.\\

    \item If $\w_1\stackrel{\mu}{\sim}\w_2$ is a weave mutation at the Lusztig cycle $\gamma_{v_k}$ of $\w_1$, then $A^{\w_1}_{v_i}=A^{\w_2}_{v_i}$ if $i\neq k$ and
    \begin{equation}\label{eq:cluster_mutation_weave_variables}
        \displaystyle A^{\w_1}_{v_k}A^{\w_2}_{v_k}=\prod_{\langle \g_i,\g_k\rangle>0} (A^{\w_1}_{v_i})^{\langle \g_i,\g_k\rangle}+\prod_{\langle \g_i,\g_k\rangle<0} (A^{\w_1}_{v_i})^{-\langle \g_i,\g_k\rangle}.
    \end{equation}

    \item There exists an isomorphism
    \begin{equation}\label{eq:braid_variety_is_cluster2}
        \Z[X(\beta)]\cong\SA(Q_\w)
    \end{equation}
    where $\SA(Q_\w)$ denotes the cluster algebra associated to the quiver $Q_\w$. In addition, the tuple $(Q_\w,T_\w,\{A^{\w}_{v_1},\ldots,A^{\w}_{v_d}\})$ is a cluster seed. In particular, any two Demazure weaves for $\beta$ provide two cluster seeds for the same cluster algebra structure in \eqref{eq:braid_variety_is_cluster2}.\\
\end{enumerate}
\end{fthm}
\vspace{0.2cm}

In particular, a Demazure weave $\wb$ provides a cluster $\SA$-seed $(T^\SA_\w,Q_\w,\{A_i^\w\})$ so that $A_i^\w\in\Z[X(\beta)]$, and $Q_\w$ encodes the intersections of dual Lusztig cycles. Similarly, a Demazure weave $\wb$ also provides a cluster $\SX$-seed $(T^\SX_\w,Q_\w,\{X_i^\w\})$ so that $X_i^\w\in\Z(X(\beta))$, $Q_\w$ encodes the intersections of Lusztig cycles, and $X(\beta)$ is the affinization of the cluster $\SX$-scheme associated to $Q_\w$. See \Cref{ssec:A_and_X_schemes_weave} for details on the relation between the cluster $\SA$ and $\SX$-schemes and $X(\beta)$.\\

\noindent We also discuss additional results, including the introduction of Reeb coordinates for $X(\beta)$. For instance, several geometric features of $X(\beta)$, including a defining ideal, a torus action and a cyclic symmetry, an algebraic closed $2$-form compatible with the cluster structures, and a natural Landau-Ginzburg potential associated to a positive braid word $\beta$, can all be expressed explicitly in terms of Reeb coordinates of $X(\beta)$. The goal of Part I of these notes is to provide the necessary details and scientific context for all these notions and results above.\\

%%%%%%%%%%%%%%%%%%%%%%%%%%%%%%%%%%%%%
%%%%%%%%%%%%%%%%%%%%%%%%%%%%%%%%%%%%%
%%%%%%%%%%%%%%%%%%%%%%%%%%%%%%%%%%%%%
\noindent {\bf Related themes for Part I.} There are different paths forward building on Theorems \ref{thm:weave_intro1} and \ref{thm:cluster_variables_weave_intro2} above. Some of them are discussed in the notes, such as studying Poisson brackets and closed 2-forms as in the work of M.~Gekhtman, M.~Shapiro and A.~Vainshtein, cf.~\cite{GSV,GSV05}, the natural superpotentials associated to a positive braid word, cf.~\Cref{thm:superpotential_Reeb_coordinates}, and finding combinatorial expansion formulas and tropical coordinates for cluster variables in terms of a given weave. The lecture notes do not discuss other equally interesting lines of continuing research. A few of these include:\\

\begin{enumerate}
    \item The combinatorial study of quivers $Q_\w$ associated to Demazure weaves $\w:\beta\lr\delta(\beta)$. It is proven in \cite{CGGLSS25} that these quivers admit reddening sequences and the associated cluster algebras are locally acyclic. But, for instance, it remains open to understand whether these quivers $Q_\w$ are either Banff or Louise.\\

    \item The additive categorification of the cluster algebra $\Z[X(\beta)]$ using weaves. This has been achieved in \cite{CasalsChrist} by constructing a perverse schober $\mathcal{F}_\w$ for each Demazure weave, and using its category of global sections $R\Gamma(\mathcal{F}_\w)$. The details of the comparison with the more classical additive categorifications in the reduced case, see e.g.~\cite{MR2242628}, are the subject of ongoing work. In particular, there is a number of interesting new structures, from the categorical and homological viewpoint, that arise when considering such categorifications, cf.~Table \ref{table:intro_homological}. A salient instance is the use of semi-twists $T^{\pm,\geq0}$ appearing as:
$$T^{+,\geq 0}_Y(X)\coloneqq \on{cof}(\tau_{\geq 0}(\on{Mor}(Y,X))\otimes Y\to X)$$
$$ T^{-,\geq0}_Y(X)\coloneqq \on{fib}(X\to \tau_{\geq0}(\on{Mor}(X,Y)^*)\otimes Y),$$
which can be understood as truncated versions of the more standard (often spherical) twists. These semi-twists preserve rigidity, behave remarkably well in the derived categorical setting, and seem to categorically capture some of the key properties of cluster algebras. The diagrammatic calculus of Demazure weaves provides an abundance of interesting examples where simple-minded collections, silting collections and their associated cluster tilting objects can be manipulated and understood in terms of weave diagrams.\\

    \item The monoidal categorification of the cluster algebra $\Z[X(\beta)]$ using weaves. In a sense, this has been partially achieved for double inductive weaves in \cite{KKOP2025MonoidalIII,Bi2026Quantization,HuhJungKimPark2026Comparison}. The general case is ongoing work, especially in relation to general Demazure weaves and $\beta$-Lusztig data.\\

    \item A systematic study of positivity aiming to address braid varieties $X(\beta)$ has been started in \cite{BaoHe2022TwistedPositivity}. A precise identification of Lusztig cycles and their duals as tropical points, building on a positive $\R_{\geq0}$-atlas still remains to be explored.\\
\end{enumerate}
In general, as hinted in the scientific context above, the established results and ongoing work seem to point towards a rather well-structured combinatorial and algebraic theory in and around braid varieties $X(\beta)$, as a non-reduced generalization of many classical Lie theoretic results. It is perhaps a sign of health for a topic in a conference titled {\it Current developments in Mathematics} that some of its natural follow-ups continue to develop, though that prevents us from providing detailed accounts of them for the time being.\\

\small
\begin{table}[htbp]
  \centering
  \caption{Homological objects featuring in the additive categorification of the cluster algebra $\Z[X(\beta)]$. Among other consequences, the existence of such structures in the derived setting hints towards a derived analogue for highest weight structures in Abelian categories.}\label{table:intro_homological}
\begin{tabular}{ |P{3cm}|P{5cm}|P{4cm}|P{5cm}|}
	\hline
    
	{\bf Highest weight abelian category} & {\bf stable $k$-linear $\infty$-category} & {\bf Expression from $(\Delta_1,\ldots,\Delta_m)$} & {\bf Categorical weave calculus for $R\Gamma(\mathcal{F}_\w)$}\\
	\hline

    standard objects & exceptional collection $\{\Delta_i\}$ & $\Delta_i$ & right weave thimbles\\
	\hline

    \rule{0pt}{0.5cm}costandard objects & coexceptional collection $\{\nabla_i\}$ & $T_{\Delta_1}^+T_{\Delta_2}^+\cdots T_{\Delta_{i-1}}^+(\Delta_i)$ & left weave thimbles\\
	\hline

    \rule{0pt}{0.5cm}simples $S_i$ & simple-minded collection & $T_{\Delta_{1}}^{+,\geq 0}(\cdots (T_{\Delta_{i-1}}^{+,\geq 0}(\Delta_i)))$ & Categorical Lusztig cycles\\
	\hline

    \rule{0pt}{0.5cm}projectives $P_i$ & first silting collection & $T_{\Delta_{m}}^{-,>0}(\cdots (T_{\Delta_{i+1}}^{-,>0}(\Delta_i)))$ & Left-dual categorical cycles\\
    \hline

    \rule{0pt}{0.5cm}injectives $I_i$ & second silting collection & $T_{\nabla_{m}}^{+,<0}(\cdots (T_{\nabla_{i+1}}^{+,<0}(\nabla_i)))$ & Right-dual categorical cycles\\
    \hline

    \rule{0pt}{0.5cm}tilting objects $T_i$ & Ringel dual silting collection & $T^{+,<0}_{\Delta_{m}^{\uparrow}}(\dots (T^{+,<0}_{\Delta_{i+1}^{\uparrow}}(\Delta_i^{\uparrow})))$ & Upward-type Lusztig cycles\\
    \hline
\end{tabular}
\end{table}
%\vspace{0.2cm}
%\noindent\rule{\linewidth}{0.4pt}
\vspace{0.5cm}
\normalsize
%%%%%%%%%%%%%%%%%%%%%%%%%%%%%%%%%%%%%
%%%%%%%%%%%%%%%%%%%%%%%%%%%%%%%%%%%%%
%%%%%%%%%%%%%%%%%%%%%%%%%%%%%%%%%%%%%
\noindent {\bf Scientific context for Part II.} For the Lie group $\G=\SL_n$, and more generally Type A, both braid varieties and Demazure weaves have a geometric model. In fact, this geometric model historically led to the systematic development of weave calculus for any simple algebraic Lie group $\G$, see \cite{CasalsZaslow,CasalsWeng22}. Specifically, for $\G=\SL_n$, a positive braid word $\beta\in\Br^+(\SL_n)$ can be used to construct a Legendrian link $\La_\beta\sse(\R^3,\xi_\std)$, see \Cref{ssec:Legendrian_link_beta}. Such a Legendrian link is a smooth link $\La_\beta\sse\R^3$ which is tangent to a certain 2-plane field, i.e.~to the contact structure $\xi_\std$.\\

A central problem in low-dimensional contact and symplectic topology is the classification of Lagrangian surfaces. In particular, given such a Legendrian link $\La_\beta\sse(\R^3,\xi_\std)$ associated to $\beta$, it is interesting and useful to try to classify all Lagrangian surfaces in $(\D^4,\la_\std)$ whose boundary is fixed to be $\La_\beta$. In precise terms, we aim at classifying all embedded exact Lagrangian surfaces $L\sse(\D^4,\la_\std)$ such that $\dd L\sse\La_\beta$, up to compactly supported Hamiltonian isotopy. These surfaces are often referred to as the Lagrangian fillings of $\La_\beta$, see \Cref{ssec:Lagrangian_fillings}. There is a rich literature constructing invariants for Lagrangian fillings, ranging from numbers to dg algebras and $A_\infty$-categories. To wit, the ruling polynomial of $\La_\beta$ \cite{Rutherford11_HOMFLY,Rutherford06_HOMFLY}, the Legendrian contact dg algebra of $\La_\beta$ \cite{Chekanov,EtnyreNg22_Survey}, similarly the partially wrapped Fukaya category stopped at $\La_\beta$ \cite{GPSCV,GPS18,Syl19}, or the dg category of sheaves with singular support on $\La_\beta$ \cite{KashiwaraSchapira_Book,GKS_Quantization,MR4612528,STZ2017}, are instances of reasonably effective invariants that aim at understanding which Lagrangian surfaces $L$ can bound a given Legendrian link $\La_\beta$. As an example of their relation to braid varieties, see \Cref{ssec:braid_varieties_Legendrian_links}, there exists an isomorphism 

\begin{equation}\label{eq:intro_braid_variety_DGA}
    \Z[X(\beta)]\cong\hH_0(\mathbb{A}^{ab}(\La_\beta,\tt))
\end{equation}

\noindent between the coordinate ring of $X(\beta)$ and the $0$th homology of the Abelianization of the Legendrian contact dg algebra of $\La_\beta$, with one basepoint in $\tt$ per strand of $\beta$.\\

\noindent Intuitively, from this symplectic geometric viewpoint, the $R$-points of the braid variety $X(\beta)$, for a commutative ring $R$, should be thought of as the coarse moduli space of Lagrangian fillings of $\La_\beta$, each endowed with a $\GL_1(R)$-local system, see e.g.~\cite{MR4575870,EtnyreNg22_Survey,EHK} and \Cref{ssec:torus_Lagrangian_filling_from_weave}. As we explain in Part II, a Demazure weave $\wb$ can be used to construct an embedded exact Lagrangian filling $L_\w$ of $\La_\beta$. In fact, the weave itself {\it is} the self-intersection locus of the front for the Legendrian lift of this Lagrangian filling $L_\w$, cf.~\Cref{ssec:Lagrangian_filling_from_AG} and \Cref{sec:weaves_Lagrangian_fillings}. The existence of such a geometric realization of $X(\beta)$ and its Demazure weaves $\wb$ leads to useful insights, among them:

\begin{enumerate}
    \item Any two weaves $\w,\w'$ for which $L_\w$ and $L_{\w'}$ are Hamiltonian isotopic must necessarily give the same seed data in $X(\beta)$. That is, they must define the same tori, quivers and cluster variables. This informs our choice when deciding when two weaves must be declared equivalent.\\

    \item There is a geometric operation on Lagrangian surfaces known as Lagrangian disk surgery, see \cite{Polterovich_Surgery} or \Cref{ssec:cluster_mutation_Lagrangian_filling_from_weave}. It has a universal local model and thus, whenever applied, it should always have the same effect on any other invariants. The choice of local model for a weave mutation $\w\stackrel{\mu}{\sim}{\w'}$ is made so as to model a Lagrangian disk surgery between two Lagrangian fillings $L_\w$ and $L_{\w'}$. At the cluster level, it induces the cluster mutation formulas.\\

    \item In this geometric realization of $\wb$ as an embedded surface $L_\w$, the Lusztig cycles and their duals are geometrically realized by 1-cycles in this surface $L_\w$. The natural $\hH_1$-lattices associated to $L_\w$ then become the $\SA$ and $\SX$ seed lattices, with their intersection form being given by the cap product on $L_\w$, see \Cref{sssec:X_lattice,ssec:torus_Lagrangian_filling_from_weave}. The algebraic duality between seed lattices is geometrically realized as Poincar\'e-Lefschetz duality.\\

    \item Since a point of $X(\beta)$ can be understood as a Lagrangian filling with a $\GL_1(R)$-local system, it is geometrically natural to consider functions on $X(\beta)$ that are extracted from the holonomies of such a local system, cf.~\Cref{ssec:torus_Lagrangian_filling_from_weave,ssec:X_coordinates_Lagrangian_filling_from_weave}. For instance, at a point of $X(\beta)$, each mutable Lusztig cycle $\g_i^\w$ naturally defines a function: the monodromy of the $\GL_1(R)$-local system on $L_\w$ parametrized by that point along the absolute 1-cycle realizing $\g_i^\w$. Furthermore, there is a geometric criterion to decide whether such functions are globally regular, see e.g.~\cite[Corollary 1.2]{CasalsLi24}.\\

    \item The introduction of the Reeb coordinates to study braid varieties $X(\beta)$. These are geometrically motivated, as they directly encode the degree 0 Floer-theoretic generators of the Legendrian contact dg algebra associated to $\La_\beta$.
\end{enumerate}

\noindent We have summarized some of these aspects above in Table \ref{table:PartII_introduction2}. A goal of Part II of these notes is to provide the necessary geometric details and context for these items above, through the lens of low-dimensional contact and symplectic topology.\\

\begin{table}
\begin{center}
\caption{A few key concepts being studied in Part II. These explain the appearance of cluster algebra structures in the study of Lagrangian surfaces in 4-dimensional symplectic topology, both from the Floer-theoretic and sheaf-theoretic viewpoints.}\label{table:PartII_introduction2}
\small
\begin{tabular}{ |P{3.8cm}|P{5.5cm}|P{4.8cm}|  }
	\hline
	\multicolumn{3}{|c|}{ {\bf First contact and symplectic aspects of Demazure weaves}} \\
	\hline

{\bf Weaves} & {\bf 4D Symplectic Topology} & {\bf Cluster theory}\\
\hline

Positive braid word $\beta$ & Legendrian link $\La_\beta$ \& basepoints $\tt$ & Cluster algebra $\hH_0(\mathbb{A}^{ab}(\La_\beta,\tt))$\\
\hline
Demazure weave $\w:\beta\to\delta(\beta)$ & Embedded Lagrangian filling $L_\w$ of Legendrian link $\partial L_\w=\La_\beta$ & Augmentation dg morphism $\tiny{\varepsilon_{L_\w}:\mathbb{A}^{ab}(\La_\beta)\to\Z[H_1(L_\w,\tt;\Z)]}$\\
\hline
Weave equivalence & Ham.~isotopy of $L_\w$ fixing $\La_\beta$ & Identity\\
\hline
Weave mutation & Lagrangian disk surgery on $L_\w$ & Cluster mutation\\
\hline
$\SA$ and $\SX$ lattices & $H_1(L_\w,\tt;\Z)$ $\&$ $H_1(L_\w\setminus\tt,\La_\beta\setminus\tt;\Z)$ & $\SA$ and $\SX$ seed lattices\\
\hline
Intersecting weave cycles and lattice duality & Algebraic intersection form on $L_\w$ and Poincar\'e duality & Quiver $Q_\w$ and duality between seed lattices\\
\hline
Lusztig cycles $\{\g_v\}$ & Basis of $H_1(L_\w,\tt;\Z)$ & Seed basis for $\SX$-lattice\\
\hline
Rational function from Lusztig cycles $\{\g_v\}$ & Holonomies along 1-cycles in $H_1(L_\w,\tt;\Z)$ & cluster $\SX$-variables $\{X_v^\w\}$ belong to $\mbox{Frac}(\hH_0(\mathbb{A}^{ab}(\La_\beta,\tt)))$\\
\hline
Dual Lusztig cycles $\{\r_v\}$ & Basis of $H_1(L_\w\setminus\tt,\La_\beta\setminus\tt;\Z)$ & Seed Basis for $\SA$-lattice\\
\hline
Regular function from dual Lusztig cycles $\{\r_v\}$ & Holonomies along 1-cycles in $H_1(L_\w\setminus\tt,\La_\beta\setminus\tt;\Z)$ & cluster $\SA$-variables $\{A_v^\w\}$ as elements in $\hH_0(\mathbb{A}^{ab}(\La_\beta,\tt))$\\
\hline
\end{tabular}
\normalsize
\end{center}
\end{table}

%%%%%%%%%%%%%%%%%%%%%%%%%%%%%%%%%%%%%
%%%%%%%%%%%%%%%%%%%%%%%%%%%%%%%%%%%%%
%%%%%%%%%%%%%%%%%%%%%%%%%%%%%%%%%%%%%
\noindent {\bf Main Results related to Part II.} The content of Part II consists of several geometric constructions, realizing many of the algebraic structures related to cluster algebras by studying Lagrangian fillings of Legendrian links. That said, we summarize some of the outcomes of the constructions as theorems below, hopefully illustrating some of the advantages of these geometric realizations more succinctly.\\

First, we can use Demazure weaves $\wb$ as a tool to construct many Lagrangian fillings of $\La_\beta$, and use the cluster algebra structure of $\Z[X(\beta)]$ as a tool to distinguish them:

\begin{fthm}\label{thm:symplectic_intro1}
Let $\beta\in\brn$ be a positive braid word, $X(\beta)$ its braid variety for $\G=\SL_n$, and $\La_\beta\sse(\R^3,\xi_\std)$ its associated Legendrian link. Suppose that $\w_1$ and $\w_2$ induce two different cluster seeds on $\Z[X(\beta)]$. Then the Lagrangian fillings $L_{\w_1}$ and $L_{\w_2}$ are not Hamiltonian isotopic.\\
\end{fthm}
\vspace{0.2cm}

\noindent In particular, if any aspect of the cluster seeds for two Demazure weaves $\w_1,\w_2:\beta\lr\delta(\beta)$ is different, then $L_{\w_1}$ must be different than $L_{\w_2}$. Note that the cluster seeds have much more structure than their underlying tori. For instance, it might be a priori hard to distinguish whether $T_{\w_1}$ and $T_{\w_2}$ coincide as subvarieties of $X(\beta)$. It is often simpler to compare the quivers $Q_{\w_1}$ and $Q_{\w_2}$ or the cluster $A$-variables $\{A_i^{\w_1}\}$ and $\{A_i^{\w_2}\}$. For instance, for the latter: one expresses them as polynomials in the Reeb variables and, if one of those polynomials differs, then $L_{\w_1}$ cannot be Hamiltonian isotopic to $L_{\w_2}$.\\

Second, at a more foundational level, we can try to understand how much of the geometric classification of Lagrangian fillings of a given Legendrian link $\La_\beta$ is captured by the cluster structure of $\Z[X(\beta)]$. A compelling fact in this direction is the following result, stating that every cluster seed of $\Z[X(\beta)]$ is inhabited by a Lagrangian filling:

\begin{fthm}[A Lagrangian filling for every seed]\label{thm:symplectic_intro2}
Let $\beta\in\brn$ be a positive braid word, $X(\beta)$ its braid variety for $\G=\SL_n$ and $\La_\beta\sse(\R^3,\xi_\std)$ its associated Legendrian link. Suppose that $\mathfrak{s}$ is a cluster seed for $\Z[X(\beta)]$ with torus $T_\mathfrak{s}\sse X(\beta)$. Then there exists an embedded exact Lagrangian filling $L\sse(\D^4,\la_\std)$ of $\La_\beta$ such that
\begin{equation}\label{eq:filling_for_seed}
\fbox{$\mbox{im}(L)=T_\mathfrak{s}$}
\end{equation}
as open toric charts in the braid variety $X(\beta)$.\\
\end{fthm}
\vspace{0.2cm}

Theorem \ref{thm:symplectic_intro2} follows from the main result of \cite{CasalsGao24}. A few comments on Theorem \ref{thm:symplectic_intro2} might be helpful:

\begin{enumerate}
    \item Theorem \ref{thm:symplectic_intro2} is a result that inputs algebraic data, that of a cluster seed for $\Z[X(\beta)]$ and outputs geometric data, that of a Lagrangian filling for $\La_\beta$. Thus far, most results in the area focus on the transition from geometric data to algebraic data, e.g.~the construction of invariants or aiming at obstructing some geometric fact. Theorem \ref{thm:symplectic_intro2} starts to indicate that certain algebraic structures in low-dimensional contact and symplectic topology are rich enough to recover some geometric data.\\

    \item The key step in the proof of Theorem \ref{thm:symplectic_intro2} is introducing a quiver with potential $(Q_\w,W_\w)$ associated to a Demazure weave and showing that its QP-mutation class algebraically captures the different geometric intersection patterns that appear between certain 1-cycles under Lagrangian disk surgery, initially being the Lusztig cycles. At its core, the result proves that these cycles have representatives for which the geometric intersection equals algebraic intersection and this persists under Lagrangian disk surgeries.\\

    \item The notation $\mbox{im}(L)$ in \eqref{eq:filling_for_seed} refers to the torus in $X(\beta)$ induced by a Lagrangian filling of $\La_\beta$ either by using Floer theory, where we use the isomorphism \eqref{eq:intro_braid_variety_DGA}
    and use the associated augmentation
    $$\varepsilon_{L}:\mathbb{A}^{ab}(\La_\beta)\to\Z[H_1(L,\tt;\Z)],$$
    or we use a sheaf quantization of $L$, and identify $X(\beta)$ as part of the moduli of microlocal rank-one sheaves with singular support in $\La_\beta$, cf.~\Cref{thm:braid_varieties_from_Legendrian}.
\end{enumerate}

\noindent A conceptually important consequence of Theorem \ref{thm:symplectic_intro2} is that Legendrian links of the form $\La=\La_\beta$ typically admit infinitely many Lagrangian fillings. Such infinitely many Lagrangian fillings are all smoothly isotopic to each other, relative to their boundaries, and thus this result is purely a symplectic phenomenon. These examples were first discovered in \cite{CasalsGao22}, see also \cite{CasalsGao24,CasalsNg22,GSW2024}. Among others, the following result, obtained by combining the ADE classification for cluster algebras from \cite{FZ2003II} and Theorem \ref{thm:symplectic_intro2}, illustrates this fact:

\begin{corollaryalpha}
Let $\La_{\beta w_0}\sse(\R^3,\xi_\std)$ be a Legendrian link. Suppose that the mutable quiver associated to a Demazure weave $\w:\beta w_0\lr w_0$ is not mutation equivalent to a disjoint union of oriented Dynkin diagrams. Then $\La_{\beta w_0}$ admits infinitely many distinct Lagrangian fillings.\sbs
\end{corollaryalpha}

\noindent An optimistic view of Theorem \ref{thm:symplectic_intro2} is that the cluster structure on $\Z[X(\beta)]$ from Theorem \ref{thm:cluster_variables_weave_intro2} not only helps us classify Lagrangian fillings of $\La_\beta$ but entirely resolves their classification. I suggest the following statement:

\begin{fconj}\label{conj:ADE_classification}
Let $\beta$ be a positive braid word, $X(\beta)$ its braid variety and $\La_\beta\sse(\R^3,\xi_\std)$ its associated Legendrian link. Then the following hold:\\

\begin{enumerate}[label=$(\arabic*)$]
    \item Any embedded exact Lagrangian filling $L\sse(\D^4,\la_\std)$ of $\La_\beta$ induces a cluster seed $\mathfrak{s}(L)$ for the cluster algebra  $\Z[X(\beta)]$.\\

    \item For any cluster seed $\mathfrak{s}_0$ of $\Z[X(\beta)]$ there exists a unique embedded exact Lagrangian filling $L$ such that $\mathfrak{s}(L)=\mathfrak{s}_0$, unique up to compactly supported Hamiltonian isotopy.\\
\end{enumerate}
 In particular, two Lagrangian fillings $L_1,L_2$ are Hamiltonian isotopic relative to $\La_\beta$ if and only if their cluster seeds coincide: $\mathfrak{s}(L_1)=\mathfrak{s}(L_2)$, up to relabeling. That is, \color{blue}{\bf the classification of Lagrangian fillings of $\La_\beta$, up to compactly supported Hamiltonian isotopy, coincides with the classification of cluster seeds in $\Z[X(\beta)]$}\color{black}.\\
\end{fconj}
\vspace{0.2cm}
A few comments on Conjecture \ref{conj:ADE_classification}. First, in view of Theorem \ref{thm:symplectic_intro2}, which is a form of surjectivity, a central remaining aspect of Conjecture \ref{conj:ADE_classification} is injectivity. That is, to establish that two Lagrangian fillings $L_1,L_2$ of $\La_\beta$ which induce the same cluster seed in $X(\beta)$ must necessarily be the same, i.e.~compactly supported Hamiltonian isotopic. The only cases in which this is known to hold are those of $\La_\beta$ being the max-tb Legendrian unknot, by \cite{EliashbergPolterovich96}, and $\La_\beta$ being the standard Legendrian Hopf link, by \cite{thomson2026legendrianhopflinkexactly}.\\

Second, the ADE classification for cluster algebras, as established in \cite{FZ2003II}, translates geometrically via Conjecture \ref{conj:ADE_classification}. Recall that the number of cluster seeds for ADE quivers is finite and computed, see e.g.~\cite[Section 5.9]{FWZ16}. Namely, if Conjecture \ref{conj:ADE_classification} holds then the following consequence must occur:

\begin{fconj}[ADE Classification of Lagrangian fillings]\label{conj:ADE_classification2}
Let $\La\sse(\R^3,\xi_\std)$ be a Legendrian link of the form $\La=\La_{\beta w_0}$ for $\beta$ a positive braid word. Suppose that the mutable quiver is connected. Then the following hold:\\

\begin{enumerate}[label=$(\arabic*)$]
		\item If $\La$ is smoothly isotopic to the link of the $A_{n}$-singularity, then $\La$ has precisely $\frac{1}{n+2}{2n+2\choose n+1}$ exact Lagrangian fillings.\\
		
		\item If $\La$ is smoothly isotopic to the link of the $D_{n}$-singularity, then $\La$ has precisely $\frac{3n-2}{n}{2n-2\choose n-1}$ exact Lagrangian fillings.\\
		
		\item If $\La$ is smoothly isotopic to the link of the $E_6$, $E_7$ or the $E_8$-singularities, then $\La$ has precisely 833, 4160, and 25080 exact Lagrangian fillings, respectively.\\
		
		\item Else, $\La$ has infinitely many exact Lagrangian fillings.\\
\end{enumerate}
\end{fconj}
\vspace{0.2cm}

\noindent Both Conjecture \ref{conj:ADE_classification} and its consequence Conjecture \ref{conj:ADE_classification2} can be understood as generalizations of the nearby Lagrangian conjecture for surfaces. Indeed, they are conjectural classifications of all embedded exact Lagrangians in a Weinstein neighborhood of a (relative) arboreal Lagrangian skeleton, consisting of the union of a Lagrangian filling and its $\L$-compressing disks, which have boundaries in the homology classes dictated by the Lusztig cycles.\\

%%%%%%%%%%%%%%%%%%%%%%%%%%%%%%%%%%%%%
%%%%%%%%%%%%%%%%%%%%%%%%%%%%%%%%%%%%%
%%%%%%%%%%%%%%%%%%%%%%%%%%%%%%%%%%%%%
\begin{center}
	\begin{figure}[h!]
		\centering
		\includegraphics[scale=0.7]{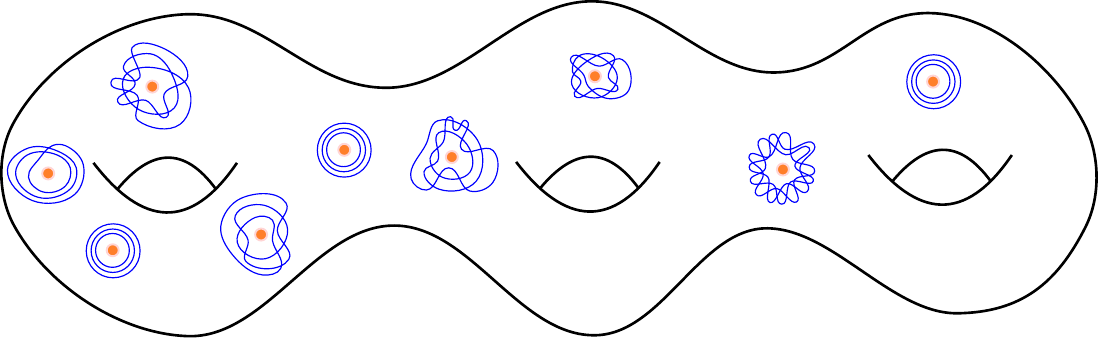}
		\caption{An instance of a Betti surface for $\G=\SL_n$: a smooth surface with marked points endowed with the data of an $n$-stranded positive braid at each puncture. The depicted case is $n=3$. The cylindrical braid diagrams associated to these positive braid words are seen as Stokes diagrams for the local irregular data at that puncture for the irregular $\G$-local system on that surface. The case of concentric circles, i.e.~the trivial braid word, corresponds to tame behavior, whereas non-trivial positive braids lead to wild behavior.}
		\label{fig:BettiSurface}
	\end{figure}
\end{center}

\noindent {\bf Related themes for Part II.} There are different, equally interesting, research directions following the developments surveyed in these notes. A few among them are:\\

\begin{enumerate}
    \item The study of Demazure weaves $\wb$ can be understood as a generalization of the study of Betti moduli spaces for the smooth surface $S^2$ with an irregular singularity at $\infty$ whose Stokes diagram is described by the front given by the cylindrical closure of $\beta$. From that viewpoint, it is natural to explore the role of weaves through the irregular Riemann-Hilbert correspondence, and through the non-Abelian Hodge correspondence. First such steps, employing Floer theory, are established in \cite{casals2025spectralnetworksbettilagrangians}, with a focus on the non-abelianization maps and the relation to spectral networks of \cite{GNMSN,GMN14_Snakes}.\\

    More generally, the diagrammatic weave calculus exists on any smooth surface, cf.~\cite{CasalsZaslow}, and it is natural to try to extend the results for Demazure weaves and their cluster algebras to the more general context of Betti moduli spaces for Betti surfaces, as defined in \cite{casals2025spectralnetworksbettilagrangians} and see \Cref{fig:BettiSurface}. The key outstanding task is a precise generalization of the Lusztig cycles and their duals for (free) weaves on Betti surfaces. In this sense, weave calculus is particularly apt to tackle the generalization of the study of V.~Fock and A.~Goncharov for Betti surfaces with regular singularities, see e.g.~\cite{FockGoncharov09_Dilog2,FG2006,FG06a,Gon17}, or that of generic irregular data, to the case of Betti surfaces with arbitrary wild singularities. In this context, the (cylindrical closure) of the positive braid word associated to each puncture is understood as the Stokes diagram associated to the irregular data for that singularity, cf.~\cite[Section 2]{casals2025spectralnetworksbettilagrangians}.\\

    \item The additive categorification of $\Z[X(\beta)]$ using Demazure weaves, as presented in \cite{CasalsChrist}, can also be geometrically realized within the context of symplectic topology, cf.~Table \ref{tab:weaves_6D_symplectic} and \Cref{eq:NCHMS}. In this case, there is a geometric picture intertwining the symplectic topology and Lagrangian skeleta of certain Weinstein 6-manifolds, whose different flavours of Fukaya categories realize the different categories featuring in the additive categorification, and the symplectic topology of Lagrangian fillings in 4-dimensions. Details of this will appear in \cite{CasalsChrist2}.\\

    \noindent Such a geometric realization leads to a few interesting consequences. To start, using the model for additive cluster categorifications that employs the relative 3-Calabi-Yau Ginzburg dg-algebra ${\bf\Gamma}_{rel}$ of the quiver with potential $(Q_\w,W_\w)$ one obtains a series of equivalences that align with non-commutative mirror symmetry, see \Cref{eq:NCHMS}. These compare different Fukaya categories, to be understood on the A-side, to categories related to the derived category of ${\bf\Gamma}_{rel}$, which can be geometrically understood as being part of the B-side, as a form of noncommutative crepant resolution.\\

\small

\begin{center}
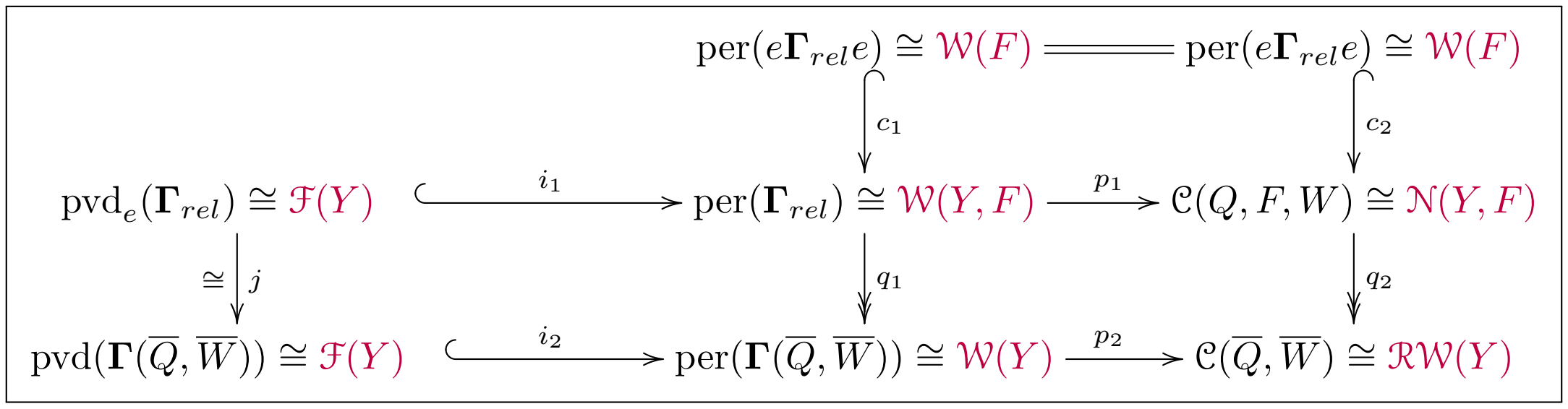

\begin{align*}
		\fbox{\xymatrix{
			&&\mathrm{per}(e\bm\Gamma_{rel}e)\cong \textcolor{purple}{\SW(F)}\ar@{=}[r]\ar@{^{(}->}[d]^{c_1}&\mathrm{per}(e\bm\Gamma_{rel}e)\cong \textcolor{purple}{\SW(F)}\ar@{^{(}->}[d]^{c_2}\\
			\mathrm{pvd}_{e}(\bm\Gamma_{rel})\cong \textcolor{purple}{\mathcal{F}(Y)\quad}\ar@{^{(}->}^{i_1}"2,3"\ar@{->}_{\cong}[d]^{j}&&\mathrm{per}(\bm\Gamma_{rel})\cong\textcolor{purple}{\SW(Y,F)}\ar[r]^{p_1}\ar@{->>}[d]^{q_1}&\mathcal{C}(Q,F,W)\cong\textcolor{purple}{\mathcal{N}(Y,F)}\ar@{->>}[d]^{q_2}\\
			\mathrm{pvd}({\bm\Gamma}(\overline{Q},\overline{W}))\cong \textcolor{purple}{\mathcal{F}(Y)\quad}\ar@{^{(}->}^{i_2}"3,3"&&\mathrm{per}({\bm\Gamma}(\overline{Q},\overline{W}))\cong\textcolor{purple}{\SW(Y)}\ar[r]^{p_2}&\mathcal{C}(\overline{Q},\overline{W})\cong \textcolor{purple}{\RW(Y)}
		}}
\end{align*}\captionof{figure}{A diagram relating categories built from the relative 3-Calabi-Yau Ginzburg dg-algebra ${\bf\Gamma}_{rel}$, in black, and their symplectic geometric realizations, employing flavours of Fukaya categories, in red.\\}
    \label{eq:NCHMS}
\end{center}
\normalsize

    \item To our knowledge, it remains an interesting question to construct a symplectic topological realization of the monoidal categorifications in \cite{KKOP2025MonoidalIII}, perhaps starting with a symplectic geometric realization of the Hernandez-Leclerc category and the associated braid group symmetries. Their direct relation to modules over KLR algebras and, in turn, the realization of such algebras in terms of endomorphisms of certain perverse sheaves, are first hints towards the fact that a symplectic geometric realization likely should exist.\\

    \item Similarly, many of the results connecting symplectic topology and cluster algebras, some illustrated in these notes, hint toward a well-structured theory of {\it positive} symplectic topology. It would be enlightening to better formalize how such positive properties of algebraic invariants arising in Floer theory (or microlocal sheaf theory) can be justified geometrically. For instance, one could seek a natural geometric explanation for the semi-twist functors.\\

    \item Finally, at least for $\G=\SL_n$, the symplectic topological constructions provide algebraic structures much richer than a commutative algebra, including noncommutative versions of braid varieties, parabolic ones, and $A_\infty$ or dg categories for which $X(\beta)$ appears as part of the moduli of pseudoperfect objects. There is not yet a completely agreed upon notion of non-commutative cluster algebra, or the notion of a cluster derived stack, or the concept of a cluster $A_\infty$ or dg category, so that their pseudoperfect moduli stack is a cluster derived stack. That said, as a personal view, any of the Legendrian dg algebra $\mathbb{A}(\La_\beta)$, or its dg category of perfect modules, or equivalently the partially wrapped Fukaya category of $\La_\beta$, or the dg category of sheaves in $\R^2$ with singular support in $\La_\beta$, should be salient examples of these to-be-defined more general notions of being cluster.\\
\end{enumerate}
\normalsize

We conclude with an explanation of the concepts featuring in \Cref{eq:NCHMS}. The black part of the diagram (the left hand side of the equivalences) focuses on the representation theory of the relative Ginzburg dg-algebra $\bm\Gamma_{rel}:=\bm\Gamma_{rel}(Q,W)$, which could be understood as a form of non-commutative crepant resolution on the $B$-side of this non-commutative form of homological mirror symmetry. Here $e\in \bm\Gamma_{rel}$ is the frozen idempotent, $\bm\Gamma(\overline{Q},\overline{W})$ is the absolute CY3 Ginzburg dg-algebra, $\mathcal{C}(Q,F,W)$ the relative cluster category and $\mathcal{C}(\overline{Q},\overline{W})$ the absolute cluster category, cf.~\cite{amiot2011generalizedclustercategories,MR4581969,KellerWu2023Higgs}. The red part of the diagram (the right hand side of the equivalences) focuses on categorical symplectic invariants of the Weinstein 6-manifold $Y=(Y_\w,\la_\w)$ and its associated Lefschetz fibration $\pi:(Y_\w,\la_\w)\lr\cC$, where $F:=F_\pi$ denotes a regular fiber, which is a 4-dimensional Weinstein subdomain. In this viewpoint, the Demazure weave $\w$ can be understood as a form of branch cut for the Lefschetz fibration $\pi$: while a rather excessive one, it is provably useful in structuring many of the symplectic objects involved in the study of $Y=(Y_\w,\la_\w)$.\\

It is important to emphasize that this 6-manifold $(Y_\w,\la_\w)$ also admits a geometric description as a $T^*S^1$-branched cover of $(\C^2,\la_\std)$ branched along the symplectic push-off of $L_\w$. Since the symplectic push-offs of any Demazure weaves $\w:\beta\lr\delta(\beta)$ are symplectically isotopic, the symplectomorphism type of $(Y_\w,\la_\w)$ is independent of $\w$. (From the Lefschetz fibration viewpoint, this can be argued via Hurwitz moves.) This geometric model factoring through 4-dimensions also explains geometrically why these Fukaya categories in 6-dimensions actually additively categorify the cluster algebra associated to the 4-dimensional framework, whose clusters are given by Lagrangian fillings. This 6-dimensional red part of the diagram can be understood as the $A$-side and, quite coarsely put, the cluster varieties extracted from the 4-dimensional setting are closely related (degenerate) to the space of stability conditions of the Fukaya-Seidel category of such $(Y_\w,\la_\w)$. About the red part of the diagram:\\

\begin{enumerate}
    \item (Middle column) The functor $c_1:\SW(F)\lr\SW(Y,F)$ is the Orlov cup functor between the wrapped Fukaya category of $F$, and the partially wrapped Fukaya category of $Y$ with a stop at $F$, which is equivalent to the Fukaya-Seidel category of $\pi$. See e.g.~\cite[Section 2.4]{Syl19}, \cite[Section 3.5]{GPSCV} and \cite[Example 1.23]{GPS18}. The functor $q_1:\SW(Y,F)\lr\SW(Y)$ is the stop removal (pushforward) functor, cf.~\cite[Theorem 1.20]{GPS18}.\\

    \item (Bottom row) The Rabinowitz wrapped Fukaya category $\RW(Y)$ is constructed in \cite[Section 4]{ganatra2022rabinowitzfukayacategoriescategorical}. By \cite[Theorem 1.1]{ganatra2022rabinowitzfukayacategoriescategorical}, it is equivalent to the formal punctured neighborhood of infinity of $\SW(Y)$, and the functor $p_2:\SW(Y)\lr \RW(Y)$ is given in \cite[Lemma 4.11]{ganatra2022rabinowitzfukayacategoriescategorical}. Note that $\RW(Y)$ is an absolute smooth 2-Calabi-Yau category, cf.~\cite[Theorem 1.2.(1)]{bae2024calabiyaustructuresrabinowitzfukaya}, as it should be according to the theory of additive cluster categorifications, cf.~e.g.~\cite[Sections 2\&3]{amiot2011generalizedclustercategories}. The natural inclusion functor $i_2:\mathcal{F}(Y)\lr\SW(Y)$ gives the proper subcategory (of closed compact Lagrangians) of the smooth category $\SW(Y)$, and exactness of the bottom row follows from \cite[Corollary 1.4]{ganatra2022rabinowitzfukayacategoriescategorical}.\\
    
    \item (Right column \& Middle row) The category $\mathcal{N}(Y,F)$ is introduced in \cite[Section 3.5]{gao2024fukayaainfinitystructurenear}, where it is referred to as the negatively wrapped Fukaya category, cf.~also \cite[Section 5]{seidel2025fukayaainftystructuresassociatedlefschetz}. Since $F$ is a full stop, \cite[Theorem 1.1.(ii)]{gao2024fukayaainfinitystructurenear} implies that $\mathcal{N}(Y,F)$ is equivalent to the formal completion of $\SW(Y,F)$ along the subcategory given by the image of $\SW(F)$ under the cup functor. The cup functor thus descends to $c_2:\SW(F)\lr\mathcal{N}(Y,F)$ and the pushforward functor $q_2:\mathcal{N}(Y,F)\lr\RW(Y)$ is given by \cite[Theorem 1.5]{gao2024fukayaainfinitystructurenear}, as $\RW(Y)\cong \mathcal{N}(Y,\emptyset)$. Exactness of the right column follows from \cite[Corollary 1.2]{gao2024fukayaainfinitystructurenear}. Similarly, for the middle row, $p_1:\SW(Y,F)\lr\mathcal{N}(Y,F)$ is given by \cite[(1.12)]{gao2024fukayaainfinitystructurenear}.
\end{enumerate}
\tiny
\begin{center}
\begin{table}
\centering
  \caption{A summary of the key objects featuring in the symplectic geometric realization of the additive cluster categorification of $\Z[X(\beta)]$ via weaves, cf.~\cite{CasalsChrist,CasalsChrist2}, and their relation to the viewpoint presented in these notes, which is focused on Lagrangian fillings in 4-dimensions.}
  \label{tab:weaves_6D_symplectic}
\qquad\begin{tabular}{ |P{3cm}|P{4cm}|P{3cm}|P{3.5cm}|P{3.5cm}|}
	\hline
    
	{\bf Demazure weaves $\w$} & {\bf Sympl.~Topology in $T^*\mathbb{R}^2$} & {\bf $\mbox{Perf}(A_\beta)$ - Rel.~CY2} & {\bf Symplectic 6-folds} & {\bf $R\Gamma(\mathcal{F}_\w)$ - Relative CY3}\\
	\hline
    
	Boundary braid $\beta=\dd\w$ ($\tt$-based) & Legendrian Link $\Lambda_\beta\subseteq T^*_\infty\mathbb{R}^2$ ($\tt$-based, e.g.~Stokes diagram) & Braid variety $X(\beta)$ from Leg.~dg algebra $\mathbb{A}_\beta$ & Liouville sector $(Y^6,F^4)$, $F\cong A_n\mbox{-plumbing of } T^*S^2$ & Relative CY3 category $R\Gamma(\mathcal{F}_\w)$ - saturated\\
	% & (with marked points $T$) & \\
	\hline
    
	Demazure weave $\w:\beta\to\delta(\beta)$ &  Exact Lagrangian filling $L_\w$ of $\La_\beta$ ($\w$ describes front) & Open toric chart $T_L=(\mathbb{C}^\times)^{d}\subseteq X(\beta)$ & Lefschetz fibration $\pi:Y\stackrel{F}{\lr}\cC$ & Exceptional collection of thimbles in $R\Gamma(\mathcal{F}_\w)$\\
	\hline

    Lusztig cycle $\g_i$ on $\w$ (tropical propagation) & Lagrangian disk $D_i\sse T^*\R^2$ with boundary on Lagrangian $L$ & Cluster $\SX$-coordinate as $\dd D_i$-monodromy & Skeleton component (e.g.~Lagrangian $S^3$) & Categorical Lusztig cycle (object of $R\Gamma(\mathcal{F}_\w)$)\\
	
	\hline
    
	  Lusztig cycles & $\lL$-compressing system  & Cluster $\SX$-coordinates & Lagrangian skeleton & Simple-minded collection\\
	
	\hline
    
	Dual Lusztig cycles & Dual $\lL$-compressing system & Cluster $\SA$-coordinates & Cocores of skeleton & Silting collection\\
	
	\hline

%	\hline

	Intersections $\g_i\cap\g_j$ & Intersections of $\dd D_i,\dd D_j$ & Quiver in seed for $X(\beta)$ & Lagr.~intersections & Ext groups \\
	\hline

    Weave potential $W_\w$, via planar regions & $A_\infty$-structure on symplectic push-off $S(L_\w)$ of $L_\w$ & Rigid QP $(Q_\w,W_\w)$ for $\Z[X(\beta)]$ & $A_\infty$-structure in $Y^6$ (e.g.~between 3-spheres) & formal completion along $\mathcal{D}$ as relative cluster category \\

    \hline
\end{tabular}
\end{table}
\end{center}
\normalsize

Similar to Part I, the study of symplectic topological objects associated to $\beta$ and its Demazure weaves seems to provide fertile ground for new insights and ideas, as hopefully illustrated by the instances above. From this viewpoint, we hope that the explicit diagrammatic calculus of weaves and its developing theory can continue to help us enjoy the wonders of contact and symplectic topology, providing new understanding and results within the field, and also bridging across to different areas of mathematics.\\

%%%%%%%%%%%%%%%%%%%%%%%%%%%%%%%%%%%%%
%%%%%%%%%%%%%%%%%%%%%%%%%%%%%%%%%%%%%
%%%%%%%%%%%%%%%%%%%%%%%%%%%%%%%%%%%%%
\noindent {\bf Acknowledgments}: I want to thank the organizers of the Current Developments in Mathematics conference for their invitation. I am also grateful to my coauthors in \cite{CasalsZaslow,CGGLSS25,CasalsChrist,casals2023demazureweavesreducedplabic,CasalsLi24} for valuable discussions. Special thanks are due to J.~Etnyre, H.~Gao, L.~Ng, Y.~J.~Nho and D.~Weng, whose savoir-faire, ideas and support are a source of professional joy. The author is supported by the National Science Foundation under the grant DMS-2505760.\sbs\\

%%%%%%%%%%%%%%%%%%%%%%%%%%%%%%%%%%%%%
%%%%%%%%%%%%%%%%%%%%%%%%%%%%%%%%%%%%%
%%%%%%%%%%%%%%%%%%%%%%%%%%%%%%%%%%%%%

\noindent {\bf Conventions}: The base ring throughout the notes is $\kk=\Z$. In Part II, we have focused on geometric realization of parts of Part I using Floer theory, specifically through the use of the Legendrian dg algebra and Reeb chords. A similar description exists from the viewpoint of the microlocal theory of sheaves, see e.g.~\cite{CasalsZaslow,CasalsWeng22,casals2023demazureweavesreducedplabic}.\sbs\\

%%%%%%%%%%%%%%%%%%%%%%%%%%%%%%%%%%%%%
%%%%%%%%%%%%%%%%%%%%%%%%%%%%%%%%%%%%%
%%%%%%%%%%%%%%%%%%%%%%%%%%%%%%%%%%%%%
\noindent {\bf A comment on references}: Part I articulates many ideas implicit in \cite{CasalsZaslow,CasalsWeng22} and \cite{CGGLSS25}. That said, in order to favor conceptual clarity, the choice of presentation differs from \cite{CGGLSS25}. In particular, I have emphasized the role and importance of dual Lusztig cycles throughout the manuscript and clarified the construction of cluster $\SA$-variables in terms of them. Many of the facts surveyed in this article can be rigorously deduced from the available literature; the few that require more serious work are currently being written up.\sbs

\part{\centering Weave diagrammatics, braid varieties, and cluster structures}

The object of this first part of the notes is to present the diagrammatic calculus of Demazure weaves and develop its applications to the study of braid varieties and their cluster structures. Specifically, a central result of this first part is the construction of a cluster algebra structure on the coordinate ring of a braid variety, cf.~\Cref{table:PartI_introduction}. The construction we present builds on the study of Demazure weaves, which can be depicted as certain diagrams on a surface, see e.g.~\Cref{fig:Example_TwoWeavesSL3SL4}. As we will explain, different moves, operations and combinatorial constructions on these diagrams lead to a fruitful interaction between combinatorics, algebra, and geometry. In particular, the relation between Demazure weaves and cluster structures leads to new results, both in the context of cluster algebras and in areas where braid varieties appear for different reasons, such as representation theory or contact and symplectic topology.

\section{The diagrammatic calculus of Demazure weaves}\label{sec:diagrammatic_weaves}

The diagrammatic calculus of weaves was initiated in \cite{CasalsZaslow}, and has subsequently been developed in a number of additional works, see e.g.~\cite{CGGJ24,CGGLSS25,CasalsWeng22,casals2023demazureweavesreducedplabic,casals2025comparingclusteralgebrasbraid,casals2025spectralnetworksbettilagrangians}. The goal of this first section is to introduce the diagrammatic calculus of Demazure weaves. For the sake of clarity, we focus on Demazure weaves on the 2-disk, and restrict to Type $A$, i.e.~to the special linear group. This particular case already illustrates many of the features and uses of weaves. There exists a more general theory of non-Demazure weaves, on an arbitrary surface, and for any simply-laced simple algebraic Lie group, cf.~\cite[Section 2]{CasalsZaslow} and \cite[Section 4]{CGGLSS25}. Two instances of Demazure weaves on the 2-disk in Type A are depicted in \Cref{fig:weave_intro_example}.\\

Let us thus fix $\mathsf{G}=\SL_n$, consider $\mbox{B}(\mathsf{G})$ the braid group associated to its $A_{n-1}$-Dynkin diagram and $\mbox{B}^+(\mathsf{G})$ the positive braid monoid. We denote by $\Br^+_n=\Br^+(\mathsf{G})$ the set of positive $n$-stranded braid words, i.e.~words on the positive Artin generators $\s_1,\ldots,\s_{n-1}$ of $\mbox{B}(\mathsf{G})$, which serve as representatives for elements in $\mbox{B}^+(\mathsf{G})$.

%%%%%%%%%%%%%%%%%%%%%%%%%%%%%%%%%%%%%%%%%%%%%%%%%%%%%%%%%%%%
%%%%%%%%%%%%%%%%%%%%%%%%%%%%%%%%%%%%%%%%%%%%%%%%%%%%%%%%%%%%
%%%%%%%%%%%%%%%%%%%%%%%%%%%%%%%%%%%%%%%%%%%%%%%%%%%%%%%%%%%%

\subsection{The definition of a Demazure weave}\label{ssec:primer_Demazure_weaves} Let us introduce Demazure weaves as follows:

\begin{definition}[Demazure weaves]\label{def:Demazure_weave}
Let $\beta,\beta'\in\Br^+_n$ be two positive braid words on $n$ strands. By definition, a Demazure weave $\w:\beta\lr\beta'$ is a sequence of positive braid words $\beta\to\beta_1\to\beta_2\to\cdots\to\beta'$, with $\beta_i\in\Br^+_n$, such that each move $\beta_j\to\beta_{j+1}$ is either a braid move in the positive braid monoid or the move $\s_i^2\to \s_i$.\sbs
\end{definition}

\noindent A Demazure weave $\w:\beta\lr\beta'$ where $\beta,\beta'\in\Br^+_n$ will be referred to as an $\SL_n$-weave if we need to emphasize the rank of $\mathsf{G}=\SL_n$. Note that both braid moves and $\s_i^2\to \s_i$ preserve the Demazure product $\delta(\beta)$ of $\beta$, as defined in \cite[Definition 3.1]{KnutsonMiller_04}, which partly justifies the choice of terminology.\\

\begin{center}
	\begin{figure}[h!]
		\centering
		\includegraphics[scale=2.1]{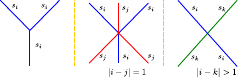}
		\caption{The three types of vertices allowed in the diagrammatic version of a weave $\w$. (Left) Trivalent vertex, where all edges are decorated with the same permutation $\s_i\in S_n$. It represents the move $\s_i^2\to \s_i$. (Center) Hexavalent vertex, with the edge decorations alternating between $\s_i,\s_j$ with $|i-j|=1$. It represents the braid move $\s_i\s_j\s_i\to \s_j\s_i\s_j$. (Right) Tetravalent vertex, with edge decorations as drawn with $|i-k|>1$. It represents the braid move $\s_i\s_k=\s_k\s_i$.}
		\label{fig:weaves_types_vertices}
	\end{figure}
\end{center}

\begin{center}
	\begin{figure}[h!]
		\centering
		\includegraphics[scale=1]{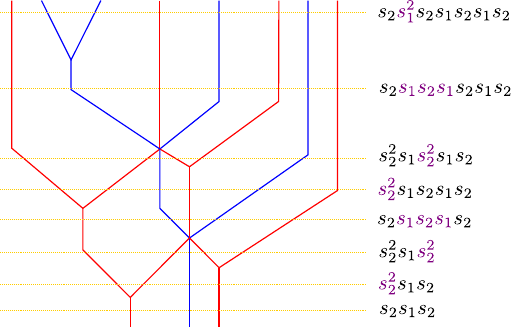}
		\caption{An instance of a Demazure weave $\w:\beta\lr\delta(\beta)$ from the 3-stranded braid word $\beta=\s_2\s_1(\s_1\s_2)^3$ to $\beta'=\s_2\s_1\s_2$, where blue weave edges are labeled by $s_1\in S_3$ and red weave edges by $s_2\in S_3$. Generic horizontal slices are depicted in yellow dashed lines, and the positive braid word associated to each such slice is written to its right. We have highlighted in purple the letters of each braid word to which the next move will be applied when scanning the Demazure weave downward.}
		\label{fig:weave_intro_example}
	\end{figure}
\end{center}

\noindent A Demazure weave, as in \Cref{def:Demazure_weave}, can be depicted diagrammatically as follows. Let us draw a positive braid word $\s_{i_1}\cdots \s_{i_k}$ as a set of parallel vertical strands on the plane, with the $j$th strand from the left labeled by $\s_{i_j}$. Then the braid moves and the move $\s_i^2\to \s_i$ can be depicted as in \cref{fig:weaves_types_vertices}. More generally, \cref{fig:weave_intro_example} depicts an example of a Demazure weave $\w$ diagrammatically. In this diagrammatic presentation, each horizontal slice of a weave yields a braid word $\s_{i_1}\cdots \s_{i_k}$, cf.~\cref{fig:weave_intro_example}. %A core aspect of weave calculus are the Lusztig cycles of a Demazure weave, their intersection theory, and their behavior under weave equivalences and mutations, which we will momentarily summarize.

\noindent When read from top to bottom, the diagram itself is often referred to as a Demazure weave. Considered as a planar graph, a Demazure weave admits the three types of vertices as in \cref{fig:weaves_types_vertices}, and the connected components of the complement of the vertices are said to be weave edges, cf.~\cite[Definition 2.1]{CasalsWeng22} or \cite[Section 2]{CasalsZaslow} for the general definition. The trivalent vertices of a weave $\w$ will be denoted $v_1,\dots,v_m$, ordered top-to-bottom, e.g.~with $v_i$ at height $i$. Each horizontal slice not intersecting a vertex determines a positive braid word by reading the weave edges left-to-right, as illustrated in \cref{fig:weave_intro_example}.\\

\begin{example}
For $\SL_2$, positive 2-stranded braid words are of the form $\beta=\sigma_1^n$ for some $n\in\N$. Thus Demazure weaves $\w:\sigma_1^n\lr\sigma_1$, up to equivalence, are in bijection with planar rooted binary trees, by reading the weave diagram bottom-to-top. For $\SL_n$ with $n\geq3$, the combinatorial study of weaves becomes significantly richer.\sbs
\end{example}

\begin{center}
	\begin{figure}[h!]
		\centering
		\includegraphics[scale=1.2]{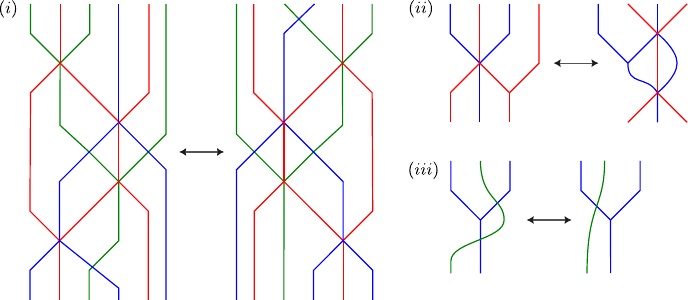}
		\caption{A few weave equivalences, in all of them blue indicates $s_1\in S_4$, red $\s_2$, and green $\s_3$. $(i)$ This weave equivalence encodes a syzygy between braid moves, it is also known as the Zamolodchikov relation. $(ii)$ A weave equivalence known as the push-through move, as it allows for pushing trivalent vertices through hexavalent vertices. $(iii)$ A weave equivalence allowing for distant colors to pass through vertices of each other.}
		\label{fig:weave_equivalences}
	\end{figure}
\end{center}

%%%%%%%%%%%%%%%%%%%%%%%%%%%%%%%%%%%%%%%%%%%%%%%%%%%%%%%%%%%%
%%%%%%%%%%%%%%%%%%%%%%%%%%%%%%%%%%%%%%%%%%%%%%%%%%%%%%%%%%%%

\subsection{Weave equivalence and weave mutation}\label{sssec:weave_equivalences_mutations} There are two important diagrammatic operations on Demazure weaves: weave equivalences and weave mutations. Equivalences between weaves were first discussed in \cite[Theorem 1.1]{CasalsZaslow}, cf.~also \cite[Section 4]{CGGJ24}, and can be summarized as follows:

\begin{definition}[Weave equivalences]\label{def:weave_equivalences} By definition, a weave equivalence is either of the following four types of equivalences:
\begin{enumerate}[label=$(\roman*)$]
    \item Two weaves $\w,\w':\beta\to\beta'$ consisting only of braid moves are declared to be equivalent, see e.g.~\cref{fig:weave_equivalences}.(i).\\
    
    \item If $|i-j|=1$, the weaves $\w:\s_i\s_j\s_i\s_j\to \s_j\s_i\s_j\s_j\to \s_j\s_i\s_j$ and $\w':\s_i\s_j\s_i\s_j\to \s_i\s_i\s_j\s_i\to \s_i\s_j\s_i\to \s_j\s_i\s_j$  are declared to be equivalent, cf.~\cref{fig:weave_equivalences}.(ii).\\
    
    \item If $|i-j|\geq2$, the weaves $\w:\s_i\s_j\s_i\to \s_i\s_i\s_j\to \s_i\s_j\to \s_j\s_i$ and $\s_i\s_j\s_i\to \s_j\s_i\s_i\to \s_j\s_i$ are declared to be equivalent, cf.~\cref{fig:weave_equivalences}.(iii).\\

    \item Two weaves $\w,\w':\beta\to\beta'$ are declared to be equivalent if there exists a compactly supported isotopy of the plane, which is the identity at and above the north boundary of $\w$ and at and below the south boundary of $\w$, which brings $\w$ to $\w'$ through Demazure weaves.\\
\end{enumerate}
\noindent By definition, two weaves $\w,\w':\beta\to\beta'$ are said to be equivalent if there is a sequence of weave equivalences, as in $(i),(ii),(iii)$ and $(iv)$ above, starting at $\w$ and ending at $\w'$. $\sbs$
\end{definition}

\begin{center}
	\begin{figure}[ht]
		\centering
		\includegraphics[scale=1.2]{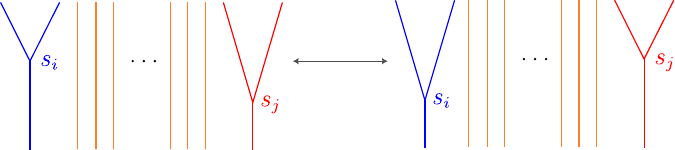}
		\caption{The local weave equivalence that exchanges the height of trivalent vertices, sometimes referred to as a pulley move. The weave lines in color orange indicate arbitrary and possibly different $\s_k$ labels, and the two labels $\s_i,\s_j$ are also arbitrary.\\}
        \label{fig:Weave_HeightChangeEquivalence}
	\end{figure}
\end{center}

\noindent Note that the equivalences in \Cref{def:weave_equivalences}.$(ii)$ and $(iii)$ are parameterized by rank $2$ subdiagrams of types $A_2$ and $A_1\times A_1$ respectively. A compactly supported isotopy as in \Cref{def:weave_equivalences}.$(iv)$ allows, for instance, to exchange the height of trivalent vertices that are far from each other, as depicted in \Cref{fig:Weave_HeightChangeEquivalence}.

\begin{comment}
\begin{example}
{\bf INSERT EXAMPLE OF TWO EQUIVALENT WEAVES HERE}
\end{example}
\end{comment}

There is a second important operation on weaves: weave mutation, which is not a weave equivalence. The notion of weave mutation was introduced in \cite[Section 4.8]{CasalsZaslow}.

\begin{center}
	\begin{figure}[h!]
		\centering
		\includegraphics[scale=1.5]{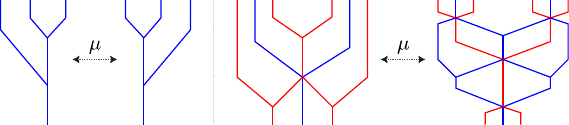}
		\caption{(Left) An elementary weave mutation between the only two Demazure weaves $\w,\w':\s_1^3\to \s_1$. (Right) A more elaborate form of weave mutation, obtained by composing weave equivalences with an elementary weave mutation.}
		\label{fig:weave_mutation}
	\end{figure}
\end{center}

\begin{definition}[Weave mutation]\label{def:weave_mutation}
By definition, the two Demazure weaves $\w,\w':\s_i^3\to \s_i$ are said to be related by an elementary weave mutation, cf.~\cref{fig:weave_mutation} (left). In general, two weaves are said to be related by a weave mutation if they can be related by a sequence of weave equivalences, an elementary weave mutation, and another sequence of weave equivalences.    
\end{definition}

\noindent Two weaves that are mutation equivalent are not typically weave equivalent. For instance, the only two Demazure weaves from $\s_1^3\lr\s_1$, as drawn in \cref{fig:weave_mutation} (left), are not weave equivalent.

\begin{center}
	\begin{figure}[ht!]
		\centering
		\includegraphics[scale=1]{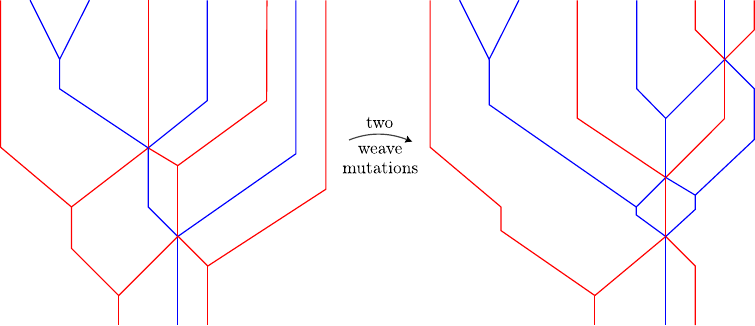}
		\caption{Two instances of Demazure weaves $\w,\w':\beta\lr\beta'$ from $\beta=\s_2\s_1(\s_1\s_2)^3\in\Br^+_3$ to $\beta'=\s_2\s_1\s_2$, for \Cref{ex:weave_mutation}. Here the blue weave edges are labeled by $\s_1$ and red weave edges by $\s_2$.}
		\label{fig:weave_mutation_example}
	\end{figure}
\end{center}

\begin{example}\label{ex:weave_mutation}
(1) \cref{fig:weave_mutation} (right) depicts two weaves $\w,\w':(\s_2\s_1\s_2)^2\to\s_2\s_1\s_2$ that are not equivalent to each other. Nevertheless, $\w,\w'$ are mutation equivalent to each other.\\

\noindent (2) \Cref{fig:weave_mutation_example} depicts two Demazure weaves $\w,\w':\s_2\s_1(\s_1\s_2)^3\lr \s_2\s_1\s_2$. These two weaves $\w,\w'$ are not weave equivalent. Instead they differ from each other by a sequence intertwining weave equivalences and two weave mutations.
\sbs
\end{example}

\noindent The following result is \cite[Theorem 4.6]{CGGJ24}, which establishes the sufficiency of weave equivalences and weave mutation as a set of moves between Demazure weaves with fixed top and bottom boundaries:

\begin{theorem}
Let $\beta$ be a positive braid word. Then any two Demazure weaves $\w,\w':\beta\to\delta(\beta)$ are related by a sequence of equivalence moves and elementary mutations.\sbs
\end{theorem}

%%%%%%%%%%%%%%%%%%%%%%%%%%%%%%%%%%%%%
%%%%%%%%%%%%%%%%%%%%%%%%%%%%%%%%%%%%%
%%%%%%%%%%%%%%%%%%%%%%%%%%%%%%%%%%%%%

\subsection{Weave cycles and weighted braid words}\label{sssec:weave_cycles} Let $\w:\beta\to\beta'$ be a Demazure weave and $E(\w)$ its set of weave edges. The geometric notion of weave cycles was studied in \cite[Section 2]{CasalsZaslow} and \cite[Section 3.2]{CasalsWeng22}, see also \cite[Section 7.4]{CGGLSS25}. The following combinatorial notion of weave cycles was introduced in \cite[Section 4.4]{CGGLSS25}:

\begin{definition}[Weave cycles]\label{def:weave_cycles}
By definition, a weave cycle on $\w$ is a function $\g: E(\w) \to \Z_{\geq 0}$ that assigns a non-negative integer to each edge of the weave. The values of $\g$ are said to be the weights of the edges for the cycle $\g$.\sbs    
\end{definition}

In the same manner that each horizontal slice of a weave yields a braid word $\s_{i_1}\cdots \s_{i_k}$, cf.~\cref{fig:weave_intro_example}, given a weave cycle $\g$ and a horizontal slice we obtain a braid word $\s_{i_1}\cdots \s_{i_k}$ and a collection of non-negative integers ${\bf a}\coloneqq (a_1,\ldots,a_k)$, where $a_j\coloneqq \g(e_{i_j})$ is the weight of the $j$th weave edge $e_{i_j}\in E(\w)$ being intersected in that horizontal slice, corresponding to the letter $\s_{i_j}$ in the braid word $\s_{i_1}\cdots \s_{i_k}$. By definition, we refer to such a pair $(\s_{i_1}\cdots \s_{i_k},{\bf a})$ as a weighted braid word.

\begin{comment}
\begin{example}
{\bf INSERT EXAMPLE OF A WEAVE CYCLE.}
\end{example}
\end{comment}

%%%%%%%%%%%%%%%%%%%%%%%%%%%%%%%%%%%%%
%%%%%%%%%%%%%%%%%%%%%%%%%%%%%%%%%%%%%
%%%%%%%%%%%%%%%%%%%%%%%%%%%%%%%%%%%%%

\subsection{Lusztig cycles}\label{sssec:Lusztigcycles} Of particular importance are the following weave cycles associated to the trivalent vertices of the weave $\w$:

\begin{definition}[Lusztig cycles]\label{def:Lusztig_cycles}
Let $\w$ be a Demazure weave and $p\in\w$ a trivalent vertex. By definition, the Lusztig cycle $\g_p$ associated to $p$ is the unique weave cycle $\g_p: E(\mathfrak{W}) \to \Z_{\ge 0}$ such that:
\begin{itemize}
    \item[$(i)$] $\gamma_p(e_\mathsf{s})=1$ for the unique south edge $e_\mathsf{s}$ at the trivalent vertex $p$, and $\gamma_p(e)=0$ if $e$ is any weave edge that contains any points above $p$.\\
    
    \item[$(ii)$] At any other trivalent vertex with top edges $e_1,e_2$ and south edge $e$, $\g_p$ satisfies
    \begin{equation}\g_p(e)=\min(\gamma_p(e_1),\g_p(e_2))\,.\label{eq:3_move}\end{equation}
    
    \item[$(iii)$] For a 4-valent vertex  with top edges $e_1,e_2$ and  bottom edges $e_1',e_2'$, $\g_p$ satisfies
    \[(\g_p(e'_1),\g_p(e'_2))=(\g_p(e_2),\g_p(e_1))\,.\]
    
    \item[$(iv)$] For a 6-valent vertex  with top edges $e_1,e_2,e_3$ and bottom edges $e_1',e_2',e'_3$, $\g_p$ satisfies
    \begin{align*}
        \g_p(e'_1) & = \g_p(e_2)+\g_p(e_3)-\min(\gamma_p(e_1),\g_p(e_3))\,,\\
        \g_p(e'_2) & = \min(\gamma_p(e_1),\g_p(e_3))\,,\\
        \g_p(e'_3) & = \g_p(e_2)+\g_p(e_1)-\min(\gamma_p(e_1),\g_p(e_3))\,.
    \end{align*}
\end{itemize}
\sbs
%\noindent \qed
\end{definition}

\begin{center}
	\begin{figure}[h!]
		\centering
        \includegraphics[scale=1.6]{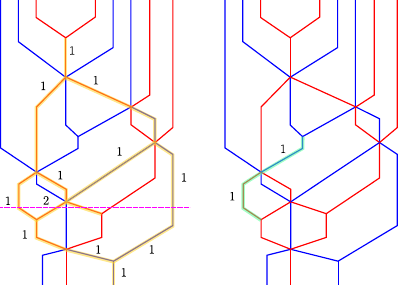}
		\caption{Two instances of Lusztig cycles on a Demazure weave $\w:(\s_1^2\s_2^2)^2\to \s_1\s_2\s_1$. (Left) The Lusztig cycle $\g_1$ associated to the first (top) trivalent vertex. The non-zero weave edges for $\g_1$ are highlighted in yellow, and the specific weight is written next to each weave edge. For instance, at the horizontal pink slice, the weighted braid word is $((\s_2^2\s_1)^2,(1,2,0,1,0,1))$. (Right) The Lusztig cycle $\g_2$ associated to the second trivalent vertex, with the non-zero weave edges for $\g_2$ highlighted in cyan.}\label{fig:weave_Lusztigcycle_example}
	\end{figure}
\end{center}

See \cref{fig:weave_Lusztigcycle_example} for two instances of Lusztig cycles; a non-highlighted edge (with no number) denotes weight 0. \cref{def:Lusztig_cycles} can also be described in terms of weighted braid words as follows. To construct $\g_p$, first consider a horizontal weave slice right below the height of $p$, which yields a certain braid word $\beta(p)\coloneqq \s_{i_1}\cdots \s_{i_{k(p)}}$ and assign the weights ${\bf a}(p)\coloneqq (0,0,\ldots,0,1,0,\ldots,0)$, where the only non-zero weight $1$ is assigned to the weave edge south of $p$. Then, propagate this weighted braid word $(\beta(p),{\bf a}(p))$ downwards by using the following rules for the weights:
\begin{multline}
\label{eq: lusztig trop}
\qquad \qquad \qquad \qquad \quad \Phi_1:(a_1,a_2)\mapsto \min(a_1,a_2), \quad \Phi_2:(a_1,a_2)\mapsto (a_2,a_1),\\ \qquad\qquad \Phi_3: (a_1,a_2,a_3)\mapsto  \left(a_2+a_3-\min(a_1,a_3),\min(a_1,a_3),a_1+a_2-\min(a_1,a_3)\right),
\end{multline}
where $\Phi_1$ is used at trivalent vertices, where the braid word changes according to $\s_i^2\lr \s_i$, $\Phi_2$ at tetravalent vertices, where the braid word changes as $\s_i\s_k=\s_k\s_i$ with $i,k$ non-adjacent, and $\Phi_3$ at hexavalent vertices, where the braid word changes via $\s_i\s_j\s_i=\s_j\s_i\s_j$ with $i,j$ adjacent.

\begin{remark} The downward propagation rules in \cref{eq: lusztig trop} above are referred to as tropical Lusztig rules, as they are obtained by tropicalizing the relations
\begin{multline*}
\label{eq: lusztig trop2}
\qquad \qquad \qquad \qquad \quad x_i(t_1)x_i(t_2)=x_i(t_1+t_2), \quad x_i(t_1)x_k(t_2)=x_k(t_2)x_i(t_1),\\ x_i(t_1)x_j(t_2)x_i(t_3)=
x_j\left(\frac{t_2t_3}{t_1+t_3}\right)x_i(t_1+t_3)x_j\left(\frac{t_1t_2}{t_1+t_3}\right), \hskip 2.8cm
\end{multline*}
\noindent where $x_i(t)\coloneqq \exp(E_{i,i+1}t)$ is the 1-parameter subgroup in $\on{SL}_n$ corresponding to the positive simple root $\alpha_i$, cf.~\cite[Proposition 2.5]{Lusztig}. Here the indices $i,j$ are adjacent, and $i,k$ non-adjacent, in the Dynkin diagram.\sbs
\end{remark}

\noindent For reasons related to cluster algebras, to become clearer later, we introduce the following notation:

\begin{definition}[Mutable and frozen Lusztig cycles]\label{def:mutable_frozen}
A Lusztig cycle $\g$ in a Demazure weave $\w$ is said to be mutable if its weights vanish, i.e.~$\g(e)=0$, on all edges $e\in\w$ of the south boundary of $\w$; it is otherwise said to be frozen.\sbs
\end{definition}
It is important to emphasize that the adjective {\it mutable} in \Cref{def:mutable_frozen} does not refer to the Lusztig cycle itself being able to undergo a weave mutation, as this might not be the case. The adjective refers to the notion of {\it mutable} vertex of a quiver and its corresponding cluster variable.

%%%%%%%%%%%%%%%%%%%%%%%%%%%%%%%%%%%%%
%%%%%%%%%%%%%%%%%%%%%%%%%%%%%%%%%%%%%
%%%%%%%%%%%%%%%%%%%%%%%%%%%%%%%%%%%%%

\subsection{Intersection form on weave cycles}\label{ssec:intersection_cycles} Let us now introduce a skew-symmetric bilinear form on the set of weave cycles. Given two weave cycles $\g,\g': E(\w)\to\Z_{\geq0}$ on a weave $\w$, not necessarily Lusztig cycles, we now define their intersection number $\sharp_{\w}(\g \cdot \g')$ as a sum of local contributions from intersections at the 3-valent vertices, at the 6-valent vertices, and an additional boundary intersection term.

\begin{center}
	\begin{figure}[h!]
		\centering
        \includegraphics[scale=2]{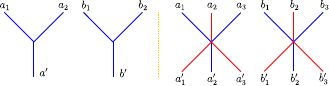}
		\caption{(Left) Two weave cycles near a trivalent vertex, with their weights indicated on each edge. Their local intersection at the trivalent vertex is computed via the determinant in \Cref{def:locint_3valent}. (Right) Two weave cycles near a hexavalent vertex, with their weights indicated on each edge. Their local intersection at this hexavalent vertex is computed via \Cref{def:locint_6valent}.}\label{fig:weaves_types_vertices_trivalent_intersection}
	\end{figure}
\end{center}

\begin{definition}[Local intersection at 3-valent vertex]\label{def:locint_3valent} Let $\w$ be a Demazure weave, $v\in\w$ a trivalent vertex, and $\g,\g': E(\w)\to\Z_{\geq0}$ two weave cycles as in \Cref{fig:weaves_types_vertices_trivalent_intersection} (left). That is, $\g$ $($resp.~$\g')$ has weights $a_1,a_2$ $($resp.~$b_1,b_2)$ on the northwest and northeast edges of the trivalent vertex $v$, respectively, and weight $a'$ $($resp.~$b')$ on the south edge. By definition, the local intersection number of $\g,\g'$ at $v$ is
$$
\sharp_{v}(\g\cdot \g'):=\left|
\begin{matrix}
1 & 1 & 1\\
a_1 & a' & a_2\\
b_1 & b' & b_2\\
\end{matrix}
\right|.
$$
\sbs
\end{definition}

%CHECK SIGNS {Jos\'e: Signs seem correct.}
%\begin{example}
%Suppose that $(a_1,a',a_2)=(0,1,0)$. Then 
%$\sharp(C\cdot \g')=b_2-b_1$.
%\end{example}

\begin{definition}[Local intersection at 6-valent vertex]\label{def:locint_6valent} Let $\w$ be a Demazure weave, $v\in\w$ a hexavalent vertex, and $\g,\g': E(\w)\to\Z_{\geq0}$ two weave cycles as in \Cref{fig:weaves_types_vertices_trivalent_intersection} (right).
That is, $\g$ $($resp.~$\g')$ has weights $a_1,a_2,a_3$ $($resp.~$b_1,b_2,b_3)$ on the northwest, north and northeast edges of the 6-valent vertex $v$, and weights $a'_1,a'_2,a'_3$ $($resp.~$b'_1,b'_2,b'_3)$ on the southwest, south and southeast edges of $v$. By definition, the local intersection number of $\g,\g'$ at $v$ is
$$
\sharp_{v}(\g\cdot \g'):=\frac{1}{2}\left(
\left|
\begin{matrix}
1 & 1 & 1\\
a_1 & a_2 & a_3\\
b_1 & b_2 & b_3\\
\end{matrix}
\right|-
\left|
\begin{matrix}
1 & 1 & 1\\
a'_1 & a'_2 & a'_3\\
b'_1 & b'_2 & b'_3\\
\end{matrix}
\right|
\right).
$$
\sbs
\end{definition}

\begin{center}
	\begin{figure}[h!]
		\centering
        \includegraphics[scale=1]{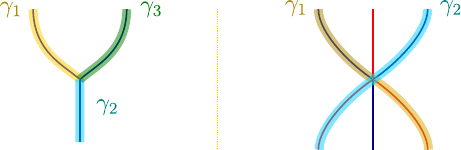}
		\caption{(Left) Three weave cycles $\g_1,\g_2,\g_3$ near a trivalent vertex: they are each declared to have weight 0 on every edge except for the edge highlighted in their corresponding color, where they have weight 1. (Right) Two weave cycles $\g_1,\g_2$ near a hexavalent vertex: each is declared to have weight 0 on every edge except for the edges highlighted in their corresponding color, where they have weight 1. Their intersections are computed in \Cref{ex:local_intersection}.}\label{fig:example_local_intersection}
	\end{figure}
\end{center}

\begin{example}\label{ex:local_intersection} (1) Consider the cycles $\g_1,\g_2,\g_3$ depicted in \Cref{fig:example_local_intersection} (left) near a trivalent vertex $v$. Specifically, $\g_1$ has weight 1 on the northwest edge of the trivalent vertex, and weight 0 elsewhere. Similarly, $\g_2$ has weight 1 on the south edge of the trivalent vertex, and weight 0 elsewhere, and $\g_3$ has weight 1 on the northeast edge of the trivalent vertex, and weight 0 elsewhere. By \cref{def:locint_3valent}, their local intersections are
$$\sharp_{v}(\g_1\cdot \g_2)=1,\quad \sharp_{v}(\g_2\cdot \g_3)=1,\quad\sharp_{v}(\g_3\cdot \g_1)=1.$$

\noindent Note that their local self-intersections $\sharp_{v}(\g_i\cdot \g_i)=0$ at $v$ all vanish, $i\in[1,3]$, as the matrix in \cref{def:locint_3valent} has two identical rows and thus the determinant vanishes.\\

\noindent (2) Consider the cycles $\g_1,\g_2$ depicted in \Cref{fig:example_local_intersection} (right) near a hexavalent vertex $v$. Specifically, $\g_1$ has weight 1 on the northwest and southeast edges of $v$, and weight 0 otherwise, and $\g_2$ has weight 1 on the northeast and southwest edges of  $v$, and weight 0 otherwise. By \Cref{def:locint_6valent}, their local intersection at $v$ is
$$\sharp_{v}(\g_1\cdot \g_2)=-1.$$
\sbs
\end{example}

\begin{definition}[Intersection form on weave cycles]\label{def:intersection_number}
Let $\w$ be a Demazure weave and $\g,\g':E(\w)\to\Z_{\geq0}$ two weave cycles whose weights at the top and bottom edges of $\w$ vanish. By definition, the intersection number $\langle\g, \g'\rangle\in\mathbb{Q}$ of $\g$ and $\g'$ is
\begin{equation}\label{eq:intersection_form}
\langle\g, \g'\rangle := \sum_{\mbox{\tiny $v$ 3-valent}}\sharp_{v}(\g \cdot \g') + \sum_{\mbox{\tiny $v$ 6-valent}}\sharp_{v}(\g\cdot \g').
\end{equation}
\sbs
\end{definition}

\noindent It follows from \Cref{def:locint_3valent} and \Cref{def:locint_6valent} that the intersection number is skew-symmetric:
$$\langle\g, \g'\rangle = -\langle\g', \g\rangle.$$
Also, the intersection number $\langle\g, \g'\rangle$ is an integer when $\g, \g'$ are both mutable Lusztig cycles, as in \Cref{def:Lusztig_cycles}. In particular, the collection of mutable Lusztig cycles endowed with this skew-symmetric intersection form determines an intersection quiver (with no frozen vertices). The intersection numbers from \Cref{def:intersection_number} define a skew-symmetric intersection form, by extending bilinearly, on the lattice of $\Z$-valued weave cycles.

\begin{center}
	\begin{figure}[ht!]
		\centering
		\includegraphics[scale=1]{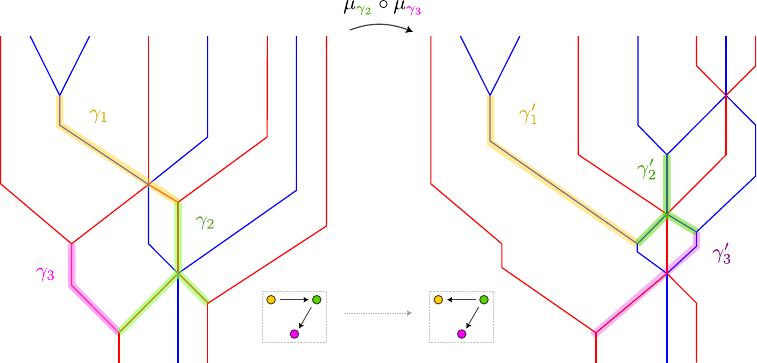}
		\caption{Two instances of Demazure weaves. (Left) The mutable Lusztig cycles $\gamma_1,\gamma_2,\gamma_3$ for $\w$ are highlighted in yellow, green and pink respectively, and their intersection quiver is displayed in the square box at the lower right of $\w$.  (Right) The weave $\w'$ obtained from $\w$ by performing weave mutations along the Lusztig cycles $\gamma_3$ and $\gamma_2$, in this order. Its corresponding Lusztig cycles $\gamma'_1,\gamma'_2,\gamma'_3$ are also highlighted and their intersection quiver is displayed to the lower left of $\w'$.}
		\label{fig:weave_intro_example2}
	\end{figure}
\end{center}

\begin{example}\label{ex:Lusztig_cycles}
Consider the two weaves in \cref{fig:weave_intro_example2}, as in \Cref{ex:weave_mutation}.(2) above. For each of the weaves, \Cref{fig:weave_intro_example2} depicts the three Lusztig cycles associated to the top vertices and their corresponding intersection quivers, computed via \Cref{def:intersection_number}. Specifically, their intersection matrices are
$$(\langle\g_i, \g_j\rangle)=\begin{pmatrix}
0 & 1 & 0\\
-1 & 0 & 1\\
0 & -1 & 0
\end{pmatrix},\qquad (\langle\g'_i, \g'_j\rangle)=\begin{pmatrix}
0 & -1 & 0\\
1 & 0 & 1\\
0 & -1 & 0
\end{pmatrix},\quad \mbox{for }i,j\in[1,3].$$
respectively for the Lusztig cycles $\gamma_1,\gamma_2,\gamma_3$ on the left weave $\w$, and the Lusztig cycles $\gamma'_1,\gamma'_2,\gamma'_3$ on the right weave $\w'$.\sbs
\end{example}
 
\begin{remark}[Boundary intersections]\label{rmk:boundary_intersections}
\Cref{def:locint_3valent} and \Cref{def:locint_6valent} suffice for computing intersections between weave cycles whose weights at the top and bottom edges of the weave are zero, as in \Cref{def:intersection_number}. This suffices for most applications, as it concerns the mutable part of the intersection quiver, cf.~\Cref{def:mutable_frozen}. Nevertheless, there is in general an intersection contribution $\sharp_{\partial}(\g\cdot \g')$ from the top and bottom boundary edges of $\w$ as well, see e.g.~\cite[Section 4.6]{CGGLSS25}. For any two weave cycles $\g,\g'$, possibly not vanishing at the boundary edges of $\w$, we can take into account this boundary intersection contribution for $\g,\g'$ and define their intersection as
\[
\langle\g, \g'\rangle := \sharp_{\partial}(\g\cdot \g')+\sum_{\mbox{\tiny $v$ 3-valent}}\sharp_{v}(\g \cdot \g') + \sum_{\mbox{\tiny $v$ 6-valent}}\sharp_{v}(\g\cdot \g').
\]
\sbs
\end{remark}

\subsection{Dual Lusztig cycles}\label{sssec:dual_Lusztigcycles} Let $\w:\beta\lr\delta(\beta)$ be a Demazure weave. The Lusztig cycles $\{\g_v\}$ of $\w$, indexed by the trivalent vertices of $\w$ as in \Cref{def:Lusztig_cycles}, are an important set of weave cycles of $\w$. Equally relevant is the dual collection $\{\rho_v\}$ of dual Lusztig cycles, dual with respect to the intersection product introduced in \Cref{ssec:intersection_cycles}. In fact, given $\{\g_v\}$ and the skew-symmetric product, it follows from induction and the normalized propagation rules that there exists a distinguished family of weave cycles $\{\rho_v\}$ that vanish at the bottom boundary of $\w$ such that
\begin{equation}\label{eq:duality}
\langle \g_v,\rho_{v'}\rangle =\delta_{v,v'}\,.
\end{equation}

\noindent Specifically, such dual Lusztig cycles $\{\rho_v\}$ for a Demazure weave, which are also indexed by the trivalent vertices of $\w$, can be constructed as follows:

\begin{definition}[Dual Lusztig cycles]\label{def:dual_Lusztig_cycles}
Let $\w$ be a Demazure weave, $\{\g_v\}$ its collection of Lusztig cycles, and $v\in\w$ a trivalent vertex. By definition, the dual Lusztig cycle $\r_v$ associated to $v$ is the unique weave cycle $\r_v: E(\w) \to \Z_{\ge 0}$ such that:\\
\begin{itemize}
    \item[$(i)$] The weights of $\r_v$ at the edges of $v$ are:
    $$\r_v(e_\mathsf{ne})=1,\quad \r_v(e_\mathsf{nw})=0,\quad \r_v(e_\mathsf{s})=0,$$
    where $e_\mathsf{ne}$ denotes the unique northeast edge at the trivalent vertex $v$, and $e_\mathsf{nw}$ and $e_\mathsf{s}$ the unique northwest and south edges at $v$. In addition,
    $$\r_v(e)=0$$
    for any weave edge $e$ that contains any points below $v$.\\
    
    \item[$(ii)$] At a trivalent vertex $p$ different from $v$, where $e_\mathsf{nw}$, $e_\mathsf{ne}$ and $e_\mathsf{s}$ are respectively the northwest, northeast and south edges of $p$, the weights of the weave cycle $
    \rho_v$ are
    \begin{equation}\label{eq:dualLusztigcycle_trivalentrule_left}
    \rho_v(e_\mathsf{nw})=\rho_v(e_\mathsf{s})+\textcolor{green}{\max(-\langle (\rho_v)_{<p},\g_p\rangle,0)},\quad \rho_v(e_\mathsf{ne})=\rho_v(e_\mathsf{s})+\textcolor{blue}{\max(\langle (\rho_v)_{<p},\g_p\rangle,0)},
    \end{equation}
    where $(\rho_v)_{<p}$ denotes the weave cycle obtained by restricting $\rho_v$ to the entire weave below the trivalent vertex $p$.\\
    
    \item[$(iii)$] For a 4-valent vertex  with top edges $e_1,e_2$ and  bottom edges $e_1',e_2'$, $\r_v$ satisfies
    \[(\r_v(e'_1),\r_v(e'_2))=(\r_v(e_2),\r_v(e_1))\,.\]
    
    \item[$(iv)$] For a 6-valent vertex  with top edges $e_1,e_2,e_3$ and bottom edges $e_1',e_2',e'_3$, $\r_v$ satisfies
    \begin{align*}
        \r_v(e'_1) & = \r_v(e_2)+\r_v(e_3)-\min(\r_v(e_1),\r_v(e_3))\,,\\
        \r_v(e'_2) & = \min(\r_v(e_1),\r_v(e_3))\,,\\
        \r_v(e'_3) & = \r_v(e_2)+\r_v(e_1)-\min(\r_v(e_1),\r_v(e_3))\,.
    \end{align*}
\end{itemize}
\sbs
%\noindent \qed
\end{definition}

\noindent The genuine difference between dual Lusztig cycles and Lusztig cycles is in \Cref{def:dual_Lusztig_cycles}.(i) and (ii). We often refer to the terms highlighted in green and blue in \eqref{eq:dualLusztigcycle_trivalentrule_left} as correction terms. Indeed, \Cref{def:dual_Lusztig_cycles}.(iii) and (iv) are identical to the rules of \Cref{def:Lusztig_cycles} at 4-valent and 6-valent vertices. A few comments on \Cref{def:dual_Lusztig_cycles}:

\begin{enumerate}
    \item For each trivalent vertex $v\in\w$, the non-zero weights of the Lusztig cycle $\g_v$ all lie below $v$. In a sense, $\g_v$ propagates downwards: it can be constructed by starting at $v$ and assigning non-zero weights to certain weave edges below $v$ according to \Cref{def:Lusztig_cycles}. In contrast, the non-zero weights of the dual Lusztig cycle $\r_v$ all lie above $v$. Similarly, $\r_v$ can be considered to propagate upwards: it can be constructed by starting at $v$ and assigning non-zero weights to certain weave edges above $v$ according to \Cref{def:dual_Lusztig_cycles}.\\

    \item By construction, specifically \eqref{eq:dualLusztigcycle_trivalentrule_left}, the Lusztig cycles $\{\g_v\}$ and the dual Lusztig cycles $\{\rho_v\}$ satisfy the duality in \eqref{eq:duality},  where we have used the intersection product from \Cref{ssec:intersection_cycles}. Indeed, suppose $v'$ is a trivalent vertex below $v$ and $e_\mathsf{s}$ is the south edge of $v$. Then
    \begin{equation*}
\begin{split}
\langle \g_v,\rho_{v'}\rangle & = \langle \g_v,(\rho_{v'})_{<v}\rangle+\sharp_{v}(\g_v \cdot \r_{v'}) \\
  & \stackrel{\eqref{eq:dualLusztigcycle_trivalentrule_left}}{=} \langle \g_v,(\rho_{v'})_{<v}\rangle+\left|
\tiny\begin{matrix}
1 & 1 & 1\\
0 & 1 & 0\\
\r(e_\mathsf{s})+\textcolor{green}{\max(\langle \g_v,(\rho_{v'})_{<v}\rangle,0)} & \r(e_\mathsf{s}) & \r(e_\mathsf{s})+\textcolor{blue}{\max(-\langle \g_v,(\rho_{v'})_{<v}\rangle,0)}\\
\end{matrix}
\right|\\
\normalsize
 & = \langle \g_v,(\rho_{v'})_{<v}\rangle+\r_{v'}(e_\mathsf{s})+\max(-\langle \g_v,(\rho_{v'})_{<v}\rangle,0) - \r_{v'}(e_\mathsf{s})-\max(\langle \g_v,(\rho_{v'})_{<v}\rangle,0)\\
 & = 0.
\end{split}
\end{equation*}
    
\noindent In the last equality we used that $a=[a]_+-[-a]_+$, where we denote $[a]_+:=\max(a,0)$. The identity \eqref{eq:duality} partially justifies the adjective {\it dual} for dual Lusztig cycles.\\

    \item Note also that \eqref{eq:dualLusztigcycle_trivalentrule_left} implies that $\rho_v$ behaves like a Lusztig cycle near any trivalent vertex $p$ above $v$. Indeed, near such $p$, with edges $e_\mathsf{nw},e_\mathsf{ne}$ and $e_\mathsf{s}$, the dual Lusztig cycle $\rho_v$ satisfies
     $$\rho_v(e_\mathsf{s})=\min(\rho_v(e_\mathsf{nw}), \rho_v(e_\mathsf{ne}))$$
     as at least one of the correction terms in \eqref{eq:dualLusztigcycle_trivalentrule_left} is zero.
\end{enumerate}

\begin{remark}
The notion of ``propagation'' of a weave cycle, either downwards or upwards, as emphasized in $(1)$ above, can be conceptualized in categorical terms by using certain left adjoint induction functors and semitwists, cf.~\cite[Sections 7.3 \& 9.4]{CasalsChrist}. In addition, there is a categorical explanation for the intersection term $\langle (\rho_v)_{<p},\g_p\rangle$ in $\eqref{eq:dualLusztigcycle_trivalentrule_left}$, in terms of Ext groups, cf.~\cite[Section 9.3]{CasalsChrist} and \eqref{eq:y} below.\sbs
\end{remark}

\noindent An interesting combinatorial fact about dual Lusztig cycles is that the northwest correction term in \eqref{eq:dualLusztigcycle_trivalentrule_left} always vanishes:

\begin{proposition}[Right-veering property of dual Lusztig cycles]\label{prop:rightveering} Let $\w:\beta\lr\delta(\beta)$ be a Demazure weave, $v\in\w_3$ a trivalent vertex and $\rho_v$ its dual Lusztig cycle. Then, for any other trivalent vertex $p\in\w_3$, 
$$\rho_v(e_{\mathsf{nw}})=\rho_v(e_{\mathsf{s}})$$
where $e_{\mathsf{nw}},e_{\mathsf{s}}$ denote the northwest and south edges of $p$.\sbs
\end{proposition}

\noindent That is, \Cref{prop:rightveering} is equivalent to the fact that the green correction term in \eqref{eq:dualLusztigcycle_trivalentrule_left} in \Cref{def:dual_Lusztig_cycles}.(ii) always vanishes. In view of \Cref{prop:rightveering}, one can in practice directly substitute \eqref{eq:dualLusztigcycle_trivalentrule_left} by
  \begin{equation}\label{eq:dualLusztigcycle_trivalentrule_left2}
    \rho_v(e_\mathsf{nw})=\rho_v(e_\mathsf{s}),\quad \rho_v(e_\mathsf{ne})=\rho_v(e_\mathsf{s})+\textcolor{blue}{\langle (\rho_v)_{<p},\g_p\rangle}
    \end{equation}
in \Cref{def:dual_Lusztig_cycles}, implicitly understanding that it is a non-trivial fact that $\langle (\rho_q)_{<v},\g_v\rangle$ is always non-negative. A proof of such non-negativity, and thus of \Cref{prop:rightveering}, follows from the equality
\begin{equation}\label{eq:y}
\langle (\rho_q)_{<v},\g_v\rangle = \on{dim}\on{Ext}^1(S_{p},\on{tot}(\mathcal{C}_{\ast\geq b})),
\end{equation}
where the terms on the right hand side of \eqref{eq:y} are explained in \cite[Section 9]{CasalsChrist}. In a nutshell, the right hand side of \eqref{eq:y} is the dimension of an $\on{Ext}^1$ group between certain dg-modules over the derived preprojective algebra, which are used in the additive categorifications introduced in \cite[Theorems A\&B]{CasalsChrist}. In particular, since dimensions are non-negative, it follows that $\langle (\rho_q)_{<v},\g_v\rangle$ is non-negative. In view of \Cref{prop:rightveering}, from now on we adopt \Cref{def:dual_Lusztig_cycles} always using \eqref{eq:dualLusztigcycle_trivalentrule_left2} instead of \eqref{eq:dualLusztigcycle_trivalentrule_left}.

\begin{remark}\label{rmk:rightveering_combinatorial} It would be wonderful to have a combinatorial proof of \Cref{prop:rightveering}. The challenge is that $\langle (\rho_q)_{<v},\g_v\rangle$ is a sum of local contributions and, typically, these contributions are both positive and negative. Thus, the difficulty lies in finding a combinatorial argument for a global cancellation, whereby all the negative terms, which come from some of the hexavalent vertices, actually must be cancelled by the positive terms coming from the trivalent and some of the other hexavalent vertices.\sbs
\end{remark}

\begin{center}
	\begin{figure}[h!]
		\centering
        \includegraphics[scale=1]{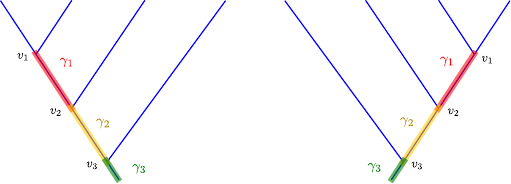}
		\caption{The Lusztig cycles $\{\g_i\}$ for the two Demazure weaves $\s_1^4\lr\s_1$ used in \Cref{ex:dual_Lusztig_cycles}.(1) and (2). Here $\gamma_i$ denotes the Lusztig cycle associated to the trivalent vertex $v_i$, and an edge highlighted with the corresponding color denotes weight 1 for the corresponding Lusztig cycle. Weave edges not highlighted in the corresponding color denote weight 0 for the corresponding Lusztig cycle.}\label{fig:Example_SL2LusztigCycles}
	\end{figure}
\end{center}

\begin{center}
	\begin{figure}[h!]
		\centering
        \includegraphics[scale=1]{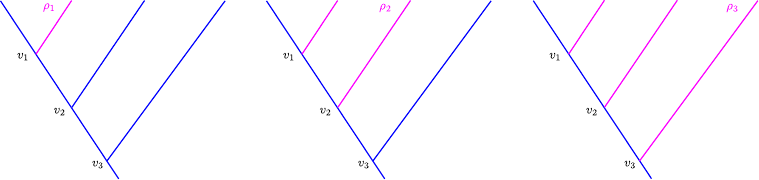}
		\caption{The dual Lusztig cycles $\{\rho_i\}$ for the weave $\w:\s_1^4\lr\s_1$ for \Cref{ex:dual_Lusztig_cycles}.(1). Here $\rho_i$ denotes the dual Lusztig cycle associated to the trivalent vertex $v_i$.}\label{fig:Example_DualLusztigCycles}
	\end{figure}
\end{center}

\begin{example}\label{ex:dual_Lusztig_cycles} (1) Consider the Demazure weave $\w:\beta\to\delta(\beta)$ for $\beta=\s_1^4$ drawn in \cref{fig:Example_DualLusztigCycles}. The Lusztig cycles $\{\g_1,\g_2,\g_3\}$ of $\w$ are depicted in \Cref{fig:Example_SL2LusztigCycles} (left), and the dual Lusztig cycles $\{\rho_1,\rho_2,\rho_3\}$ are depicted in pink in \cref{fig:Example_DualLusztigCycles}, where a pink edge indicates weight 1, and a blue edge indicates weight 0.

\begin{center}
	\begin{figure}[h!]
		\centering
		\includegraphics[scale=1]{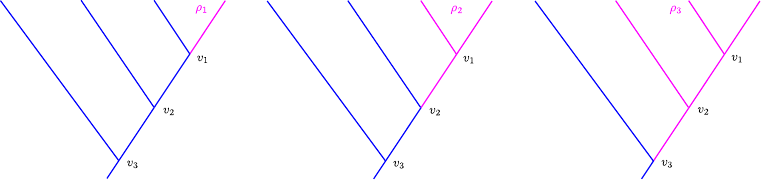}
		\caption{The dual Lusztig weave cycles $\{\rho_i\}$ from \cref{ex:dual_Lusztig_cycles}.(2), highlighted in pink.}\label{fig:Example_DualLusztigCycles2}
	\end{figure}
\end{center}

Suppose that $\{\g_1,\g_2,\g_3\}$ are given, as in \Cref{fig:Example_SL2LusztigCycles} (left). Verifying that $\rho_1$ is indeed the correct dual Lusztig cycle at $v_1$ is immediate from the definition, as only the clause \Cref{def:dual_Lusztig_cycles}.$(i)$ needs to be applied. To justify $\rho_2$ and $\rho_3$ we also need  \Cref{def:dual_Lusztig_cycles}.$(ii)$, and thus we must first compute certain intersections with the Lusztig cycles in order to apply \eqref{eq:dualLusztigcycle_trivalentrule_left}. Specifically, we can construct $\r_2$ by propagating upwards from the northeast edge of $v_2$, with weight 1. Then we scan the weave upwards until we reach the next vertex, which is the trivalent vertex $v_1$. By \eqref{eq:dualLusztigcycle_trivalentrule_left}, the northwest edge of $v_1$ automatically gets assigned weight 0, as this is the weight that was propagating upwards from $\r_2$ on the south edge $e_\mathsf{s}$ of $v_1$ when we reached $v_1$. By \eqref{eq:dualLusztigcycle_trivalentrule_left}, the weight of $\r_2(e_\mathsf{ne})$ assigned to the northeast edge $e_\mathsf{ne}$ of $v_1$ is given by
$$\rho_2(e_\mathsf{ne})=\rho_2(e_\mathsf{s})+\textcolor{blue}{\langle (\rho_2)_{<v_1},\g_1\rangle}=0+1.$$
Here we have used that $\g_1$ and $(\rho_2)_{<v_1}$ intersect only at the trivalent vertex $v_2$, and with contribution +1, as determined in \Cref{ex:local_intersection}.(1). This shows that $\r_2$ is drawn correctly in \cref{fig:Example_DualLusztigCycles}. Iterating this argument also verifies that $\r_3$ is correct as drawn in \cref{fig:Example_DualLusztigCycles}.\\

\noindent (2) Consider the weave $\w:\beta\to\delta(\beta)$ for $\beta=\s_1^4$ drawn in \cref{fig:Example_DualLusztigCycles2}. The Lusztig cycles are depicted in \Cref{fig:Example_SL2LusztigCycles} (right), and the dual Lusztig cycles in \Cref{fig:Example_DualLusztigCycles2}. In this case, all the correction terms $\langle (\rho_q)_{<v},\g_v\rangle$ appearing in \eqref{eq:dualLusztigcycle_trivalentrule_left} vanish. For instance, $\r_2$ starts at $v_2$ and propagates upwards with weight 1 on the south edge $e_\mathsf{s}$ of $v_1$ when we reach $v_1$. By \eqref{eq:dualLusztigcycle_trivalentrule_left}, the northwest edge of $v_1$ is assigned weight 1, and the weight of $\r_2(e_\mathsf{ne})$ assigned to the northeast edge $e_\mathsf{ne}$ of $v_1$ is given by
$$\rho_2(e_\mathsf{ne})=\rho_2(e_\mathsf{s})+\textcolor{blue}{\langle (\rho_2)_{<v_1},\g_1\rangle}=1+0.$$
Here we have used that $\g_1$ and $(\rho_2)_{<v_1}$ intersect only at the trivalent vertex $v_2$ and with local intersection 0 at $v_2$, cf.~ \Cref{ex:local_intersection}.(1).\\

\noindent (3) More generally, consider an arbitrary $\mbox{SL}_2$ Demazure weave, i.e.~a weave $\w:\sigma_1^n\lr\sigma_1$ consisting of one color $\sigma_1$, always depicted in blue in our diagrams. Suppose that we are propagating a dual Lusztig cycle $\r$ upwards and we reach a trivalent vertex. Let $(\b_\downarrow,{\bf a}_\downarrow)$ be the weighted braid word given by $\r$ on a horizontal slice below (and arbitrarily close to) a trivalent vertex $v$. Let us denote by $(\b_\uparrow,{\bf a}_\uparrow)$ the weighted braid word given by $\r$ on a horizontal slice above (and arbitrarily close to) the trivalent vertex $v$. Since $\r$ starts at some trivalent vertex below $v$, we are interested in propagating $\r$ upwards past $v$, using \eqref{eq:dualLusztigcycle_trivalentrule_left}: that is, we seek a formula for $(\b_\uparrow,{\bf a}_\uparrow)$ in terms of $(\b_\downarrow,{\bf a}_\downarrow)$. For the braid words this is just going from $\b_\downarrow=\ldots\s_{i_k}\ldots$ to $\b_\uparrow=\ldots\s_{i_k}\s_{i_k}\ldots$ for the index $i_k$ corresponding to the trivalent $v$, and in particular the lengths satisfy $\ell(\b_{\uparrow})=\ell(\b_{\downarrow})+1$. In this $\mbox{SL}_2$ case, there is also a simple formula near a trivalent vertex $v$ for the weights. Namely, if ${\bf a}_\downarrow=(\ldots,a_{k-1},a_k,a_{k+1},\ldots)\in\N^{\ell(\b_{\downarrow})}$, where $a_k$ is the weight of $\r$ at the south edge $e_\mathsf{s}$ of $v$, then
\begin{equation}\label{eq:upwards_propagation_SL2}
{\bf a}_\uparrow=(\ldots,a_{k-1},a_k,\textcolor{purple}{\max(a_k,a_{k+1})},a_{k+1},\ldots)\in\N^{\ell(\b_{\downarrow})+1}
\end{equation}
where the newly added entry is highlighted in purple. Here we declare $a_{k+1}=0$ if $i_k$ is the rightmost crossing of $\b_\downarrow$, i.e.~if $k=\ell(\b_\downarrow)$. By comparing to \eqref{eq:dualLusztigcycle_trivalentrule_left}, the term $\max(a_k,a_{k+1})$ in \eqref{eq:upwards_propagation_SL2} is obtained via the sum $\max(a_k,a_{k+1})=a_k+\max(a_{k+1}-a_k,0)$, where $a_k=\r(e_\mathsf{s})$ and $\max(a_{k+1}-a_k,0)=\langle \rho_{<v},\g_v\rangle$ is the correction term. The latter equality is obtained via \Cref{def:locint_3valent} and using that $\g_v$ must end at the other trivalent boundary vertex of $e_\mathsf{s}$, which lies below $v$.
\color{black}

\begin{center}
	\begin{figure}[h!]
		\centering
		\includegraphics[scale=0.8]{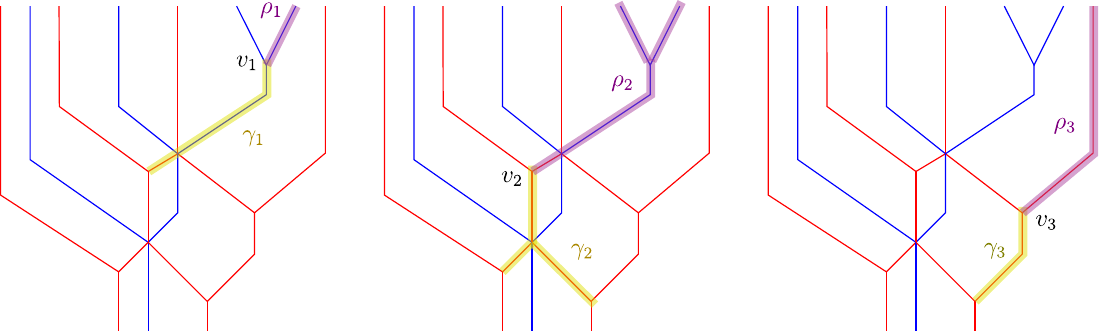}
		\caption{The dual Lusztig cycles $\{\rho_1,\r_2,\r_3\}$ from \cref{ex:dual_Lusztig_cycles}.(4), highlighted in pink, with the corresponding Lusztig cycles $\{\g_1,\g_2,\g_3\}$ highlighted in yellow.}\label{fig:Example_DualLusztigCycles3}
	\end{figure}
\end{center}

\noindent (4) Consider the Demazure weave $\w:\beta\to\delta(\beta)$ depicted in \cref{fig:Example_DualLusztigCycles3}, where $\beta=\s_2\s_1\s_2\s_1\s_2\s_1^2\s_2$. The Lusztig cycles $\g_1,\g_2,\g_3$ associated to the trivalent vertices $v_1,v_2,v_3$ are highlighted in yellow. Specifically, these are all binary weave cycles, with the only values being $0$ and $1$, so we directly highlight the edges where the cycle evaluates to 1. The upwards dual Lusztig cycles $\rho_1,\rho_2,\rho_3$
 are highlighted in purple. Note that there are five Lusztig cycles but $\g_4,\g_5$ have non-zero weights in the bottom boundary. Direct computation also verifies that the corresponding intersections are
$$(\langle\g_i, \rho_j\rangle)=\begin{pmatrix}
1 & 0 & 0\\
0 & 1 & 0\\
0 & 0 & 1
\end{pmatrix},\quad \mbox{for }i,j\in[1,3],$$
with $i=j=1$ being the top left entry and $i=j=3$ the bottom right entry. The only subtle entry is $\langle \g_1,\r_2\rangle=0$, where the intersection is 0 despite the support of the weave cycles geometrically overlapping. Specifically, if $v$ denotes the highest hexavalent vertex, all the three contributions to the intersection number vanish:
$$\langle \g_1,\r_2\rangle=\sharp_{v_1}(\g_1 \cdot \r_2)+\sharp_{v_2}(\g_1 \cdot \r_2)+\sharp_{v}(\g_1 \cdot \r_2)=0+0+0=0.$$
Hence the duality \eqref{eq:duality} holds.\sbs
\end{example}

%%%%%%%%%%%%%%%%%%%%%%%%%%%%%%%%%%%%%
%%%%%%%%%%%%%%%%%%%%%%%%%%%%%%%%%%%%%
%%%%%%%%%%%%%%%%%%%%%%%%%%%%%%%%%%%%%

\subsection{Intersection quiver under weave equivalences and weave mutations}\label{ssec:intersection_Lusztig_cycles_weave_operations} Let $\w,\w':\beta\lr\eta$ be two Demazure weaves, for some braid words $\beta,\eta$. \Cref{def:Lusztig_cycles} provides their collections of Lusztig weave cycles $\{\g_i\}$ and $\{\g'_i\}$, indexed by the corresponding set of trivalent vertices of $\w$ and $\w'$ respectively. \Cref{def:intersection_number}, and see also \Cref{rmk:boundary_intersections}, provides the intersection matrices
$$B_\w=(\langle\g_i,\g_j\rangle),\qquad B_{\w'}=(\langle\g'_i,\g'_j\rangle).$$
The following result explains how $B_\w$ and $B_{\w'}$ are related to each other:

\begin{proposition}\label{prop:intersection_Lusztig_cycles_weave_operations}
Let $\beta,\eta$ be two positive braid words, and $\w,\w':\beta\lr\eta$ two Demazure weaves. Consider the intersection matrices $B_\w$ and $B_{\w'}$ associated to the Lusztig cycles of $\w$ and $\w'$. Then

\begin{enumerate}
    \item If $\w\simeq\w'$ are weave equivalent, there exists a permutation matrix $P$ such that
    \begin{equation}
        B_{\w'}=PB_{\w}P^{-1}.
    \end{equation}
    That is, $B_\w$ and $B_{\w'}$ coincide up to permutation of the indexing of the trivalent vertices.\\

    \item If $\w\stackrel{\mu}{\sim}\w'$ differ by a weave mutation at the $k$th Lusztig cycle, then
\begin{equation}
\langle \gamma'_i,\gamma'_j\rangle =
\begin{cases}
-\langle \gamma_i,\gamma_j\rangle & \text{if } i=k \text{ or } j=k; \\
\langle \gamma_i,\gamma_j\rangle + \langle \gamma_i,\gamma_k\rangle \langle \gamma_k,\gamma_j\rangle & \text{if } \langle \gamma_i,\gamma_k\rangle > 0 \text{ and } \langle \gamma_k,\gamma_j\rangle > 0; \\
\langle \gamma_i,\gamma_j\rangle - \langle \gamma_i,\gamma_k\rangle \langle \gamma_k,\gamma_j\rangle & \text{if } \langle \gamma_i,\gamma_k\rangle < 0 \text{ and } \langle \gamma_k,\gamma_j\rangle < 0; \\
\langle \gamma_i,\gamma_j\rangle & \text{otherwise.}
\end{cases}
\end{equation}
That is, $B_{\w'}$ is obtained by matrix mutation of $B_{\w}$ at the $k$th row and column.
\end{enumerate}
\end{proposition}

\noindent \Cref{prop:intersection_Lusztig_cycles_weave_operations}.(1) is \cite[Corollary 4.31]{CGGLSS25} and \Cref{prop:intersection_Lusztig_cycles_weave_operations}.(2) follows from \cite[Lemma 4.33]{CGGLSS25}. The concept of matrix mutation referred to in \Cref{prop:intersection_Lusztig_cycles_weave_operations}.(2) is explained in \cite[Lemma 2.7.3]{FWZ16}. In the language of quivers, \Cref{prop:intersection_Lusztig_cycles_weave_operations} implies:

\begin{fcor}\label{cor:quiver_cycles_weave_operations}
Let $\beta,\eta$ be two positive braid words, $\w,\w':\beta\lr\eta$ two Demazure weaves. Let $Q_\w$ and $Q_{\w'}$ be the intersection quivers of the Lusztig cycles of $\w$ and $\w'$. Then

\begin{enumerate}
    \item If $\w\simeq\w'$ are weave equivalent, up to vertex relabeling, $Q_{\w'}=Q_\w$.
    
    \item If $\w\stackrel{\mu}{\sim}\w'$ differ by a weave mutation at the $k$th Lusztig cycle, then $Q_{\w'}=\mu_k(Q_\w)$.
\end{enumerate}
\end{fcor}
\vspace{0.2cm}
\noindent In \Cref{cor:quiver_cycles_weave_operations}.(2) we denoted by $\mu_k(Q_\w)$ the mutation of the (ice) quiver $Q_\w$ at its $k$th vertex, which corresponds to the mutable Lusztig cycle $\g_k$. Specifically, $\mu_k(Q_\w)$ is obtained from $Q_\w$ by the following process:
\begin{enumerate}
    \item For each oriented two-arrow $i \to k \to j$ in $Q_\w$ with $i,j$ not both frozen, we add a new arrow $i \to j$.\footnote{If both $i$ and $j$ are frozen, then do nothing.}
    \item Reverse the direction of all arrows in $Q_\w$ incident to the vertex $k$.
    \item Repeatedly remove oriented $2$-cycles, until unable to do so.
\end{enumerate}
\noindent See \cite[Section 2]{FWZ16} for more details on quiver mutation.\\

\noindent The quivers $Q_\w$ obtained as intersection quivers of Lusztig cycles in a Demazure weave $\w$ are reasonably well-behaved. For instance, \cite[Proposition 7.10]{CGGLSS25} implies that, for any Demazure weave $\w$, the quiver $Q_\w$ admits a reddening sequence. The quivers $Q_\w$ are typically of infinite quiver mutation type if $\w$ is an $\SL_n$-weave, $n\geq3$, i.e.~$\beta$ has at least three strands.

%%%%%%%%%%%%%%%%%%%%%%%%%%%%%%%%%%%%%%%%%%%%%%%%%%%%%%%%%%%%
%%%%%%%%%%%%%%%%%%%%%%%%%%%%%%%%%%%%%%%%%%%%%%%%%%%%%%%%%%%%
%%%%%%%%%%%%%%%%%%%%%%%%%%%%%%%%%%%%%%%%%%%%%%%%%%%%%%%%%%%%

\subsection{Lattices}\label{ssec:lattices_seed_data} The $\SA$ and $\SX$-lattices associated to a Demazure weave $\wb$ can be readily named once we introduce braid varieties. Indeed, each Demazure weave $\w$ produces tori $T_\w^\SA$ and $T_\w^\SX$, each embedded in $X(\beta)$, depending on whether it is seen as an $\SA$-scheme or $\SX$-scheme. Then the $\SA$ and $\SX$-lattices for the seed data associated to $\w$ are the character lattices of the respective tori. For now, we emphasize that \textcolor{blue}{the Lusztig cycles span the $\SX$-lattice and the dual Lusztig cycles span the $\SA$-lattice.} The perfect pairing between them, which in particular realizes one as the dual of the other, is given by the intersection form discussed in \Cref{ssec:intersection_cycles}.

%%%%%%%%%%%%%%%%%%%%%%%%%%%%%%%%%%%%%%%%%%%%%%%%%%%%%%%%%%%%
%%%%%%%%%%%%%%%%%%%%%%%%%%%%%%%%%%%%%%%%%%%%%%%%%%%%%%%%%%%%
%%%%%%%%%%%%%%%%%%%%%%%%%%%%%%%%%%%%%%%%%%%%%%%%%%%%%%%%%%%%

\subsection{Demazure weaves from known combinatorics}\label{ssec:Demazure_weaves_from_combinatorics} The literature on algebraic combinatorics and cluster algebras contains many interesting combinatorial models that help manipulate and better understand cluster seeds, e.g.~double wiring diagrams \cite{FZ2000Tests} or $N$-triangulations of surfaces \cite{FG2006}. Central among them are plabic graphs and their alternating strand diagrams, cf.~\cite[Section 11]{Pos}. For instance, the geometry of positroid varieties is intimately connected to the study of reduced plabic graphs, see e.g.~\cite{GL2023,Scott06,SSBW2019}, \cite[Chapter 7]{FWZ16} and references therein. Note also that one can transition from an $N$-triangulation to a plabic graph as in \cite[Section 1.1]{Gon17}.\\

\begin{center}
	\begin{figure}
		\centering
		\includegraphics[scale=1]{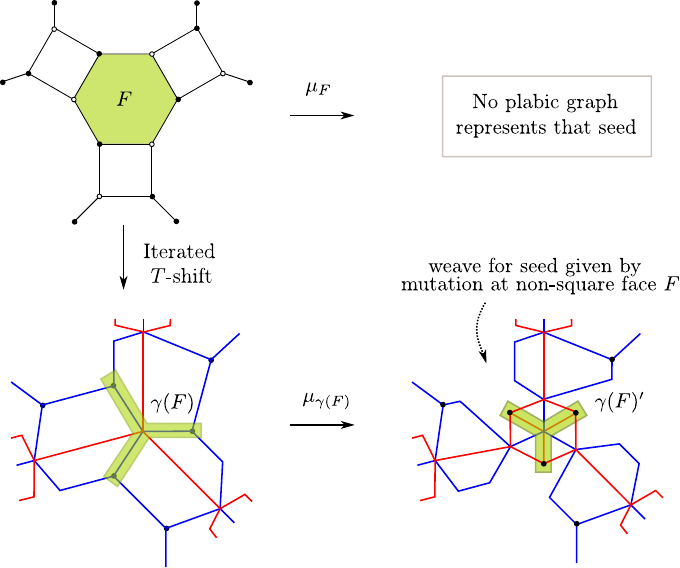}
		\caption{(Upper left) A reduced plabic graph $\bG$ for the top positroid variety in $\Gr(3,6)$. If we label the boundary vertices $1, 2, \dots, 6$ clockwise, with the upper right boundary vertex labeled 1, then the cluster variable for the hexagonal face $F$ is $\Delta_{135}$. (Lower left) The weave $\w(\bG)$ corresponding to this plabic graph under iterated $T$-shifts, cf.~\cite{casals2023demazureweavesreducedplabic}. The highlighted Lusztig cycle $\gamma(F)$ corresponds to the face $F$ of $\bG$, and the seed associated to $\w(\bG)$ and its Lusztig cycles is the target seed of that plabic graph. (Lower right) The weave $\mu_{\gamma(F)}(\w(\bG))$ obtained from mutating the weave $\w(\bG)$ at $\gamma(F)$, which corresponds to a mutation of $\Sigma_T(\bG)$ \emph{not} represented by a plabic graph. The new cycle $\gamma(F)'$, corresponding to the vertex where we have mutated, is highlighted. The cluster variable associated to $\gamma(F)'$ is $\Delta_{136}\Delta_{245} - \Delta_{126} \Delta_{345}$, which is not a Pl\"ucker coordinate.}
		\label{fig:ConjugateSurface10}
	\end{figure}
\end{center}

In general, the idea is that one can transition from a plabic graph to a weave by iteratively applying the $T$-shift construction for plabic graphs, cf.~\cite[Section 8.2]{PSBW}. An instance of this is depicted in \Cref{fig:ConjugateSurface10}. This is explained in detail for reduced plabic graphs in \cite[Section 3]{casals2023demazureweavesreducedplabic}, for plabic fences in \cite[Section 3]{casals2022conjugate} and for the more general grid plabic graphs in \cite[Section 3]{CasalsWeng22}. See also \cite[Section 6.1]{casals2025comparingclusteralgebrasbraid}. Intuitively, the weave is the superposition of the plabic graph $\bG$ and all its $T$-iterates. In all the instances above, this construction naturally bijects the bounded faces of $\bG$ with the mutable Lusztig cycles of the corresponding weave $\w(\bG)$. After appropriately choosing base points, the unbounded faces naturally biject to the frozen Lusztig cycles.

\noindent An advantage of using the weave $\w(\bG)$, even when studying $\bG$, is that it is sometimes possible to mutate at the Lusztig cycle $\g(F)$ associated to a face $F$ of $\bG$ and obtain a new weave. This is the case even if the face $F$ was not square. In particular, several seeds in positroids that are {\it not} Pl\"ucker seeds and do not have a description in terms of plabic graphs, can be diagrammatically described via weaves.

\begin{example}
\Cref{fig:ConjugateSurface10} is an illustration of such a phenomenon, where we are effectively mutating at the hexagonal face of a plabic graph $\bG$ for the top positroid variety in $\Gr(3,6)$ and describing the resulting seed in terms of a weave obtained from a weave mutation applied to $\w(\bG)$. Specifically, the initial seed given by the plabic graph $\bG$ at the top left is a Pl\"ucker seed, and the associated weave $\w(\bG)$ obtained by iterated $T$-shifts is depicted at the bottom left. In this case, the mutation at the quiver vertex associated to the hexagonal face results in a non-Pl\"ucker seed for which no plabic graph exists. Nevertheless, the weave $\w(\bG)$ can be mutated within the diagrammatics of Demazure weaves and we obtain a new weave $\mu_{\gamma(F)}(\w)$ describing this non-Pl\"ucker seed. In fact, in this case of the top positroid variety of $\Gr(3,6)$, all cluster seeds can be represented by weaves, cf.~\cite{AnBaeLee21ADE,Hughes21D}.\sbs
\end{example}

%%%%%%%%%%%%%%%%%%%%%%%%%%%%%%%%%%%%%%%%%%%%%%%%%%%%%%%%%%%%
%%%%%%%%%%%%%%%%%%%%%%%%%%%%%%%%%%%%%%%%%%%%%%%%%%%%%%%%%%%%
%%%%%%%%%%%%%%%%%%%%%%%%%%%%%%%%%%%%%%%%%%%%%%%%%%%%%%%%%%%%

\subsection{Summary of section}\label{ssec:summary_Demazure_weaves} Let $\beta,\eta\in\Br_n^+$ be positive braid words. We have introduced Demazure weaves $\w:\beta\lr\eta$ and their diagrammatic calculus. Such calculus involves:

\begin{enumerate}
    \item The notion of weave equivalences $\w\simeq\w'$ between two weaves $\w,\w':\beta\lr\eta$, cf.~\Cref{def:weave_equivalences}.\\
    \item The notion of weave mutation $\w\stackrel{\mu}{\sim}\w'$ between two weaves $\w,\w':\beta\lr\eta$, cf.~\Cref{def:weave_mutation}.\\
    \item The construction of the Lusztig cycles $\{\g_v\}$ of $\w$, indexed by the trivalent vertices of $\w$, cf.~\Cref{def:Lusztig_cycles}. Similarly, we introduced the dual Lusztig cycles in \Cref{def:dual_Lusztig_cycles}.\\
    \item The construction of a skew-symmetric intersection product for weave cycles, cf.~\Cref{def:intersection_number}. Thus, in particular, the construction of a quiver $Q_\w$ encoding the intersections between the Lusztig cycles $\{\g_v\}$ of $\w$. In addition, \Cref{cor:quiver_cycles_weave_operations} explained how the quiver $Q_\w$ behaves when $\w$ undergoes weave equivalences,  where $Q_\w$ just undergoes a vertex relabeling, and under weave mutation, where $Q_\w$ undergoes a quiver mutation.
\end{enumerate}

\section{Demazure weaves and braid varieties}\label{sec:weaves_braid_varieties}

Let $\beta\in\Br_n^+$ be a positive braid word in $n$-strands. In \Cref{sec:diagrammatic_weaves} we introduced the diagrammatic calculus of Demazure weaves $\w:\beta\to\delta(\beta)$, cf.~\Cref{ssec:summary_Demazure_weaves}. In this section we start to explore a use of such diagrammatic weave calculus: studying the geometry of a certain smooth affine variety $X(\beta)$ associated to $\beta$, known as its braid variety. Braid varieties include and generalize many interesting instances of varieties appearing in Lie theory and parts of algebraic geometry, see \Cref{ssec:braid_varieties_examples}. A central result on such braid varieties, which we approach and establish via Demazure weaves, is that their coordinate ring of regular functions admits the structure of a cluster algebra. The quiver and cluster variables of many cluster seeds for these cluster algebras can be described via the diagrammatic calculus of Demazure weaves.

\begin{comment}
\color{purple} {\bf TO BE EDITED}
\noindent The object of this section is to define braid varieties, provide examples and discuss some properties. The key definition is \cref{def:braid_variety} in \cref{ssec:braidvarieties_def}. Examples are provided in \cref{ssec:braid_varieties_examples}. Invariance properties, including independence of braid word and a torus action, are presented in \cref{ssec:braid_varieties_invariance}. See \cref{lem:braid_invariance} and \cref{lem:torusaction_coordinates}. A useful coordinate system to describe braid varieties is then introduced in \cref{ssec:braid_varieties_coordinates}, cf.~\cref{lem:equations_braidvar}.
\color{black}
\end{comment}
%%%%%%%%%%%%%%%%%%%%%%%%%%%%%%%%%%%%%%%%%%%%%
%%%%%%%%%%%%%%%%%%%%%%%%%%%%%%%%%%%%%%%%%%%%%

\subsection{Definition of braid varieties}\label{ssec:braidvarieties_def} Let $\G$ be a simple algebraic group with Weyl group $W(\G)$. We fix a Borel subgroup  $\borel\sse\G$, a Cartan subgroup $\T\subset \borel$ and consider the flag variety $\G/\borel$. As before, for the sake of the exposition, we focus on $\G=\SL_n$ and thus $W(\SL_n)\simeq S_n$ is the permutation group on $n$ elements. We can also fix $\borel$ to be the subgroup of special upper-triangular matrices and $\T$ the subgroup of special diagonal matrices. In this case, the flag variety $\SL_n/\borel$ has dimension $n(n-1)/2$.\\

Pairs of flags $\borel_1,\borel_2\in \G/\borel$ in relative position $w\in W(\G)$, i.e.~$\borel_1=[g_1\borel]$ and $\borel_2=[g_2\borel]$ such that $g_1^{-1}g_2\in \borel w\borel$, are denoted by $\borel_1 \buildrel w \over \longrightarrow \borel_2$. Braid varieties are defined as follows:\\

\begin{fdef}\label{def:braid_variety}
Let $\beta = \sigma_{i_1}\cdots \sigma_{i_{\ell}}$ be a positive braid word and $\delta(\beta) \in W(\G)$ its Demazure product. The {\it braid variety} associated with $\beta$ is
\begin{equation}\label{eq:braid_variety}
X(\beta) := \{(\borel_1, \dots, \borel_{\ell+1}) \in (\G/\borel)^{\ell + 1} \mid \borel_1 = \borel, \borel_{k} \buildrel s_{i_{k}} \over \longrightarrow \borel_{k+1}, \borel_{\ell + 1} = \delta(\beta) \borel \},
\end{equation}
\noindent where $\delta(\beta) \in W(\G) \cong N_{\G}(\T)/\T$ has been lifted to $N_{\G}(\T)$; this is well-defined since the flag $\delta(\beta) \borel$ does not depend on such a lift.
\end{fdef}
\vspace{0.2cm}
\noindent That is, \textcolor{blue}{$X(\beta)$ parametrizes linear configurations of flags whose consecutive relative positions are dictated by the word $\beta$ and their ends are fixed}: they start with the standard $\borel$ and end with $\delta(\beta)\borel$. A few comments on Definition \ref{def:braid_variety}:
\begin{itemize}
    \item[(i)] The flag variety $\G/\borel$ is a projective variety, and thus so is the Cartesian product $(\G/\borel)^{\ell + 1}$. By choosing appropriate coordinates, we can show that the subset $X(\beta)\sse (\G/\borel)^{\ell + 1}$ is itself an affine variety. It is smooth, irreducible, of dimension $\ell(\beta)-\ell(\delta(\beta))$ and defined over $\Z$, cf.~\Cref{lem:equations_braidvar}. Note that its dimension is independent of the rank of $\G$, e.g.~independent of the number $n$ of strands in $[\beta]\in\Br(\G)$ for $\G=\SL_n$.\\

    \item[(ii)] A braid word $\beta$ is used to define $X(\beta)$, but the isomorphism type of $X(\beta)$ as an affine variety only depends on the braid $[\beta]\in\Br(\G)$. In fact, we will explain that each braid move in $\Br(\G)$, from a braid word $\beta_1$ to another $\beta_2$, with $[\beta_1]=[\beta_2]$, gives an explicit isomorphism $X(\beta_1)\cong X(\beta_2)$. In addition, if two braid words $\beta_1,\beta_2$ differ by cyclic rotation $\s_i\beta\to\beta\s_{n-i}$ and $\delta(\beta_1)=w_0$, we still have an explicit isomorphism $X(\beta_1)\cong X(\beta_2)$, despite the braids $[\beta_1]\neq[\beta_2]$ being different. See \cref{ssec:braid_varieties_invariance} below for details.\\

    \item[(iii)] At this stage, $X(\beta)$ is defined only using the information of $\beta$. As we move forward, we study the geometry of $X(\beta)$ by using different Demazure weaves $\w:\beta\lr\delta(\beta)$. So, to emphasize, Demazure weaves are {\it not} used in the definition of $X(\beta)$.
\end{itemize}

\begin{center}
	\begin{figure}
		\centering
		\includegraphics[width=\textwidth]{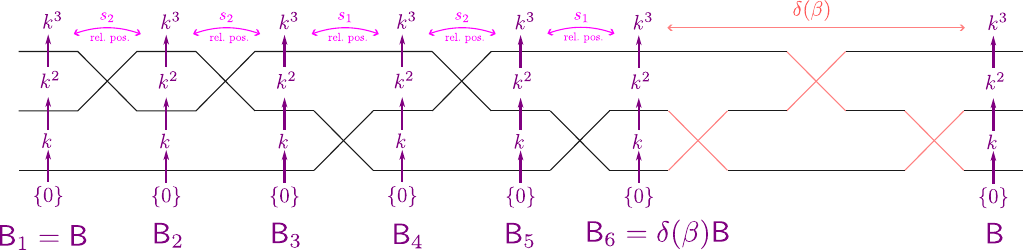}
		\caption{A visualization of the linear configuration of flags parametrized by $X(\beta)$, with $\beta=\sigma_2^2\s_1\s_2\s_1$. The crossings corresponding to $\delta(\beta)=w_0=\s_1\s_2\s_1$ are drawn in orange.}
		\label{fig:BraidVariefyDef1}
	\end{figure}
\end{center}

\begin{remark}\label{rmk:borels_flags} Since we focus on $\G=\SL_n$, we can switch between Borel subgroups of $\SL_n$ and full flags in $k^n$, understood as nested sequences
$$\{0\}\subset k\subset k^2\subset\ldots \subset k^n$$
of linear subspaces of $k^n$.\footnote{Technically this discussion also applies to $\Z$-modules and we use $\Z$ as the ground ring throughout. We focus here on a field $k$ and $k$-vector spaces for clarity, which we can do by base change.} The Borel subgroup that corresponds to a given full flag is given by all the elements $g\in\SL_n$ that leave that particular flag invariant. That is, given the flag $g\borel$, the stabilizer is the Borel subgroup of matrices of the form $g\borel g^{-1}$, i.e.~$\G$ acting by conjugation. For instance, given the standard basis $\{e_i\}$ of $k^n$, the Borel subgroup corresponding to the standard flag
$$\{0\}\subset \langle e_1\rangle\subset \langle e_1,e_2\rangle\subset\ldots \subset \langle e_1,e_2,\ldots,e_n\rangle,$$
is the subgroup $\borel\sse\SL_n$ of upper-triangular matrices. That corresponding to the antistandard flag
$$\{0\}\subset \langle e_n\rangle\subset \langle e_{n-1},e_n\rangle\subset\ldots \subset \langle e_1,e_2,\ldots,e_n\rangle,$$
is the opposite subgroup $\borel_-$ of lower-triangular matrices. In terms of matrices, this corresponds to the fact that $w_0Bw_0$ is a lower-triangular matrix if $B\in\borel$ is upper-triangular, where we used a permutation matrix to lift $w_0$ from $W(\G)$ to $\G$.
\begin{comment}
In terms of matrices, this corresponds to the fact that $w_0Bw_0$ is a lower-triangular matrix if $B\in\borel$ is upper-triangular, where we use the permutation matrix for $w_0$:
$$w_0=\left(\begin{matrix}
0 & 0 & \ldots & 0 & 1\\ 
0 & 0 &  \ldots & 1 & 0 \\
0 &  0 & \ldots & 0 &0 \\ 
\vdots & \vdots & \ddots & \vdots & \vdots \\
1 & 0 & \ldots & 0 & 0
\end{matrix}\right)\in \SL_n.$$
\end{comment}
Under the identification between Borel subgroups of $\SL_n$ and full flags in $k^n$, the $i$-dimensional subspace of the full flag associated to a Borel $\borel\sse\G$ is denoted by $\borel^i$. Then $\borel_1\stackrel{s_i}{\to}\borel_2$ if and only if the associated flags differ {\it exactly} at the $i$ dimensional subspace, i.e.~ $\borel_1^k=\borel_2^k$ if $k\neq i$ and $\borel_1^i\neq\borel_2^i$. The standard and antistandard flags above are in $w_0$-relative position, as transverse as possible: $\borel^j\cap(\borel_-)^{n-j}=\{0\}$ for all $j\in[0,n]$.\sbs
\end{remark}

In this $\G=\SL_n$ case, and with \Cref{rmk:borels_flags}  in mind, it can be helpful to visualize Definition \ref{def:braid_variety} as follows. First, draw horizontally the braid diagram for $\beta$ in the plane $\R^2$ and assign vector spaces in the open regions $\R^2\setminus\beta$ according to the following rule: the unbounded region below $\beta$ has $\{0\}$, the region above has $k^n$ and, in general, a bounded region is assigned the vector space $k^i$ if and only if it has $i$ strands of $\beta$ below it; cf.~Figure \ref{fig:BraidVariefyDef1}. Second, the tuple of flags $(\borel_1, \dots, \borel_{r+1})$ can then be thought of as drawing a vertical segment in between each pair of crossings of $\beta$ and choosing injective linear maps $k^i\lr\ k^{i+1}$ going upwards. This is drawn in Figure \ref{fig:BraidVariefyDef1} in an example, where the flags are emphasized in purple. Third, the two conditions cutting out $X(\beta)\sse (\G/\borel)^{\ell + 1}$ in Definition \ref{def:braid_variety} are of two types:

\begin{enumerate}
    \item {\it Transversality conditions}: for each crossing $\sigma_{i_k}$ of $\beta$, the flags at the left and right of the crossing are in $s_{i_k}$-relative position. That is, all subspaces of the flags at the left and right of the crossing coincide except the $i_k$-dimensional subspace, which is forced to be different, cf.~\cref{rmk:borels_flags}. This transversality condition is indicated in pink in Figure \ref{fig:BraidVariefyDef1}.\\

    \item {\it Boundary conditions}: the first flag is forced to be standard, so $\borel_1=\borel$, and the last one ``as antistandard as possible'', so $\borel_{r + 1} = \delta(\beta) \borel$. If $\delta(\beta)=w_0$, which is an important case, then it is literally the antistandard flag. From a visual viewpoint, these boundary conditions allow us to essentially think of the tuple as being cyclic, so the right part of $\beta$ connects back to its left part in Figure \ref{fig:BraidVariefyDef1}.
\end{enumerate}

\begin{example}\label{ex:permutation_vs_braid} $(1)$ Consider the case where $\beta=s$ is the lift of a permutation $s\in W(\G)$, i.e.~each pair of strands intersects at most once. \cite[Lemma 3.1]{CGGLSS25} shows that given two flags $\FF_1,\FF_2$ in relative position $\FF_1\stackrel{s}{\to}\FF_2$, then the $(\ell(s)-1)$ intermediate flags between $\FF_1$ and $\FF_2$ are uniquely determined by $\FF_1$ and $\FF_2$.\\

\noindent $(2)$ For a general braid word $\beta$, the fact in (1) is not true. Namely, for almost any braid diagram for a $\beta$, with two fixed flags $\FF_1,\FF_2$ left and right of this diagram, the space of intermediate flags between $\FF_1$ and $\FF_2$ is positive dimensional. In fact, $X(\beta)$ is zero-dimensional if and only if $\beta$ is a permutation braid. In this case where $\beta$ is a permutation, $X(\beta)$ is a point. This is a starting point for braid varieties: the transversality conditions for $X(\beta)$, which can be abbreviated $\borel\stackrel{\beta}{\to}\delta(\beta)\borel$, typically define a moduli of tuples of flags of non-zero dimension with rather rich geometry. For more details on relative position of flags, see e.g.~\cite[Section 3.2]{CGGLSS25}.\sbs
\end{example}

\begin{center}
	\begin{figure}
		\centering
		\includegraphics[scale=1]{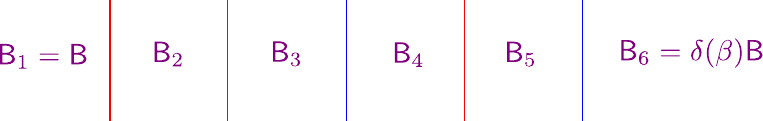}
		\caption{A visualization of the linear configuration of flags parametrized by $X(\beta)$, with $\beta=\sigma_2^2\s_1\s_2\s_1$, which is closer to the diagrammatics of weaves. Here a blue vertical line indicates a $\s_1$-label and a red vertical line indicates a $\s_2$-label. Since we label blue and red weave edges by $\s_1$ and $\s_2$ respectively, two flags which are to the left and right of a vertical blue strand must be in $s_1$-relative position, whereas two flags which are to the left and right of a vertical red strand must be in $s_2$-relative position.}
		\label{fig:BraidVariety_DiagramExample}
	\end{figure}
\end{center}
%\noindent \cref{ex:permutation_vs_braid} is a starting point for braid varieties: the transversality conditions for $X(\beta)$, which can be abbreviated $\borel\stackrel{\beta}{\to}w_0\borel$, typically define a moduli of tuples of flags of non-zero dimension with rather rich geometry. For more details on relative position of flags, see e.g.~\cite[Section 3.2]{CGGLSS25}.\\
\subsubsection{Braid varieties via weave diagrams}\label{sssec:first_connection_weave_diagrams} Given $\beta\in\mbox{Br}_n^+$, and to start connecting with the diagrammatic calculus of Demazure weaves in \Cref{sec:diagrammatic_weaves}, we can also represent points in the braid variety $X(\beta)$ as follows:

\begin{enumerate}
    \item As in \Cref{ssec:primer_Demazure_weaves}, we draw the positive braid word $\beta=\s_{i_1}\cdots \s_{i_\ell}\in \mbox{Br}_n^+$ as a set of $\ell$ parallel vertical strands on the plane, with the $j$th strand from the left labeled by $\s_{i_j}$.\\

    \item In each space in between two such parallel vertical strands we assign a flag in $\G/\borel$.\\

    \item In such a diagrammatic representation, the two flags $\borel_l,\borel_r$ respectively to the left and to the right of a vertical strand corresponding to $\s_{i_j}$ (labeled by $s_{i_j}$) are forced to be in $s_{i_j}$-relative position, i.e.~$\borel_l\stackrel{s_{i_j}}{\to}\borel_r$.\\

    \item The leftmost and rightmost flags are forced to be the standard flag $\borel$ and the $\delta(\beta)\borel$, respectively.
\end{enumerate}

\noindent For instance, \Cref{fig:BraidVariety_DiagramExample} illustrates the configuration of flags in \Cref{fig:BraidVariefyDef1} using these diagrammatic rules. This hints towards a first important relation between braid varieties and Demazure weaves: a Demazure weave $\w:\beta\to\delta(\beta)$ will allow us to study the braid variety $X(\beta)$, where $X(\beta)$ will appear by considering flags along the north boundary of $\w$, as in \Cref{fig:BraidVariety_DiagramExample}.\\

\noindent In general, given a Demazure weave $\w:\beta\to\delta(\beta)$, a horizontal slice $h$ away from a weave vertex determines a positive braid word $\beta_h$. In particular, we can assign the braid variety $X(\beta_h)$ to a horizontal slice $h$. For instance, the braid variety $X(\beta)$ corresponds to a horizontal slice at the north boundary of $\w$. The braid variety corresponding to the horizontal slice at the south boundary is $X(\delta(\beta))\cong\{\mbox{pt}\}$, which is a point by \cref{ex:permutation_vs_braid}.(2).

\begin{remark}
In this sense, if a Demazure weave $\w:\beta\to\delta(\beta)$ can be understood as a way of simplifying $\beta$ down to the permutation $\delta(\beta)$, via braid moves and $\s_i^2\to\s_i$ moves, then the same principle applies to the braid variety $X(\beta)$. Namely, a Demazure weave $\w:\beta\to\delta(\beta)$ provides a way to find smaller and smaller open subsets of $X(\beta)$, which we will see ultimately yield an open toric embedding $\iota_\w:(\bG_m)^d\lr X(\beta)$, with $d={\dim X(\beta)}$, see \Cref{cor:weave_to_toric_embedding} below.\sbs
\end{remark}

%%%%%%%%%%%%%%%%%%%%%%%%%%%%%%%%%%%%%%%%%%%%%
%%%%%%%%%%%%%%%%%%%%%%%%%%%%%%%%%%%%%%%%%%%%%
%%%%%%%%%%%%%%%%%%%%%%%%%%%%%%%%%%%%%%%%%%%%%

\subsection{Examples of braid varieties}\label{ssec:braid_varieties_examples} Certain varieties studied in the context of representation theory, algebraic geometry and algebraic combinatorics, are \color{blue}instances of braid varieties\color{black}:\\

{\bf 1.} The \color{blue}{\it open Richardson variety} \color{black} $\SR_{u,v}\sse \G/\B$ is defined as the intersection
        $$\SR_{u,v}:=(BvB/B)\cap (B_-uB/B)$$
        of the Schubert (and opposite Schubert) cells associated to two permutations $u,v\in W(\G)$ with $u\leq v$. Then, the open Richardson variety $\SR_{u,v}\sse \G/\B$ is isomorphic to the braid variety $X(\beta_{u,v})$, with
        $$\beta_{u,v}=\beta(v)\beta(u^{-1}w_0),\quad u,v\in W(\G),$$
        where $\beta(s)$ denotes the braid lift of a reduced permutation $s\in W(\G)$ to the positive braid monoid, and $w_0\in W(\G)$ is the longest element, a.k.a.~the half-twist. This is proven in \cite[Section 4]{CGGS2} and \cite[Section 3.6]{CGGLSS25}.\\

{\bf 2.} The \color{blue}{\it open positroid varieties} \color{black} $\Pi_{u,v}\sse\Gr(k,n)$ are smooth subvarieties of the Grassmannian $\Gr(k,n)$, indexed by two permutations $u,v\in W(\G)$ with $v$ satisfying a certain combinatorial property (being $k$-Grassmannian). Following the Richardson case $(1)$ above, each positroid $\Pi_{u,v}\sse\Gr(k,n)$  is isomorphic to the braid variety $X(\beta_{u,v})$, with
        $$\beta_{u,v}=\beta(v)\beta(u^{-1}w_0),\quad u,v\in W(\G).$$
        This is proven in \cite[Theorem 1.3]{CGGS2}. For instance, the braid variety for the top-dimensional (open) positroid $\Pi_{u,v}\sse\Gr(k,n)$ is given by the braid (typically) associated to the $(k,n+k)$-torus link.\\
        
        The relation between open positroid varieties and open Richardson varieties is explained in \cite{KLS,KLS2}; see also the survey \cite{speyer2024richardsonvarietiesprojectedrichardson}. At its core, the former are projections of the latter, under the natural projection from the (complete) flag variety $\G/\borel$ to the Grassmannian $\Gr(k,n)$, seen as a partial flag variety. Notably, unlike a general Richardson variety, positroid varieties are also described by many other combinatorial gadgets, different from the pair $(u,v)$. These include $k$-bounded affine permutations $f$, cyclic rank matrices $r$, and Le diagrams $\Le$. It is explained in \cite{CGGS2} how to extract braids $J_k(f),M_k(r),D_{k}(\Le)$ associated to each of these combinatorial objects so that the associated braid varieties are isomorphic to the given positroid variety.\\

{\bf 3.} As precursors of $(1)$ and $(2)$ above, \color{blue}{\it double Bruhat cells} \color{black} were one of the first classes of varieties to be studied cluster-theoretically, cf.~\cite{FZ1999}. The double Bruhat cells $\G_{u,v}$ are intersections of cells from the Bruhat decomposition of $\G$ with respect to opposite Borel subgroups, thus indexed by pairs $(u,v)\in W(\G)$. The quotient of the double Bruhat cell $\G_{u,v}$ by the corresponding torus is isomorphic to a braid variety $X(\beta_{u,v})$, e.g.~choosing the positive braid word $\beta_{u,v}=uvw_0$. More generally, the study of double Bott-Samelson varieties, cf.~\cite{ShenWeng2021DoubleBottSamelson}, can be subsumed within that of braid varieties. For instance, the (half-decorated) double Bott-Samelson variety $\mbox{Conf}^b_d$, with $b,d\in \Br^+_n$ positive braid words, is isomorphic to $X(\beta_{b,d})$, where $\beta_{b,d}=bdw_0$; this is explained in \cite[Section 3.7]{CGGLSS25}.\\

\noindent Another example from algebraic combinatorics and Lie-theory comes from brick manifolds. Brick manifolds were introduced and studied in \cite{Escobar}, building on \cite{MR2864447,MR3327085} and in part focusing on the relation to toric geometry; see also \cite[Section 3.4]{speyer2024richardsonvarietiesprojectedrichardson}. It is an exercise to show that the \color{blue} maximal open stratum of a {\it brick manifold} \color{black} is a braid variety, cf.~\cite[Theorem 1.4.(i)]{CGGS2}.\\

\begin{comment}
\color{orange}
Theorem \ref{thm:main} proves that the rings of regular functions of each such types of varieties are cluster algebras. Even in cases where Theorem \ref{thm:main} was known, e.g.~double Bruhat cells and positroids, the proof we give -- following \cite{CGGLSS25} -- is different from the original sources. In these cases, some comparisons with the previously known cluster algebras are discussed in \cite[Section 10]{CGGLSS25}, \cite[]{CLSW23}, \cite[Section 10.5]{GLSBS} and \cite{clustercomparison25}.
\color{black}
\end{comment}

\begin{remark}
$(i)$ Braid varieties also appear in the study of affine Springer fibers, non-abelian Hodge theory and their connections to knot invariants. Specifically, part of the homology of $X(\beta)$ can be identified with lowest $a$-degree Khovanov–Rozansky link homology. See \cite[Theorem 4.20]{CGGS2}, and also \cite[Section 4]{gorsky23_braidlectures} and \cite[Appendix B]{trinh2021heckecategoryunipotentlocus}. See \cite[Appendix B]{trinh2021heckecategoryunipotentlocus} for a summary of braid varieties in the literature.\\

\noindent $(ii)$ In an appropriate sense, braid varieties can be seen as a generalization of the Betti side, under the irregular Riemann-Hilbert correspondence for curves, for a De Rham moduli of meromorphic flat connections on the 2-sphere with one irregular singularity at infinity. This ties braid varieties closely with the work of P.~Boalch on Stokes local systems, cf.~\cite{Boalch21_TopologyStokes}. Specifically, braid varieties for a positive braid word $\beta$ are (a form of) moduli spaces of filtered Stokes local systems on the 2-sphere.\sbs
\end{remark}
\begin{center}
	\begin{figure}
		\centering
		\includegraphics[scale=1]{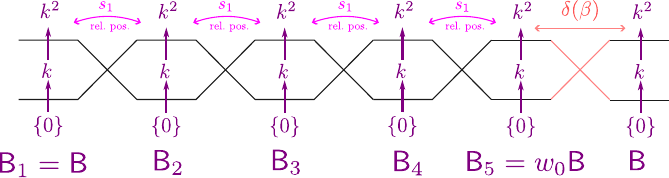}
		\caption{The linear configurations of flags parametrized by $X(\s_1^4)$, as in \cref{ex:simple_braidvarieties}.(1). The first flag $\borel_1$ on the left is the standard flag $\{0\}\sse\langle e_1\rangle\sse k^2$, the last flag $\borel_5$ is the antistandard flag $\{0\}\sse\langle e_2\rangle\sse k^2$.}
		\label{fig:BraidVariety_Example_2Strands}
	\end{figure}
\end{center}

\noindent If the reader is not familiar with the types of varieties discussed above, the following explicit examples might be helpful:

\begin{example}\label{ex:simple_braidvarieties}
(1) Consider the braid $\beta=\s_1^n\in\Br(\SL_2)$. The Weyl group $S_2=W(\SL_2)=\langle s_1\rangle$ is generated by a single order-2 transposition, a complete flag in $k^2$ is given by a line $k\sse k^2$, and two lines $\ell_1\stackrel{s_1}{\to}\ell_2$ are in $s_1$-relative position if and only if $\ell_1\neq\ell_2$. The standard flag associated to the upper-triangular Borel $\borel\sse \SL_2$ is given by the line $\langle e_1\rangle\sse k^2$, spanned by $e_1=(1,0)\in k^2$, and the antistandard flag is given by $\langle e_2\rangle\sse k^2$, spanned by $e_2=(0,1)$. The associated braid variety can be described as:
$$X(\s_1^n)\cong\{(\ell_1,\ldots,\ell_{n+1}):\ell_i\subset k^2\mbox{ lines such that }\ell_i\neq\ell_{i+1}, \ell_1=\langle e_1\rangle,\ell_{n+1}=\langle e_2\rangle\}.$$
See Figure \ref{fig:BraidVariety_Example_2Strands} for a visualization. In particular $\dim(X(s_1^n))=n-1$ and it is visibly an affine variety. In simple cases, $X(s_1)=\{\mbox{pt}\}$ is a point and $X(s_1^2)\cong k^\times$, since the space of lines in $k^2$ different from both $\langle(1,0)\rangle$ and $\langle(0,1)\rangle$ is given by the graphs of non-zero linear maps from $k$ to $k$. In general, the varieties $X(\s_1^n)$ are rather interesting for arbitrary $n\in\N$.\\

\noindent (2) Consider the braid $\beta=c^nw_0\in\Br(\SL_{m})$, where $c=\s_1\s_2\cdots \s_{m-1}$ is a Coxeter element, $m\leq n$, and $\ell=\ell(\beta)=n(m-1)+{m\choose 2}$. A point in the braid variety $X(\beta)$ is given by a tuple $(\borel_1,\ldots,\borel_{\ell+1})$ of $(\ell+1)$ full flags in $k^{m}$ satisfying certain conditions. From the piece of $\beta$ corresponding to $c^n$, let $w_1,\ldots,w_n\in k^m$ be vectors such that $\langle w_i\rangle$ coincides with the 1-dimensional subspace of the flag $V_{1+i(m-1)}$, $i\in[1,n]$. Note that in this case the flag $V_{1+j(m-1)}$ is completely determined by $\langle v_j\rangle$ and the flag $V_{1+(j+1)(m-1)}$, $j\in[1,n-1]$, as $V_{1+j(m-1)}$ and $V_{1+(j+1)(m-1)}$ are in $c$-relative position. Using the piece corresponding to $w_0\delta(\b)=w_0^2$, we can construct a set of vectors $w_{n+1},\ldots, w_{m+n}\in k^m$ such that the collection of 1-dimensional vector subspaces $\langle w_i\rangle$ for $i\in[1,n+m]$ uniquely determines a point in $X(\beta)$. Indeed, these vectors $w_j$ are chosen exactly as the previous $w_i$-vectors upon using the identification $c^m=w_0^2$, which holds because $m$ is the Coxeter number. That is, $w_j$ are vectors generating the lines in the corresponding tuples of flags that appears after moving the piece of the braid $w_0^2$ to $c^m$. In terms of such collections of 1-dimensional subspaces, the conditions for the tuple of flags $(\borel_1,\ldots,\borel_{\ell+1})\in X(\beta)$ to belong to the braid variety become
$$w_1\wedge\ldots\wedge w_m\neq 0,\quad w_2\wedge\ldots\wedge w_{m+1}\neq 0,\quad\ldots,\quad w_{n+m}\wedge w_1\wedge\ldots\wedge w_{m-1}\neq 0.$$
Therefore $X(\beta)$ is isomorphic to the locus of the Grassmannian $\Gr(m,n+m)$ where the cyclically consecutive Pl\"ucker coordinates do not vanish quotiented by the corresponding torus.\footnote{The torus quotient is needed, as the dimensions of $X(c^nw_0)$ and $\Gr(m,n+m)$ are $n(m-1)$ and $nm$.} This is known as the top-dimensional open positroid variety, modulo its torus action. Note that the $(-1)$-framed satellite of such $\beta$ along the unknot in $\R^3$ is the $(m,n)$-torus link.\footnote{The 0-framed satellite is the $(m,n+m)$-torus link.} Understood as a max-tb Legendrian link, this explains geometrically the combinatorial relations between invariants of torus links and certain finite point counts on this positroid in $\Gr(m,n+m)$, see e.g.~\cite{LamGalashin_Catalan,GL2023}.\sbs
\end{example}

\begin{comment}
\color{orange}
\noindent It can be enlightening to understand how a braid variety $X(\beta)$ changes when we insert a new crossing to a given braid word $\beta$. The following result, which follows from \cite[Lemma 3.4 \& Section 5.5]{CGGLSS25}, can help build intuition:

\begin{lemma}[Braid varieties under adding crossings]\label{lemma:add_crossing} Let $\beta,\beta'\in\Br_n^+$ be positive braid words such that $\beta'$ is obtained from $\beta$ by inserting one new crossing $\s_i$ somewhere in $\beta$. Then
\begin{enumerate}
\item If $\delta(\beta)=\delta(\b')$, then $X(\beta)\cong X(\beta')$.\\
\item If $\delta(\beta)=\delta(\b')$, then $X(\beta)$ is isomorphic to a codimension-one locally closed subvariety $X(\beta')$.
\end{enumerate}
\end{lemma}
\noindent It might be good to contrast \Cref{lemma:add_crossing} with \Cref{ex:permutation_vs_braid}.
\color{black}
\end{comment}
%%%%%%%%%%%%%%%%%%%%%%%%%%%%%%%%%%%%%%%%%%%%%
%%%%%%%%%%%%%%%%%%%%%%%%%%%%%%%%%%%%%%%%%%%%%
%%%%%%%%%%%%%%%%%%%%%%%%%%%%%%%%%%%%%%%%%%%%%

\subsection{Invariance under braid moves}\label{ssec:braid_varieties_invariance} Let us consider the two braid moves for $\G=\SL_n$:
$$\s_i\s_{i+1}\s_i\to\s_{i+1}\s_i\s_{i+1},\qquad \s_i\s_j\to\s_j\s_i,$$
where $|i-j|\geq2$. In topology, the first move is often referred to as a Reidemeister III (R3) move. The second move can be referred to as a commutation. In terms of the diagrammatics of weaves, these two local braid moves are respectively encoded by the hexavalent and tetravalent vertices, see \Cref{fig:weaves_types_vertices} or \Cref{fig:BraidVariety_DiagramMoves}. Let us now insert flags to these pieces of the weave, building on the diagrammatics from \Cref{sssec:first_connection_weave_diagrams}, see e.g.~\Cref{fig:BraidVariety_DiagramExample}. For that we insert a (horizontal) linear configuration of flags at the horizontal slices above and below the vertex, as in \Cref{fig:BraidVariety_DiagramMoves}. The following result explains how to propagate downwards a tuple of flags in the horizontal slice above a hexavalent or tetravalent vertex:

\begin{center}
	\begin{figure}
		\centering
		\includegraphics[scale=0.65]{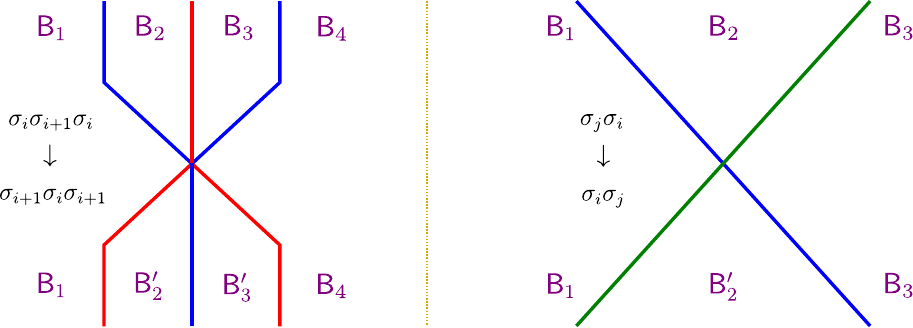}
		\caption{(Left) A hexavalent weave vertex, with the tuple of flags $(\borel_1,\borel_2,\borel_3,\borel_4)$ along a horizontal slice above the vertex, and tuple of flags $(\borel_1,\borel'_2,\borel'_3,\borel_4)$ below the vertex. \Cref{lem:braid_invariance}.(1) shows that either tuple uniquely determines the other. (Right) A tetravalent weave vertex, with the tuple of flags $(\borel_1,\borel_2,\borel_3)$ above, and $(\borel_1,\borel'_2,\borel_3)$ below the vertex. \Cref{lem:braid_invariance}.(2) shows that either tuple uniquely determines the other.}
		\label{fig:BraidVariety_DiagramMoves}
	\end{figure}
\end{center}

\begin{center}
	\begin{figure}
		\centering
		\includegraphics[width=\textwidth]{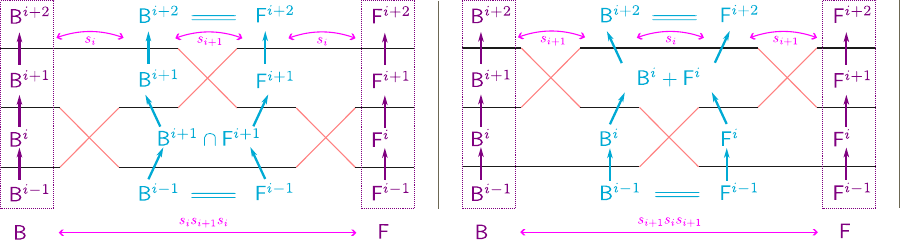}
		\caption{The proof of isomorphism invariance of $X(\beta)$ under the Reidemeister 3 ($R3_i$) braid move $\s_i\s_{i+1}\s_i\leftrightarrow\s_{i+1}\s_i\s_{i+1}$, applied to a piece of the braid word $\beta$. (Left) The two flags $\borel,\FF$ at the left and right of $\s_i\s_{i+1}\s_i$, drawn in purple, determine the two intermediate flags in between, drawn in blue. (Right) The two flags $\borel,\FF$ at the left and right of $\s_{i+1}\s_i\s_{i+1}$, in purple, determine the two flags in between, in blue.}
		\label{fig:BraidVariety_R3_Invariance}
	\end{figure}
\end{center}

%Braid varieties $X(\beta)$, as in \cref{def:braid_variety}, have a number of interesting symmetries. Some are discussed as follows:
\begin{lemma}[Invariance of $X(\beta)$ under braid moves]\label{lem:braid_invariance}
Let $\beta_1,\beta_2\in\Br(G)$ be two positive $n$-stranded braid words of lengths $\ell_1,\ell_2$. Then:\\

\noindent $(1)$ \underline{R3 move}. For any $i\in[1,n-2]$, consider the regular map
    $$R3_i:X(\beta_1\s_i\s_{i+1}\s_i\beta_2)\lr X(\beta_1\s_{i+1}\s_i\s_{i+1}\beta_2)$$
    \noindent that sends a tuple $(\borel_1,\ldots,\borel_{\ell_1+1},\color{blue}\borel_{\ell_1+2}\color{black},\color{blue}\borel_{\ell_1+3}\color{black},\borel_{\ell_1+4},\ldots,\borel_{\ell_1+\ell_2+4})$ to the tuple
    $$(\borel_1,\ldots,\borel_{\ell_1+1},\color{blue}\borel'_{\ell_1+2}\color{black},\color{blue}\borel'_{\ell_1+3}\color{black},\borel_{\ell_1+4},\ldots,\borel_{\ell_1+\ell_2+4}),$$
    \noindent where the two flags $\color{blue}\borel'_{\ell_1+2}\color{black},\color{blue}\borel'_{\ell_1+3}\color{black}$ are given by

    $$\color{blue}(\borel'_{\ell_1+2})^k\color{black}:=\begin{cases}
			(\borel_{\ell_1+1})^i, & \text{if $k=i$,}\\
            (\borel_{\ell_1+1})^i+(\borel_{\ell_1+4})^i, & \text{if $k=i+1$,}\\
            (\borel_{\ell_1+1})^k, & \text{otherwise,}		 \end{cases}
          \quad \color{blue}(\borel'_{\ell_1+3})^k\color{black}:=\begin{cases}
			(\borel_{\ell_1+4})^i, & \text{if $k=i$,}\\
            (\borel_{\ell_1+1})^i+(\borel_{\ell_1+4})^i, & \text{if $k=i+1$,}\\
            (\borel_{\ell_1+1})^k, & \text{otherwise.}		 \end{cases}$$
    
    \noindent Then $R3_i$ is an isomorphism of algebraic varieties over $\Z$. See Figure \ref{fig:BraidVariety_R3_Invariance} for a visual depiction of the isomorphism $R3_i$.\\

\noindent $(2)$ \underline{Commutation}. For any $i,j\in[1,n-1]$, $2\leq|i-j|$, consider the regular map
    $$C:X(\beta_1\s_i\s_j\beta_2)\lr X(\beta_1\s_j\s_i\beta_2)$$
    \noindent that sends a tuple $(\borel_1,\ldots,\borel_{\ell_1+1},\color{blue}\borel_{\ell_1+2}\color{black},\borel_{\ell_1+3}\color{black},\borel_{\ell_1+4},\ldots,\borel_{\ell_1+\ell_2+3})$ to the tuple
    $$(\borel_1,\ldots,\borel_{\ell_1+1},\color{blue}\borel'_{\ell_1+2}\color{black},\borel_{\ell_1+3}\color{black},\borel_{\ell_1+4},\ldots,\borel_{\ell_1+\ell_2+3}),$$
    where
    $$\color{blue}(\borel'_{\ell_1+2})^k\color{black}:=\begin{cases}
			(\borel_{\ell_1+1})^k, & \text{if $k=i$,}\\
            (\borel_{\ell_1+3})^k, & \text{if $k=j$,}\\
            (\borel_{\ell_1+2})^k, & \text{otherwise.}		 \end{cases}
    $$
    Then $C$ is an isomorphism of algebraic varieties over $\Z$.\\
\end{lemma}

\noindent \cref{lem:braid_invariance} is the reason we use the terminology {\it braid variety} for $X(\beta)$, as it follows that the isomorphism type of $X(\beta)$, as an algebraic variety over $\Z$, only depends on the braid $[\beta]\in\Br(\G)$, rather than the specific choice of word. The proof of \Cref{lem:braid_invariance}.(1) is essentially contained in \Cref{fig:BraidVariety_R3_Invariance}, and Part (2) follows from similar considerations. See \cite[Lemma 3.5]{CGGLSS25} for an alternative proof via coordinates. Another conceptual, more geometric, explanation of the invariance properties in \cref{lem:braid_invariance} is provided by symplectic topology, cf.~\cref{thm:braid_varieties_from_Legendrian}. Note also that \Cref{lem:braid_invariance} applies both ways, since these moves are invertible.\\

From the viewpoint of weaves, \Cref{lem:braid_invariance} has the following consequence. Let $\w$ be a Demazure weave and $v\in\w$ a hexavalent or tetravalent vertex. If $(\borel_1,\ldots,\borel_k)$ is a tuple of flags assigned to a horizontal slice of a Demazure weave $\w$ exactly above $v$, then \Cref{lem:braid_invariance} implies that there exists a unique allowable tuple of flags $(\borel'_1,\ldots,\borel'_k)$ to be assigned to the horizontal slice of $\w$ exactly below $v$. Conversely, a tuple of flags below such a vertex $v$ uniquely determines the tuple above such a vertex $v$. Therefore, we can uniquely propagate a given tuple of flags assigned to a horizontal slice of a Demazure weave up or down past a hexavalent or tetravalent vertex. \Cref{fig:BraidVariety_DiagramMoves} illustrates this propagation locally, with the formulas for $\borel'_2$ and $\borel'_3$ in terms of the flags above given as in \Cref{lem:braid_invariance}.\\

\begin{center}
	\begin{figure}
		\centering
		\includegraphics[scale=0.65]{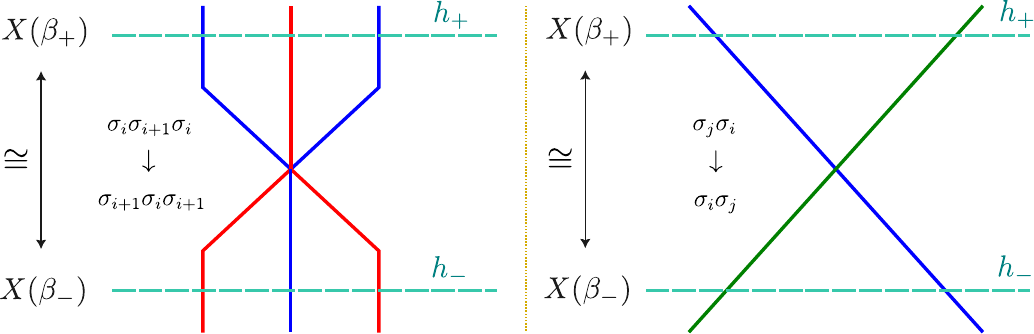}
		\caption{(Left) The isomorphism $X(\beta_+)\stackrel{\cong}{\longleftrightarrow}X(\beta_-)$ between the braid varieties given by the horizontal slices above and below a hexavalent vertex is provided in \Cref{lem:braid_invariance}.(1).  (Right) The isomorphism $X(\beta_+)\stackrel{\cong}{\longleftrightarrow}X(\beta_-)$ between the braid varieties given by the horizontal slices above and below a tetravalent vertex is provided in \Cref{lem:braid_invariance}.(2).}
		\label{fig:BraidVariety_DiagramMoves2}
	\end{figure}
\end{center}

Therefore, at this stage, we can enhance the diagrams for Demazure weaves by inserting compatible tuples of flags at each horizontal slice, and we understand how these propagate past hexavalent and tetravalent vertices. Specifically, if $h_\pm$ are two horizontal slices exactly above and below a hexavalent or tetravalent vertex, as in \Cref{fig:BraidVariety_DiagramMoves2}, and we let $\beta_\pm:=\beta_{h_\pm}$ denote their corresponding braid words, then \Cref{lem:braid_invariance} provides an explicit isomorphism $X(\beta_+)\stackrel{\cong}{\longleftrightarrow}X(\beta_-)$.

A remaining key aspect is to understand what happens to these tuples of flags when we scan downwards past a {\it trivalent vertex}, which we momentarily address in \Cref{ssec:braid_varieties_trivalent}. 

\subsubsection{Further symmetries of braid varieties}\label{sssec:braid_variety_further_symmetries} In addition to \cref{lem:braid_invariance}, there is a cyclic symmetry for braid varieties if $\delta(\beta)=w_0$, which we now assume. Specifically, for any $i\in[1,n-1]$, consider the index $i^*$ such that $s_{i^*}=w_0s_iw_0$, choose $\mathrm{w}\in W(\G)$ with $w_0=\mathrm{w}s_{i}=s_{i^*}\mathrm{w}$ and set $\ell=\ell(\beta)$. Then there exists an isomorphism $R_g^\circlearrowright:X(\s_i\beta)\lr X(\beta \s^*_i)$ of algebraic varieties over $\Z$, see \cite[Section 5.5]{CGGLSS25}. For $\G=\SL_n$, the index is $i^*=n-1$, so that $\s_{i^*}=\s_{n-i}$. Beyond \cref{lem:braid_invariance}, it is also possible to relate certain braid varieties for different ranks $n$ via a Markov-type stabilization process, cf.~\cite[Section 3.2]{CGGS2}. That is, certain braid varieties for $\SL_n$ are isomorphic to those for $\SL_{n+1}$, up to a trivial $\bG_m$-factor. For instance, such isomorphisms capture the $n\leftrightarrow m$ symmetry of $(n,m)$-torus links.\\

Braid varieties have another type of symmetry: they are acted on by a Cartan torus $\T\sse\G$. Since $\T$ fixes the standard and antistandard flags $\borel$ and $w_0\borel$, it follows that the boundary conditions for $X(\beta)$ are preserved, cf.~\cref{eq:braid_variety}. Since $\T$ acts on the flag variety $\G/\borel$ preserving the relative position of pairs of flags, the torus action
\begin{equation}\label{eq:torus_action}
\T\times X(\beta)\lr X(\beta),\quad (\tau;\borel_1,\borel_2,\dots, \borel_{\ell+1})\mapsto (\tau\borel_1,\tau\borel_2,\dots, \tau\borel_{\ell+1})
\end{equation}
is well-defined. In fact, the isomorphisms in \cref{lem:braid_invariance} are $\T$-equivariant, see e.g.~\cite[Section 2.2]{CGGJ24} or \cite[Section 4.5]{CGSS_clusterdeeploci}. The $\T$-action is, in general, not free. For instance, it is not free for $\beta=\s_1^3\in\mbox{Br}_2^+$ as the point $(\borel,\borel_-,\borel,\borel_-)\in X(\s_1^3)$ is fixed. That said, the corresponding $(\T/Z(\G))$-action from \cref{eq:torus_action} is free if $\beta\delta(\beta)$ defines a knot, i.e.~if the decomposition of the permutation $\pi_{\beta\delta(\beta)}\in S_n$ given by the Coxeter projection of $\beta\delta(\beta)$ has a unique factor in its cycle decomposition. In this case, $X(\beta)/(\T/Z(\G))$ is an affine variety and there is an isomorphism $X(\beta)\cong \left(X(\beta)/(\T/Z(\G))\right )\times(\T/Z(\G))$ of algebraic varieties over $\Z$.

\begin{remark}\label{rmk:torus_quasimoment_map} This torus action is part of a quasi-Hamiltonian structure on $X(\beta)$, which further exemplifies the rich geometry of braid varieties. In particular, there is a $\T$-valued moment map $\mu:X(\beta)\lr\T$, where $\T$ acts on the codomain via the $\pi_\beta$-twisted conjugation $(\tau;t)\mapsto (\pi_\beta^{-1}\tau\pi_\b)t\tau^{-1}$.\sbs
\end{remark}
%%%%%%%%%%%%%%%%%%%%%%%%%%%%%%%%%%%%%%%%%%%%%
%%%%%%%%%%%%%%%%%%%%%%%%%%%%%%%%%%%%%%%%%%%%%
%%%%%%%%%%%%%%%%%%%%%%%%%%%%%%%%%%%%%%%%%%%%%

\begin{center}
	\begin{figure}
		\centering
		\includegraphics[scale=2.65]{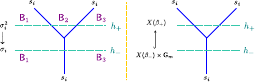}
		\caption{(Left) A trivalent vertex of weave edges labeled by $s_i$, with two horizontal slices $h_\pm$ and the corresponding tuples of flags. In this case, along $h_+$, the relative positions of $(\borel_1,\borel_2,\borel_3)$ are forced to be $\borel_1\stackrel{s_i}{\to}\borel_2$ and $\borel_2\stackrel{s_i}{\to}\borel_3$. Along $h_-$, the relative position of $(\borel_1,\borel_3)$ is $\borel_1\stackrel{s_i}{\to}\borel_3$. (Right) The injective map between the braid varieties, where the braid word $\beta_-$ along $h_-$ is obtained from the word $\beta$ along $h_+$ by applying a $\s_i^2\to\s_i$ move, cf.~\Cref{lem:braidvarieties_trivalent}.}
		\label{fig:BraidVariety_Trivalent2}
	\end{figure}
\end{center}

\subsection{Trivalent weave vertices and braid varieties}\label{ssec:braid_varieties_trivalent} Let us start by studying the local situation of flags near a weave trivalent vertex $v$ for the $\SL_2$-weave $\w:\s_1^2\to\s_1$, as depicted in \Cref{fig:BraidVariety_Trivalent2} (left). At the horizontal slice $h_+$ above $v$ the tuple of flags is of the form $(\borel_1,\borel_2,\borel_3)$, and the relative positions are $\borel_1\stackrel{s_i}{\to}\borel_2$ and $\borel_2\stackrel{s_i}{\to}\borel_3$.\footnote{For now there is merit in not imposing $\borel_1=\borel$ and $\borel_3=s_1\borel$, as maintaining the freedom at the boundary produces a local model that we can insert anywhere in the middle of a braid word to model $\s_i^2\to\s_i$.} At the horizontal slice $h_-$ below $v$ the tuple of flags is of the form $(\borel_1,\borel_3)$, and the relative positions are $\borel_1\stackrel{s_i}{\to}\borel_3$. Therefore, if we are scanning top-to-bottom, implementing the move $\s_i^2\to\s_i$ via a trivalent vertex of a weave, \color{blue}the effect of a trivalent vertex on a tuple of flags as we scan downwards is to impose an additional open condition \color{black} on the tuple of flags above the vertex. In this local model, this is the additional condition $\borel_1\stackrel{s_i}{\to}\borel_3$, which was not necessarily satisfied by all the tuples $(\borel_1,\borel_2,\borel_3)$ along $h_+$.\\

If we read this situation in \Cref{fig:BraidVariety_Trivalent2} (left) from bottom-to-top, the tuple $(\borel_1,\borel_3)$ below determines only two of the flags in the tuple $(\borel_1,\borel_2,\borel_3)$ above. The flag $\borel_2$ can be chosen: since we are in $\SL_2$, and the flag variety is $\SL_2/\borel\cong\PP_k^1$, the choice of $\borel_2$ is given by any point in $(\SL_2/\borel\setminus\{\borel_1,\borel_3\})\cong(\PP_k^1\setminus\{0,\infty\})$. By coordinatizing $(\PP_k^1\setminus\{0,\infty\})$ via $z\in\bG_m$, we conclude that the choices of $\borel_2$ are given by an invertible value $z\in\bG_m$. More generally, by inserting this local model in a positive braid word we obtain the following result:

\begin{lemma}\label{lem:braidvarieties_trivalent}
Let $\beta_1,\beta_2$ be two positive braid words and consider the length $\ell_1:=\ell(\beta_1)-1$. Then there exists a natural open embedding $X(\beta_1\s_i\beta_2)\times\bG_m\hookrightarrow X(\beta_1\s_i^2\beta_2)$ of algebraic varieties over $\Z$ given by
\begin{equation}\label{eq:braidvarieties_trivalent}
(\ldots,\borel_{\ell_1},\borel_{\ell_1+1},\ldots;z)\mapsto (\ldots,\borel_{\ell_1},\color{blue}\borel_{z}\color{black},\borel_{\ell_1+1},\ldots)
\end{equation}
where $\color{blue}\borel_{z}\color{black}$ denotes a flag parametrized by $z\in\bG_m\cong ((\SL_2/\borel)\setminus\{\borel_{\ell_1},\borel_{\ell_1+1}\})$.\sbs
\end{lemma}

\noindent See \Cref{fig:BraidVariety_Trivalent2} (right) for a depiction of \Cref{lem:braidvarieties_trivalent}. The map in \eqref{eq:braidvarieties_trivalent} is typically not surjective. Indeed, any tuple in its image must satisfy $\borel_{\ell_1+1}\stackrel{s_i}{\to}\borel_{\ell_1+2}$, since the domain is $X(\beta_1\s_i\beta_2)$, whereas there is typically no such constraint $\borel_{\ell_1+1}\stackrel{s_i}{\to}\borel_{\ell_1+3}$ in the target $X(\beta_1\s_i^2\beta_2)$.

%%%%%%%%%%%%%%%%%%%%%%%%%%%%%%%%%%%%%%%%%%%%%
%%%%%%%%%%%%%%%%%%%%%%%%%%%%%%%%%%%%%%%%%%%%%
%%%%%%%%%%%%%%%%%%%%%%%%%%%%%%%%%%%%%%%%%%%%%
\subsection{Weave tori in braid varieties}\label{ssec:weave_tori_braid_varieties} Now, given a Demazure weave $\w:\beta\to\delta(\beta)$, we can consider the braid variety $X(\beta)$ as the space of tuples of flags which can be compatibly assigned to the horizontal slice of the north boundary of $\w$. By iteratively scanning top-to-bottom we encounter weave vertices:
\begin{enumerate}
    \item At hexavalent and tetravalent vertices, \Cref{lem:braid_invariance} implies that the tuple of flags can be propagated downwards uniquely {\it without} any loss of information. That is, the braid varieties that arise above and below such vertices are isomorphic, and the particular isomorphism from \Cref{lem:braid_invariance} is chosen to propagate the tuples of flags, cf.~\Cref{fig:BraidVariety_DiagramMoves2}.\\

    \item At a trivalent vertex, \Cref{lem:braidvarieties_trivalent} implies that the tuple of flags can be propagated downwards past the trivalent vertex {\it iff} we impose an additional relative position constraint to the above tuple of flags. \Cref{lem:braidvarieties_trivalent} shows that this determines a locus inside the braid variety above which splits up a $\bG_m$ factor, cf.~\Cref{fig:BraidVariety_Trivalent2} (right).\\
\end{enumerate}

In particular, if we consider two horizontal slices $h,h'$ of the weave $\w$, with $h$ below $h'$ and $q$ trivalent vertices in between, then there exists an open embedding $X(\beta_h)\times(\bG_m)^q\hookrightarrow X(\beta)$, cf.~\cite[Section 5]{CGGJ24}. In particular, if we consider $h$ to be a horizontal slice of the south boundary, so $X(\beta_h)$ is a point and $h'$ a horizontal slice of the north boundary, so $X(\beta_{h'})=X(\beta)$, we obtain an open embedding $\iota_\w:(\bG_m)^d\hookrightarrow X(\beta)$ of an algebraic torus, where $d=\dim(X)$. In conclusion, \Cref{lem:braid_invariance} and \Cref{lem:braidvarieties_trivalent} lead to the following:

\begin{fcor}[Weave tori $T_\w$ in braid varieties]\label{cor:weave_to_toric_embedding}
Let $\beta\in\Br_n^+$ be a positive braid word, $X(\beta)$ its braid variety and $d:=\dim(X)$ its dimension. Consider a Demazure weave $\w:\beta\lr\delta(\beta)$ with $d$ trivalent vertices. Then:
\begin{enumerate}[label=$(\roman*)$]
    \item There exists an open toric embedding
\begin{equation}\label{eq:weave_tori_map}
\iota_\w:(\bG_m)^d\hookrightarrow X(\beta)
\end{equation}
whose image is uniquely determined by assigning flags in $\SL_n/\borel$ to each face of $\w$ such that:\\
\begin{enumerate}[label=$(\alph*)$]
    \item the relative position of such flag is determined by the weave edges, so that two flags $\borel_l,\borel_r$ to the left and right of a weave $s_i$-edge are forced to be in $s_i$-relative position.\\
    \item the leftmost flag must be $\borel$, and the rightmost flag must be $\delta(\beta)\borel$.\\
\end{enumerate}

\item If two Demazure weaves $\w,\w':\beta\lr\delta(\beta)$ are weave equivalent $\w\simeq\w'$, then the images of the corresponding open embeddings $\iota_\w$ and $\iota_{\w'}$ coincide.\\
\item If two Demazure weaves $\w,\w':\beta\lr\delta(\beta)$ are not weave equivalent $\w\not\simeq\w'$, then the images of the corresponding open embeddings $\iota_\w$ and $\iota_{\w'}$ do not coincide.\sbs
\end{enumerate}
\end{fcor}
\vspace{0.2cm}
\noindent \Cref{cor:weave_to_toric_embedding}.$(ii)$ requires a few computations beyond the discussion above, cf.~\cite[Section 5.3]{CGGLSS25}, and \Cref{cor:weave_to_toric_embedding}.$(iii)$ follows from \cite[Theorem 1.1]{CGGLSS25} and \cite[Prop.~ 5.3]{CaoKellerQin22}. We denote the image of $\iota_\w$ in $X(\beta)$ by $T_\w\sse X(\beta)$ and refer to $T_\w$ as a weave torus. A fundamental non-trivial aspect of Demazure weaves is that \color{blue}there are enough Demazure weaves $\w:\beta\to\delta(\beta)$ so that the corresponding collection of weave tori $T_\w$ covers $X(\beta)$ up to codimension two\color{black}. In particular, by the algebraic Hartogs' lemma, the coordinate ring $k[X(\beta)]$ can be entirely studied by using such open tori $T_\w$.\footnote{In fact, $\dim(X)+1$ suffice for such a covering.}

\begin{remark}
More generally, given a weave $\w\sse S$ in a smooth surface $S$, we can consider the flag moduli $\M_\w$ associated to $\w$, as in \cite[Section 5.2]{CasalsZaslow} and see also \cite[Section 5.1]{CGGJ24}. In the case of Demazure weave $\w$, its flag moduli $\M_\w$ is closely related to the weave tori $T_\w$.\sbs
\end{remark}

\subsection{Summary of section}\label{ssec:summary_braid varieties} At this stage we have the following summary of mathematical objects and their relations:
\begin{enumerate}
    \item Given a positive braid word $\beta$, we can consider its braid variety $X(\beta)$, as in \Cref{def:braid_variety}. As illustrated in \Cref{ssec:braid_varieties_examples}, braid varieties generalize a number of important varieties considered in the literature.\\
    
    \item Given a Demazure weave $\w:\beta\to\delta(\b)$, \Cref{cor:weave_to_toric_embedding} provides an open torus $T_\w\sse X(\beta)$. By construction, this open torus $T_\w$ is obtained by imposing further relative position constraints to the tuples of flags giving the points in $X(\beta)$. These additional relative position constraints are dictated by the diagram of $\w$.\\

    \item Given a Demazure weave $\w:\beta\to\delta(\b)$, we can consider the Lusztig cycles $\{\g_v\}$, as in \Cref{def:Lusztig_cycles}, and \Cref{def:intersection_number} provides a quiver $Q_\w$ encoding their intersections.\\
\end{enumerate} 
\Cref{cor:weave_to_toric_embedding} and \Cref{cor:quiver_cycles_weave_operations} imply that both $T_\w$ and $Q_\w$ remain the same under weave equivalences applied to $\w$, and they necessarily change under a weave mutation.\footnote{The resulting quiver $Q_\w$ might, nevertheless, be abstractly isomorphic to $Q_\w$, but $T_\w$ will be genuinely different.}\\

Let $\beta\in\Br_n^+$ be a positive braid word, $\w:\beta\to\delta(\beta)$ a Demazure weave and $X(\beta)$ its braid variety.\\

\noindent Our next task will be to explain how each Lusztig cycle $\g_v$ and each dual Lusztig cycle $\rho_v$ respectively define a rational function $X_{\g_v}\in\Z(X(\beta))$ and a regular function $A_{\r_v}\in\Z[X(\beta)]$. These functions will be the cluster $X$-variables and cluster $A$-variables associated to the quiver $Q_\w$, geometrically realized as elements of the rational field of functions $\Z(X(\beta))$ and the coordinate ring of functions $\Z[X(\beta)]$ of the braid variety $X(\beta)$.

\section{Cluster structures on braid varieties}\label{sec:braid_variety_cluster}

The object of this section is to explain the construction of a cluster algebra structure on the commutative algebra $\Z[X(\beta)]$ given by the coordinate ring of regular functions on $X(\beta)$. In order to do that, we first review framed flags in \Cref{ssec:framed_rep} and discuss certain quantitative measurements of transversality between them in \Cref{ssec:rational_functions}, which help towards better understanding the tori underlying the cluster seeds for $\Z[X(\beta)]$. The cluster $\SA$-variables are constructed in \Cref{ssec:cluster_Avariables_weave} and \Cref{ssec:cluster_A_variables_propagation_vertices} explains how to obtain them algorithmically by scanning the weave bottom-to-top. Similarly, \Cref{ssec:cluster_Xvariables_weave} constructs the dual cluster $\SX$-variables and presents their key properties.\\

In \Cref{ssec:A_and_X_schemes_weave} we explain the relation between the braid variety $X(\beta)$ and the cluster $\SA$ and $\SX$-schemes, and \Cref{ssec:brick_compactifications} concludes by discussing certain projective compactifications of braid varieties that can be used to better understand them, cf.~\Cref{sssec:applications_compactifications}.\\

\begin{comment}
Theorem \ref{thm:main} proves that the rings of regular functions of each such types of varieties are cluster algebras. Even in cases where Theorem \ref{thm:main} was known, e.g.~double Bruhat cells and positroids, the proof we give -- following \cite{CGGLSS25} -- is different from the original sources. In these cases, some comparisons with the previously known cluster algebras are discussed in \cite[Section 10]{CGGLSS25}, \cite[]{CLSW23}, \cite[Section 10.5]{GLSBS} and \cite{clustercomparison25}.
\end{comment}

%%%%%%%%%%%%%%%%%%%%%%%%%%%%%%%%%%%%%%%%%%%%%
%%%%%%%%%%%%%%%%%%%%%%%%%%%%%%%%%%%%%%%%%%%%%
\subsection{Framed representatives for points in $X(\beta)$}\label{ssec:framed_rep} Let $\beta\in\brn$ be a positive braid word. The braid variety $X(\beta)$ and the weave tori $T_\w\sse X(\beta)$ are both defined by using relative positions between flags in $\G/\borel$. For a pair $\borel_1,\borel_2\in\G/\borel$ we can decide their relative position, given by a permutation in $W(\SL_n)$, but there is not enough information in $\G/\borel$ to measure this relative position in a quantifiable manner. For instance, if $\borel_1\stackrel{s_i}{\to}\borel_2$, we would want to be able to measure how (or how much) $\borel^i_1$ and $\borel^i_2$ differ. This leads to the notion of framed flags. To define them, let $\uni\sse \borel$ be a fixed unipotent radical, e.g.~the subgroup of upper unitriangular matrices for $\G=\SL_n$.

\begin{definition}[Framed flags]
$\G/\uni$ is said to be the principal flag variety, and a framed flag is a point in $\G/\uni$. Given $w\in W(G)$, two framed flags $x\uni,y\uni\in\G/\uni$ are said to be in relative $w$-position, denoted $x\uni\stackrel{w}{\to}y\uni$, if $x^{-1}y\in \uni w\uni$.\hfill$\Box$
\end{definition}
Since $\T\cong\borel/\uni$, the projection map $\pi:\G/\uni\lr\G/\borel$ induced by the inclusion $\uni\sse \borel$ is a principal $\T$-bundle; in particular its fibers are isomorphic to $\T$. In the correspondence between Borel subgroups and flags for $\G=\SL_n$, as in \cref{rmk:borels_flags}, framed flags in $\G/\uni$ can be interpreted as pairs $(\FF,\{f_1,\ldots,f_n\})$ consisting of a full flag
$$\FF=(\FF^0\subset \FF^1\subset \FF^2\subset\ldots \subset \FF^n)\in\G/\borel,\quad \FF^i\cong k^i,$$
equipped with a set of generators $\{f_1,\ldots,f_n\}$ of the 1-dimensional vector spaces $\FF^i/\FF^{i-1}=\langle f_i\rangle$ given by the consecutive quotients of the full flag $\FF$ such that these generators lift to vectors $\tilde{f}_i\in k^n$ with volume one, i.e.~$\tilde{f}_1\wedge\ldots\wedge \tilde{f}_n=\mbox{vol}_n$, where $\mbox{vol}_n$ is the volume element intrinsically defining $\SL_n\sse\on{GL}_n$. Note that the projection $\pi(\FF,\{f_1,\ldots,f_n\})=\FF$ above is forgetting the data of the choice of generators and just remembering the underlying flag $\FF$. Each flag in $\G/\borel$ can be represented by many framed flags, and such framed flags are a $\T$-torsor: this simply corresponds to the choices of generators $\{f_1,\ldots,f_n\}$.

\begin{example}[Framed flags in $\SL_2$]\label{ex:framed_flags_SL2}
Consider $\G=\SL_2$, its framed flag variety is $\SL_2/\uni\cong\kk^2 \setminus \{(0,0)\}$. For instance, if we represent $g\in\SL_2$ via $g=(v_1\quad v_2)$ where $v_1,v_2\in\kk^2$ denote the columns of $g$, then an isomorphism is given by $g=(v_1\quad v_2)\uni\mapsto v_1$. The forgetful projection $\pi:\SL_2/\uni\lr\SL_2/\borel$ is the quotient of $\kk^2 \setminus \{(0,0)\}$ by the $\bG_m$-action diagonally scaling the entries of $v_1\in \kk^2 \setminus \{(0,0)\}$. In particular, the fibers of $\pi$, which parametrize the possible framings of a given flag $g\borel\in\SL_2/\borel$, are the $\bG_m$-orbits and thus isomorphic to $\T=\bG_m$, since the action is free. In the coordinates above, the fiber over $g\borel$ for $g=(v_1\quad v_2)$ consists of those framed flags $\tilde{g}\uni$ such that
\begin{equation}\label{eq:framed_flags_sl2}
\tilde{g}=(av_1\quad a^{-1}v_2+bv_1),\quad \mbox{for some }(a,b)\in\kk^\times\kk.
\end{equation}
Varying $a\in\kk^\times$ provides the $\G_m$-fiber, while varying $b$ does not change the unipotent coset. As in \cref{rmk:borels_flags},  the flag $g\borel$ is represented by the sequence of nested subspaces $\{0\}\subset \langle v_1\rangle \subset \kk^2\cong\langle v_1,v_2\rangle$. Since all framed flags $g'\uni$ with $\pi(g'\uni)=g\borel$ are parametrized by \eqref{eq:framed_flags_sl2}, the framing $\{f_1,f_2\}$ of such a framed flag $g'\uni$ is given by the image $f_1=av_1$ of $v_1$ in $\langle v_1\rangle/\{0\}$ and $f_2=a^{-1}v_2$ in $\kk^2/\langle v_1\rangle$. Note that $av_1\wedge(a^{-1}v_2+\lambda v_1)=v_1\wedge v_2=\on{vol}_{\kk^2}$ since $v_1,v_2$ are the columns of a matrix in $\SL_2$. In general, given a framed flag $g'\uni\in\SL_n/\uni$, the framing $f_i$ at the $i$th quotient of the flag $\pi(g'\uni)$ is given by the image of the $i$th column of $g'$.\sbs
\end{example}

For us, it is useful that {\it every} point $(\borel_1,\borel_2,\ldots,\borel_{\ell+1})\in X(\beta)$ in the braid variety can be represented by a tuple of framed flags and, moreover, uniquely so if we require the framed lift of $\borel_1$ to be the specific lift $\uni$ of $\borel_1=\borel$. To show that, it suffices to note that given any framed flag $x\uni\in\G/\uni$ and (unframed) flag $\FF\in\G/\borel$ with $x\borel\stackrel{s_i}{\to}\FF$, there exists a {\it unique} framed flag $y\uni\in\G/\uni$ such that $x\uni\stackrel{s_i}{\to}y\uni$. Indeed, at an $s_i$-crossing with $x\borel$ on its left and $\FF$ on its right, there are canonical isomorphisms
$$\frac{(x\borel)^k}{(x\borel)^{k-1}}\cong \frac{\FF^k}{\FF^{k-1}}\mbox{ if }k\neq i,i+1,\quad \frac{(x\borel)^i}{(x\borel)^{i-1}}\cong \frac{\FF^{i+1}}{\FF^{i}}\mbox{ and }\quad \frac{(x\borel)^{i+1}}{(x\borel)^{i}}\cong \frac{\FF^{i}}{\FF^{i-1}}$$
and thus a framing of $x\borel$, lifting it to a framed flag, is uniquely transferred to a framing for $\FF$. Therefore, given $(\borel_1,\borel_2,\ldots,\borel_{\ell+1})\in X(\beta)$ we can consider the {\it unique} tuple
$$(\uni_1,\uni_2,\ldots,\uni_{\ell+1}),\quad \uni_i\in\G/\uni,\quad \mbox{s.t. }\uni_j\stackrel{s_{i_j}}{\uni_{j+1}}\mbox{ and }\uni_1=\uni.$$
\noindent By construction, $(\pi(\uni_1),\pi(\uni_2),\ldots,\pi(\uni_{\ell+1}))\in X(\beta)$.

\begin{remark}\label{rmk:unipotent_representative_lastflag}
This lifting procedure frames each flag in $(\borel_1,\borel_2,\ldots,\borel_{\ell+1})\in X(\beta)$ by first framing $\borel$ by $\uni$ and propagating the framing to the right. For instance if $\delta(\beta)=w_0$, this process frames the last flag $\borel_{\ell+1}=w_0\borel$, but the framing is {\it not} necessarily $w_0\uni$. Rather, the last framing is given by $tw\delta(\beta)\uni$ for some element $t\in\T$, but not necessarily the identity. This explains geometrically the map $\mu:X(\beta)\lr\T$ in \Cref{rmk:torus_quasimoment_map}.\sbs
\end{remark}

%%%%%%%%%%%%%%%%%%%%%%%%%%%%%%%%%%%%%%%%%%%%%
%%%%%%%%%%%%%%%%%%%%%%%%%%%%%%%%%%%%%%%%%%%%%
\subsection{Rational functions on $X(\beta)$ via weave edges}\label{ssec:rational_functions} Let $\beta\in\brn$ be a positive braid, $X(\beta)$ its braid variety, $\w:\beta\lr\delta(\beta)$ a Demazure weave and $T_\w\sse X(\beta)$ its weave torus. The Demazure weave $\w$ not only provides $T_\w$, it also provides a set of coordinates on $T_\w$. In order to explicitly describe them, we fix a pinning $(x_{i}, y_{i})_{i \in \dynkin}$ of the group $\G$, i.e.~ for every simple root indexed by $i \in \dynkin$ in its Dynkin diagram we choose two isomorphisms $x_{i}: \kk \to \uni_{i}^{+}$ and $y_{i}: \kk \to \uni_{i}^{-}$, where $\uni_{i}^{+}$ and $\uni_{i}^{-}$ are the corresponding root subgroups of $i \in \dynkin$, such that the assignment
\[
\left(\begin{matrix} 1 & z \\ 0 & 1 \end{matrix}\right) \mapsto x_{i}(z), \quad \left(\begin{matrix} 1 & 0 \\ z & 1 \end{matrix}\right) \mapsto y_{i}(z),\quad \left(\begin{matrix} t & 0 \\ 0 & t^{-1} \end{matrix}\right) \mapsto \chi_{i}(t)
\]
gives an algebraic group morphism $\varphi_{i}: \SL_{2}(k) \to \G$, with $\chi_{i}: \bG_m \to T$ the simple coroot corresponding to $i \in \dynkin$.\footnote{Recall that a simple algebraic group $\G$ admits a pinning and any two pinnings are conjugate.} Fixing a pinning is useful in our context because it allows us to explicitly measure flag transversality and framing differences in a quantifiable manner. Specifically, we have the following result:

\begin{lemma}\label{lem:framing_variables}
Let $g\uni,g'\uni\in\G/\uni$ be a pair of framed flags whose underlying flags are in $s_i$-relative position. If a representative $g$ of $g\uni$ is chosen, then there exist unique values $(\textcolor{purple}{\tilde{z}},\textcolor{blue}{u})\in\kk\times\kk^\times$ such that
$$g'\uni=gx_{i}(\textcolor{purple}{\tilde{z}})s_{i}\chi_i(\textcolor{blue}{u})\uni$$
where $s_{i} := \varphi_{i}\left(\begin{matrix}0 & -1 \\ 1 & 0 \end{matrix} \right)\in \G$, so $s_{i} \in N_{\G}(T)$ lifts the simple reflection corresponding to $i \in \dynkin$.\qed
\end{lemma}

\noindent The importance of \Cref{lem:framing_variables} is that it provides us with explicit values $\textcolor{purple}{\tilde{z}},\textcolor{blue}{u}$ that measure the transversality of two framed flags. That is, given a representative for a framed flag $g\uni$, the variable $\textcolor{purple}{\tilde{z}}$ parametrizes the choices of (unframed) flags which are $s_i$-transverse to the underlying flag $g\borel$, and the variable $\textcolor{blue}{u}$ parametrizes the choices of framing, cf.~\Cref{ex:framed_flags_SL2}.\\

Let us fix a representation of the points of $X(\beta)$ as tuples of framed flags, as in \Cref{ssec:framed_rep}, by choosing the constant functions $u_e=1$ at the top edges of the weave. Now we can use \Cref{lem:framing_variables} to construct more interesting rational functions, as follows. For any point in a weave edge $e\in\w$ labeled by $s_i$, we can consider a horizontal slice $h$ through that point. Given a point $(\borel_1,\ldots,\borel_{\ell_h+1})\in X(\beta_h)$ in the corresponding flag variety, the two flags $\borel_j,\borel_{j+1}$ to the left and right of the weave edge $e$ must be in $s_{i_j}$-relative position, where $\s_{i_j}$ is the $j$th crossing of $\beta_h$ reading left-to-right. Since points $X(\beta)$ are represented as tuples of framed flags, we can use the raking procedure in \cite[Section 5]{CGGLSS25} and obtain a morphism $X(\beta_h)\lr X(\beta)$, which is oftentimes a subvariety, and in this case $(\borel_1,\ldots,\borel_{\ell_h+1})\in X(\beta_h)$ is already lifted to a tuple of framed flags $(\uni_1,\ldots,\uni_{\ell_h+1})\in X(\beta_h)$. By writing $\uni_j=g\uni$ and $\uni_{j+1}=g'\uni$, and choosing a representative $g\in\G$, \Cref{lem:framing_variables} implies that 
there exist unique values $(\textcolor{purple}{\tilde{z_e}},\textcolor{blue}{u_e})\in\kk\times\kk^\times$ such that
$$g'\uni=gx_{i_j}(\textcolor{purple}{\tilde{z_e}})s_{i_j}\chi_{i_j}(\textcolor{blue}{u_e})\uni$$
The value of $\textcolor{blue}{u_e}$, as $(\uni_1,\ldots,\uni_{\ell_h})$ ranges over the allowed tuples of flags in $T_\w$, defines a rational function on the open torus $T_\w$. Thus, by implicitly including the raking corrections at the virtual vertices so as to always compare to $X(\beta)$, as in \cite[Section 5]{CGGLSS25}, we can write:

\begin{definition}[Rational function of a weave edge]\label{def:rational_function_weave_edge}
Let $\beta\in\brn$ be a positive braid word, $\w:\beta\lr\delta(\beta)$ a Demazure weave and $e\in E(\w)$ an edge labeled by $s_i$. By definition, and in the notation above, the rational function $\textcolor{blue}{u_e}:X(\beta)\lr\kk$ associated to $e$ is the unique rational function on $X(\beta)$ such that
\begin{equation}
g'\uni=gx_{i_j}(\tilde{z_e})s_{i_j}\chi_{i_j}(\textcolor{blue}{u_e})\uni
\end{equation}
for all pairs of flags $g\uni,g'\uni$ which appear left and right of $e$ in a configuration of flags that belongs to $T_\w$, where a representative $g$ of $g\uni$ is fixed, for some varying value of $\tilde{z_e}$.\sbs
\end{definition}

\noindent The rational functions in \Cref{def:rational_function_weave_edge} are typically not regular functions on $X(\beta)$, i.e.~they do not belong to $\kk[X(\beta)]$. That said, they already provide a coordinate system in $T_\w$:

\begin{lemma}[Coordinates in $T_\w$ from weave vertices]\label{lem:rational_functions_from_weave_vertices}
Let $\beta\in\brn$ be a positive braid and $\w:\beta\lr\delta(\beta)$ a Demazure weave with $d$ trivalent vertices. Consider the set of edges $e_1,\ldots,e_d$ with $e_i$ the south edge of the $i$th trivalent vertex of $\w$ and write $u_i:=u_{e_i}$. Then the set of rational functions $\{u_1,\ldots,u_d\}\in\kk(X(\beta))$ satisfies
$$T_\w\cong \Spec\kk[u_1^{\pm1},\ldots,u_d^{\pm1}].$$
That is, $u_1,\ldots,u_d$ is a coordinate system on the weave torus $T_\w$.
\end{lemma}

\begin{remark}
In \Cref{def:rational_function_weave_edge} and \Cref{lem:rational_functions_from_weave_vertices} there is a subtlety with regards to $u_e$ being well-defined on an entire edge $e$, rather than depending on the chosen point in $e$. Technically, $u_e$ is only well-defined on the edges of the refined framed weave. This is related to the discussion on raking from \cite[Section 5]{CGGLSS25}, as is the choice of $g$ representing $g\uni$, see specifically \cite[Definition 5.8.(5)]{CGGLSS25}. As said above, we treat this implicitly in these notes and refer to \cite{CGGLSS25} for details.\sbs
\end{remark}

%%%%%%%%%%%%%%%%%%%%%%%%%%%%%%%%%%%%%%%%%%%%%
%%%%%%%%%%%%%%%%%%%%%%%%%%%%%%%%%%%%%%%%%%%%%
\subsection{The cluster $A$-variables of a Demazure weave}\label{ssec:cluster_Avariables_weave} The object of this subsection is to define regular functions $A^\w_v\in\kk[X(\beta)]$ associated to the trivalent vertices of any Demazure weave $\w:\beta\lr\delta(\beta)$ via the dual Lusztig cycles introduced in \Cref{def:dual_Lusztig_cycles}. For that, let $\sigma$ be a trivalent or hexavalent vertex of $\w$ and denote by $e_{c}(\sigma)$, the weave edge incident to $\sigma$ in the cardinal directions $c\in\{\mathsf{nw},\mathsf{ne},\mathsf{s}\}$ of northwest, northeast and south respectively. Similarly denote by $u_{c}(\sigma)$ the rational function associated to $e_{c}(\sigma)$ via \Cref{def:rational_function_weave_edge}.

\begin{fdef}[Cluster $A$-variables of $\w$]\label{def:cluster_Avariables_weave}
Let $\beta\in\brn$ be a positive braid word, $\w:\beta\to\delta(\beta)$ a Demazure weave, $v\in\w$ a trivalent vertex and $\rho_v$ its dual Lusztig cycle. By definition, the rational function $A^\w_v:X(\beta)\lr\PP_k^1$ associated to $v$ is the product
\begin{equation}\label{eq:cluster_A_variables}
\boxed{
A^\w_v:=\prod_{\sigma\in \w_3\cup \w_6}L_\sigma(\rho_v).
}
\end{equation}
where we have denoted by $\w_3$ and $\w_6$ the sets of 3-valent and 6-valent vertices of $\w$ and each factor $L_\sigma(\rho_v)$ is given by
$$L_\sigma(\rho_v):=
\begin{cases}
\displaystyle
\left(\frac{u_{\rm s}}{u_{\rm nw}}\right)^{\rho_v(e_{\rm ne})-\rho_v(e_{\rm nw})},
&\mbox{ if }\rho_v(e_{\rm ne})\ge \rho_v(e_{\rm nw}),\\[1.2em]
\displaystyle
\left(\frac{u_{\rm ne}}{u_{\rm s}}\right)^{\rho_v(e_{\rm nw})-\rho_v(e_{\rm ne})},
&\mbox{ if }\rho_v(e_{\rm nw})\ge \rho_v(e_{\rm ne}).
\end{cases}$$
where we abbreviated $e_{\rm c}:=e_{\rm c}(\sigma)$ and $u_{c}:=u_c(\sigma)$ for $c\in\{\mathsf{nw},\mathsf{ne},\mathsf{s}\}$.\sbs
\end{fdef}
\vspace{0.2cm}
\Cref{eq:cluster_A_variables} makes it apparent that $A^\w_v$ is constructed by using solely the information of the dual Lusztig cycles $\rho_v$ and the rational functions $u_e$ associated to the weave edges $e\in\w$. If $\w$ is implicitly understood by context or fixed, we often abbreviate $A_v:=A^\w_v$. Given a Demazure weave $\w$, it can be useful to notice that it is possible to write each such $u_e$ in terms of the collection $\{A_v\}$ and the Lusztig cycles $\g_v$ of $\w$. Indeed, it follows from \cite[Theorem 5.12]{CGGLSS25} that
\begin{equation}\label{eq:u_via_cluster_A_variables}
\boxed{
u_e=\prod_{v\in\w_3}A_v^{\g_v(e)}.
}
\end{equation}

\begin{remark}\label{rmk:Avars_different_defs} The presentation of the cluster $A$-variables $\{A_v\}$ in \Cref{def:cluster_Avariables_weave} differs from the original presentation $A_v$ in \cite[Section 5]{CGGLSS25}. They can be shown to be equivalent.\sbs
\end{remark}

\Cref{eq:cluster_A_variables} does not make it apparent that $A_v$ is actually a regular function on $X(\beta)$. It is also important to understand their behavior under weave equivalence and weave mutation. The following facts are proven in \cite[Section 5.2 \& 5.3]{CGGLSS25}:

\begin{fthm}\label{thm:cluster_variables_weave}
Let $\beta\in\brn$ be a positive braid word, $X(\beta)$ its braid variety and $\w:\beta\lr\delta(\beta)$ a Demazure weave with $d$ trivalent vertices $v_1,\ldots,v_d$. Then the following hold:

\begin{enumerate}
    \item For every trivalent vertex $v\in\w$, the rational function $A_v\in\kk(X(\beta))$ is a regular function. That is, $A_v\in\kk[X(\beta)]$, so $A_v:X(\beta)\lr\mathbb{A}^1_\kk$ has no poles.\\

    \item For every trivalent vertex $v\in\w$, the rational function $A_v\in\kk[X(\beta)]$ is non-vanishing if and only if at least one of the weights of the Lusztig cycle $\g_v$ at the south boundary of $\w$ does not vanish.\\

    \item The restrictions of $A_{v_1},\ldots,A_{v_d}$ to the weave torus $T_\w\sse X(\beta)$ are a coordinate system, i.e.~
    $$T_\w\cong \Spec\kk[A_{v_1}^{\pm1},\ldots,A_{v_d}^{\pm1}].$$
    Thus $T_\w\sse X(\beta)$ is cut-out as an open subvariety as the non-vanishing of all the regular functions $A_{v_i}$.\\

    \item If $\w_1\simeq\w_2$ is a weave equivalence, the two sets of regular functions $\{A^{\w_1}_{v_1},\ldots,A^{\w_1}_{v_d}\}$ and $\{A^{\w_2}_{v_1},\ldots,A^{\w_2}_{v_d}\}$ coincide.\\

    \item If $\w_1\stackrel{\mu}{\sim}\w_2$ is a weave mutation at the Lusztig cycle $\gamma_{v_k}$ of $\w_1$, then $A^{\w_1}_{v_i}=A^{\w_2}_{v_i}$ if $i\neq k$ and
    \begin{equation}\label{eq:cluster_mutation_weave_variables}
        \displaystyle A^{\w_1}_{v_k}A^{\w_2}_{v_k}=\prod_{\langle \g_i,\g_k\rangle>0} (A^{\w_1}_{v_i})^{\langle \g_i,\g_k\rangle}+\prod_{\langle \g_i,\g_k\rangle<0} (A^{\w_1}_{v_i})^{-\langle \g_i,\g_k\rangle}.
    \end{equation}
\end{enumerate}
\end{fthm}
\vspace{0.2cm}
The identity \eqref{eq:cluster_mutation_weave_variables} in Theorem \ref{thm:cluster_variables_weave}.(5) is known as the cluster exchange relation, cf.~\cite[Definition 3.1.3]{FWZ16}. In particular, if we denote by $\SA(Q)$ the commutative cluster algebra associated to a quiver $Q$, cf.~\cite[Def.~3.1.9]{FWZ16}, we obtain the inclusion
\begin{equation}\label{eq:inclusion_cluster_algebra}
\SA(Q_\w)\sse \kk[X(\beta)]
\end{equation}
where $Q_\w$ is the intersection quiver of Lusztig cycles of any Demazure weave $\w:\beta\to\delta(\beta)$. Specifically, any Demazure weave $\w:\beta\to\delta(\beta)$ provides a cluster seed for $\SA(Q_\w)$ as a subalgebra of $\kk[X(\beta)]$: the cluster seed is given by the torus $T_\w\sse X(\beta)$ endowed with the quiver $Q_\w$ and the set of regular functions $\{A^\w_v\}$, indexed by the trivalent vertices $v\in\w_3$ of $\w$. We will soon see that \eqref{eq:inclusion_cluster_algebra} is an equality.

%%%%%%%%%%%%%%%%%%%%%%%%%%%%%%%%%%%%%%%%%%%%%
%%%%%%%%%%%%%%%%%%%%%%%%%%%%%%%%%%%%%%%%%%%%%
%%%%%%%%%%%%%%%%%%%%%%%%%%%%%%%%%%%%%%%%%%%%%
\subsection{Propagation of weave cluster $A$-variables after adding weave vertex on top}\label{ssec:cluster_A_variables_propagation_vertices} Let $\beta\in\brn$ be a positive braid word of length $\ell=\ell(\beta)$, and $\w:\beta\lr\delta(\beta)$ a Demazure weave. In \Cref{def:cluster_Avariables_weave} we constructed the cluster $A$-variables $A_v^\w\in\kk[X(\beta)]$ in terms of the dual Lusztig cycles $\{\rho_v\}$ of the weave $\w$. Our task now is to analyze how $\{A_v^\w\}$ changes as we modify the braid word $\beta$ according to a braid move and $\s_i^2\lr\s_i$, i.e.~using a weave vertex, cf \Cref{fig:weaves_types_vertices}.\\

For that, let $\w^+$ be a weave obtained from $\w$ by adding exactly one weave vertex $p\in(\w^+\setminus\w)$ to the north boundary of $\w$. Let $h_\pm$ be horizontal slices exactly above and below this new vertex $p$, cf.~\Cref{fig:BraidVariety_DiagramMoves2} and \Cref{fig:BraidVariety_Trivalent2}, and $\beta_\pm$ the corresponding braid words. By construction, $\beta_-=\beta$ and so we denote by $\{A_v\}$ the cluster variables of $\w$, as they are elements of $\kk[X(\beta)]$. The cluster variables of $\w^+$ are denoted by $\{A_v^+\}$ and, if $p$ is a trivalent vertex the additional cluster variable is denoted by $A^+_p$. This notation is consistent with the fact that $\w^+$ coincides with $\w$ below $p$, and thus the trivalent vertices of $\w$ and $\w^+$ below $p$ are identical, and can be indistinguishably referred to by using the letter $v$. We then have a relation between $X(\beta)=X(\beta_-)$ and $X(\beta_+)$ depending on the type of weave vertex $p$:

\begin{enumerate}
\item {\bf Trivalent.} If $p$ is trivalent, then \Cref{lem:braidvarieties_trivalent}, or see \Cref{fig:BraidVariety_Trivalent2} (right), provides an embedding
$$\iota:X(\beta)\times\bG_m\hookrightarrow X(\beta_+).$$
\Cref{prop:Avariables_up_new_trivalent} below describes how to relate $\{A_v\}$ to the restriction of $\{A^+_v\}$ and $A^+_p$ to the image $\mbox{im}(\iota)$.\\

\item {\bf Hexavalent.} If $p$ is hexavalent, then \Cref{lem:braid_invariance}.(1) provides an isomorphism
$$R3:X(\beta)\stackrel{\cong}{\longrightarrow} X(\beta_+).$$
It is a cluster isomorphism but, importantly, it is {\it not} the case that the pull-back of $\{A_v^+\}$ under $R3$ yields the same initial cluster seed from $\w$. Rather, the pull-back of $(T_{\w^+},Q_{\w^+},\{A_v^+\})$ under $R3$ gives a cluster seed in $\kk[X(\beta)]$ obtained by applying certain mutations to the cluster seed $(T_\w,Q_\w,\{A_v\})$ associated to $\w$, cf.~\cite[Remark 2.7]{casals2025comparingclusteralgebrasbraid}. The local model for this sequence of mutations is already present in the theory of wiring diagrams, see e.g.~\cite[Prop.~3.2.3]{FWZ16}.\\

\item {\bf Tetravalent.} If $p$ is tetravalent, then \Cref{lem:braid_invariance}.(2) provides an isomorphism
$$C:X(\beta)\stackrel{\cong}{\longrightarrow} X(\beta_+).$$
It is a cluster isomorphism that satisfies $C^*A_v^+=A_v$.\\
\end{enumerate}

\noindent In conclusion, the genuinely new phenomenon is at a trivalent vertex, so let us now assume $p\in\w^+$ is a trivalent vertex. By \Cref{lem:braidvarieties_trivalent} and \Cref{def:cluster_Avariables_weave}, the image of $\iota:X(\beta)\times\bG_m\hookrightarrow X(\beta_+)$ coincides with
$$\mbox{im}(\iota)=\{A^+_p\neq0\}$$
as subvarieties of $X(\beta_+)$, with the second factor identified as $\bG_m=\Spec[(A_p^+)^{\pm1}]$. Therefore, we can consider the projection $\pi:X(\beta)\times\bG_m\lr X(\beta)$ onto the first factor and ask how $\iota^*A^+_v$ relate to $\pi^*A_v$, as a comparison of regular functions on $X(\beta)\times\bG_m$.\\

\noindent For specificity, let us write $\beta=\beta_1\s_{i_k}\beta_2$ and $\beta_+=\beta_1\s_{i_k}\s_{i_k}\beta_2$, so that the south edge of the trivalent vertex $p$ corresponds to the $k$th crossing $\s_{i_k}$ of $\beta$. In line with above, we denote the trivalent vertices of $\w$ with the letter $v$, and use $p$ for the top trivalent vertex of $\w^+$, which is the only new trivalent of $\w^+$ (not in $\w$). Denote the north boundary edges of $\w$, corresponding to crossings of $\beta_1\s_{i_k}\beta_2$, by $e_1^{top},\ldots,e_\ell^{top}$. By definition, we let $\gamma_p|_\w$ be the restriction of the Lusztig cycle of $p$ in $\w^+$ to the subweave $\w$, i.e.~ $\gamma_p|_\w$ is the unique weave cycle on $\w$ which has weights $\gamma_p|_\w(e_i^{top})=\delta_{ik}$ and otherwise satisfies the tropical Lusztig rules at all the vertices of $\w$. Let us refer to the cluster variables of $\w^+$ associated to trivalent vertices $v$ also in $\w$ as $\{A_v^+\}$.\\

Finally, let $\{\alpha_i\}$ and $\{\alpha_i^\vee\}$ be the set of simple roots and simple coroots of $\G$, indexed by the vertices in $\dynkin$. Given the braid word $\beta=\sigma_{i_1}\cdots\sigma_{i_\ell}$, we need the sequences of root and coroots from \cite[Equation (16)]{CGGLSS25}, namely:
\begin{equation}\label{eq:roots_coroots}
 \vartheta_j:=s_{i_1}\cdots s_{i_{j-1}}(\alpha_{i_j}),
 \qquad
 \vartheta_j^\vee:=s_{i_1}\cdots s_{i_{j-1}}(\alpha_{i_j}^\vee).
\end{equation}

\begin{proposition}\label{prop:Avariables_up_new_trivalent} In the notation above, the cluster $A$-variables of $\w$ and $\w^+$ satisfy
\begin{equation}\label{eq:Avariables_up_new_trivalent}
\iota^*A^+_v=(\iota^*A_p^+)^{c_v}\pi^*A_v,
\end{equation}
\noindent where the exponent in \eqref{eq:Avariables_up_new_trivalent} is
\begin{equation}\label{eq:c_exponent}
 \boxed{
 c_v:=\rho_v(e_k^{top})
 +\langle\rho_v,\gamma_p|_\w\rangle
 +\sum_{k<j}(\vartheta_j,\vartheta_{k}^\vee)\rho_v(e_j^{top})}
\end{equation}
\sbs
\end{proposition}

\noindent The following clarifications on \Cref{prop:Avariables_up_new_trivalent} might be helpful:\\

\begin{enumerate}
    \item The exponent $c_v$ might be negative, zero, or positive. This is not an issue as $A_p^+$, which is typically not a unit in $\kk[X(\beta_+)]$, restricts to a unit in $\mbox{im}(\iota)$ which is where the identity \eqref{eq:Avariables_up_new_trivalent} takes place.\footnote{In particular, it is not a contradiction with the fact that $A_v^+$ is an irreducible function on $X(\beta_+)$.}\\

    \item The expression $\langle\rho_v,\gamma_p|_\w\rangle$ in \eqref{eq:c_exponent} refers to the intersection number as introduced in \Cref{def:intersection_number}, without any boundary intersection terms.\\

    \item The sum
    $$\sum_{k<j}(\vartheta_j,\vartheta_{k}^\vee)\rho_v(e_j^{top})$$
    in \eqref{eq:c_exponent} appears because of the raking procedure in \cite[Section 5.1]{CGGLSS25}, which also explains why only contributions to the right of $p$ are relevant. It can be conceptually understood as a boundary intersection contribution, cf.~\Cref{rmk:boundary_intersections} and \cite[Definition 4.24]{CGGLSS25}.\\
    
    \item The index $k$ in \eqref{eq:c_exponent} is fixed and independent of $v$. It is fixed by the location of the vertex $p$, as we have written $\beta=\beta_1\s_{i_k}\beta_2$ and $\beta_+=\beta_1\s_{i_k}^2\beta_2$ above.
\end{enumerate}

%%%%%%%%%%%%%%%%%%%%%%%%%%%%%%%%%%%%%%%%%%%%%
%%%%%%%%%%%%%%%%%%%%%%%%%%%%%%%%%%%%%%%%%%%%%
\subsection{The cluster $X$-variables of a Demazure weave}\label{ssec:cluster_Xvariables_weave} The object of this subsection is to define rational functions $X^\w_v\in\kk(X(\beta))$ associated to the trivalent vertices of any Demazure weave $\w:\beta\lr\delta(\beta)$ via the Lusztig cycles introduced in \Cref{def:Lusztig_cycles}, as a dual version of \Cref{ssec:cluster_Avariables_weave}. For that, let us denote by $e_1^{b},\ldots,e_r^{b}$ the south edges of $\w:\beta\to\delta(\beta)$, where $r=\ell(\delta(\beta))$ and consider the sequence of roots $\vartheta_j$ in \eqref{eq:roots_coroots} but associated to the permutation $\delta(\beta)$. Given a Lusztig cycle $\g_v$ of $\w$, we consider the following rational function on the $u$-variables:

\begin{equation}\label{eq:boundary_term_Xvariables}
B_\dd(\g_v):=
\displaystyle
\prod_{j=1}^{\ell(\delta(\beta))}u(e_j^{b})^{-\g_v(e_j^{b})-R_v^{(j)}},\qquad\mbox{where }R_v^{(j)}:=\displaystyle\sum_{j<k}\g_v(e_k^{b})(\vartheta_k,\vartheta_j).
\end{equation}
\noindent The term \eqref{eq:boundary_term_Xvariables} is a boundary contribution, in line with \Cref{rmk:boundary_intersections}. It is only possibly non-trivial for frozen Lusztig cycles, that is, $B_\dd(\g_v)=1$ for any mutable Lusztig cycle $\g_v$.\\

\begin{fdef}[Cluster $X$-variables of $\w$]\label{def:cluster_Xvariables_weave}
Let $\beta\in\brn$ be a positive braid word, $\w:\beta\to\delta(\beta)$ a Demazure weave, $v\in\w$ a trivalent vertex and $\g_v$ its Lusztig cycle. By definition, the rational function $X^\w_v:X(\beta)\lr\PP_k^1$ associated to $v$ is the product
\begin{equation}\label{eq:cluster_X_variables}
\boxed{
X^\w_v:=B_\dd(\g_v)\prod_{\sigma\in \w_3\cup \w_6}\mathcal{L}_\sigma(\g_v).
}
\end{equation}
where we have denoted by $\w_3$ and $\w_6$ the sets of 3-valent and 6-valent vertices of $\w$ and each factor $\mathcal{L}_\sigma(\g_v)$ is given by
$$\mathcal{L}_\sigma(\g_v):=
\displaystyle
u_{\rm nw}^{\g_v(e_{\rm ne})-\g_v(e_{\rm s})}u_{\rm s}^{\g_v(e_{\rm nw})-\g_v(e_{\rm ne})}u_{\rm ne}^{\g_v(e_{\rm s})-\g_v(e_{\rm nw})},$$
where we abbreviated $e_{\rm c}:=e_{\rm c}(\sigma)$ and $u_{c}:=u_c(\sigma)$ for $c\in\{\mathsf{nw},\mathsf{ne},\mathsf{s}\}$.\sbs
\end{fdef}
\vspace{0.2cm}

In order to relate the cluster $A$-variables from \Cref{def:cluster_Avariables_weave} to the cluster $X$-variables in \Cref{def:cluster_Xvariables_weave}, the exchange matrix $B=(B_{ij})$ associated to $Q_\w$ can be modified to the $p$-map, as in \cite[Section 8.2]{CGGLSS25}, between the cluster $A$-variety and the cluster $X$-variety. In short, we can construct a matrix $(p_{ij})$ such that:

\begin{enumerate}
\item The entries $B_{ij}=p_{ij}$ coincide if one of the vertices associated to $i$ and $j$ are mutable.
    \item The entries $p_{ij}\in\Z$ are integers.
    \item $(p_{ij})$ is unimodular, i.e.~$\det((p_{ij}))=\pm1$.
\end{enumerate}

\noindent Note nevertheless that $p$ might have non-zero diagonal entries due to the frozen vertices. We refer to \cite[Section 8.2]{CGGLSS25} for the precise expression of the correct $(p_{ij})$ as $B_{ij}$ plus a (half) correction summand. \\

\begin{lemma}\label{lem:X_vars_via_A_vars} Let $\w:\beta\to\delta(\beta)$ be a Demazure weave and $v_i\in\w$ its $i$th trivalent vertex. Then
\begin{equation}
\label{eq:X_vars_via_A_vars}
\boxed{\pi^*X^\w_{v_i}= \prod_{v_j\in\w_3} (A^\w_{v_j}) ^{p_{ij}}}
\end{equation}
where $(p_{ij})$ is a unimodular integral matrix with $p_{ij}=\langle\g_i,\g_j\rangle$ if $\g_i$ or $\g_j$ are mutable Lusztig cycles.\sbs
\end{lemma}

\noindent At this stage, we read \eqref{eq:X_vars_via_A_vars} as an equality of rational functions on $X(\beta)$ but see \Cref{ex:X_in_terms_of_A} below. The analogue of \Cref{ssec:cluster_A_variables_propagation_vertices} and \Cref{prop:Avariables_up_new_trivalent} for cluster $X$-variables is simpler. The behavior under adding a hexavalent and tetravalent vertex on top, as in \Cref{ssec:cluster_A_variables_propagation_vertices}.(2) and (3) is identical to that of the cluster $A$-variables. In contrast, since cluster $X$-variables are constructed from Lusztig cycles $\{\g_v\}$, rather than dual Lusztig cycles $\{\rho_v\}$, the cluster $X$-variables under adding a trivalent vertex on top do not significantly change. Specifically, the analogue of \Cref{eq:Avariables_up_new_trivalent} for the cluster $X$-variables is
\begin{equation}\label{eq:Xvariables_up_new_trivalent}
 \boxed{\iota^*X^+_v=X_v}\quad\mbox{if }v\neq p,
\end{equation}
and one can express the pull-back under $\iota$ of the unique new cluster $X$-variable $X^+_p$ in $\w^+$ as
$$\boxed{\iota^*X^+_p=(\iota^*A_p^+)^{-1}\prod_{v\in\w_3}(\pi^*X_v)^{q_v}}$$
for some exponents $q_v$ that depend on the extended $p$-map $(p_{ij})$. Namely, if we extend the $p$-map $(p_{ij})$ of $\w$ to a $p$-map for $\w^+$ by adding\footnote{Let us use the convention here that the new vertex is added with the last index.} a row to the bottom of $(p_{ij})$ and a column to its right, then $q=(q_v)$ is the vector obtained by multiplying the bottom row of the extended $p$-map for $\w^+$ by the inverse of the $p$-map matrix $(p_{ij})$ of $\w$. Here, as above, all equalities are understood in the locus $\{A_p^+\neq0\}\sse X(\beta_+)$.

%%%%%%%%%%%%%%%%%%%%%%%%%%%%%%%%%%%%%%%%%%%%%
%%%%%%%%%%%%%%%%%%%%%%%%%%%%%%%%%%%%%%%%%%%%%
\subsection{The cluster schemes of a Demazure weave}\label{ssec:A_and_X_schemes_weave} Building on \cite{FG2006,FG09}, the $\SA$-scheme and $\SX$-schemes are constructed from seed data in \cite[Section 2]{GHK15}. Given a Demazure weave $\w:\beta\to\delta(\beta)$, we denote by $\SA_\beta$ and $\SX_\beta$ the $\SA$-scheme and $\SX$-schemes associated to the seed data associated to $\w$. This seed data is typically encoded by the intersection quiver $Q_\w$ of Lusztig cycles. Intuitively, these schemes are obtained by gluing tori, see e.g.~\cite[Appendix A]{GHKK18}, indexed by the cluster $\SA$-seeds and the cluster $\SX$-seeds respectively:

$$\displaystyle\SA_\beta:=\bigcup_{{\bf s}\in S(\SA,\w)} T_{{\bf s}}^\SA,\qquad\qquad\displaystyle\SX_\beta:=\bigcup_{{\bf s}\in S(\SX,\w)} T^\SX_{{\bf s}}$$

\noindent where we have denoted by $S(\SA,\w)$ and $S(\SX,\w)$ the sets of all seeds obtained from the seed associated to $\w$ via a (possibly empty) sequence of cluster $\SA$-mutations and $\SX$-mutations, respectively. We refer to \cite[Prop.~2.4]{GHK15} for details on the gluing. As in \cite[Remark 2.7.$(iii)$]{GHK15}, there is a map $p:\SA_\beta\lr\SX_\beta$ which is compatible with the local maps between the respective tori. For the cluster structures associated to Demazure weaves, it follows from \cite[Section 8]{CGGLSS25} that there are scheme isomorphisms
\begin{equation}\label{eq:scheme_iso}
\Spec \Gamma(\SA_\beta,\mathcal{O}_{\SA_\beta}) \cong X(\beta),\qquad \Spec \Gamma(\SX_\beta,\mathcal{O}_{\SX_\beta}) \cong X(\beta).
\end{equation}
\noindent Thus, at the same time, the braid variety $X(\beta)$ can be considered as the affinization of the cluster $\SA$-scheme $\SA_\beta$, or as the affinization of the cluster $\SX$-scheme $\SX_\beta$. In the former viewpoint, the emphasis is on the fact that the coordinate ring $\kk[X(\beta)]$ is the cluster algebra $\SA(Q_\w)$, which in this case happens to be isomorphic to the upper cluster algebra $\Gamma(\SA_\beta,\mathcal{O}_{\SA_\beta})$. In the latter sense, the emphasis is on the fact that $X(\beta)$ admits a cluster Poisson structure, i.e.~a Poisson bracket which is given on each seed by the given skew-symmetric bilinear form. When we want to emphasize $X(\beta)$ as the affinization of the cluster $\SA$-scheme $\SA_\beta$, respectively the cluster $\SX$-scheme $\SX_\beta$, we will denote it by $\SA(\beta)$, respectively $\SX(\beta)$.\\

\begin{example}\label{ex:X_in_terms_of_A} For instance, given the map $p:\SA(\beta)\lr\SX(\beta)$, identity \eqref{eq:X_vars_via_A_vars} should be understood as expressing the pull-back of the $\SX$-variables under $p$ and expressing the resulting rational function in terms of the cluster $\SA$-variables:
$$p_\w^*(X^\w_{v_i}) := \prod_{v_j\in\w_3} (A^\w_{v_j}) ^{p_{ij}},$$
where $p_\w:T_\w^\SA\lr T_\w^\SX$ is the restriction of $p$ to the corresponding seed torus.\sbs
\end{example}

More generally, to ease notation, if we are discussing cluster $\SA$-variables, then $X(\beta)$ will implicitly mean $\SA(\beta)$, and if we discuss cluster $\SX$-variables, then $X(\beta)$ will refer to $\SX(\beta)$. If we discuss both, we will explicitly use $\SA(\beta)$ and $\SX(\beta)$. 

\begin{remark}\label{rmk:embedding_not_iso}
By \eqref{eq:scheme_iso}, there are natural embeddings
\begin{equation}\label{eq:embedding_not_iso}
\SA_\beta\hookrightarrow \SA(\beta),\qquad \SX_\beta\hookrightarrow \SX(\beta),
\end{equation}
but they are {\bf not} typically isomorphisms. Indeed, there are points in $X(\beta)$ that are not contained in any cluster $\SA$-seed or $\SX$-seed torus for many interesting positive braid words $\beta$. For instance, consider $\beta=\s_1^3$, then the point $(\borel,\borel_-,\borel,\borel_-)\in X(\beta)$ does not belong to any cluster seed, and so it is not in the image of either map in \eqref{eq:embedding_not_iso}. See \cite{CGSS_clusterdeeploci} for further discussions on this matter.\sbs
\end{remark}

In addition to the $\SA$ versus $\SX$ distinction, we recall that a Lie group $\G$ has been fixed in order to be able to consider the braid variety $X(\beta)$. Let us denote the braid variety $X(\beta)$ by $X_\G(\beta)$ when we want to emphasize the choice of $\G$. Following the convention above, $\SA_\G(\beta)$ denotes $X_\G(\beta)$ when seen as the affinization of the cluster $\SA$-scheme, and $\SX_\G(\beta)$ denotes $X_\G(\beta)$ when seen as the affinization of the cluster $\SX$-scheme. The notational emphasis on $\G$ is relevant in the context of cluster ensembles and Fock-Goncharov duality, cf.~\cite[Sections 1.2.9 \& 1.2.10]{FG09}, \cite[Definition A.4]{GHKK18} and \cite[Section 7]{GHK15}. Specifically, given $\G$, we denote by $\G^\vee$ the Langlands dual group, see e.g.~\cite[Appendix A]{GHKK18} or \cite[Section 1.17]{FG2006}. For instance if $\G=\SL_n$ then $\G^\vee=\operatorname{PSL}_n$. Then \cite[Theorem 8.8]{CGGLSS25} proves the following duality result on braid varieties:

\begin{theorem}[Fock-Goncharov duality for braid varieties]\label{thm:fock_goncharov_duality} Let $\beta$ be a positive braid word, $\G$ a simply-laced simple algebraic Lie group and $\G^\vee$ its Langlands dual. Then both
\begin{equation}\label{eq:cluster_ensembles}
(\SA_\G(\beta),\SX_{\G^\vee}(\beta))\qquad\qquad (\SA_{\G^\vee}(\beta),\SX_{\G}(\beta))
\end{equation}
are cluster ensembles, and Fock-Goncharov cluster duality holds for each of them.\sbs
\end{theorem}

\noindent For instance, \Cref{thm:fock_goncharov_duality} implies that the cluster algebra $\kk[\SA_\G(\beta)]$ admits a canonical basis of theta
functions parameterized by the integral tropical points of the dual braid variety $\SX_{\G^\vee}(\beta)$. The key aspect of the proof of \Cref{thm:fock_goncharov_duality} is showing that there exists a Donaldson-Thomas transformation\footnote{Combinatorially, that there exists a reddening sequence for the quiver, cf.~\cite{Keller2017_DT,Keller12_DT2}.} and that the exchange matrices from Demazure weaves are full rank (for the extended exchange matrix with mutable columns). In fact, there is a rather explicit map $D_\beta:X(\beta)\lr X(\beta)$, generalizing the usual twist automorphisms, that realizes the Donaldson-Thomas transformation, cf.~\cite[Theorem 8.7]{CGGLSS25}.

%%%%%%%%%%%%%%%%%%%%%%%%%%%%%%%%%%%%%%%%%%%%%
%%%%%%%%%%%%%%%%%%%%%%%%%%%%%%%%%%%%%%%%%%%%%
%%%%%%%%%%%%%%%%%%%%%%%%%%%%%%%%%%%%%%%%%%%%%
%%%%%%%%%%%%%%%%%%%%%%%%%%%%%%%%%%%%%%%%%%%%%
\subsection{Brick compactifications of a braid variety}\label{ssec:brick_compactifications} Let $\beta=\s_{i_1}\cdots\s_{i_\ell}$ be a positive braid word and $X(\beta)$ its braid variety. Following \cite[Section 4]{CGGS2}, we now construct a smooth projective compactification $Y(\beta)$ of $X(\beta)$. As we will explain, such compactification $Y(\beta)$ has a number of desired properties both from the combinatorial viewpoint and as it relates to aspects of mirror symmetry, cf.~\Cref{thm:superpotential_Reeb_coordinates}. Importantly, $Y(\beta)$ depends on the choice of positive braid word $\beta$ and not just on its underlying braid $[\beta]$. Therefore, since the isomorphism type of $X(\beta)$ only depends on the underlying braid $[\beta]$, we can construct several well-behaved compactifications of $X(\beta)$ by using different positive braid words for $[\beta]$.\\

\noindent Given a positive braid word $\beta$, the intuition for compactifying $X(\beta)$ is simple: we allow pairs of consecutive flags in \Cref{def:braid_variety} to be equal, not just in a certain relative position. This intuition is formalized as follows.

%%%%%%%%%%%%%%%%%%%%%%%%%%%%%%%%%%%%%%%%%%%%%
%%%%%%%%%%%%%%%%%%%%%%%%%%%%%%%%%%%%%%%%%%%%%
\subsubsection{The brick compactification $Y(\beta)$}\label{sssec:brick_compactification} Let \(P_i\) denote the minimal parabolic subgroup associated with
\(s_i\) which contains the fixed Borel $\borel\sse\G$. For instance, if $\G=\SL_n$ and $\borel$ is the subgroup of invertible upper-triangular matrices, then
$$P_i\cong\{(M_{jk})\in\SL_n: M_{jk}=0\mbox{ if }k<j\mbox{ and }(j,k)\neq(i+1,i)\}.$$
That is, the subgroup of invertible upper-triangular matrices where, in addition, the below-diagonal entry $(i+1,i)$ is allowed to be non-zero.\\

By definition, the closed Bott--Samelson variety associated to the braid word $\beta=\s_{i_1}\ldots\s_{i_{\ell}}$ is
\begin{equation}\label{eq:BSquotient}
 \BS(\beta)=(P_{i_1}\times_\borel P_{i_2}\times_\borel\cdots\times_\borel P_{i_{\ell}})/\borel^{\ell},
\end{equation}
where $\borel^{\ell+1}$ acts on the right by
$$(p_1,\ldots,p_{\ell})\cdot (b_1,\ldots,b_{\ell}):=(p_1b_1,b_1^{-1}p_2b_2,\ldots,b_{\ell-1}^{-1}p_{\ell} b_{\ell}).$$
That is, $\BS(\beta)$ parametrizes the following tuples of flags:
\begin{equation}\label{eq:BSgallery}
 \BS(\beta)=\Bigl\{(\borel_1,\ldots,\borel_{\ell+1}):
 \borel_1=\borel,\quad \borel_{k} \buildrel s_{i_{k}}\mbox{\tiny or }e \over \longrightarrow \borel_{k}\}.
\end{equation}
For any braid word $\beta$, $\BS(\beta)$ is a smooth projective variety\footnote{In fact, it is an iterated $\PP^1$-bundle, and thus smooth and projective.} of dimension $\dim_\kk\BS(\beta)=\ell$. Consider the endpoint map
\begin{equation}\label{eq:endpoint}
  m_\beta:\BS(\beta)\longrightarrow \G/\borel,
  \qquad (\borel_1,\ldots,\borel_{\ell+1})\longmapsto \borel_{\ell+1}.
\end{equation}
For instance, the image of $m_\beta$ is the classical Schubert variety $\overline{B\dot uB/B}$ with $u=\delta(\beta)$.

\begin{definition}[The compactification $Y(\beta)$]\label{def:brick}
Let $\beta=\s_{i_1}\cdots\s_{i_\ell}$ be a positive braid word. By definition, the brick compactification $Y(\beta)$ associated to $\beta$ is
\begin{equation}\label{eq:brickdef}
  Y(\beta):=m_\beta^{-1}(\delta(\beta)\borel).
\end{equation}
where $\delta(\beta)\in W(\G)$ has been lifted to $N_\G(\T)$. Equivalently,
\begin{equation}\label{eq:brickgallery}
Y(\beta):=\Bigl\{(\borel_1,\ldots,\borel_{\ell+1}):
 \borel_1=\borel, \borel_{k} \buildrel s_{i_{k}}\mbox{\tiny or }e \over \longrightarrow \borel_{k+1},\borel_{\ell+1}=\delta(\beta)\borel\}.
\end{equation}
\sbs
\end{definition}

These varieties were introduced and studied in \cite{Escobar}. It is proven in \cite[Theorem 20]{Escobar} that $Y(\beta)$ is a smooth irreducible projective variety of dimension $\dim_\kk Y(\beta)=\ell(\beta)-\ell(\delta(\beta))$, which coincides with the dimension of $X(\beta)$. By comparing \Cref{def:brick}, especially \eqref{eq:brickgallery}, with \Cref{def:braid_variety} it thus follows that $X(\beta)\sse Y(\beta)$ is an open (and dense) subvariety.

%%%%%%%%%%%%%%%%%%%%%%%%%%%%%%%%%%%%%%%%%%%%%
%%%%%%%%%%%%%%%%%%%%%%%%%%%%%%%%%%%%%%%%%%%%%
\subsubsection{The boundary divisors}\label{sssec:dual_boundary_complex}
The complement $Y(\beta)\setminus X(\beta)$ is a union of irreducible boundary divisors $\{D_i\}$. Specifically, we can consider
\begin{equation}\label{eq:Djdef}
  D_j:=\{(\borel_1,\ldots,\borel_{\ell+1})\in Y(\beta):\borel_{j}=\borel_{j+1}\},\quad j\in[1,\ell].
\end{equation}
Let us denote by $\beta\setminus \s_{i_j}$ the positive braid word resulting from deleting the $j$th crossing $\s_{i_j}$ from $\beta$. More generally, if $S$ is a subset of crossings of $\beta$, we denote by $\beta\setminus S$ the positive braid word resulting from deleting the crossings in $S$ from $\beta$. The following holds:

\begin{lemma}[Boundary divisors]\label{lem:boundary_div} Let $\beta$ be a positive braid word, $X(\beta)$ its braid variety and $Y(\beta)$ its brick compactification. Then:
\begin{enumerate}
    \item $D_j\ne\emptyset
  \quad\Longleftrightarrow\quad
  \delta(\beta\setminus \s_{i_j})=\delta(\beta)$\\

  \item If $\delta(\beta\setminus \s_{i_j})=\delta(\beta)$, then $D_j\sse Y(\beta)$ is a smooth irreducible divisor and $D_j\cong Y(\beta\setminus \s_{i_j})$.\\

  \item More generally, $\displaystyle\bigcap_{\s_{i_j}\in S}D_j\ne\emptyset\iff \delta(\beta\setminus S)=\delta(\beta)$ and, in such a case,
  \begin{equation}\label{eq:multiisomorphism}
 \bigcap_{\s_{i_j}\in S}D_j\cong Y(\beta\setminus S).
\end{equation}
In addition, such intersection in \eqref{eq:multiisomorphism} has codimension $|S|$ in $Y(\beta)$.\\
 \item The union of divisors
 $$\bigcup_{j}D_j\sse Y(\beta),$$
 where the union is over indices $j\in[1,\ell]$ satisfying $\delta(\beta\setminus \s_{i_j})=\delta(\beta)$, is a simple normal crossings divisor in $Y(\beta)$.\\

  \item The dual boundary complex of $\{D_j\}$ is isomorphic to the subword complex of $(\beta,\delta(\beta))$. In addition, it is homotopy equivalent to a sphere.
  \sbs
\end{enumerate}
\end{lemma}
Parts \Cref{lem:boundary_div}.(1) and (2) follow from direct computation, \Cref{lem:boundary_div}.(3) follows from iterating \Cref{lem:boundary_div}.(2). \Cref{lem:boundary_div}.(4) follows from \Cref{lem:boundary_div}.(3), as the divisors intersect transversely, since all the intersections are smooth, cf.~\cite[Theorem 4.11]{CGGS2} for more details. The subword complexes mentioned in \Cref{lem:boundary_div}.(5) were studied in \cite{KnutsonMiller_04} and sphericity follows from \cite[Corollary 3.8]{KnutsonMiller_04}, see also \cite[Section 4]{CGGS2}.\footnote{The assertion that the dual boundary complex is homotopy equivalent to a sphere is an instance of the geometric $P=W$ conjecture.}\\

\noindent Let $D_\beta:=\sum D_j$ denote the linear sum of the boundary divisors of the brick compactification $Y(\beta)$ of $X(\beta)$. Beyond \Cref{lem:boundary_div}, the following property of $D_\beta$ is significant from the viewpoint of mirror symmetry and the study of cluster algebras via algebraic geometry:

\begin{proposition}\label{prop:logCY}
Let $\beta$ be a positive braid word, $X(\beta)$ its braid variety, $Y(\beta)$ its brick compactification and $D_\beta$ the linear sum of its compactifying divisors. Then:
\begin{enumerate}
    \item $(Y(\beta),D_\beta)$ is a log Calabi-Yau pair.
    \item $(Y(\beta),D_\beta)$ has maximal boundary and $D_\beta$ is the support of an ample divisor.\sbs
\end{enumerate}

\end{proposition}

\noindent Equivalently, after base change\footnote{If we are to remain over $\Z$, then the resulting object is a relative regular form.} to $\C$, \Cref{prop:logCY}.(1) states that there exists a nowhere vanishing holomorphic top form on $X(\beta)$ which acquires only simple poles along $D_\beta$ when it compactifies to $Y(\beta)$, see \cite{HackingKeel_ICM18} for details on log Calabi-Yau pairs as they relate to cluster structures. \Cref{prop:logCY} follows from studying the divisors and canonical bundle of $\BS(\beta)$, see e.g.~\cite{Anderson,LauTho} or \cite[Lemma 4.9]{CasalsWeng22}. For braid varieties $X(\beta)$ such that the closed 2-form $\omega_\beta$ in \Cref{sssec:symplectic_form_coordinates} is non-degenerate, we can directly take the top power of $\omega_\beta$.\\

\noindent The maximal-boundary statement in \Cref{prop:logCY}.(2), i.e.~that there is a 0-dimensional stratum in the intersection complex of $D_\beta$, follows from \Cref{lem:boundary_div}. Note that \Cref{prop:logCY}.(2) only states that $D_\beta$ is the support of an ample divisor, but $D_\beta$ itself might not be ample. Indeed, by \Cref{prop:logCY}.(1) $D_\beta$ is itself ample iff $Y(\beta)$ is a Fano variety, which is not always the case.%, e.g.~for $\beta=c\cdot S_c(w_0)$ where $c$ represents the Coxeter element and $S_c(w_0)$ is the $c$-sorting reduced word for $w_0$, $Y(\beta)$ is a toric variety obtained from the associahedron, as in \cite[Corollary 2.10]{Escobar}, and it is not Fano for $\G=\SL_n$ if $n\geq3$.

%%%%%%%%%%%%%%%%%%%%%%%%%%%%%%%%%%%%%%%%%%%%%
%%%%%%%%%%%%%%%%%%%%%%%%%%%%%%%%%%%%%%%%%%%%%
\subsubsection{Two applications of brick compactifications}\label{sssec:applications_compactifications}

First, the cohomology of the variety $X(\beta)$ can be computed by using the fact that the cohomology of $Y(\beta)$ is rather well-understood and applying the inclusion-exclusion principle. Specifically, consider the $\T$-action on $X(\beta)$ introduced in \cref{sssec:braid_variety_further_symmetries}, cf.~\eqref{eq:torus_action}, and let $H_{c,T}^*(X(\beta),\Z)$ denote the compactly supported $\T$-equivariant cohomology of $X(\beta)$. Since $X(\beta)$ is a smooth affine variety, $H_{c,T}^*(X(\beta),\Z)$ is equipped with a weight filtration (after base change from $\Z$ to $\C$) and we denote by $\operatorname{gr}^W H_{c,T}^*(X(\beta),\Z)$ the associated graded quotients of this weight filtration.\\

A neat result about braid varieties is that such graded pieces of the equivariant cohomology of $X(\beta)$ recover a link invariant if $\G=\SL_n$. Indeed, let $\operatorname{HHH}$ be the triply-graded Khovanov-Rozansky link homology, as in \cite{KhS,KR}, and we denote by $\operatorname{HHH}(\beta w_0)$ the triply-graded Khovanov-Rozansky link homology of the smooth link obtained by the closure of $\beta w_0$, where $w_0$ denotes a braid lift of a reduced expression for $w_0$. Recall that one of the three degrees is the $a$-degree, which ranges from $0$ to the number of strands of $\beta$. Let us denote by $\operatorname{HHH}^{a=n}$ the highest $a$-degree piece and we use $\ateb$ to denote $\beta$ read right-to-left. The following result is then proven in \cite[Theorem 4.22]{CGGS2}:

\begin{theorem}\label{thm:cohomology_and_HHH}
Let $\beta$ be a positive $n$-stranded braid word with $\delta(\beta)=w_0$. Then
\begin{equation}\label{eq:cohomology_and_HHH}
\boxed{\operatorname{gr}^W H_{c,T}^*(X(\beta),\Z)\cong \operatorname{HHH}^{a=n}(\ateb w_0)}
\end{equation}
\sbs
\end{theorem}

\noindent \Cref{thm:cohomology_and_HHH} is a cohomological enhancement of the fact that the ruling polynomial of a Legendrian link appears as the appropriate $a$-coefficient  of the HOMFLY polynomial of the underlying smooth link, see e.g.~\cite{Rutherford06_HOMFLY,Rutherford11_HOMFLY}. It can also be understood as a cohomological enhancement in characteristic zero of certain identities between finite points counts, cf.~\cite{trinh2021heckecategoryunipotentlocus,trinh2022unipotentelementstwistinglink}.\\

Second, it is possible to use $Y(\beta)$ as a compactification of $X(\beta)$ in order to establish a form of mirror symmetry for the appropriate cluster $\SA$ and $\SX$-schemes associated to $\beta$. Specifically, in the context of Fock-Goncharov duality, as in \Cref{thm:fock_goncharov_duality}, the cluster ensemble
$$(\SA_\G(\beta),\SX_{\G^\vee}(\beta))$$
can be refined by compactifying $\SX_{\G^\vee}(\beta)$ and, dually, adding a potential $W_\beta:\SA_\G(\beta)\lr\mathbb{A}^1_\kk$ to $\SA_\G(\beta)$. This turns $(\SA_\G(\beta),W_\beta)$ into a Landau-Ginzburg model from which it is possible to compute certain invariants. These potentials $W_\beta$ have clean expressions in terms of Reeb chord asymptotics, as we explain in \Cref{thm:superpotential_Reeb_coordinates} below.

\begin{remark} We use the notation $W_\beta$ because the potential depends on the positive braid word $\beta$, as does the compactification $Y(\beta)$ of $\SX_{\G^\vee}(\beta)$. Thus, as above, one can change the positive braid word $\beta$ by braid moves and obtain different Landau-Ginzburg models for the same variety $\SA_\G(\beta)$.\sbs
\end{remark}

%%%%%%%%%%%%%%%%%%%%%%%%%%%%%%%%%%%%%%%%%%%%%
%%%%%%%%%%%%%%%%%%%%%%%%%%%%%%%%%%%%%%%%%%%%%
%%%%%%%%%%%%%%%%%%%%%%%%%%%%%%%%%%%%%%%%%%%%%
%%%%%%%%%%%%%%%%%%%%%%%%%%%%%%%%%%%%%%%%%%%%%
\section{The Reeb coordinates of a braid variety}\label{sec:Reeb_variables} Let $\beta\in\brn$ be a positive braid word with length $\ell:=\ell(\beta)$. By \Cref{def:braid_variety}, the braid variety $X(\beta)$ is defined as a subvariety of the projective variety $(\G/\borel)^{\ell(\beta)+1}$. The description of $X(\beta)$ in \cref{def:braid_variety} and discussions above are coordinate-free. There is nevertheless some merit in introducing explicit coordinates to study $X(\beta)$. In this section we introduce a particular set of coordinates in $\mathbb{A}_\kk^{\ell(\beta)}$ and describe a particular embedding $X(\beta)\sse\mathbb{A}_\kk^{\ell(\beta)}$ as an affine subvariety using such coordinates. We refer to such coordinates as Reeb coordinates or $z$-variables, as they are geometrically motivated by the study of Reeb chords in contact topology, see \cite[Section 5.1]{CasalsNg22}.\\

Because of their geometric origin and several facts stemming from contact and symplectic topology, these Reeb coordinates have remarkable properties to study braid varieties. A few advantages that Reeb coordinates have are:

\begin{enumerate}

    \item All the cluster variables $\{A^\w_v\}$ associated to any weave $\w:\beta\lr\delta(\beta)$ will be polynomials in the Reeb $z$-coordinates $(z_1,\ldots,z_\ell)$, $\ell=\ell(\beta)$, cf.~\Cref{thm:polynomiality}. In addition, the exponents in the leading monomial of such a polynomial $A^\w_v(z_1,\ldots,z_\ell)$, in the lexicographic order given by $z_1\prec z_2\prec \ldots\prec z_\ell$, provide the crystal vertex in the $\beta$-Lusztig parametrization of a global basis of $\kk[X(\beta)]$, see e.g.~\Cref{thm:leading_monomial_cluster_A_variables}. Plus, there is a combinatorial expression for such polynomials $A^\w_v(z_1,\ldots,z_\ell)$ in terms of a hook content formula, cf.~\Cref{thm:hook_content_cluster_A_variables}.\\
    
    \item The proof that the inclusion \cref{eq:inclusion_cluster_algebra} is an equality provided in \cite[Section 5]{CGGLSS25} uses Reeb coordinates. In addition, proofs obtained by computations in Reeb coordinates work identically for all $\G$. For instance, some invariance properties directly transfer from $\SL_3$ to $\G$ via pinnings, and having explicit coordinate-dependent equations for $X(\beta)$ as an algebraic subvariety of $\mathbb{A}_\kk^{\ell(\beta)}$ can help study properties of $X(\beta)$ such as smoothness, dimension, group actions, and Gr\"obner bases for $\kk[X(\beta)]$, among others. In particular, the torus action on braid varieties has a neat expression in terms of Reeb coordinates, cf.~\Cref{lem:torusaction_coordinates}, as does the cluster 2-form, cf.~\Cref{prop:2forms}. Similarly, the superpotential of a natural log Calabi-Yau compactification of $X(\beta)$ has a remarkably simple expression in the Reeb coordinates, cf.~\Cref{thm:superpotential_Reeb_coordinates}.\\
    
    \item For everyday use, it is also often useful in applications to have explicit formulas for the cluster variables $\{A^\w_v\}$ in terms of a specific set of generators for $\kk[X(\beta)]$. In particular, there is computer software available to compute cluster variables from Demazure weaves which uses this explicit Reeb coordinate system.\footnote{A Javascript applet is available at \href{https://daping-weng.github.io/math/demazure_weave.html}{Daping Weng's website here}.}
\end{enumerate}

%%%%%%%%%%%%%%%%%%%%%%%%%%%%%%%%%%%%%%%%%%%%%
%%%%%%%%%%%%%%%%%%%%%%%%%%%%%%%%%%%%%%%%%%%%%
%%%%%%%%%%%%%%%%%%%%%%%%%%%%%%%%%%%%%%%%%%%%%
\subsection{Reeb coordinates for braid varieties}\label{ssec:Reeb_coordinates_braid_variety} We will now consider ambient affine coordinates in $\mathbb{A}_\kk^{\ell}=\Spec k[z_1,\ldots,z_\ell]$ and an embedding $X(\beta)\sse \mathbb{A}_\kk^{\ell}$ with explicit equations cutting out $X(\beta)$ in terms of these $z$-variables. These ambient $z$-variables will also be referred to as Reeb coordinates, and they are closely related to the $\widetilde{z}$-variables in \Cref{lem:framing_variables}. As before, we fix a pinning of the group $\G$ as in \Cref{ssec:rational_functions} and consider
\[
s_{i} := \varphi_{i}\left(\begin{matrix}0 & -1 \\ 1 & 0 \end{matrix} \right)\in \G,
\]
so $s_{i} \in N_{\G}(T)$ lifts the simple reflection corresponding to $i \in \dynkin$. Given a permutation $v\in W(\G)$, we define its lift to $\G$ by choosing an arbitrary reduced expression and multiplying these elements $s_i$ accordingly. In order to discuss tuples of flags, we introduce a variable $z \in k$ and define
\[
B_{i}(z):= x_{i}(z)s_{i} = \varphi_{i}\left(\begin{matrix} z & -1 \\ 1 & 0\end{matrix}\right) \in \G.
\]

\noindent {\it What is this variable $z$ parametrizing in terms of flags?} Suppose we fix a flag $x\borel \in \G/\borel$ and we want to parametrize all flags in $\G/\borel$ which are in relative $s_i$-position with $x\borel$. The claim is that this set is isomorphic to $\mathbb{A}_k^1$ and $z$ provides a coordinate for it. Explicitly, we can parametrize this set via
\begin{equation}\label{eq:Biz_parametrizing_flags_rel_pos}
\{y\borel \in \G/\borel \mid x\borel \buildrel s_{i} \over \longrightarrow y\borel\}\cong \{xB_{i}(z)\borel | z \in \kk\},
\end{equation}
and so two flags $xB_i(z)\borel=xB_i(z')\borel$ are equal if and only if $z=z'$. This is essentially \Cref{lem:framing_variables} for unframed flags.

\begin{example} (1). For $\G=\SL_n$, we choose the following pinning and $B_i(z)$ matrices:
\begin{equation}\label{eq:pinning1}
x_{i}(z) = \begin{pmatrix}
1 & \cdots  & & & \cdots & 0\\
\vdots & \ddots & & & & \vdots\\
0 & \cdots & 1 & z & \cdots & 0\\
0 & \cdots & 0 & 1 & \cdots & 0\\
\vdots &  & & &\ddots & \vdots\\
0 & \cdots & & & \cdots & 1\\
\end{pmatrix},\qquad
y_{i}(z) = \begin{pmatrix}
1 & \cdots  & & & \cdots & 0\\
\vdots & \ddots & & & & \vdots\\
0 & \cdots & 1 & 0 & \cdots & 0\\
0 & \cdots & z & 1 & \cdots & 0\\
\vdots &  & & &\ddots & \vdots\\
0 & \cdots & & & \cdots & 1\\
\end{pmatrix},
\end{equation}
\begin{equation}\label{eq:pinning1}
\chi_{i}(t) = \begin{pmatrix}
1 & \cdots  & & & \cdots & 0\\
\vdots & \ddots & & & & \vdots\\
0 & \cdots & t & 0 & \cdots & 0\\
0 & \cdots & 0 & t^{-1} & \cdots & 0\\
\vdots &  & & &\ddots & \vdots\\
0 & \cdots & & & \cdots & 1\\
\end{pmatrix},\qquad
B_{i}(z) = \begin{pmatrix}
1 & \cdots  & & & \cdots & 0\\
\vdots & \ddots & & & & \vdots\\
0 & \cdots & z & -1 & \cdots & 0\\
0 & \cdots & 1 & 0 & \cdots & 0\\
\vdots &  & & &\ddots & \vdots\\
0 & \cdots & & & \cdots & 1\\
\end{pmatrix},
\end{equation}
where the $(2\times 2)$-submatrices sit at the $i$th and $(i + 1)$st rows and columns. The variable $z$ is thus the $(i,i+1)$-entry of $x_i(z)$, the $(i+1,i)$-entry of $y_i(z)$, and the $(i,i)$-entry of $B_i(z)$. The variable $t$ is the $(i,i)$-entry of $\chi_i(t)$. Note that $B_i(z)^{-1}=s_iB_i(-z)s_i$.\\

\noindent (2) Consider $\G=\SL_2$, and the standard flag $\borel$, which is determined by the line $\borel^1=\langle e_1\rangle\in k^2$. Then the Borel subgroup
$$B_1(z)\borel=\left\{g\in\SL_2: g=\begin{pmatrix}
z & -1\\
1 & 0
\end{pmatrix}B\begin{pmatrix}
0 & 1\\
-1 & z
\end{pmatrix},\mbox{ for some }B\mbox{ upper-triangular}\right\}$$
corresponds (i.e.~it is the stabilizer of) the unique full flag in $k^2$ whose 1-dimensional subspace is $\langle ze_1+e_2\rangle$. Thus $(B_1(z)\borel)^1=\langle ze_1+e_2\rangle$. In particular, since $e_1$ is never parallel to $ze_1+e_2$ for $z\in k$, these two flags are different, i.e.~$\borel\stackrel{s_1}{\to}B_1(z)\borel$. Conversely, any line in $k^2$ which is different from $\langle e_1\rangle$ has a unique generator of the form $ze_1+e_2\in k^2$, and thus $z\in k$ bijectively parametrizes such flags.\sbs
\end{example}

%%%%%%%%%%%%%%%%%%%%%%%%%%%%%%%%%%%%%%%%%%%%%
%%%%%%%%%%%%%%%%%%%%%%%%%%%%%%%%%%%%%%%%%%%%%
%%%%%%%%%%%%%%%%%%%%%%%%%%%%%%%%%%%%%%%%%%%%%
\subsection{Equations for braid varieties via the Reeb coordinates}\label{ssec:equation_braid_varieties_Reeb_coordinates} Given a positive braid word $\beta = \sigma_{i_{1}}\cdots \sigma_{i_{\ell}}$, we define $B_{\beta}(z)\in\G$ as
\begin{equation}\label{eq:braid_matrix}
B_{\beta}(z) := B_{i_{1}}(z_1)\cdots B_{i_{\ell}}(z_\ell) \in \G,\quad z=(z_1,\ldots,z_\ell)\in \mathbb{A}_\kk^{\ell}.
\end{equation}

\noindent As above, these variables $(z_1,\ldots,z_\ell)\in \mathbb{A}_\kk^{\ell}$ associated to the crossings of $\beta$ will be referred to as the $z$-variables or Reeb coordinates for $\beta$. It follows from \cref{eq:Biz_parametrizing_flags_rel_pos} that, in these $z$-coordinates, the braid variety $X(\beta)$ can be described as follows:

\begin{lemma}[\color{blue}Equations for braid varieties in Reeb coordinates\color{black}]\label{lem:equations_braidvar}
Let $\beta = \sigma_{i_{1}}\cdots \sigma_{i_{\ell}}$ be a positive braid word, with Demazure product $\delta(\b)$, and $\delta(\b)^{-1}$ a lift of the inverse of $\delta(\b)$ to $\G$ using $s_i$. Then there is an isomorphism
\begin{equation}\label{eq:braid_var_eq}
X(\beta) \cong \{(z_1, \dots, z_\ell) \in \mathbb{A}_\kk^{\ell} \mid \delta(\b)^{-1}B_{i_{1}}(z_1)\cdots B_{i_{\ell}}(z_\ell) \in \borel\},
\end{equation}
$X(\beta)$ is a smooth irreducible affine variety, and $\dim_\kk(X(\beta))=\ell(\b)-\ell(\delta(\b))$.
\sbs
\end{lemma}

\noindent See \cite[Corollary 3.7]{CGGLSS25} for a detailed proof. The equation $\delta^{-1}B_{i_{1}}(z_1)\cdots B_{i_{r}}(z_r) \in \borel$ is given by the vanishing of some generalized minors, and thus $X(\beta)$ is indeed an affine variety defined over $\Z$. In these $z$-coordinates, the braid invariance properties from \cref{lem:braid_invariance} follow from the identities
$$B_{i}(z)B_{j}(w) = B_{j}(w)B_{i}(z),\quad i,j\in\dynkin\mbox{ not adjacent}$$
$$B_{i}(z_1)B_{j}(z_2)B_{i}(z_3) = B_j(z_3)B_i(z_1z_3 - z_2)B_j(z_1),\quad i,j\in\dynkin\mbox{ adjacent}.$$
The latter is proven by directly transporting the $\SL_3$-identity
$$B_{1}(z_1)B_{2}(z_2)B_{1}(z_3) = B_2(z_3)B_1(z_1z_3 - z_2)B_2(z_1)$$
to $\G$ via pinnings. This $\SL_3$-identity is checked by hand:
$$\begin{pmatrix}
z_1 & -1 & 0\\
1 & 0 & 0\\
0 & 0 & 1\\
\end{pmatrix}
\begin{pmatrix}
1 & 0 & 0\\
0 & z_2 & -1\\
0 & 1 & 0\\
\end{pmatrix}
\begin{pmatrix}
z_3 & -1 & 0\\
1 & 0 & 0\\
0 & 0 & 1\\
\end{pmatrix}=\begin{pmatrix}
z_1z_3-z_2 & -z_1 & 1\\
z_3 & -1 & 0\\
1 & 0 & 0\\
\end{pmatrix},$$
$$\begin{pmatrix}
1 & 0 & 0\\
0 & z_3 & -1\\
0 & 1 & 0\\
\end{pmatrix}
\begin{pmatrix}
z_1z_3-z_2 & -1 & 0\\
1 & 0 & 0\\
0 & 0 & 1\\
\end{pmatrix}
\begin{pmatrix}
1 & 0 & 0\\
0 & z_1 & -1\\
0 & 1 & 0\\
\end{pmatrix}=\begin{pmatrix}
z_1z_3-z_2 & -z_1 & 1\\
z_3 & -1 & 0\\
1 & 0 & 0\\
\end{pmatrix}.
$$
\normalsize
\noindent The cyclic symmetry mentioned in \Cref{sssec:braid_variety_further_symmetries} can also be verified through matrix computation, cf.~\cite[Lemma 3.10]{CGGLSS25}.

\begin{example}\label{ex:braid_matrices} (1) Let $\beta=\s_{i_1}\ldots \s_{i_\ell}$ be a positive braid for $\G=\SL_n$ and consider $(\borel_1,\ldots,\borel_{\ell+1})\in X(\beta)$ a point in the braid variety. Then the $i$-dimensional subspace of the $j$th flag $\borel_j^i$ is spanned by the first $i$ columns of the $(n\times n)$-matrix
$$B_{i_{1}}(z_1)\cdots B_{i_{j-1}}(z_{j-1}) \in \G, \quad z=(z_1,\ldots,z_{j-1})\in X(\beta),\quad j\in[1,\ell+1],$$
\noindent where the empty product for $j=1$ is understood as the identity matrix, so $\borel_1=\borel$ is the standard flag, and for $j=\ell+1$ the condition $\borel_{\ell+1}=\delta(\beta)\borel$ is the condition \cref{eq:braid_var_eq} for $(z_1,\ldots,z_{j-1})\in X(\beta)$.\\

\noindent (2) For instance, consider the braid $\beta=\s_1^3\in\Br(\SL_2)$. Since
$$s_1B_1(z_1)B_1(z_2)B_1(z_3)=\begin{pmatrix}
1-z_2z_3 & z_2\\
z_1z_2z_3-z_3-z_1 & 1-z_1z_2
\end{pmatrix},$$
the braid variety $X(\beta)$ in the $z$-coordinates is given by
\begin{align*}
X(\beta) & \cong \{(z_1,z_2,z_3)\in k^3:s_1B_1(z_1)B_1(z_2)B_1(z_3)\in\borel\}\\
         & \cong \left\{(z_1,z_2,z_3)\in k^3:z_1z_2z_3-z_3-z_1=0,\quad (1-z_1z_2)(1-z_2z_3)=1\right\}\\
         & \cong \left\{(z_1,z_2)\in k^2:(1-z_1z_2)\neq0\right\},
\end{align*}
where in the last isomorphism we used that $(1-z_1z_2)(1-z_2z_3)=1$ implies $(1-z_1z_2)\neq0$ and $z_1z_2z_3-z_3-z_1=0$ allows us to write $z_3=z_3(z_1,z_2)$ restricted to $X(\beta)$ as a function of $(z_1,z_2)$ via
$$z_3=\frac{-z_1}{1-z_1z_2}.$$
Under this isomorphism, the $\T$-action mentioned in \Cref{sssec:braid_variety_further_symmetries}, $\bG_m\cong\T\subset\SL_2$, is given by
$$\bG_m\times X(\beta)\lr X(\beta),\quad (\tau;z_1,z_2)\lr (\tau^2 z_1,\tau^{-2}z_2),$$
which indeed leaves the quantity $z_1z_2$ invariant, and thus the condition $(1-z_1z_2)\neq0$ is preserved. The fixed point $(\borel,\borel_-,\borel,\borel_-)\in X(\s_1^3)$ is the point $(z_1,z_2)=(0,0)$ under these isomorphisms, which is indeed fixed by the $\T$-action above.\\

\noindent (3) More generally, braid varieties of the form $X(\gamma w_0)$, where $\beta=\gamma w_0$ contains a consecutive subword $w_0$, have simpler and fewer equations. These braid varieties can be described as an open set in $\mathbb{A}_\kk^{\ell(\gamma)}$ instead of a subvariety of $k^{\ell(\gamma w_0)}$. Indeed, this follows by noting that $\delta(\gamma w_0)=w_0$ and that the matrix $B_{w_0}^{-1}(w)\in\SL_n$ has the form
$$B_{w_0}^{-1}(w)=\left(\begin{matrix} 0 & 0 & \ldots & 0 & \pm1\\ 
0 & 0 &  \ldots & \pm1 & \pm w_{\binom{n}{2}} \\ 0 &  0 & \ldots & \pm w_{\binom{n}{2}-(n-3)} & \pm w_{\binom{n}{2}-1} \\ 
\vdots & \vdots & \ddots & \vdots & \vdots \\ \pm1 & \pm w_{1} & \ldots & \pm w_{\binom{n}{2}-2(n-3)} & \pm w_{\binom{n}{2}-(n-2)}
\end{matrix}\right)=w_0N\mbox{ with } N\in\borel,$$
where $w=(w_1,w_2,\ldots,w_{\ell(w_0)})\in \mathbb{A}_\kk^{\ell(w_0)}$ and $\pm$ indicates a possible sign. The exact signs are irrelevant, what matters is that $B_{w_0}^{-1}(w)$ is of the form $w_0N$ for a $N\in\borel$ and the variables $w\in\C^{\ell(w_0)}$ appear each alone in exactly one entry. Therefore \cref{eq:braid_var_eq} reads
\begin{align*}
X(\gamma w_0) & \cong \{(z,w) \in \mathbb{A}_\kk^{\ell(\gamma)+\ell(w_0)} \mid w_0^{-1}B_\gamma(z)B_{w_0}(w) \in \borel\}\\
         & \cong \{z \in \mathbb{A}_\kk^{\ell(\gamma)} \mid w_0B_\gamma(z) \in \borel w_0\borel\}.
\end{align*}
\noindent The condition $w_0 M\in\borel w_0\borel$ is equivalent to $M\in (w_0\borel w_0)\borel=\borel_-\borel$. This is equivalent to $M$ admitting an $LU$-decomposition, which holds if and only if all the $n$ leading principal minors\footnote{The leading principal minor of $M$ of order $k$ is the minor of order $k$ obtained by deleting the last
$(n-k)$ rows and columns of $M\in\SL_n$. There are $n$ of those, for $k\in[1,n]$.} $\Delta_k(M)$ of $M$ do not vanish for any $k\in[1,n]$. In consequence,
$$X(\gamma w_0)\cong\{z\in \mathbb{A}_\kk^{\ell(\gamma)}: \Delta_1(B_\gamma(z))\cdot\ldots\cdot\Delta_{n-1}(B_\gamma(z))\cdot\Delta_{n}(B_\gamma(z))\neq0\}\sse \mathbb{A}_\kk^{\ell(\gamma)}.$$
\noindent In $\G=\SL_n$, the determinant is $\Delta_n(B_\gamma(z))=1$, and thus this condition for $X(\gamma w_0)$ is equivalent to the $(n-1)$ polynomials $\Delta_1(B_\gamma(z)),\ldots,\Delta_{n-1}(B_\gamma(z))$ being non-zero.
\sbs
\end{example}

%%%%%%%%%%%%%%%%%%%%%%%%%%%%%%%%%%%%%%%%%%%%%
%%%%%%%%%%%%%%%%%%%%%%%%%%%%%%%%%%%%%%%%%%%%%
%%%%%%%%%%%%%%%%%%%%%%%%%%%%%%%%%%%%%%%%%%%%%
\subsection{Cluster $A$-variables and Reeb coordinates}\label{ssec:polynomiality_A_variables_Reeb_coordinates}
There are a number of useful applications of the Reeb $z$-variables. First, by using the embedding $\iota:X(\beta)\hookrightarrow\mathbb{A}_\kk^{\ell(\beta)}$ from \Cref{lem:equations_braidvar}, the following is a consequence of \cite[Theorem 7.6]{CGGLSS25}:

\begin{theorem}[Weave cluster variables are polynomials in the Reeb coordinates]\label{thm:polynomiality} Let $\beta\in\brn$ be a positive braid word, $\iota:X(\beta)\hookrightarrow\mathbb{A}_\kk^{\ell(\beta)}$ its braid variety with the ambient Reeb $z$-variables on $\mathbb{A}_\kk^{\ell(\beta)}$. Then, for any Demazure weave $\w:\beta\to\delta(\beta)$ and any trivalent vertex $v$, there exists a polynomial ${\bf A}^\w_v\in \kk[z_1,\ldots,z_\ell]$ such that
$$A^\w_v=\iota^*{\bf A}^\w_v.$$
Informally, all weave cluster variables are polynomials in the $z$-variables.\sbs
\end{theorem}

\noindent In particular, \Cref{thm:polynomiality} makes it apparent that the cluster variables $\{A^\w_v\}\in\kk[X(\beta)]$ are indeed regular functions on $X(\beta)$, as they are the restriction of a regular function on $\mathbb{A}_\kk^{\ell(\beta)}$ to $X(\beta)$. Second, the $z$-variables can be used to show that the inclusion \eqref{eq:inclusion_cluster_algebra} is an equality. Recall that we denote by $\SA(Q)$ the commutative cluster algebra associated to a quiver $Q$, cf.~\cite[Def.~3.1.9]{FWZ16}, which we take to be defined over $\kk$. The precise statement reads as follows:

\begin{fthm}[Cluster structures on braid varieties]\label{thm:cluster_structures_braid_varieties} Let $\beta\in\brn$ be a positive braid word and $X(\beta)$ its braid variety. Then there exists an isomorphism of commutative $\kk$-algebras
\begin{equation}\label{eq:cluster_algebra_braid_variety}
\SA(Q_\w)\cong \kk[X(\beta)]
\end{equation}
where $\w:\beta\lr\delta(\beta)$ is an arbitrary Demazure weave. In addition, under this isomorphism the following holds:
\begin{enumerate}
    \item Each Demazure weave $\w:\beta\lr\delta(\beta)$ provides a cluster seed $(T_\w,Q_\w,\{A^\w_v\})$ for $\SA(Q_\w)$.\\
    \item If $\w_1\simeq\w_2$ are weave equivalent, then their respective cluster seeds coincide.\\
    \item If $\w_1\stackrel{\mu}{\sim}\w_2$ are related by one weave mutation, then their respective cluster seeds are related by one cluster mutation.\\
    \item If $\iota:X(\beta)\hookrightarrow \mathbb{A}_\kk^{\ell(\beta)}$ denotes the embedding given by the Reeb $z$-coordinates, then all weave cluster variables $\{A^\w_v\}$, for any Demazure weave $\w:\beta\lr\delta(\beta)$, are restrictions of polynomials in the Reeb $z$-variables.
\end{enumerate}
\end{fthm}
\vspace{0.2cm}
\noindent Most of Theorem \ref{thm:cluster_structures_braid_varieties} is Theorem \ref{thm:cluster_variables_weave} and \Cref{thm:polynomiality}, we only restate it for clarity. The only crucial fact remaining is the inclusion $\kk[X(\beta)]\sse \SA(Q_\w)$. To prove that inclusion, we note that the Reeb $z$-variables generate $\kk[X(\beta)]$, as they generate $\kk[z_1,\ldots,z_{\ell}]$ and we have the embedding $\iota:X(\beta)\hookrightarrow\mathbb{A}_\kk^{\ell(\beta)}$ from \Cref{lem:equations_braidvar}. The inclusion $\kk[X(\beta)]\sse \SA(Q_\w)$ then follows from the fact that, for any non-zero $z_i$-variable, there exists a Demazure weave $\w_{z_i}:\beta\lr\delta(\beta)$ such that its associated cluster seed has $z_i$ as one of its cluster variables $\{A_v^\w\}$, cf.~\cite[Section 4.3]{CGGLSS25} and \Cref{ex:z_cluster_variables} below. It is {\it not} the case that all $z$-variables are cluster variables in the same cluster seed, so the choice of Demazure weave $\w_{z_i}$ must typically depend on the particular $z_i$-variable.

\begin{example}\label{ex:z_cluster_variables} Suppose that we have a positive braid word of the form $\beta=\beta_1\s_i^2\beta_2\in\brn$ and set $\ell_1=\ell(\beta_1)$. Then, any Demazure weave $\w:\beta\lr\delta(\beta)$ that starts with the move $\beta_1\s_i^2\beta_2\to\beta_1\s_i\beta_2$ will be such that its top trivalent vertex $v\in\w$ satisfies $A_v^\w=z_{\ell_1+2}$. This shows that such a Reeb $z$-variable $z_{\ell_1+2}$ is a cluster variable. In order to show that the remaining $z$-variables are cluster variables, one can exploit the cyclic symmetry and bring any given crossing of $\beta$ to the end of the braid word and use an inductive weave, as defined in \cite[Section 4.3]{CGGLSS25}.\sbs
\end{example}

%%%%%%%%%%%%%%%%%%%%%%%%%%%%%%%%%%%%%%%%%%%%%
%%%%%%%%%%%%%%%%%%%%%%%%%%%%%%%%%%%%%%%%%%%%%
%%%%%%%%%%%%%%%%%%%%%%%%%%%%%%%%%%%%%%%%%%%%%
\subsection{Reeb coordinates after adding a vertex}\label{ssec:Reeb_coordinates_adding_vertex} Let us consider the situation as in \Cref{ssec:cluster_A_variables_propagation_vertices}, where we are given a Demazure weave $\w:\beta\to\delta(\beta)$ and we add a weave vertex $p$ on top of it, creating a new weave $\w^+$. Recall that $\ell=\ell(\beta)$ and consider the following labeling of the Reeb coordinates for $\beta_+$:

\begin{enumerate}
    \item If $p$ is a hexavalent or tetravalent vertex, then $\ell(\beta_+)=\ell$ and the Reeb coordinates of $\beta_+$ are labeled $(\zeta_1,\ldots,\zeta_\ell)$, corresponding to the crossings of $\beta_+$ read left-to-right.\\

    \item If $p$ is a trivalent vertex, then $\ell(\beta_+)=\ell+1$. If the south edge of the trivalent vertex $p$ corresponds to the $k$th crossing of $\beta$, counting from the left, then the Reeb coordinates of $\beta_+$ are labeled $(\zeta_1,\ldots,\zeta_k,\textcolor{blue}{\zeta},\zeta_{k+1},\ldots,\zeta_\ell)$, which again correspond to the crossings of $\beta_+$ read left-to-right. Thus, $\textcolor{blue}{\zeta}$ corresponds to the northeast edge of $p$. This notation is chosen, in part, so that the top weave edge of $\w$ carrying the variable $z_i$, for any $i\neq k$, is now extended to $\w^+$ as a straight line up carrying the variable $\zeta_i$.
\end{enumerate}

\begin{figure}[ht!]
\centering
    \includegraphics[scale=2]{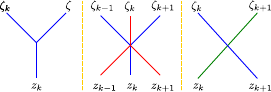}
\caption{The Reeb coordinates above and below each of the weave vertices. \Cref{prop:change_Reebcoordinates_add_vertex_on_top} explains how to relate the bottom Reeb $z$-coordinates and the top Reeb $\zeta$-coordinates.}
    \label{fig:weaves_types_vertices_ReebCoordinates}
\end{figure}

\begin{proposition}\label{prop:change_Reebcoordinates_add_vertex_on_top} In the notation above, the following holds:
\begin{enumerate}[label=$(\roman*)$]
    \item If $p$ is a hexavalent vertex inserted on top of the crossings $\s_{i_{k-1}}\s_{i_{k}}\s_{i_{k+1}}$ of $\beta$, then
    $$z_{k-1}=\zeta_{k+1},\quad z_{k}=\zeta_{k-1}\zeta_{k+1}-\zeta_{k},\quad z_{k+1}=\zeta_{k-1},\quad z_{j}=\zeta_{j}\quad \mbox{if }j\not\in\{k-1,k,k+1\}.$$
    Equivalently, the inverse transformation is
    $$\zeta_{k-1}=z_{k+1},\quad \zeta_{k}=z_{k-1}z_{k+1}-z_{k},\quad \zeta_{k+1}=z_{k-1},\quad \zeta_{j}=z_{j}\quad \mbox{if }j\not\in\{k-1,k,k+1\}.$$

    \item If $p$ is a tetravalent vertex inserted on top of the crossings $\s_{i_{k}}\s_{i_{k+1}}$ of $\beta$, then
    $$\zeta_{k}=z_{k+1},\quad \zeta_{k+1}=z_{k},\quad \zeta_{j}=z_{j}\quad \mbox{if }j\not\in\{k,k+1\}.$$

    \item If $p$ is a trivalent vertex inserted on top of the crossings $\s_{i_{k}}$ of $\beta$, then

$$
\boxed{
\begin{aligned}
\iota^*\left(\zeta_j|_{\{\zeta\neq0\}}\right) & =  z_j \qquad\qquad\quad\mbox{ if } j < k, \\
\iota^*\left(\zeta_k|_{\{\zeta\neq0\}}\right) &=  z_k + t^{-1}, \\
\iota^*(\zeta)&=t,\\
\iota^*\left(\zeta_j|_{\{\zeta\neq0\}}\right) & =  \eta_j(z_1,\ldots,z_{j};t) \mbox{ if } k < j,
\end{aligned}
}
$$
\noindent where $\iota:X(\beta)\times\bG_m\lr X(\beta_+)$ is the embedding induced by the trivalent vertex $p$, with $t\in\bG_m$, the equalities hold in the preimage of the locus $\mbox{im}(\iota)\cong\{\zeta\neq0\}\sse X(\beta_+)$, and we have recursively defined
$$\eta_j(z_1,\ldots,z_{j};t) := \frac{(C_{j-1})_{i_j+1,i_j+1}\cdot z_j - (C_{j-1})_{i_j,i_j+1}}{(C_{j-1})_{i_j,i_j}}, \qquad C_j := B_{i_j}(z_j)^{-1}C_{j-1}B_{i_j}(\eta_j).$$
for $j\in(k,\ell]$ with initial value
$$C_k(t):= \varphi_{i_k} \begin{pmatrix} t & -1 \\ 0 & t^{-1} \end{pmatrix} = \chi_{i_k}(t)x_{i_k}(-t^{-1}) \in \borel.$$
\sbs
\end{enumerate}
\end{proposition}

\noindent \Cref{prop:change_Reebcoordinates_add_vertex_on_top}.$(i)$ and $(ii)$ follow from the $\SL_3$ matrix identities
\begin{equation}
\label{eq: braid moves}
B_i(z_{k-1})B_j(z_k)B_i(z_{k+1}) = B_j(z_{k+1})B_i(z_{k-1}z_{k+1} - z_k)B_j(z_{k-1}),\mbox{ }B_i(z_{k})B_k(z_{k+1}) = B_k(z_{k+1})B_i(z_k).
\end{equation}
in combination with \eqref{eq:braid_var_eq}.
\noindent \Cref{prop:change_Reebcoordinates_add_vertex_on_top}.$(iii)$ is not as local, as \textcolor{blue}{the new trivalent vertex affects all the Reeb coordinates to the right of where it is inserted}, cf.~\cite[Section 5.1]{CGGLSS25}.\footnote{Geometrically, this is due to the nature of the Ng resolution, see \cite[Section 2.1]{Ng03_computable}.} The core $\SL_2$ identity in that case is
\begin{equation}
\label{eq: trivalent}
B_{i_k}(\zeta_k)B_{i_k}(\zeta) = B_{i_k}(\zeta_k - \zeta^{-1})U,\ \quad  U := \varphi_{i_k}\left(\begin{matrix} \zeta & -1 \\ 0 & \zeta^{-1} \end{matrix} \right).
\end{equation}

\noindent The equation \eqref{eq: trivalent} is defined only when $\zeta \neq 0$ and can be inserted in the middle of a product of several braid matrices. Thus, from the viewpoint of matrix computations, the reason the Reeb coordinates to the right of $\zeta$ are affected is because we must move the element $U \in \borel$ to the right of all the braid matrices $B_{i_j}(z)$ appearing to the right of $B_{i_k}(\zeta)$ in \eqref{eq:braid_var_eq}. This implies that at every trivalent vertex we must modify all the variables appearing to the right of this vertex.\\

\begin{remark}
In the notation of \cite{CGGJ24,CGGLSS25}, the second equation in \Cref{prop:change_Reebcoordinates_add_vertex_on_top}.$(iii)$ would be written as $z_k=\zeta_k-\zeta^{-1}$.\sbs
\end{remark}

\noindent A conceptually important result, which can be proven by combining \Cref{prop:Avariables_up_new_trivalent} and \Cref{prop:change_Reebcoordinates_add_vertex_on_top}.$(iii)$ is the following identity:

\begin{theorem}[in progress]\label{thm:leading_monomial_cluster_A_variables}
Let $\w:\beta\to\delta(\beta)$ be a Demazure weave, $v\in\w$ a trivalent vertex and $\top(\rho_v)\in\N^\ell$ the top weights of the dual Lusztig cycle $\rho_v$. Then the leading monomial of $A^\w_v(z_1,\ldots,z_\ell)$ in the lexicographic order $z_1\prec z_2\prec \ldots\prec z_\ell$ is ${\bf z}^{\top(\rho_v)}$.\sbs
\end{theorem}

\noindent In \Cref{thm:leading_monomial_cluster_A_variables} we used the notation ${\bf z}^{\vec{a}}:=z_1^{a_1}\ldots z_\ell^{a_\ell}$ if $\vec{a}:=(a_1,\ldots,a_\ell)\in\N^\ell$. In the statement, we also implicitly lifted a regular function $A^\w_v(z_1,\ldots,z_\ell)\in\Z[X(\beta)]$ to an actual polynomial, as a regular function in an ambient affine space. There are choices of such lifts, but the scanning of Demazure weaves gives a natural one, which is the one we work with, cf.~\cite[Section 7.2]{CGGLSS25}.

%%%%%%%%%%%%%%%%%%%%%%%%%%%%%%%%%%%%%%%%%%%%%
%%%%%%%%%%%%%%%%%%%%%%%%%%%%%%%%%%%%%%%%%%%%%
\subsection{Hook content formula for cluster $A$-variables}\label{ssec:hook_content_formula} Let $\w:\beta\to\delta(\beta)$ be a Demazure weave and $v\in\w_3$ a trivalent vertex. By \Cref{thm:polynomiality}, $A_v(z)$ is a polynomial in the Reeb $z$-coordinates and, by \Cref{thm:leading_monomial_cluster_A_variables}, the exponents of its leading term are captured by $\top(\rho_v)$. It could be desirable to have a way to compute the entirety of the polynomial $A_v(z)$, with all its monomials and their coefficients, in terms of the combinatorics of the weave $\w$.\\

In order to establish such formula, one must work with raked weaves as in \cite[Section 5.1]{CGGLSS25} and introduce some variations on the dual Lusztig cycle $\rho_v$, which are indexed by bonsai forests (\begin{CJK*}{UTF8}{min}
寄せ植え
\end{CJK*}, yose--ue). A bonsai forest $\mathsf{y}$ involves a choice of certain rooted trees, the bonsais, along each rake of the weave. Such a choice of bonsais along each rake dictates an upwards propagation from $v$, leading to a Lusztig-like weave cycle $\rho_{\mathsf{y}}$ on $\w$. Each bonsai forest $\mathsf{y}$ has a hook content $h(\mathsf{y})$ obtained by multiplying the hooks of certain vertices in its bonsai trees. Let us denote by $\mathsf{F}_v(\w)$ the set of bonsai forests associated to a trivalent vertex $v\in\w_3$. The cluster $A$-variables of a Demazure weave can be expressed combinatorially in terms of the Reeb coordinates as follows:

\begin{theorem}[Hook content formula]\label{thm:hook_content_cluster_A_variables}
Let $\w:\beta\to\delta(\beta)$ be a Demazure weave, $v\in\w$ a trivalent vertex and $A^\w_v(z_1,\ldots,z_\ell)$ the cluster $A$-variable of $v$ in terms of Reeb coordinates. Then
\begin{equation}\label{eq:hook_formula}
\boxed{A^\w_v(z_1,\ldots,z_\ell)=\sum_{\mathsf{y}\in\mathsf{F}_v(\w)} \frac{z^{\on{top}(\rho_{\mathsf{y}})}}{h(\mathsf{y})} }
\end{equation}
where the sum is taken over the collection of bonsai forests starting at $v$.\sbs
\end{theorem}

\noindent There is currently no proof of \Cref{thm:hook_content_cluster_A_variables} printed in the literature, though it can be deduced from the relation between spectral networks and Floer theory as developed in \cite{casals2025spectralnetworksbettilagrangians}. Hopefully a more direct proof of \eqref{eq:hook_formula}, without using \cite{casals2025spectralnetworksbettilagrangians}, will appear within the next year. In terms of \eqref{eq:hook_formula}, it is not clear that the polynomial $A^\w_v(z_1,\ldots,z_\ell)$ has integral coefficients, as the sum occurs with rational coefficients. That said, it follows from the construction of \cite[Section 7.2]{CGGLSS25}, which is well-defined over $\Z$, that $A^\w_v(z_1,\ldots,z_\ell)$ must have integral coefficients. It would be enlightening to have a combinatorial proof of this fact.

%%%%%%%%%%%%%%%%%%%%%%%%%%%%%%%%%%%%%%%%%%%%%
%%%%%%%%%%%%%%%%%%%%%%%%%%%%%%%%%%%%%%%%%%%%%
\subsection{The torus action in the Reeb coordinates}\label{sssec:torus_action_coordinates} For $\G=\SL_n$, the torus action \eqref{eq:torus_action} from \Cref{sssec:braid_variety_further_symmetries} can be written in the Reeb $z$-coordinates explicitly, as established in \cite[Section 2.2]{CGGJ24}. The statement reads as follows:

\begin{lemma}[$\T$-action in Reeb coordinates]\label{lem:torusaction_coordinates} The $\T$-action on $X(\beta)$ from \cref{eq:torus_action}
$$\T\times X(\beta)\lr X(\beta),\quad (\tau;\borel_1,\borel_2,\dots, \borel_{\ell+1})\mapsto (\tau\borel_1,\tau\borel_2,\dots, \tau\borel_{\ell+1})$$
is expressed in the $z$-coordinates for $X(\beta)$ as
$$
(t;z_1,\ldots,z_{\ell}):=(\color{blue}t_{w_1(i_{1})}\color{black}z_1\color{blue}t^{-1}_{w_1(i_{1}+1)}\color{black},\color{blue}t_{w_2(i_{2})}\color{black}z_2\color{blue}t^{-1}_{w_2(i_{2}+1)}\color{black},\ldots,\color{blue}t_{w_{\ell}(i_{\ell})}\color{black}z_{\ell}\color{blue}t^{-1}_{w_{\ell}(i_{\ell}+1)}\color{black}),
$$
%where $c_k(t)$ are given by the formula
%$$c_k=t_{w_k(i_{k}+1)}t^{-1}_{w_k(i_{k})},\quad k\in[1,\ell]$$
\noindent where $w_j:= s_{i_1}\cdots s_{i_{j-1}}\in S_n$ denote the permutations associated to the partial words for $\beta$, $j\in[1,\ell+1]$, and $t\in\T$ and $(z_1,\ldots,z_{\ell})\in \mathbb{A}_\kk^{\ell}$.
%and $w=w_{r+1}$ is the permutation corresponding to $\beta$.
\end{lemma}
\begin{comment}
\begin{proof}
Let us abbreviate $c_k=t_{w_k(i_{k}+1)}t^{-1}_{w_k(i_{k})}$ and consider the diagonal matrix $D(t)$ representing the point $t=(t_1,\ldots,t_n)\in\T$, so its $(i,i)$-entry is $t_i$, $i\in[1,n]$. Then the statement follows from the identity
\begin{equation*}
%\label{eq:slide diagonal 2}
D(t)B_{i_1}(z_1)\cdots B_{i_r}(z_r)=B_{i_1}(c_1z_1)\cdots B_{i_r}(c_rz_r)D(\pi_\beta(t)),
\end{equation*}
where $\pi_\beta=w_{\ell+1}$ is the permutation associated to $\beta$. This identity is implied by the following equality in $\SL_2$:
$$\displaystyle\left(\begin{matrix} t_1 & 0\\ 0 & t_2\\ \end{matrix}\right)\left(\begin{matrix} z & -1\\ 1 & 0\end{matrix}\right)=\left(\begin{matrix}t_1zt_2^{-1} & -1\\
1 & 0\\ \end{matrix}\right)\left(\begin{matrix} t_2 & 0\\ 0 & t_1\\ \end{matrix}\right).$$
\end{proof}
\end{comment}
\noindent Visually, the indices $i_k,j_k\in[1,n]$ of $t_{i_k}z_kt_{j_k}^{-1}$ can be read from the braid diagram as follows, cf.~\cite[Section 2.2]{CGGJ24}. Consider the $k$th crossing in the braid diagram, which has the variable $z_k$ associated to it, and follow the top left strand at that crossing to the left until reaching the beginning of the braid diagram for $\beta$, at its leftmost part. This selects a strand at the leftmost part of the braid diagram $\beta$: the index of that strand is then $j_k$, indexing strands counting from the bottom. Similarly, the index $i_k$ is that of the strand at the leftmost part of the braid diagram for $\beta$ which is reached by following the bottom left strand at the $k$th crossing to the left.

%%%%%%%%%%%%%%%%%%%%%%%%%%%%%%%%%%%%%%%%%%%%%
%%%%%%%%%%%%%%%%%%%%%%%%%%%%%%%%%%%%%%%%%%%%%
\subsection{The cluster 2-form in the Reeb coordinates}\label{sssec:symplectic_form_coordinates} Given a cluster seed $\{A_i\}$ with a skew-symmetric exchange matrix $B=(b_{ij})$ we can consider the following cluster 2-form
\begin{equation}\label{eq:cluster_2form}
\omega_{cl}:=\sum_{i,j}b_{ij}\frac{dA_i\wedge dA_j}{A_iA_j}=\sum_{i,j}b_{ij}\cdot d\log(A_i)\wedge d\log(A_j)
\end{equation}
\noindent in the torus $\Spec \kk[A_1^{\pm1},\ldots,A_\ell^{\pm1}]$ underlying the seed. By the mutation formula for cluster $A$-variables, the 2-forms associated to all the different cluster seeds glue together to give a well-defined global 2-form on (the spectrum of) the cluster algebra. Such cluster 2-form is a central ingredient in the study of symplectic and Poisson structures in relation to cluster algebras, see e.g.~\cite{FG2006,FG09,GSV05,GSV}. The goal of this subsection is to present an expression of the cluster 2-form \eqref{eq:cluster_2form} in terms of the Reeb coordinates.\\

For simplicity, we assume that $\delta(\beta) = w_0$ and $\kk=\C$. Let $\theta^L$ and $\theta^{R}$ denote the left and right invariant $\mathfrak{g}$-valued forms on $\G$, respectively, a.k.a.~the Maurer-Cartan forms. For instance, if $\G = \SL_n$ then we can write $\theta^L(f) = f^{-1}df$ and $\theta^{R}(g) = dg\cdot g^{-1}$ explicitly. Let $\kappa: \mathfrak{g} \otimes \mathfrak{g} \to \kk$ be the (normalized) Killing form on the Lie algebra $\mathfrak{g}$ of $\G$. Consider the following $2$-form on $\G \times \G$:
\begin{equation}\label{eq:fg_pairing}
(f|g) := \kappa(\theta^L(f)\wedge\theta^{R}(g))
\end{equation}
\noindent where $f,g\in\G$. Here \eqref{eq:fg_pairing} defines a skew-symmetric bilinear form at $T_{(f,g)}\G\times\G$ which, when varying the point $(f,g)\in\G\times\G$ defines an element of $\Omega^2(\G\times\G)$. Note that $(\cdot|\cdot)$ as defined in \eqref{eq:fg_pairing} satisfies a cocycle condition:
\begin{equation}
\label{eq: cocycle}
(f|g)+(fg|h)=(f|gh)+(g|h).
\end{equation}
Given a collection of $\G$--valued functions $f_1,\ldots,f_{\ell}$, we define
\begin{equation}
    \label{eq: 2 form}
(f_1|\cdots |f_{\ell}):=(f_1|f_2)+(f_1f_2|f_3)+\ldots+(f_1\cdots f_{\ell-1}|f_{\ell}).
\end{equation}
\noindent By \eqref{eq: cocycle} this definition is associative in the entries $f_i$. Let $\beta=\sigma_{i_1}\cdots\sigma_{i_\ell}$, we relate this construction to the braid variety $X(\beta)$ by applying \eqref{eq: 2 form} with $f_j:=B_{i_j}(z_j)$. Specifically, by using the embedding $X(\beta)\sse\mathbb{A}_\kk^\ell$ given by the Reeb coordinates in \Cref{lem:equations_braidvar} we define:
\begin{definition}\label{def:braid_variety_2_form} Let $\beta=\sigma_{i_1}\cdots\sigma_{i_\ell}$ be a positive braid word and $X(\beta)\sse\mathbb{A}_\kk^\ell$ its braid variety embedded via the Reeb coordinates. By definition, the 2-form $\omega_\beta\in\Omega^2X(\beta)$ is the restriction
\begin{equation}\label{eq:braid2form}
\omega_\beta:=(B_{i_1}(z_1)|\cdots |B_{i_{\ell}}(z_{\ell}))|_{X(\beta)},
\end{equation}
where $(B_{i_1}(z_1)|\cdots |B_{i_{\ell}}(z_{\ell}))\in\Omega^2\mathbb{A}_\kk^\ell$ is a 2-form in the Reeb coordinates.\sbs
\end{definition}
By \Cref{def:braid_variety}, the image of the map $B_{\beta}: \mathbb{A}_\kk^\ell \to \G$ restricted to the braid variety must have its image contained in $\delta(\beta)\borel$. Thus $\omega_\beta$ is a closed 2-form on $X(\beta)$. The following result is \cite[Theorem 9.5]{CGGLSS25}:

\begin{proposition}\label{prop:2forms}
Let $\beta$ be a positive braid word, $X(\beta)$ its braid variety and $\omega_{cl}\in\Omega^2 X(\beta)$ the cluster 2-form defined by any cluster seed coming from a Demazure weave $\w:\beta\to\delta(\beta)$. Then
\begin{equation}\label{eq:equality_2_forms}
\boxed{\omega_\beta=\omega_{cl}}
\end{equation}
\sbs
\end{proposition}

\noindent In particular, since $\omega_\beta$ is defined entirely in terms of the Reeb coordinates, \Cref{prop:2forms} expresses the cluster 2-form $\omega_{cl}$ in terms of the Reeb coordinates.

%%%%%%%%%%%%%%%%%%%%%%%%%%%%%%%%%%%%%%%%%%%%%
%%%%%%%%%%%%%%%%%%%%%%%%%%%%%%%%%%%%%%%%%%%%%
\subsection{The Landau-Ginzburg brick potential in the Reeb coordinates}\label{thm:superpotential_Reeb_coordinates} Building on the brick compactifications $(Y(\beta),D_\beta)$ introduced in \Cref{ssec:brick_compactifications}, we now proceed to computing the potentials $W_\beta$ associated to such compactifications, cf.~\Cref{sssec:applications_compactifications}. Specifically, here the term {\it potential} is understood within the context of Landau-Ginzburg models, see e.g.~\cite[Sections 3\& 11]{GoncharovShen15}, \cite[Sections 3\&4]{GS19} and \cite[Section 0.4]{GHKK18}.\footnote{This use should not be confused with the notion of potentials in the context of {\it quivers with potential}, which is also an important but different aspect of the study of cluster structures in braid varieties, cf.~\cite{CasalsGao24}.} For specificity, let us assume $\delta(\beta)=w_0$.

\begin{remark} As the upcoming discussion involves a form of mirror symmetry, it is paramount to distinguish whether we are working in the affinization of the $\SA$-scheme or that of the $\SX$-scheme, especially because both happen to be isomorphic to the braid variety $X(\beta)$. In addition, the choice of Lie group also needs to be emphasized. Thus, as in \Cref{sssec:applications_compactifications}, we adopt the notation of $\SA_\G(\beta)$ and $\SX_\G(\beta)$ as presented in \Cref{ssec:A_and_X_schemes_weave}.\sbs
\end{remark}

Let us consider the braid variety $\SX_{\G^\vee}(\beta)$ and use the compactification $Y_{\G^\vee}(\beta)=Y(\beta)$ from \Cref{def:brick}, which yields the log Calabi-Yau pair $(Y_{\G^\vee}(\beta),D_\beta)$, cf.~\Cref{prop:logCY}. Following a tenet in mirror symmetry, cf.~\cite[Section 7.3]{HoriVafa} and \cite[Conjecture 1.1]{Auroux07_Tduality}, compactifying $\SX_{\G^\vee}(\beta)$ should correspond to adding the data of a potential in $\SA_{\G}(\beta)$, see also \cite[Section 1]{GHK15} or \cite[Section 5]{GHKK18}. In our case, we consider the opposite cluster ensemble $(\SA_{\G}(\beta),\SX_{\G^\vee}(\beta))$ as in \cite[Section 1.2]{FockGoncharov09_Dilog2}. This opposite ensemble has the negated $p$-map, so that the $\SX$-variables in the cluster $\SX$-seeds are the inverses, and the exchange matrix also acquires an overall negative sign.\footnote{Note that the cluster $\SA$-mutation rule is invariant under such sign change.} This is done because, in general, each of the boundary divisors in $D_\beta$ gives a discrete valuation $-\ord_D$, in the $\max$-tropical convention, but in this case we want to measure $\ord_D$, and thus it suffices to consider the opposite cluster $\SX$-structure. Thus, given an irreducible divisor component $D$ of the compactifying divisor $D_\beta$, we consider the valuation vector
\begin{equation}\label{eq:tropicallabel}
        q_D:=\bigl(-\ord_D(X_v^{-1}\bigr)_v
             =\bigl(\ord_D(X_v)\bigr)_v,
\end{equation}
where the index for the entries runs over the trivalent vertices of $\w$. Here we understand the valuation $\ord_D:\kk(\SX_{\G^\vee}(\beta))\lr\Z\cup\{\infty\}$ as defining $q_D$, whose entries are seen as tropical points in $\SX_{\G^\vee}(\beta)$, i.e.~as elements of $\SX_{\G^\vee}(\beta)(\Z^T)$. Following \cite[Corollary 11.9]{GHKK18}, we define the potential
\begin{equation}\label{eq:superpotential}
W_{Y(\beta)}:\SA_\G(\beta)\lr\mathbb{A}^1_\kk, \qquad W_{Y(\beta)}:=\sum_{j}\tht_{q_{D_j}},
\end{equation}
where the index $j$ runs over the positions $j\in[1,\ell]$ such that $\delta(\beta\setminus \sigma_{i_j})=\delta(\beta)$, cf.~\Cref{lem:boundary_div}.(1), and the theta functions $\vartheta_{q_D}$ in \eqref{eq:superpotential}, associated to a divisorial valuation $q_D$, are defined as in \cite[Definition 3.3]{GHKK18}. The following result expresses the potential \eqref{eq:superpotential} in terms of the Reeb coordinates:

\begin{theorem}\label{thm:superpotential_Reeb_coordinates} Let $\beta\in\brn$ be a positive braid word, $(\SA_\G(\beta),\SX_{\G^\vee}(\beta))$ its opposite cluster ensemble and $(Y_{\G^\vee}(\beta),D_\beta)$ the brick compactification of $\SX_{\G^\vee}(\beta)$. Then:
\begin{enumerate}
    \item For every non-empty divisor $D_j$ in $D_\beta$,
    \begin{equation}\label{eq:individual}
       \tht_{q_{D_j}}=z_j,
    \end{equation}
where $z_j$ is the $j$th Reeb coordinate of $\SA_\G(\beta)$.\\

\item The potential associated to the compactification $(Y_{\G^\vee}(\beta),D_\beta)$ satisfies
\begin{equation}\label{eq:superpotential_Reeb_coordinates}
 \boxed{\displaystyle
 W_{Y(\beta)}=\sum_{j}z_j}
\end{equation}
where the index runs over $j\in[1,\ell]$ such that $\delta(\beta\setminus \s_{i_j})=\delta(\beta)$.\sbs
\end{enumerate}
\end{theorem}

Note that \eqref{eq:superpotential_Reeb_coordinates} is independent of any choice of the Demazure weave. With regards to \Cref{thm:superpotential_Reeb_coordinates}.(1), there in fact exists a Demazure weave $\w$ whose associated cluster $\SX$-seed is such that the vector $(\ord_{D_j}X_v)_v$ is integral and primitive, with non-negative entries at the mutable vertices, and the Reeb coordinate $z_j$ is the cluster monomial with this exponent vector. About \cref{thm:superpotential_Reeb_coordinates}.(2), note that the coordinate $z_j$ vanishes identically on $X(\beta)$ if $\delta(\beta\setminus\s_{i_j})\neq \delta(\beta)$. Thus \eqref{eq:superpotential_Reeb_coordinates} in \Cref{thm:superpotential_Reeb_coordinates} can be written as
\begin{equation}\label{eq:superpotential_Reeb_coordinates2}
\displaystyle
 W_{Y(\beta)}=\sum_{i=1}^{\ell(\beta)}z_i.
\end{equation}
\Cref{thm:superpotential_Reeb_coordinates} can be proven with the same strategy as \cite[Section 5.5]{CGGLSS25}. Namely, one first computes directly with a weave as in \cite[Figure 16 (left)]{CGGLSS25} in order to deduce \Cref{thm:superpotential_Reeb_coordinates}.(1) for the first Reeb coordinate $z_1$, and then applies cyclic rotations to deduce the statement for the remaining $z_i$ coordinates.

\begin{remark} For reference, a version of \Cref{thm:superpotential_Reeb_coordinates}.(2) was discussed by one of the research groups at the AIM Workshop ``Cluster algebras and braid varieties'' (2023). I recently talked with D.~Weng, one of the active participants in that group, about whether the result would be written anywhere by that group. After that conversation, there seems to be some merit in at least recording a precise statement, for reference, as done in \Cref{thm:superpotential_Reeb_coordinates} above.\sbs
\end{remark}

\part{\centering Demazure weaves in contact and symplectic topology}

\begin{center}
	\begin{figure}[h!]
		\centering
		\includegraphics[scale=0.75]{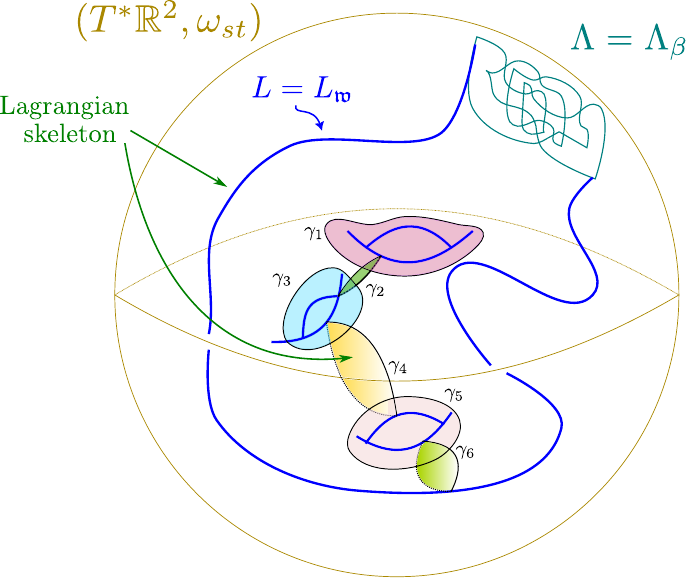}
		\caption{A depiction of the Lagrangian surface $L_\w\sse\D^4$, constructed from $\wb$ bounding the Legendrian link $\La_\beta\sse S^3$, constructed from $\beta$. The mutable Lusztig cycles $\g_i^\w$ are geometrically realized by 1-cycles in $L_\w$ that bound embedded Lagrangian disks in the complement of $L_\w$, as it is hopefully illustrated.}
        \label{fig:Filling_Merodromy2}
	\end{figure}
\end{center}

The object of this second part of the article is to present the diagrammatic calculus of Demazure weaves from the viewpoint of low-dimensional contact and symplectic topology. Specifically, we focus on the case of $\G=\SL_n$ throughout Part II and we aim at explaining and constructing the concepts in the middle column of \Cref{table:PartII_introduction2}. In a nutshell, we consider a Legendrian link $\La_\beta\sse\R^3$ associated to a positive braid word $\beta$ and $X(\beta)$ appears as a sort of moduli for Lagrangian surfaces in $\D^4$ that bound such a Legendrian link, endowed with a $\GL_1$-local system, cf.~\Cref{fig:Filling_Merodromy2}. The key concepts featuring in the cluster algebra for $\Z[X(\beta)]$ are realized through this geometric lens, including the $\SA$ and $\SX$-lattices, the quiver, the cluster $\SA$ and $\SX$-variables and the cluster mutation formula, see \Cref{table:PartII_introduction2}. Geometrically, a core ingredient is the construction of an embedded exact Lagrangian filling $L_\w$ of $\La_\beta$ for any Demazure weave $\wb$. In many key cases, the Demazure weave also endows $L_\w$ with an $\L$-compressing system, which allows us to explain cluster mutations from the viewpoint of Lagrangian disk surgeries.\\

\noindent Beyond its aesthetic value, such geometric realization of the cluster algebras $\Z[X(\beta)]$ and the weaves $\wb$ both provide new results in low-dimensional contact and symplectic topology, and informs the construction and discovery of new algebraic structures motivated by existing geometric constructions on the symplectic side. Among other instances: the introduction of weave calculus itself was originally motivated by symplectic topology in \cite{CasalsZaslow}, and the construction of cluster algebra structures on the coordinate rings of many interesting varieties, including open Richardson varieties, as illustrated in \cite{CGGLSS25}, was informed by existing results in symplectic topology.

%%%%%%%%%%%%%%%%%%%%%%%%%%%%%%%%%%%%%%%%%%%%%%%%%%%%%%%%
%%%%%%%%%%%%%%%%%%%%%%%%%%%%%%%%%%%%%%%%%%%%%%%%%%%%%%%%
%%%%%%%%%%%%%%%%%%%%%%%%%%%%%%%%%%%%%%%%%%%%%%%%%%%%%%%%
\section{Legendrian links and their Lagrangian fillings}\label{sec:Legendrian_links_Lagrangian_fillings}

The object of this section is to explain how to transition from the algebraic and diagrammatic aspects of Demazure weaves, as presented in Part I, to the geometric context of contact and symplectic topology. We focus on two constructions:\\

\begin{enumerate}
    \item Given a positive braid word $\beta$, we obtain a Legendrian link $\La_\beta\sse(S^3,\xi_\std)$ in the standard contact 3-sphere, cf.~\Cref{ssec:Legendrian_link_beta}, and explain how the braid variety $X(\beta)$ can be understood in terms of such a Legendrian link, cf.~\Cref{ssec:braid_varieties_Legendrian_links}.\\
    \item After reviewing Lagrangian fillings in \Cref{ssec:Lagrangian_fillings}, we then construct in \Cref{ssec:Lagrangian_filling_from_AG} embedded exact Lagrangian surfaces $L\sse(\D^4,\la_\std)$ by studying certain algebraic curves in $\C^2$ with purely imaginary periods. Such a class of Lagrangian surfaces, and the way in which they can be plotted, provide a compelling motivation to introduce weave diagrams $\wb$ and later construct Lagrangian surfaces $L_\w$ from them.
\end{enumerate}

\begin{comment}

\noindent A central geometric goal in low-dimensional contact and symplectic topology is to \textcolor{blue}{classify embedded exact Lagrangian surfaces with a fixed Legendrian boundary link, up to compactly supported Hamiltonian isotopy}. Part of the motivation to study Demazure weaves and cluster structures on braid varieties is that new results in symplectic topology can be obtained towards that classification.\\

\noindent In fact, from a hopeful viewpoint, I conjecture that the classification of such Lagrangian surfaces bounding $\La_\beta$ is given by the cluster seeds of the cluster algebra structure on the coordinate ring of functions of the braid variety $X(\beta)$. We will present and discuss results providing some evidence towards that conjecture. For now, we start building the geometric objects from the ground up and relate them to the diagrammatic calculus of Demazure weaves from \Cref{sec:diagrammatic_weaves}.
\end{comment}

%%%%%%%%%%%%%%%%%%%%%%%%%%%%%%%%%%%%%%%%%%%%%%%%%%%%%%%%
%%%%%%%%%%%%%%%%%%%%%%%%%%%%%%%%%%%%%%%%%%%%%%%%%%%%%%%%
%%%%%%%%%%%%%%%%%%%%%%%%%%%%%%%%%%%%%%%%%%%%%%%%%%%%%%%%
\subsection{The Legendrian link of a positive braid word}\label{ssec:Legendrian_link_beta}

Let $\beta$ be a positive braid word. The first step towards geometry is to consider a Legendrian link $\La_\beta\sse(S^3,\xi_\std)$ associated to $\beta$. A Legendrian link embedded in the 3-sphere $S^3$ is a particular type of smooth link. Specifically, it is a link constrained to a tangency condition with respect to the 2-plane distribution $\xi_\std$ given by the complex tangencies of the unit sphere
$$S^3=\{(z_1,z_2)\in\C^2: |z_1|^2+|z_2|^2=1\}\sse\C^2,\qquad \xi_\std:=TS^3\cap i(TS^3).$$

\begin{definition}[Legendrian links]\label{def:legendrian_link} Let $\La\sse S^3$ be a smooth link embedded in the 3-sphere. By definition, $\La$ is said to be a Legendrian link if $T\La\sse\xi_\std$, i.e.~if the tangent line $T_p\La\sse T_pS^3$ at every point $p\in\La$ belongs to the 2-plane of complex tangencies $T_pS^3\cap i(T_pS^3)$.\\

\noindent By definition, two Legendrian links $\La_0,\La_1\sse (S^3,\xi_\std)$ are said to be Legendrian isotopic if there exists a 1-parametric family $\La_t\sse(S^3,\xi_\std)$ of Legendrian links with $t\in[0,1]$.\sbs
\end{definition}

\begin{center}
	\begin{figure}[h!]
		\centering
		\includegraphics[scale=0.85]{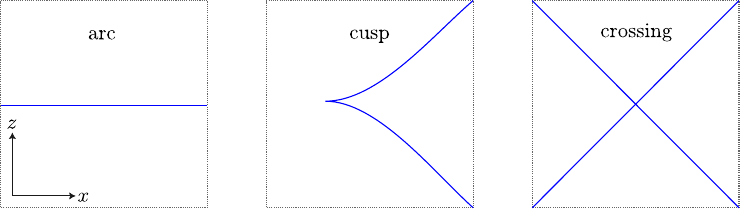}
		\caption{The generic 3 planar pieces of a front $\ff\sse\R^2_{x,z}$ of a Legendrian link $\La\sse(\R^3,\xi_\std)$.}
		\label{fig:Front_Arc_Cusp_Crossing_R3}
	\end{figure}
\end{center}

In the same manner that smooth links in the 3-sphere can be studied by using link diagrams in the plane, Legendrian links can be studied by using front diagrams. In contrast to link diagrams, fronts do not need to specify the overcrossing and undercrossing information. For our purposes, a front for a Legendrian link can be defined as an immersed curve $\mathfrak{f}$ in the plane $\R^2_{x,z}$, except for finitely many simple cusp singularities, having no vertical tangencies. See \Cref{fig:Front_Arc_Cusp_Crossing_R3_Example1} for two instances of fronts. Thus a front is locally built of the three types of pieces depicted in \Cref{fig:Front_Arc_Cusp_Crossing_R3}, consisting of arcs, simple cusps and crossings, none of which ever acquires a vertical tangent line. These are locally modeled as
$$\mbox{arc}:=\{z=0\},\quad \mbox{cusp}:=\{x^3=z^2\},\quad\mbox{crossing}:=\{(x-z)(x+z)=0\}.$$

\noindent From a front $\mathfrak{f}\sse \R^2_{x,z}$ we can construct a Legendrian link:

\begin{definition}[Legendrian link from a front]\label{def:legendrian_link_from_front} Let $\mathfrak{f}\sse \R^2_{x,z}$ be a front, then the Legendrian link $\La_\mathfrak{f}$ is
\begin{equation}\label{def:leglink_from_front}
\La_\mathfrak{f}:=\{(x,y,z)\sse\R^3: y=\dd_x z(x)\}
\end{equation}
where at each point we have locally expressed $\mathfrak{f}$ as the graph of a function $z(x)$.\sbs
\end{definition}

\noindent A few comments on \Cref{def:legendrian_link_from_front}:

\begin{enumerate}
    \item The condition $y=\dd_x z(x)$ guarantees that $\La_\mathfrak{f}$ is tangent to the 2-plane distribution $\xi:=\ker\{dz-ydx\}\sse T\R^3$ in $\R^3$. This is relevant to \Cref{def:legendrian_link} because the complement of a point in $(S^3,\xi_\std)$ is isomorphic\footnote{The precise term is contactomorphic, see \cite[Chapter 2]{Geiges08}.} to $(\R^3,\ker\{dz-ydx\})$. In particular, by considering the one-point compactification $\R^3\sse S^3$, we can always think of $\La_\mathfrak{f}$ as a Legendrian link in $(S^3,\xi_\std)$, which we will still denote by $\La_\ff$. In a sense, \Cref{def:legendrian_link_from_front} is the analogue in contact topology to the fact in smooth topology that a link diagram in the plane determines a unique link in $\R^3$, up to smooth isotopy.\\

    \item For \eqref{def:leglink_from_front}, it is possible to always locally express the restriction of the $z$-coordinate in $\R^2_{x,z}$ to $\mathfrak{f}$ as a function $z|_{\ff}=z(x)$ because of the implicit function theorem, as there are no vertical tangencies in $\mathfrak{f}$. In addition, the Legendrian link $\La_\ff\sse\R^3$ is a smooth link. Indeed, at a cusp singularity $\{x^3=z^2\}$ of $\ff$, as in \Cref{fig:Front_Arc_Cusp_Crossing_R3} (center), the lift can be smoothly parametrized by $\gamma(t)=(t^2,3t/2,t^3)$. Similarly, at a crossing $\{(x-z)(x+z)=0\}$ of $\ff$ the crossing point lifts to two distinct points $(0,\pm 1,0)$ in $\R^3$. These points are distinguished in $\R^3$ by the $y$-coordinates, which measure the slopes of the branches at the crossing point in $\ff\sse\R^2_{x,z}$.\\
\end{enumerate}

\begin{center}
	\begin{figure}[h!]
		\centering
		\includegraphics[width=\textwidth]{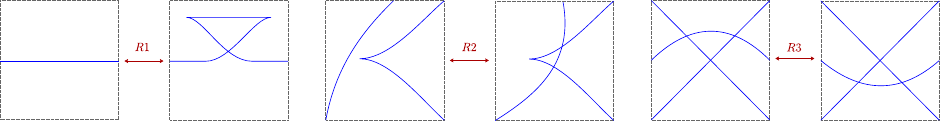}
		\caption{Legendrian Reidemeister moves in the front projection. Each of these moves can be realized by a Legendrian isotopy. Conversely, in $(\R^3,\xi_\std)$ every Legendrian isotopy is itself isotopic to a Legendrian isotopy whose front projections are exactly a sequence of these three Reidemeister moves. There are symmetric versions of these moves which we always implicitly include: R1 upside-down, R2 with the cusp pointing east.}
		\label{fig:Front_Arc_Cusp_Crossing_R3_Reidemeister}
	\end{figure}
\end{center}

In smooth topology, the Reidemeister theorem states that two link diagrams give rise to two smoothly isotopic links in $\R^3$ if and only if the link diagrams differ by certain moves, known as the Reidemeister moves. Similarly, a Legendrian enhancement of the Reidemeister theorem states that given two fronts $\ff_0,\ff_1$, then $\La_{\ff_0}\simeq \La_{\ff_1}$ are Legendrian isotopic if and only if $\ff_0$ and $\ff_1$ can be connected through the three Legendrian Reidemeister moves, as depicted in \Cref{fig:Front_Arc_Cusp_Crossing_R3_Reidemeister}, and a compactly supported ambient smooth isotopy where the fronts do not acquire vertical tangencies.

\noindent Remarkably, there are typically many Legendrian links that are not Legendrian isotopic even if their underlying smooth links are smoothly isotopic. For instance, \Cref{fig:Front_Arc_Cusp_Crossing_R3_Example1} depicts two Legendrian fronts $\ff_0,\ff_1$ such that $\La_{\ff_0}$ and $\La_{\ff_1}$ are smoothly isotopic but not Legendrian isotopic.

\begin{center}
	\begin{figure}[h!]
		\centering
		\includegraphics[scale=0.5]{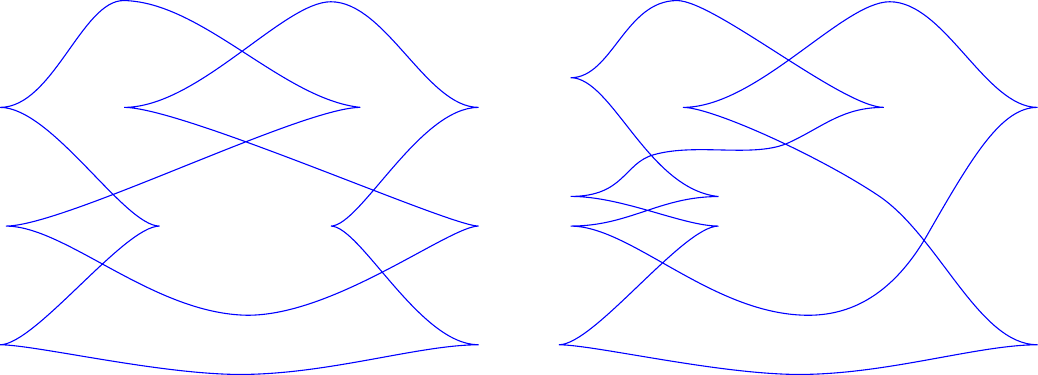}
		\caption{Two fronts for two Legendrian knots $\La\sse(\R^3,\xi_\std)$, fronts drawn in $\R^2_{x,z}$. The smooth type of both these knots in $\R^3$ is that of $m(5_2)$, the mirror of the three twist knot: the only 5-crossing knot apart from the $(2,5)$-torus knot. These Legendrian knots are smoothly isotopic but not Legendrian isotopic, as first proven in \cite[Example 4.4]{Chekanov}.}
		\label{fig:Front_Arc_Cusp_Crossing_R3_Example1}
	\end{figure}
\end{center}

\begin{center}
	\begin{figure}[h!]
		\centering
		\includegraphics[scale=0.9]{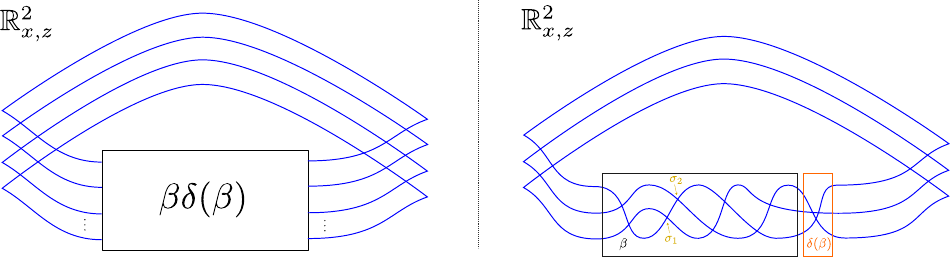}
		\caption{(Left) The front in $\R^2_{x,z}$ for $\La_\beta\sse(\R^3,\ker\{dz-ydx\})$ referred to as the $(-1)$-closure of $\beta\delta(\beta)$. (Right) The front of the Legendrian link $\La_\beta$ for the 3-stranded positive braid word $\beta=\s_2\s_1(\s_1\s_2)^3$.}
		\label{fig:FrontsBeta}
	\end{figure}
\end{center}

For our purposes, it will suffice to focus on the following class of Legendrian links:

\begin{fdef}
[Legendrian link associated to $\beta$]\label{def:Legendrian_from_braid} Let $\beta\in \brn$ be a positive braid word and $\delta(\beta)$ its Demazure product. By definition, $\La_\b\sse(\R^3,\xi_\std)$ is the Legendrian link defined by the front in \cref{fig:FrontsBeta} (left). It is said to be the $(-1)$-closure of $\beta\delta(\beta)$.\sbs
\end{fdef}
\vspace{0.2cm}

See \Cref{fig:FrontsBeta} (right) for an instance of \Cref{def:Legendrian_from_braid}. A key connection between braid varieties and symplectic topology is that the braid variety $X(\beta)$, as discussed in \Cref{sec:weaves_braid_varieties}, is intimately related to the Legendrian link $\La_\beta$, as we will momentarily explain.

\begin{center}
	\begin{figure}[h!]
		\centering
		\includegraphics[scale=0.9]{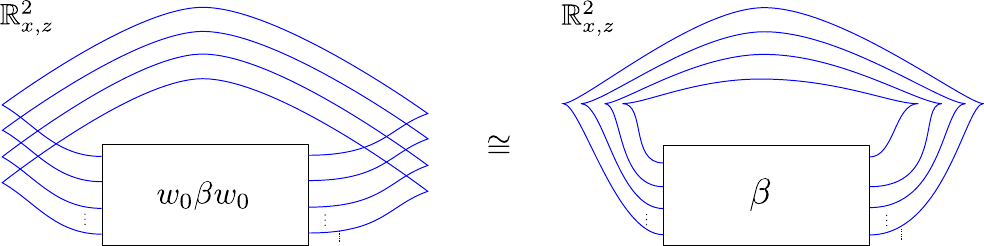}
		\caption{(Left) The $(-1)$-closure of $w_0\beta w_0$, where $w_0$ is any positive braid word representing the half-twist. (Right) The rainbow closure of $\beta$, which is Legendrian isotopic to the $(-1)$-closure of $w_0\beta w_0$.}
		\label{fig:FrontsBeta2}
	\end{figure}
\end{center}

\noindent A few comments on \Cref{def:Legendrian_from_braid} are in order:

\begin{enumerate}
    \item If two braid words $\beta_1,\beta_2$ differ by braid moves $\s_i\s_{i+1}\s_i\to\s_{i+1}\s_i\s_{i+1}$, then $\La_{\beta_1}\simeq\La_{\beta_2}$ are Legendrian isotopic by the third Legendrian Reidemeister move, as in \Cref{fig:Front_Arc_Cusp_Crossing_R3_Reidemeister} (right). Similarly, if they differ by a commutation $\s_i\s_j\to\s_j\s_i$, $|i-j|\geq2$, then $\La_{\beta_1}\simeq\La_{\beta_2}$ are Legendrian isotopic, as there is an ambient planar isotopy that relates these fronts. Therefore, the Legendrian isotopy type of $\La_\beta$ only depends on the underlying braid element $[\beta]$, i.e.~it is independent of the choice of braid word $\beta$ for $[\beta]$.\\

    \item The literature has also featured the rainbow closure of a positive braid word $\beta$, as depicted in \cref{fig:FrontsBeta2} (right), cf.~\cite[Section 6.5]{STZ2017}. It is an exercise to show that this is an instance of a $(-1)$-closure, specifically the $(-1)$-closure of $w_0\beta w_0$, cf.~\cref{fig:FrontsBeta2}. Also, the smooth type of $\La_\beta$ is that of the smooth link typically associated in smooth topology to the braid $w_0^{-1}\beta\delta(\beta) w_0^{-1}$, by closing the braid smoothly without any framing, i.e.~the 0-closure, with no twists. Briefly put, the crossings in each of the two sides of the $\beta\delta(\beta)$-box in \cref{fig:FrontsBeta} (left) are each smoothly contributing $w_0^{-1}$.\\

    \item From the viewpoint of smooth knot theory, the class of Legendrians $\La_\beta\sse(\R^3,\xi_\std)$ is reasonably broad, as it includes smooth knots in each of the three classes: torus knots, satellite knots and hyperbolic knots. Thurston's results imply that any smooth knot must be in one of these three classes, cf.~\cite[Cor.~2.5]{Thurston82}. In fact, for each of these smooth classes, there exist infinitely many Legendrian knots of the form $\La_\beta$ in such class, each not smoothly isotopic to one another, cf.~\cite[Section 6.2]{CasalsGao22}.
\end{enumerate}

\begin{example} (1) Consider $c:=\sigma_1\ldots\sigma_{n-1}\in\Br_n^+$ and let $w_0$ be any positive braid word representing the half-twist, e.g.~$w_0=\sigma_1(\sigma_2\sigma_1)\cdots(\sigma_{n-1}\cdots\sigma_2\s_1)$. Then the Legendrian link $\La_\b$ associated to $\beta=w_0c^m$ is a Legendrian representative of the $(n,m)$-torus link. If $\gcd(n,m)=1$, then $\La_\beta$ is a knot and it is the unique max-tb Legendrian representative of the $(n,m)$-torus link. In the context of cluster algebras, these Legendrian links will be closely related to the top-dimensional positroid variety in the Grassmannian $\mbox{Gr}(n,n+m)$, cf.~\Cref{ssec:braid_varieties_examples}.(2) and \Cref{ex:simple_braidvarieties}.(2). More generally, the max-tb Legendrian representative of an algebraic knot is unique and it is of the form $\La_{w_0}\beta$ for some $\beta\in\Br_n^+$, cf.~\cite{CasalsLagSkel}.\\

\noindent (2) If $\beta=\s_1^2\s_2^2\s_1^2\s_2^2\s_1\s_2\in\mbox{Br}_3^+$, then $\La_\b$ is Legendrian isotopic to the Legendrian representative of $m(5_2)$ depicted in \cref{fig:Front_Arc_Cusp_Crossing_R3_Example1} (left). The Legendrian representative of $m(5_2)$ depicted in \cref{fig:Front_Arc_Cusp_Crossing_R3_Example1} (right) is {\it not} of the form $\La_\beta$ for any $n\in\N$ and $\beta\in\mbox{Br}_n^+$, see e.g.~\cite[Section 7.2.2]{STZ2017}.\sbs
\end{example}

%%%%%%%%%%%%%%%%%%%%%%%%%%%%%%%%%%%%%%%%%%%%%%%%%%%%%%%%
%%%%%%%%%%%%%%%%%%%%%%%%%%%%%%%%%%%%%%%%%%%%%%%%%%%%%%%%
%%%%%%%%%%%%%%%%%%%%%%%%%%%%%%%%%%%%%%%%%%%%%%%%%%%%%%%%
\subsection{Braid varieties from Legendrian links}\label{ssec:braid_varieties_Legendrian_links} A driving goal in low-dimensional contact topology is the classification of Legendrian links, up to Legendrian isotopy. Part of such classification involves showing that some Legendrian links are {\it not} Legendrian isotopic to each other. This is achieved by developing Legendrian invariants. There is now a range of such Legendrian invariants, from numerical ones, such as the Thurston-Bennequin number and the rotation number, to polynomial invariants, such as the ruling polynomial, and algebraically richer ones, including certain dg-categories and $A_\infty$-categories. Let us focus on two such invariants:\\

\begin{enumerate}
    \item {\bf The Legendrian contact dg-algebra}. For a Legendrian link $\La\sse(\R^3,\xi_\std)$ and a choice of basepoints $\tt$, this is a non-commutative dg algebra $\mathbb{A}(\La,\tt)$ whose quasi-isomorphism type only depends on the Legendrian isotopy class of $\La$ and the choice of basepoints $\tt$, cf.~\cite{Chekanov,EtnyreNg22_Survey}. For context, the dg-category of perfect dg-modules over $\mathbb{A}(\La,\tt)$ is an explicit presentation of the partially wrapped Fukaya category of $(\D^4,\la_\std)$ stopped at $\La$, cf.~\cite{EkholmLekili2023}.\\
    
    For the Legendrian links of the form $\La=\La_\beta$, as in \Cref{def:Legendrian_from_braid}, this dg algebra $\mathbb{A}(\La_\beta,\tt)$ is defined over $\Z$ and it is $\Z$-graded. An explicit presentation for $\mathbb{A}(\La_\beta,\tt)$ is given in \cite[Section 5.1]{CasalsNg22}, and its homology is reasonably understood, see e.g.~\cite[Corollary 1.1]{casals2026vanishinghigherlegendrianhomology}. In this case of $\La=\La_\beta$ we always choose the basepoints $\tt$ to be one per strand of $\beta$, and located to the right of all crossings of $\beta$, cf.~\cite[Section 5]{CasalsNg22}. To ease notation, we often abbreviate $\mathbb{A}(\La_\beta,\tt)$ to $\mathbb{A}(\La_\beta)$ since the basepoints are understood to be chosen according to this convention.\\

    \item {\bf The dg category of sheaves with singular support on $\La$}. Building on the notion of singular support of a sheaf introduced in \cite{KashiwaraSchapira_Book}, the dg category $\operatorname{Sh}_\La(\R^2)$ of (categorically) compact complexes of sheaves of $\kk$-modules on $\R^2_{x,z}$ with singular support contained in the (cone of the) Legendrian link $\La$ is a Legendrian invariant of $\La$, cf.~\cite{GKS_Quantization}. It is conjectured that $\operatorname{Sh}_\La(\R^2)$ is equivalent to the dg-category of perfect dg-modules over $\mathbb{A}(\La)$.\\
\end{enumerate}

A central connection between braid varieties and Legendrian links is as follows:\\

\begin{center}
  \fbox{%
    \begin{tabular}{l}
      \itshape The braid variety $X(\beta)$ appears naturally as a Legendrian invariant of $\La_\beta$, \\
      \itshape either via its Legendrian contact dg algebra or its dg category of sheaves.
    \end{tabular}%
  }
\end{center}
\vspace{0.1cm}
\noindent Recall that the braid variety $X(\beta)$ admits a $\T$-action as described in \Cref{sssec:braid_variety_further_symmetries}, see also \Cref{lem:torusaction_coordinates}. Specifically, the tenet above can be formalized in the following result:

\begin{theorem}[Braid varieties as Legendrian invariants]\label{thm:braid_varieties_from_Legendrian}
Let $\beta\in\brn$ be a positive braid word, $\La_\beta\sse(\R^3,\xi_\std)$ its associated Legendrian link, and $X(\beta)$ its braid variety. Then:

\begin{enumerate}[label=$(\roman*)$]
    \item There exists an isomorphism of algebraic varieties over $\Z$
    \begin{equation}\label{eq:braid_variety_dga}
    X(\beta)\cong\operatorname{Spec }\operatorname{H}_0(\mathbb{A}^{ab}(\La_\beta)) 
    \end{equation}
    
    \noindent where the right hand side denotes the $0$th homology of the abelianization $\mathbb{A}^{ab}(\La_\beta)$ of the Legendrian dg algebra $\mathbb{A}(\La_\beta)$, and it is computed with one basepoint per strand of $\beta$.\\

    \item There exists an isomorphism of algebraic stacks
    \begin{equation}\label{eq:braid_variety_sheaf}
        [X(\beta)/\T]\cong \M^{\dagger}(\operatorname{Sh}_\La(\R^2))
    \end{equation}

    \noindent where $\M^{\dagger}(\operatorname{Sh}_\La(\R^2))$ denotes the component of $\GL_1$-local systems in the 1-truncation $t_1(\M(\La_\beta))$ of the derived moduli stack $\M(\La_\beta)$ of pseudoperfect objects of $\operatorname{Sh}_\La(\R^2)$.\\

    \item For a Legendrian knot $\La=\La_\beta$ and $\G=\GL_n$, there exists an isomorphism of algebraic stacks
    \begin{equation}\label{eq:braid_variety_sheaf2}
        X(\beta)\times B\bG_m\cong \M^{\dagger}_{\mbox{\tiny fr}}(\operatorname{Sh}_\La(\R^2))
    \end{equation}

    \noindent where $\M^{\dagger}_{\mbox{\tiny fr}}(\operatorname{Sh}_\La(\R^2))$ denotes the component of $\GL_1$-local systems in the 1-truncation $t_1(\M(\La_\beta))$ of the derived moduli stack $\M(\La_\beta)$ of pseudoperfect objects of the dg category $\operatorname{Sh}^{\mbox{\tiny fr}}_\La(\R^2)$ of sheaves commonly framed at the marked points of $\Lambda$.
\sbs
\end{enumerate}
\end{theorem}

\noindent Note that \eqref{eq:braid_variety_dga} in \Cref{thm:braid_varieties_from_Legendrian}.(i) can be equivalently written as an isomorphism of commutative $\Z$-algebras
\begin{equation}\label{eq:braid_variety_dga2}
\Z[X(\beta)]\cong \mbox{H}_0(\mathbb{A}^{ab}(\La_\beta))
\end{equation}
and that is the commutative algebra which is shown, via Demazure weaves, to be isomorphic to a cluster algebra in \cite[Theorem 1.1]{CGGLSS25}, cf.~Theorem \ref{thm:cluster_structures_braid_varieties}.\\

The isomorphism \eqref{eq:braid_variety_dga} follows from comparing the description of its right hand side in \cite[Prop.~5.2]{CasalsNg22} and the description of its left hand side via \Cref{lem:equations_braidvar}. The key point is that the product of braid matrices in \Cref{lem:equations_braidvar} counts the appropriate once-punctured holomorphic strips that enter the differential of the Legendrian contact dg-algebra in the case of $\La=\La_\beta$. Thus \cite[Prop.~5.2]{CasalsNg22} and  \Cref{lem:equations_braidvar} together imply \eqref{eq:braid_variety_dga}.

\begin{remark}\label{rmk:Perf_dga} A more intrinsic description of the right hand side of \eqref{eq:braid_variety_dga} arises by directly considering the dg category of perfect dg modules over the Legendrian contact dg algebra, which is homologically smooth, and taking the corresponding moduli of pseudoperfect objects, cf.~\cite{ToenVaquie07}. One of the simplest components of that moduli is of the form $X(\beta)\times B\bG_m$.\sbs
\end{remark}

For \eqref{eq:braid_variety_sheaf} and its variant \eqref{eq:braid_variety_sheaf2} in \Cref{thm:braid_varieties_from_Legendrian}, see \cite[Section 2]{CasalsLi24}. It is perhaps more enlightening to emphasize that sheaves with singular support on a Legendrian link $\La$ are a subcategory of the category of sheaves constructible with respect to the stratification of $\R^2_{x,z}$ given by a front of $\Lambda$. By the exodromy theorem, such constructible sheaves can be studied as (categorical) modules over an exit path category, cf.~\cite[Theorem A.9.3]{HA}, cf.~also \cite[Theorems 1.2 \& 6.10]{Treumann09_ExitPaths} and \cite[Chapters 0\&1]{barwick2020exodromy}. This translates the description of this type of constructible sheaves into a description of the representations of the poset of strata. By further imposing the singular support conditions that cut $\operatorname{Sh}_\La(\R^2)$ and its framed variant $\operatorname{Sh}^{\mbox{\tiny fr}}_\La(\R^2)$ inside of the category of all compact constructible sheaves, one obtains a rather explicit description of them in terms of flags. Such description allows us to deduce \eqref{eq:braid_variety_sheaf} and \eqref{eq:braid_variety_sheaf2}, see e.g.~\cite[Section 6]{STZ2017} and \cite[Section 4.1]{CasalsWeng22}.

\begin{remark}
The quasi-isomorphism type of the Legendrian contact dg algebra $\mathbb{A}^{ab}(\La_\beta)$ is independent of the braid word for $[\beta]$ and, since cyclic rotation is a Legendrian isotopy, it is also independent of the braid word $\beta$ up to cyclic rotation. Thus, \Cref{thm:braid_varieties_from_Legendrian} implies that for $\G=\SL_n$ the braid varieties $X(\beta)$ must be independent of the choice of braid word for $[\beta]$ and the isomorphism type of $X(\beta)$ must be preserved under the appropriate cyclic rotation of the crossings for the positive braid word $\beta$. For Type A, this gives a purely Floer theoretic proof of \Cref{lem:braid_invariance} and the cyclic symmetry mentioned in \Cref{sssec:braid_variety_further_symmetries}.\sbs
\end{remark}

At this stage, \Cref{thm:braid_varieties_from_Legendrian} illustrates that the braid varieties $X(\beta)$ appear as Legendrian invariants of $\La_\beta$. Nevertheless, it does not provide any geometric description of $X(\beta)$. Even if certain counts of pseudoholomorphic disks, which lead to $\mathbb{A}(\La_\beta)$, and studying certain sheaves, leading to $\operatorname{Sh}_\La(\R^2)$, are reasonably geometric operations, they are still far from providing geometric objects that relate directly to the Legendrian link $\La_\beta$. The next step is to explain how a Demazure weave  $\w:\beta\to\delta(\beta)$ can be used to construct an embedded exact Lagrangian surface $L_\w\sse(\D^4,\la_\std)$ whose boundary is precisely $\La_\beta$, directly bringing the diagrammatic calculus of Demazure weaves to contact and symplectic topology. This will start to give a precise geometric meaning to the points of $X(\beta)$: as we will explain that they parametrize certain types of Lagrangian surfaces bounding $\La_\beta$ endowed with a $\operatorname{GL}_1$-local system.

%%%%%%%%%%%%%%%%%%%%%%%%%%%%%%%%%%%%%%%%%%%%%%%%%%%%%%%%
%%%%%%%%%%%%%%%%%%%%%%%%%%%%%%%%%%%%%%%%%%%%%%%%%%%%%%%%
%%%%%%%%%%%%%%%%%%%%%%%%%%%%%%%%%%%%%%%%%%%%%%%%%%%%%%%%
\subsection{Lagrangian fillings of Legendrian links}\label{ssec:Lagrangian_fillings} Let $\La\sse(S^3,\xi_\std)$ be a Legendrian link. For instance, one arising as $\La=\La_\beta$ for a positive braid word $\beta\in\brn$ as in \Cref{def:Legendrian_from_braid}. The goal of this subsection is to discuss the notion of a Lagrangian filling of $\La$, which is a certain type of surface in the 4-ball bounding $\La$. A motivating reason is that, at least intuitively, the braid variety $X(\beta)$ can be understood as a candidate moduli space for the Lagrangian surfaces bounding $\La_\beta$, once endowed with certain local systems.\\

\noindent Legendrian links, as in \Cref{def:legendrian_link}, provide suitable boundary conditions for Lagrangian surfaces in order to study their contact and symplectic behavior.\footnote{E.g.~they are well-behaved boundary conditions to apply Floer-theoretic or microlocal sheaf-theoretic methods.} We focus on Lagrangian surfaces in the standard symplectic 4-ball $(\D^4,\la_\std)$, where
\begin{equation}\label{eq:Liouville_form}
2\la_\std:=(x_1 dy_1 - y_1 dx_1) + (x_2 dy_2 - y_2 dx_2)
\end{equation}
is a 1-form such that $d\la_\std=\omega_\std=dx_1dy_1+dx_2dy_2$, i.e.~a primitive for the standard symplectic 2-form $\omega_\std$ in $\D^4$. Here we have used Cartesian coordinates $(x_1,y_1,x_2,y_2)\in \D^4$. For us, the following definition of a Lagrangian surface bounding a Legendrian link suffices:

\begin{definition}[Lagrangian fillings]\label{def:Lagrangian_filling} Let $\La\sse (S^3,\xi_\std)$ be a Legendrian link. By definition, a Lagrangian filling $L$ of $\La$ is an embedded exact Lagrangian submanifold $L\sse(\D^4,\la_\std)$ such that:
\begin{enumerate}
    \item The restriction of $\la_\std$ to $L$ admits a smooth primitive $f:L\lr\R$, i.e.~$\la_\std|_L=df$,
    \item $L$ coincides with the cone of $\La$ outside a compact set $K\subset(\operatorname{int} \D^4,\la_\std)$,
    \item Any primitive of $\la|_L$ must be constant on $L\setminus (L\cap K)$.
\end{enumerate}
In this case we write $\dd L=\La$. Two Lagrangian fillings $L_0,L_1\sse(\D^4,\la_\std)$ are said to be exact Lagrangian isotopic if there exists a 1-parametric family $\{L_t\}_{t\in[0,1]}$ of exact Lagrangians in $(\D^4,\la)$ which is constant near the boundary $(S^3,\xi_\std)$.\sbs
\end{definition}

\Cref{def:Lagrangian_filling}.(1) is equivalent to the De Rham cohomology class $[\la_\std|_L]\in\operatorname{H}^1(L,\R)$ vanishing. This condition is known as being an exact Lagrangian and it is the main point of \Cref{def:Lagrangian_filling}, with $(2)$ and $(3)$ there to control asymptotic behavior. Thus the condition of being an exact Lagrangian can be analytically translated into a series of integrals of $\la_\std$ vanishing, e.g.~all the periods of $\la_\std$ along $L$ must vanish, among others.

\begin{remark}\label{rmk:Lagrangians_as_graphs}
Intuitively, exact Lagrangian surfaces in $(\D^4,\la_\std)$ are generalizations of the graph
\begin{equation}\label{eq:graph_df}
\operatorname{gr}(df):=\{y_1=\dd_{x_1}f,y_2=\dd_{x_2}f\}\sse \D^4
\end{equation}
of a smooth function $f:\R^2_{x_1,x_2}\lr\R$, where $f$ is generalized to be possibly multivalued (and non-graphical at certain points) and the Lagrangian surface $L$ itself no longer needs to have the topology of the plane $\R^2$.\sbs
\end{remark}

\noindent As said above, a current goal of low-dimensional contact and symplectic topology is the classification of Lagrangian fillings of a given Legendrian link $\La\sse(S^3,\xi_\std)$, up to exact Lagrangian isotopy. Note that two exact Lagrangians are exact Lagrangian isotopic if and only if they are Hamiltonian isotopic, see e.g.~\cite[Section 9.3]{McDuffSalamon98} and \cite[Theorem 3.6.7]{Oh15_SymplecticTopologFloerHomology}.\\

\noindent Here are a few starting facts on Lagrangian fillings:\\

\begin{enumerate}
    \item A given Legendrian link $\La\sse(S^3,\xi_\std)$ might or might not have Lagrangian fillings. This dichotomy is a Legendrian invariant: if $\La_0$ has a Lagrangian filling and there is a Legendrian isotopy from $\La_0$ to $\La_1$, then the trace of the Legendrian isotopy gives a Lagrangian concordance which concatenates with the filling of $\La_0$ to give a Lagrangian filling of $\La_1$.\\
    
    \noindent Since we focus on Legendrian links of the form $\La_\beta$, we emphasize that $\La_\beta$ always admits at least one Lagrangian filling. It can be shown that Legendrian knots admitting Lagrangian fillings must be quasi-positive, cf.~\cite{HaydenSabloff2015}, and thus, in a sense, focusing on the class of Legendrian links $\La_\beta$ is already quite illustrative among all Legendrian links that admit Lagrangian fillings.\\
    
    \item There are examples of pairs of Legendrian knots $\La_1,\La_2$ which have the same underlying smooth knot but $\La_1$ admits a Lagrangian filling and $\La_2$ does not, even if they are both max-tb. Thus this is a contact and symplectic property, not detected by smooth topology.\\

    \item The genus of a Lagrangian filling of a Legendrian knot $\La$ is determined by $\La$ itself. Specifically, the 4-ball genus $g_4(\La)$ of the smooth knot underlying $\La$ satisfies
                \begin{equation}\label{eq:filling_genus}
                    g(L)=\frac{\operatorname{tb}(\La)+1}{2}=g_4(\La)
                \end{equation}
    for any Lagrangian filling $L$, cf.~\cite{Chantraine10}. This is in contrast to the smooth case: if a smooth knot $K\sse S^3$ bounds an embedded genus $g$ smooth surface in the 4-ball $\D^4$, then it also bounds an embedded genus $g'$ smooth surface in the 4-ball $\D^4$ for any $g\leq g'$. Note that the braid variety $X(\beta)$ has dimension
    \begin{equation}\label{eq:dimension_braid_variety_genus}
        \dim_\kk X(\beta)=\textcolor{blue}{2g_4(\La)-1}+\textcolor{red}{n}\stackrel{\eqref{eq:filling_genus}}{=}-\textcolor{blue}{\chi(L)}+\textcolor{red}{n}
    \end{equation}
    where $n$ is the number of strands in $\beta$, and we supposed that $\delta(\beta)=w_0$ and $\La_\beta$ is a knot. The right hand sides of \eqref{eq:dimension_braid_variety_genus} are written so as to emphasize which part comes from the algebraic topology of a Lagrangian filling $L$, which is the term $-\chi(L)={2g_4(\La)-1}$ highlighted in blue, and which part comes from the marked points of $\La$, highlighted in red. In the case that $\beta$ contains $w_0$ as a suffix, then the number of mutable vertices in the cluster structure for $\Z[X(\beta)]$ from \Cref{sec:braid_variety_cluster} is $2g_4(\La)$, i.e.~two mutable vertices per genus of $L$, and the number of frozen vertices is $n-1$.\\

    \item The same Legendrian knot $\La\sse(S^3,\xi_\std)$ might bound several non-isotopic Lagrangian fillings, even infinitely many, cf.~\cite{CasalsGao22}. There are many instances where these infinitely many Lagrangian fillings are all {\it smoothly} isotopic to each other. Thus, as above, this is a contact and symplectic phenomenon, not detected by smooth low-dimensional topology.\\
\end{enumerate}

\noindent Both the construction and obstruction of embedded exact Lagrangian surfaces in symplectic 4-manifolds, including the study of Lagrangian fillings of Legendrian links, are active areas of research in low-dimensional contact and symplectic topology. While they inform each other and are naturally intertwined, the constructive techniques tend to be independent from the obstructing ones. The former often uses methods closer to geometric topology, whereas the latter is often focused on developing Lagrangian and Legendrian invariants, e.g.~via Floer theory, the microlocal theory of sheaves, or generating functions, to name some.\\

For specificity, let us consider the Legendrian contact dg algebra $\mathbb{A}(\La)$ as a Legendrian invariant of a Legendrian link $\La$, decorated with a set of basepoints $\tt$. An important property of a $\Z$-graded orientable Lagrangian filling $L$ of $\dd L=\La$ is that there is an induced morphism of dg algebras
\begin{equation}\label{eq:augmentation_map_filling}
\varepsilon_L:\mathbb{A}^{ab}(\La_\beta)\lr\Z[\hH_1(L,\tt;\Z)]
\end{equation}
\noindent where the target is understood as having trivial differential and concentrated in degree 0. See \cite[Theorem 6.8]{EtnyreNg22_Survey} for the construction of \eqref{eq:augmentation_map_filling}. By \Cref{thm:braid_varieties_from_Legendrian}.(1), the zeroth homology group of the commutative Legendrian dg algebra of $\La=\La_\beta$ is precisely the coordinate ring $\Z[X(\beta)]$ of the braid variety. In particular, we can specialize \eqref{eq:augmentation_map_filling} for a Lagrangian filling $L=L_\w$ associated to a Demazure weave $\w:\beta\to\delta(\beta)$, with $\La=\La_\beta$, and apply $\hH_0$ and the spectrum functor to the corresponding map $\varepsilon_{L_\w}$ as in \eqref{eq:augmentation_map_filling}. This leads to the following fact:

\begin{proposition}\label{prop:torus_from_filling} Let $\beta$ be a positive braid word, $(\La_\beta,\tt)\sse(\R^3,\xi_\std)$ its associated Legendrian link with basepoints $\tt$, one per strand, and $\mathbb{A}(\La_\beta,\tt)$ its Legendrian dg algebra over a field $k$. Suppose that $L\sse(\D^4,\la_\std)$ is an orientable embedded exact Lagrangian filling of $\La_\beta$ and set $d=\mbox{rk}(\hH_1(L,\tt;\Z))$. Then there exists an open inclusion
    \begin{equation}\label{eq:torus_from_filling}
        \Spec(\hH_0(\varepsilon_{L})):(\bG_m)^d\hookrightarrow X(\beta)
    \end{equation}
    which is invariant under compactly supported Hamiltonian isotopies of $(\D^4,\la_\std)$.\sbs
\end{proposition}

In particular, \Cref{prop:torus_from_filling} implies that the image of \eqref{eq:torus_from_filling} is an invariant of the embedded exact Lagrangian filling $L$, up to compactly supported exact Lagrangian isotopies. That is, if $L_1$ and $L_2$ are two Lagrangian fillings of $\La_\beta$ and the images of $\Spec(\hH_0(\varepsilon_{L_1}))$ and $\Spec(\hH_0(\varepsilon_{L_2}))$ do not coincide as open subsets of $X(\beta)$, then $L_1$ and $L_2$ are not compactly supported Hamiltonian isotopic. See also \cite{LiCobordism1,li2024lagrangiancobordismfunctormicrolocal} for a sheaf-theoretic approach.\\

\begin{example}\label{ex:Lagrangian_fillings} (1) Let $\beta=\s_1^{n+1}\in\mbox{Br}_2^+$ be a 2-stranded positive braid word. Then $\La_\beta$ is the unique max-tb Legendrian representative of the $(2,n)$-torus link. It is known to admit at least a Catalan number $$C_n:=\frac{1}{n+1}{2n \choose n}$$
of distinct Lagrangian fillings, up to compactly supported Hamiltonian isotopy. This was first proven via Floer theory in \cite{Pan} using the maps from \Cref{prop:torus_from_filling}. For $n=1,2$ we also know these are all Lagrangian fillings, cf.~ \cite{EliashbergPolterovich96,thomson2026legendrianhopflinkexactly}. It is conjectured that this is also the case for all $n\in\N$. More generally, see e.g.~\cite[Conjecture 5.1]{CasalsLagSkel}. Since the cluster algebra $\Z[X(\beta)]$ is of cluster type $A_{n-1}$, this can be seen as a first hint towards a potential relation between Lagrangian fillings and clusters.\\

\noindent (2) Consider the braid word $\beta=(\s_1\s_2)^8\s_1$, then $\La_\beta$ is the max-tb Legendrian $(3,7)$-torus knot. There is a natural Legendrian isotopy $\{\La_t\}$ that cyclically rotates the leftmost crossing of $\beta$ to its rightmost vertex by passing through the upper strands in \Cref{fig:FrontsBeta}. In particular, applying this Legendrian isotopy leads to a Lagrangian concordance $L_{c}$, via the trace, from $\La_\beta$ to itself. In this case, there exists a Lagrangian filling $L$ of $\La_\beta$ such that the concatenation $L_k$ of $L$ with $k$ stacked copies of $L_c$ is a Lagrangian filling of $\La_\beta$. It can be shown that $\{L_k\}_{k\in\N}$ are an infinite family of pairwise non-isotopic Lagrangian fillings by studying the maps from \Cref{prop:torus_from_filling}, see e.g.~\cite{CasalsNg22} and \cite{Casals_NewtonPolytopes}.\\

\noindent This construction of building Lagrangian fillings by stacking Lagrangian cylinders on top of them, which is the key idea behind \cite{CasalsGao22}, is rather general and typically produces infinitely many Lagrangian fillings. This is tightly related to the fact that the (square of the) Donaldson-Thomas transformation for the associated $X(\beta)$ has infinite order, cf.~\cite{GSW2024}.\\

\noindent (3) For an isolated plane curve singularity $f:\C^2\lr\C$, any real Morsification $\tilde{f}$ of $f$ provides a Lagrangian filling in the Milnor ball $(\D^4,\la_\std)$ with Legendrian boundary given by the unique max-tb representative of the link of the singularity, cf.~\cite[Section 2]{CasalsLagSkel}. Smoothly, this Lagrangian filling is isotopic to the Milnor fiber. Already in these instances, there are typically many distinct Lagrangian surfaces whose symplectic push-offs give symplectically isotopic symplectic Milnor fibers. The symplectic Milnor fiber is expected\footnote{This is a relative version of the symplectic isotopy problem.} to be the unique symplectic surface bounding the transverse link of the singularity. In contrast, Lagrangian surfaces bounding its unique Legendrian representatives behave rather differently. In a vague informal sense, studying Lagrangian surfaces versus symplectic surfaces is akin to studying totally real surfaces versus complex surfaces.\sbs
\end{example}

%%%%%%%%%%%%%%%%%%%%%%%%%%%%%%%%%%%%%%%%%%%%%%%%%%%%%%%%
%%%%%%%%%%%%%%%%%%%%%%%%%%%%%%%%%%%%%%%%%%%%%%%%%%%%%%%%
%%%%%%%%%%%%%%%%%%%%%%%%%%%%%%%%%%%%%%%%%%%%%%%%%%%%%%%%
\subsection{A construction of Lagrangian fillings from algebraic geometry}\label{ssec:Lagrangian_filling_from_AG} Let us present a source of Lagrangian fillings in $(\D^4,\la_\std)$, which historically motivated the introduction of weaves in \cite{CasalsZaslow}, see also \cite[Section 2]{TreuZas16}. To ease notation, we will indistinguishably work with $\D^4$ and its completion $(\R^4,\la_\std)\cong (T^*\R^2,\la_\std)$, and often identify projections to two of the Cartesian coordinates as maps to either $\D^2$ or $\R^2$. It is provided by the following lemma:

\begin{lemma}\label{lem:complex_curve_real_Lagrangian}
Let $C:=\{F(z,w)=0\}\sse\C^2_{z,w}$ be a smooth complex curve. Then:
\begin{enumerate}[label=$(\roman*)$]
    \item $C$ is a Lagrangian surface in $(\R^4,\omega_\std)$ under the identification
$$z=x_1+ix_2,\qquad w=y_1-iy_2.$$
That is, $\omega_\std=dx_1dy_1+dx_2dy_2=\Re(dz\wedge dw)$ vanishes along $C$.\\

\item $C$ is an exact Lagrangian surface in $(\R^4,\la_\std)$ if and only if 
   \begin{equation}\label{eq:lambda_periods}
   \Re\left(\oint_\gamma zdw-wdz\right)=0,\quad \forall\gamma\in \operatorname{H}_1(C,\R)
   \end{equation}
\end{enumerate}
\sbs
\end{lemma}

\noindent For \Cref{lem:complex_curve_real_Lagrangian}.(2), just note that $2\la_\std=\Re\left(zdw-wdz\right)$ where $\la_\std$ as in \eqref{eq:Liouville_form}. It can be helpful to observe that \eqref{eq:lambda_periods} is equivalent to
\begin{equation}\label{eq:lambda_periods2}
   \Re\left(\oint_\gamma zdw\right)=0,\quad \forall\gamma\in \operatorname{H}_1(C,\R)
   \end{equation}

The existence of smooth algebraic curves in $\C^2$ with purely imaginary periods for $zdw$ is studied in \cite{eynard2024existenceboutrouxcurvesgfunctions}. The intersection of such Lagrangian surfaces, as constructed in \Cref{lem:complex_curve_real_Lagrangian}, with the unit sphere $S^3=\{|z|^2+|w|^2=1\}\sse\C^2_{z,w}$ is typically not a Legendrian link, but this can be fixed by an appropriate perturbation, cf.~\cite[Section 5.1]{Chantraine10}, as these Lagrangians are conically asymptotical to a Legendrian link away from a compact subset. The Lagrangian surfaces in $(\D^4,\la_\std)$ that arise from \Cref{lem:complex_curve_real_Lagrangian} are special, cf.~\cite[Section 2.1]{Joyce_LecturesSLags}. As above, we often interchange the inclusions of a Lagrangian surface $L\sse(\R^4,\la_\std)$ and its restriction $L\sse(\D^4,\la_\std)$ to a large enough open ball.\\

\begin{example}\label{ex:special_fillings} (1) Let $f:\C^2\lr\C$ be an irreducible isolated plane curve singularity. By perturbing $f$ to $\tilde{f}$, e.g.~by morsifying, $\{\tilde{f}=0\}\sse\C^2$ is a smooth complex curve (e.g.~a Milnor fiber) and, if the parameters of the perturbation are chosen so that their $zdw$-periods are purely imaginary, then \Cref{lem:complex_curve_real_Lagrangian} provides a Lagrangian filling $L_{\tilde{f}}$ of the max-tb Legendrian representative of the link of the singularity $f$. See also \Cref{ex:Lagrangian_fillings}.(3) and \cite{CasalsLagSkel}.\\

\noindent (2) Here is an explicit instance of (1). Consider a simple cusp singularity $\{f(z,w)=0\}\sse\C^2$, where $f(z,w)=z^2-w^3$. Consider a deformation of the form
\begin{equation}\label{eq:A2_deformation}
\tilde{f}(z,w):=w^2-z^3+\textcolor{blue}{a}z+\textcolor{blue}{b},\quad\mbox{with } 
4\textcolor{blue}{a}^3-27\textcolor{blue}{b}^2\neq0.\end{equation}
\noindent The inequality for \textcolor{blue}{a}, \textcolor{blue}{b} in \eqref{eq:A2_deformation} is to ensure that $\{\tilde{f}(z,w)=0\}$ is smooth. In order to impose the condition of purely imaginary $wdz$-periods, we can use the Weierstrass elliptic function $\varrho$ to uniformize an elliptic curve, so that we can write
$$z=\varrho_{\omega_1,\omega_2}(u),\qquad 2w=\varrho'_{\omega_1,\omega_2}(u),\qquad $$
for some $\omega_1,\omega_2\in\C$ and complex parameter $u\in\C$, so the lattice $\La$ for $\{\tilde{f}(z,w)=0\}\cong(\C/\La)\setminus\{\operatorname{pt}\}$ is $2\omega_1\Z\oplus 2\omega_2\Z$. Let us assume $\omega_1\in\R_{>0}$, let $\gamma_i$ be the homology 1-cycle associated to the $\omega_i\Z$ direction, $i\in[1,2]$, and write $\tau:=\omega_2/\omega_1\in\mathbb{H}$. Then the $wdz$-periods of $\{\tilde{f}(z,w)=0\}$ are generated by

$$\oint_{\gamma_1} wdz=\frac{\pi^6}{360\omega_1^5}(E_2(\tau)E_4(\tau)-E_6(\tau)),\qquad \oint_{\gamma_2} wdz=\tau\cdot\left(\oint_{\gamma_1} wdz\right)-\frac{\pi^5E_4(\tau)}{60\omega_1^5}i$$

\noindent where we used the normalized Eisenstein series
$$E_2(\tau):=1-24\sum_{n\geq1}\frac{nq^n}{1-q^n},\quad 
 E_4(\tau)
 :=1+240\sum_{n\geq1}\frac{n^3q^n}{1-q^n},\quad E_6(\tau)
 :=1-504\sum_{n\geq1}\frac{n^5q^n}{1-q^n},$$
after abbreviating $q:=e^{2\pi i\tau}$, which are real when restricting to the line with real part $0.5$ for $\tau$, even if not necessarily on all of $\mathbb{H}$. Indeed, the formula for the periods above follows from the identifications
\begin{equation}\label{eq:a_b_tau_value_A2case}
a = \frac{\pi^4}{48\omega_1^4} E_4(\tau), \qquad b = \frac{\pi^6}{864\omega_1^6} E_6(\tau)
\end{equation}
under the Weierstrass uniformization. Therefore, if we choose $\tau=\tau_\ast\in\mathbb{H}$ such that
\begin{equation}\label{eq:tau_value_A2case}
E_2(\tau_\ast)E_4(\tau_\ast)-E_6(\tau_\ast)=0
\end{equation}
the period along $\gamma_1$ vanishes and that along $\gamma_2$ becomes purely imaginary. Recall that we have the restriction $\Re(\tau)=0.5$. Since the Eisenstein series are real, Bolzano's theorem shows that there must be at least one solution $\tau_\ast$ to \eqref{eq:tau_value_A2case}. Specifically, \eqref{eq:tau_value_A2case} applied to $\tau_+=0.5+i\cdot(0.4)$ is positive, and applied to $\tau_-=0.5+i\cdot(0.41)$ is negative. For specificity, numerical approximation gives
$$\tau_\ast\approx\frac{1}{2}+i\cdot(0.40868185996886)$$
close to a solution of \eqref{eq:tau_value_A2case}. Thus the deformation
\begin{equation}\label{eq:Boutroux_A2_case}
\tilde{f}(z,w):=w^2-z^3+\textcolor{blue}{\left(\frac{\pi^4}{48\omega_1^4} E_4(\tau_\ast)\right)}\cdot z+\textcolor{blue}{\frac{\pi^6}{864\omega_1^6} E_6(\tau_\ast)}
\end{equation}
is such that $\{\tilde{f}(z,w)=0\}$ is an explicit smoothly embedded exact Lagrangian filling of the max-tb Legendrian trefoil knot $\La=\La_\beta$, $\beta=\s_1^4$. The deformation in \eqref{eq:Boutroux_A2_case} admits the following variant:

\begin{equation}\label{eq:Boutroux_A2_case_5symmetry}
\tilde{f}_k(z,w):=w^2-z^3+\textcolor{purple}{\zeta_5^{-k}}\textcolor{blue}{\left(\frac{\pi^4}{48\omega_1^4} E_4(\tau_\ast)\right)}\cdot z+\textcolor{purple}{\zeta_5^{k}}\textcolor{blue}{\frac{\pi^6}{864\omega_1^6} E_6(\tau_\ast)},\quad k\in[0,4],
\end{equation}
where $\zeta_5:=e^{2\pi i/5}$ is a fifth root of unity. Then the smooth algebraic curves $\{\tilde{f}_k(z,w)=0\}$, $k\in[0,4]$, provide five embedded exact Lagrangian fillings of the max-tb Legendrian trefoil knot $\La_{\s_1^4}$, which can be shown to be pairwise nonisotopic. This aligns with \Cref{ex:Lagrangian_fillings}.(1), as the corresponding Catalan number is $C_3=5$.\\

\noindent (3) As illustrated in (2), finding Lagrangian fillings by using \Cref{lem:complex_curve_real_Lagrangian} requires an understanding of periods, which often involves a combination of complex analysis and number theory. This leads to interesting explicit examples, often involving special functions. For instance, consider $f(z,w)=z^3-w^4$ and let us try to find an embedded exact Lagrangian filling of the max-tb Legendrian $(3,4)$-torus knot via \Cref{lem:complex_curve_real_Lagrangian}.(2). For that we deform via
\begin{equation}\label{eq:e6-fb}
\tilde{f}(z,w)=w^3-z^4+w+\textcolor{blue}{\alpha}z^2-\frac{\textcolor{blue}{\alpha}^2}{8},
\end{equation}
for some real parameter $\textcolor{blue}{\alpha}\geq0$. The goal is to choose $\textcolor{blue}{\alpha}$ so that the $wdz$-periods of $\{\tilde{f}(z,w)=0\}$ are purely imaginary. For that we use the Gamma function $\Gamma$ and Gauss' hypergeometric series ${}_2F_1$, with continuation from its Taylor series
along the negative real axis.\footnote{\noindent A conceptual reason for the appearance of hypergeometric series ${}_2F_1$ is that the 2nd order differential equation that it solves is closely related to the Picard-Fuchs equation that appears for the periods of $\{\tilde{f}(z,w)=0\}$.} A more elaborate version of the computation in (2), as now $\{\tilde{f}(z,w)=0\}$ has genus 3, shows that the required $wdz$-periods are purely imaginary if $\alpha$ satisfies the equation
\begin{equation}\label{eq:alpha_parameter_for_34}
\frac{\Gamma(5/8)\Gamma(5/4)}{\Gamma(15/8)}\cdot {}_2F_1\left(-\tfrac7{24},\tfrac1{24};\tfrac34;
                         -\tfrac{27\textcolor{blue}{\alpha}^4}{256}\right) = \frac{\textcolor{blue}{\alpha}\Gamma(3/8)\Gamma(3/4)}{4\Gamma(9/8)}\cdot{}_2F_1\left(-\tfrac1{24},\tfrac7{24};\tfrac54;
                         -\tfrac{27\textcolor{blue}{\alpha}^4}{256}\right)
\end{equation}
\noindent Let us refer to the left- and right-hand sides of \eqref{eq:alpha_parameter_for_34} as $L(\alpha)$ and $R(\alpha)$ respectively. Then both sides of \eqref{eq:alpha_parameter_for_34} asymptotically behave like $L(\alpha)\sim R(\alpha)\sim \alpha^{7/6}$, up to leading constants, and one directly computes that
$$\displaystyle\lim_{\alpha\to\infty}\frac{R(\alpha)}{L(\alpha)}>1\quad\mbox { and }\quad \frac{R(0)}{L(0)}=0.$$
Since the quotient is strictly monotone as a function of $\alpha$, there exists a unique positive value $\alpha_\ast$ solving \eqref{eq:alpha_parameter_for_34}. For specificity, numerical approximation gives
$$\alpha_\ast\approx1.779407345096091\ldots.$$
Thus, via \Cref{lem:complex_curve_real_Lagrangian}, the algebraic curve
$$\{w^3-z^4+w+\textcolor{blue}{\alpha_\ast}z^2-\frac{\textcolor{blue}{\alpha_\ast}^2}{8}=0\}\sse\C^2$$
provides a smoothly embedded exact Lagrangian filling in $(\D^4,\la_\std)$ of the max-tb Legendrian $(3,4)$-torus knot $\La_{\beta}$, where $\beta=(\s_1\s_2)^5\s_1$. It is known that there are at least 833 exact Lagrangian fillings of $\La_\beta$, cf.~\cite{CasalsGao22}. It is conjectured that such $\La_{\beta}$ admits precisely 833 exact Lagrangian fillings, cf.~\cite[Section 5]{CasalsLagSkel}, which is the number of cluster seeds of cluster type $E_6$.\sbs
\end{example}

\begin{center}
	\begin{figure}[h!]
		\centering
		\includegraphics[width=\textwidth]{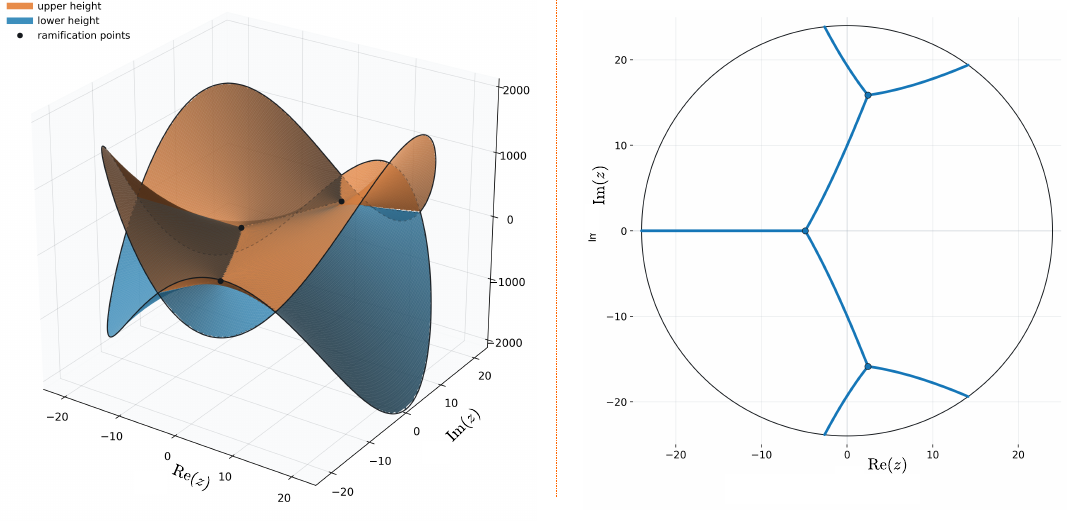}
		\caption{(Left) The graph of the map $z\mapsto (z,\Re\left(\int w(z)dz\right))$, where $\tilde{f}(z,w(z))=0$ and $\tilde{f}(z,w)$ is the deformation of $w^2-z^3$ defined in \eqref{eq:Boutroux_A2_case}, as in \Cref{ex:special_fillings}.(2). (Right) The locus $\w\sse\R^2_{x_1,x_2}$ of points where $\Re\left(\int w(z)dz\right)$ assigns one value, instead of the generic two values that are assigned to the points $(x_1,x_2)$ on the complement of $\w$. Here we identify $\R^2\cong\C_z$ via $z=x_1+ix_2$.}
		\label{fig:Numerical_Trefoil_2weave}
	\end{figure}
\end{center}

\begin{remark}\label{rmk:Boutroux_via_period} As pointed out to me by Y.J.~Nho, for singularities of the form $z^2=w^{n+1}$, including \Cref{ex:special_fillings}.(2) above, one can argue the existence of a Morsification with purely imaginary $wdz$-periods by \cite[Theorem 1]{Varchenko1988}, though it does not provide explicit coefficients.\sbs
\end{remark}

\noindent As hopefully illustrated by \Cref{ex:special_fillings}, it is a rather subtle matter to use \Cref{lem:complex_curve_real_Lagrangian} to construct and easily manipulate embedded exact Lagrangian fillings coming from algebraic curves. From a geometric and topological viewpoint, there are two downsides:\\

\begin{enumerate}
    \item It can be challenging to geometrically visualize these surfaces in 4-dimensions.\\

    \item The rigidity of complex analysis and algebraic geometry makes it challenging to manipulate by hand these Lagrangian fillings. Specifically, in the spirit of geometric topology, with which contact and symplectic topology shares many features, it would be ideal to be able to construct and manipulate these Lagrangian fillings via local models and cut-and-paste arguments.\\
\end{enumerate}

\noindent We will now explain a solution to (1), which leads to considering the notion of a weave. Once weaves are defined, they immediately resolve (2) above. The idea to address (1) above is \Cref{rmk:Lagrangians_as_graphs} above. Namely, we insist on thinking of exact Lagrangians as graphs $\gr(df)\sse\R^4_{x_1,x_2,y_1,y_2}$, where $f:\R^2_{x_1,x_2}\lr\R$ is possibly multivalued. Such a primitive $f$ can be extracted from the exact Lagrangian itself by restricting $\la_\std$ (or a fixed ambient primitive) and integrating. In order to visualize exact Lagrangians, we plot the graph $\gr(f)\sse\R^3$, which is a singular surface in 3-dimensional space: this at least we can directly visualize with our eyes.

\begin{example}\label{ex:first_weave_plots} (1) Consider a smooth algebraic curve $\{(z,w)\in\C^2:P(z,w)=0\}\sse\C^2$. Then we can express $w=w(z)$ as a multivalued function of $z$ by requiring $P(z,w(z))=0$. Suppose the curve has purely imaginary $zdw$-periods. Then the associated exact Lagrangian surface in $(\D^4,\la_\std)$ can be understood as $\gr(df)$ where the multivalued function $f$ is
\begin{equation}\label{eq:graph_Re_integral}
f:\R^2_{x_1,x_2}\lr\R,\qquad f(x_1,x_2):=\Re\left(\int w(z)dz\right)
\end{equation}
where we write $z:=x_1+ix_2$ as above. Therefore, we can directly plot the graph of $z\mapsto (z,\Re(\int w(z)dz))$, which is a subset of $\R^3$.\footnote{Similarly, we could have instead expressed $z=z(w)$ as a multivalued function of $w$ by requiring $P(z(w),w)=0$ and plot the graph of $w\mapsto (w,\Re(\int z(w)dw))$.}\\

\noindent (2) Let us consider the $P(z,w):=\tilde{f}(z,w)$ as in \eqref{eq:Boutroux_A2_case} from \Cref{ex:special_fillings}.(2), which gives a Lagrangian filling of the Legendrian trefoil $\La_{\beta}$, with $\beta=\s_1^4$. In this case the graph of the multivalued function \eqref{eq:graph_Re_integral} is depicted in \Cref{fig:Numerical_Trefoil_2weave} (left). It can be described as follows:
\begin{enumerate}[label=$(\roman*)$]
    \item $f$ generically assigns 2 values per choice of $z=(x_1,x_2)\in\R^2$.\\
    \item There is a subset $\w\sse\R^2$ of points $(x_1,x_2)\in\R^2$ where $f$ assigns one value. This subset is depicted in \Cref{fig:Numerical_Trefoil_2weave} (right), as a subset of the $z$-plane, where $z=x_1+ix_2$.\\
    \item Within $\w$, there are 3 points where $\w$ is singular, each giving a trivalent vertex in $\w$. These are the three trivalent vertices depicted in \Cref{fig:Numerical_Trefoil_2weave} (right).
\end{enumerate}
\noindent The choice of notation $\w$ is no mistake. In summary, $\gr(f)\sse\R^3$ is a singular surface obtained by weaving two sheets in $\R^3$. The locus $\w$ records where these two sheets intersect. Generically they intersect at a codimension-one locus of immersed points, and within that locus there are three points with a more elaborate singularity. The name of this singularity is $D_4^-$, see \Cref{sssec:trivalent_model} below or \cite[Section 2.2.3]{CasalsWeng22}. The graphs for the variants \eqref{eq:Boutroux_A2_case_5symmetry} are obtained by rotating the plot by an angle of $2\pi/5$ in the base $z$-plane.\\

\begin{center}
	\begin{figure}[h!]
		\centering
		\includegraphics[width=\textwidth]{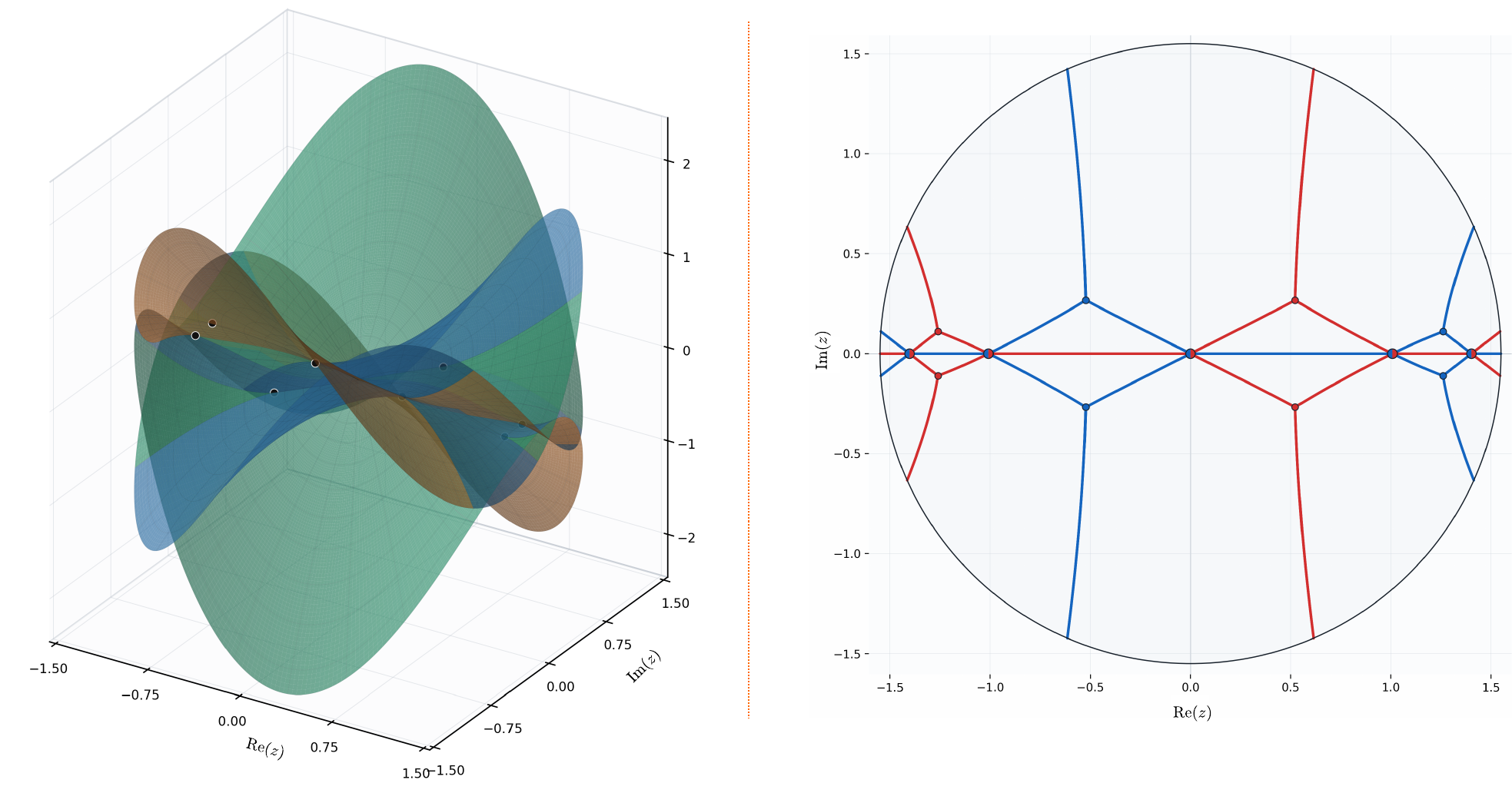}
		\caption{(Left) The graph of the map $z\mapsto (z,\Re\left(\int w(z)dz\right))$, where $\tilde{f}(z,w(z))=0$ and the deformation $\tilde{f}(z,w)$ of $w^3-z^4$ is defined as in \eqref{eq:e6-fb}, cf.~\Cref{ex:special_fillings}.(3). (Right) The locus $\w\sse\R^2_{x_1,x_2}$ of points where $\Re\left(\int w(z)dz\right)$ assigns strictly fewer than 3 values. At the hexavalent vertices the function is univalued, at the rest of $\w$ it is bi-valued, and on the complement of $\w$ it is tri-valued. We have depicted the points in $\w$ where sheets one and two intersect in blue, while the points in $\w$ where sheets two and three intersect are drawn in red.}
		\label{fig:Numerical_34_3weave}
	\end{figure}
\end{center}

\noindent (3) Let us consider the $P(z,w):=\tilde{f}(z,w)$ as in \eqref{eq:e6-fb} from \Cref{ex:special_fillings}.(3), which gives a Lagrangian filling of the Legendrian $(3,4)$-torus knot $\La_{\beta}$, with $\beta=(\s_1\s_2)^5\s_1$. The graph of the corresponding multivalued function, as in \eqref{eq:graph_Re_integral}, is depicted in \Cref{fig:Numerical_34_3weave} (left). It is described as follows:
\begin{enumerate}[label=$(\roman*)$]
    \item $f$ generically assigns 3 values per choice of $z=(x_1,x_2)\in\R^2$.\\
    \item There is a subset $\w\sse\R^2$ of points $(x_1,x_2)\in\R^2$ where $f$ assigns one or two values (instead of three). This subset $\w$ is depicted in \Cref{fig:Numerical_34_3weave} (right) as a subset of the $(x_1,x_2)$-plane, equivalently the $z$-plane with $z=x_1+ix_2$.\\

    Above a point in $\R^2\setminus\w$, we order the sheets from top to bottom. At a smooth point of $\w$ two sheets intersect: either the bottom two sheets or the top two sheets. For clarity, we have drawn the points in $\w$ where the bottom two sheets intersect in blue, and the points in $\w$ where the top two sheets intersect in red.\\
    
    \item Within $\w$, there are 13 points where the locus $\w$ is singular. Of these points, 8 give trivalent vertices in $\w$, and 5 of them give hexavalent vertices in $\w$. These are illustrated in \Cref{fig:Numerical_34_3weave} (right) as well.\\
\end{enumerate}
\noindent In summary, $\gr(f)\sse\R^3$ is a singular surface obtained by weaving three sheets in $\R^3$. The locus $\w$ records where these three sheets intersect. As before, these sheets generically intersect at a codimension-one locus of immersed points, and within that locus there are points with two types of singularities. The first type is a $D_4^-$-singularity, as above, and gives the trivalent vertices of $\w$. In these $D_4^-$-singularities only two consecutive sheets interact, the third sheet is simply entirely above or below them. The second type is an $A_1^3$-singularity, where the three sheets together actually meet at a point: it is exactly the type of singularity formed by union of the coordinate planes in $\R^3$, see \Cref{sssec:hexavalent_local_model} below.\sbs
\end{example}

\noindent \Cref{ex:first_weave_plots} should give a strong indication of the geometric relation between weaves and Lagrangian fillings. Indeed, \Cref{fig:Numerical_Trefoil_2weave} (right) is an $\SL_2$ weave, which is actually a Demazure weave $\w:\beta\to\delta(\beta)$ for $\beta=\s_1^4$. \Cref{fig:Numerical_34_3weave} (right) depicts an $\SL_3$ weave, which can be turned into a Demazure weave $\w:\beta\to\delta(\beta)$ for $\beta=(\s_1\s_2)^5\s_1$. The next subsection will add the necessary definitions and details.\\

\begin{remark}
Informally said, \Cref{fig:Numerical_Trefoil_2weave} (right) and \Cref{fig:Numerical_34_3weave} (right) each record the graph on their left by looking at it straight from above and only drawing where the sheets intersect, using the color to indicate which pair of consecutive sheets is intersecting.\sbs
\end{remark}

For now, we note that one can obtain a positive braid word $\eta_\w$, up to cyclic reordering, by reading the edges of $\w$ that intersect a large enough circle in the $z$-plane and recording which sheets intersect. Specifically, if we number the sheets bottom to top in increasing order of height, the intersection of sheets $i$ and $i+1$ is labeled by $\s_i$. This produces a positive braid word $\eta_\w$, up to cyclic reordering. For instance, if we start at $(20,0)$ in \Cref{fig:Numerical_Trefoil_2weave} (right) and go around clockwise we obtain $\eta_\w=\s_1^5$. Note that this is exactly $\beta\delta(\beta)=\s_1^5$ for $\beta=\s_1^4$, as in \Cref{ex:special_fillings}.(2). Similarly, if we start at $(0,2)$ in \Cref{fig:Numerical_34_3weave} (right) and go around clockwise we obtain $\eta_\w=(\s_1\s_2)^7=\beta\delta(\beta)$ where $\beta=(\s_1\s_2)^5\s_1$, as in \Cref{ex:special_fillings}.(3). Note that the same procedure applied to a Demazure weave $\w:\beta\to\delta(\beta)$ indeed always gives $\eta_\w=\beta\delta(\beta)$.

\begin{remark}\label{rmk:topology_Lagrangian_filling}
The topology of these embedded exact Lagrangian fillings $L\sse(\D^4,\la_\std)$ that we are visualizing by plotting the graph of $\Re\left(\int w(z)dz\right)$ can be readily deduced as follows. The restriction to $L$ of the projection $\C^2_{z,w}\lr\C_z$ onto the $z$-plane is a simple branched cover, of order the number of sheets above a generic point of $z$. The trivalent vertices of $\w$ are precisely the branch points of such projection. By Riemann-Hurwitz, all that remains to determine the diffeomorphism type of $L$ is to compute the number of boundary components: this is given by the number of cycles in which the braid word $\eta_\w$ decomposes when projected to the permutation group. Note that this is independent of the cyclic representative. Thus one can readily obtain the genus and number of boundary components from the set $\w$ if we label its edges by which sheets are joined right above them.\sbs
\end{remark}

\begin{remark}\label{rmk:front_of_lift}
For experts, \Cref{fig:Numerical_Trefoil_2weave} (left) and \Cref{fig:Numerical_34_3weave} (left) are fronts for Legendrian lifts of the corresponding exact Lagrangian surfaces. A Legendrian front for a Legendrian surface generically has singularities in real codimension one, e.g.~edges of cusps, and swallowtail singularities. (Here $A_1^2$ is not technically counted as a singularity, as it is an immersed locus.) The conceptual reason none of these fronts coming from algebraic curves have these cusp edges or swallowtails is that they are the real part of the complex front for a complex Legendrian curve. Since the front singularities of a complex Legendrian curve must be in complex codimension one, they must be points. Therefore, their real parts will have fronts with singularities in real codimension two. By the classification of Legendrian front singularities for surfaces, cf.~\cite[Chapter 21]{ArnoldGuseinZadeVarchenko85_Vol1}, they must be $A_1^3$ and $D_4^-$ singularities.\sbs
\end{remark}

%%%%%%%%%%%%%%%%%%%%%%%%%%%%%%%%%%%%%%%%%%%%%%%%%%%%%%%%
%%%%%%%%%%%%%%%%%%%%%%%%%%%%%%%%%%%%%%%%%%%%%%%%%%%%%%%%
%%%%%%%%%%%%%%%%%%%%%%%%%%%%%%%%%%%%%%%%%%%%%%%%%%%%%%%%
\section{Demazure weaves as Lagrangian fillings}\label{sec:weaves_Lagrangian_fillings}

The object of this section is to explain how to construct an exact embedded Lagrangian filling $L_\w\sse(\D^4,\la_\std)$ of $\La_\beta\sse(S^3,\xi_\std)$ from a Demazure weave. The local models are introduced in \Cref{ssec:Lagrangian_filling_from_weave_local_models}, motivated by \Cref{ssec:Lagrangian_filling_from_AG}. The global construction is explained in \Cref{ssec:Lagrangian_filling_from_weave}. The smooth topology of such embedded surfaces is discussed in \Cref{ssec:smooth_topology_Lagrangian_filling_from_weave}. Finally, \Cref{ssec:Lagrangian_coboridsms} provides a refinement of the notion of a Lagrangian filling $L_\w$ from a weave $\wb$ to the notion of a Lagrangian cobordism $L_\w:\La_{\beta_1}\prec\La_{\beta_2}$ from a weave $\w:\beta_2\to\beta_1$, where $\beta_1$ is not necessarily yet equal to $\delta(\beta_2)$.

%%%%%%%%%%%%%%%%%%%%%%%%%%%%%%%%%%%%%%%%%%%%%%%%%%%%%%%%
%%%%%%%%%%%%%%%%%%%%%%%%%%%%%%%%%%%%%%%%%%%%%%%%%%%%%%%%
%%%%%%%%%%%%%%%%%%%%%%%%%%%%%%%%%%%%%%%%%%%%%%%%%%%%%%%%
\subsection{Local models to build a Lagrangian from a weave}\label{ssec:Lagrangian_filling_from_weave_local_models} Let $\beta\in\brn$ be a positive braid word and $\La_\beta$ its associated Legendrian link, cf.~\Cref{def:Legendrian_from_braid}. Consider a Demazure weave $\w:\beta\to\delta(\beta)$, as in \Cref{sec:diagrammatic_weaves}. The goal of this subsection is to construct an embedded exact Lagrangian surface $L_\w\sse(\D^4,\la_\std)$ from $\w$ with the property that $\dd L_\w=\La_\beta$. The equality $\dd L_\w=\La_\beta$ is understood at the boundary $(S^3,\xi_\std)\cong (\dd \D^4,\ker\{\la_\std|_{\dd \D^4}\})$, which is a contact 3-manifold, cf.~\Cref{ssec:Lagrangian_fillings}. The motivating examples that lead to this construction are explained in \Cref{ssec:Lagrangian_filling_from_AG} above.\\

\noindent The embedded exact Lagrangian surface $L_\w\sse(\D^4,\la_\std)$ is constructed by gluing the following four local models, corresponding to the different parts of a weave diagram:\\

\begin{center}
	\begin{figure}[h!]
		\centering\includegraphics[width=\textwidth]{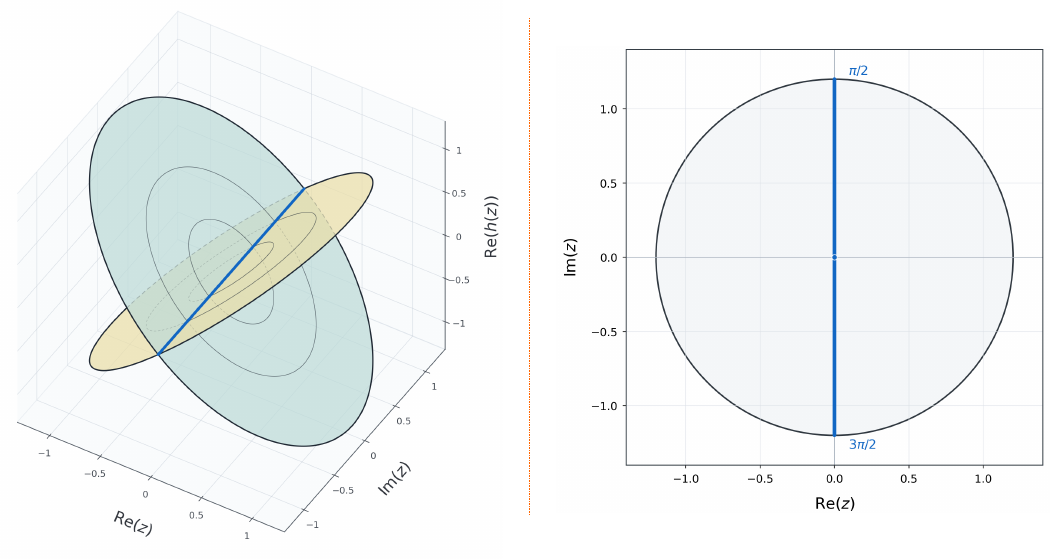}
		\caption{The local model for a weave edge, as discussed in \Cref{sssec:weave_edge_model}. (Left) The graph of the real part of the multivalued function $h(z)=\pm z$. (Right) The points in the $z$-plane where this real part is univalued.}
		\label{fig:Numerical_A12_Boutroux}
	\end{figure}
\end{center}

\begin{enumerate}
    \item {\bf Weave edge.} Two Lagrangian 2-disks associated to the neighborhood of a weave edge, as in \Cref{fig:Numerical_A12_Boutroux}.\\
    \item {\bf Trivalent vertex.} One Lagrangian 2-disk associated to the neighborhood of a trivalent vertex, as in \Cref{fig:Numerical_D4minus}.\\
    \item {\bf Hexavalent vertex.} Three Lagrangian 2-disks associated to the neighborhood of a hexavalent vertex, as in \Cref{fig:Numerical_A13minus}.\\
    \item {\bf Tetravalent vertex.} Four Lagrangian 2-disks associated to the neighborhood of a tetravalent vertex, as in \Cref{fig:Numerical_Tetra_Boutroux}.\\
\end{enumerate}

\noindent These models are precisely constructed as follows, using the same technique as in \Cref{ssec:Lagrangian_filling_from_AG}. Namely, for each such model we consider a multivalued function $f:\R^2_{x_1,x_2}\lr\R$ and consider the multigraph $\gr(df)\sse(\D^4,\la_\std)$ as in \eqref{eq:graph_df}. In these models $f(x_1,x_2)=\Re(h(z))$ for a multivalued holomorphic function $h:\C_z\lr\C$, where we identify $z=x_1+ix_2$ as in \Cref{lem:complex_curve_real_Lagrangian}.

\subsubsection{Weave edge.}\label{sssec:weave_edge_model} Consider the multivalued function $h(z)$ defined by
\begin{equation}
\boxed{\{h(z)^2=z^2\}}
\end{equation}
\noindent so that the multivalued graph of $h(z)$ is that of the functions $\pm z$. The multigraph of its real part is depicted in \Cref{fig:Numerical_A12_Boutroux} (left), and the locus of points $z\in\C$ where the sheets of $\Re(h(z))$ intersect is drawn in \Cref{fig:Numerical_A12_Boutroux} (right). Geometrically, the Lagrangian surface $\gr(d\Re(h(z)))\sse(\D^4,\la_\std)$ associated to this model consists of two disjoint embedded Lagrangian 2-disks. Note that these two 2-disks are indeed disjoint: the intersection of the two sheets in \Cref{fig:Numerical_A12_Boutroux} (right) disappears in $\gr(d\Re(h(z)))$ because the slopes of the two sheets at the intersections are different, and thus above the locus $\{\Re(z)=0\}$ the differential $dh(z)$ actually has two distinct values as well.

\subsubsection{Trivalent vertex.}\label{sssec:trivalent_model} 
\begin{center}
	\begin{figure}[h!]
		\centering		\includegraphics[width=\textwidth]{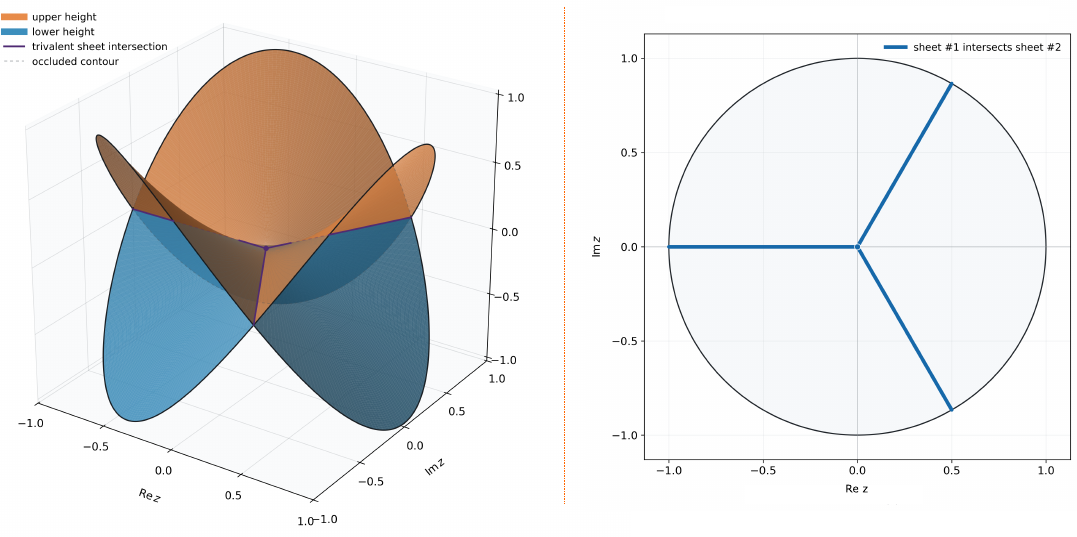}
		\caption{The local model for a trivalent vertex of a weave, as discussed in \Cref{sssec:trivalent_model}. (Left) The graph of the real part of the multivalued function $z^{3/2}$. (Right) The locus of points in the $z$-plane where the function is univalued.}
		\label{fig:Numerical_D4minus}
	\end{figure}
\end{center}

Consider the multivalued function $h(z)$ defined by
\begin{equation}
\boxed{\{h(z)^2=z^3\}}
\end{equation}
\noindent so that the multivalued graph of $h(z)$ is that of the function $z^{3/2}$. The multigraph of $\Re(h(z))$ is depicted in \Cref{fig:Numerical_D4minus} (left), and the locus of points $z\in\C$ where the sheets of $\Re(h(z))$ intersect is drawn in \Cref{fig:Numerical_D4minus} (right). Such a locus consists of a trivalent vertex, which occurs at $z=0$. There are two points above a point in the $z$-plane which is not in this locus, and one point above the points in this locus. In coordinates, the locus of points in \Cref{fig:Numerical_D4minus} (right), where the graph has at most one value, consists of the union of the three rays $(-1)\R_{\geq0}\cup \left(e^{\pi i/3}\R_{\geq0}\right)\cup \left(e^{-\pi i/3}\R_{\geq0}\right)$.\\

\noindent The projection of the Lagrangian $\gr(d(\Re h(z)))$ onto the $z$-plane, which is identified with the first two coordinates of $\R^4$, is a double branched cover with a unique branch point. Its ramification value is $z=0$, corresponding to the trivalent vertex. In particular the Lagrangian $\gr(d(\Re h(z) ))$ is a 2-plane and it can be verified directly that it is embedded.

\subsubsection{Hexavalent vertex.}\label{sssec:hexavalent_local_model} 
\begin{center}
	\begin{figure}[h!]
		\centering		\includegraphics[width=\textwidth]{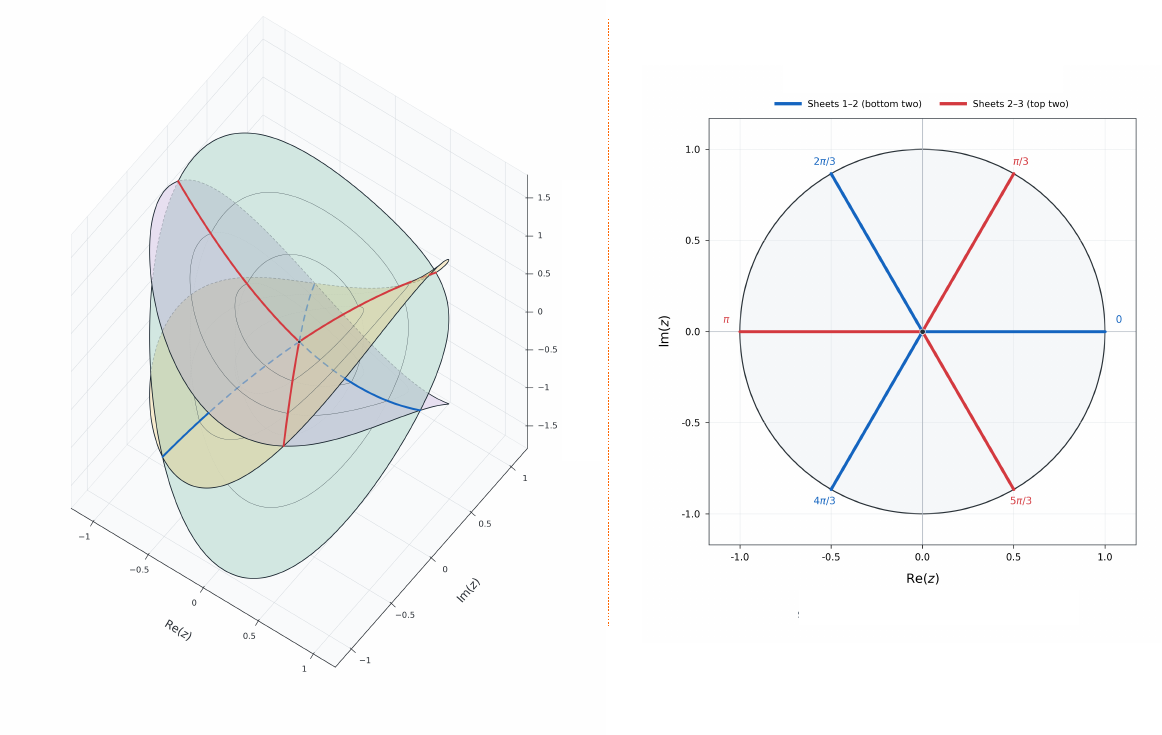}
		\caption{The local model near a hexavalent vertex of a weave, as in \Cref{sssec:hexavalent_local_model}. (Left) The graph of the real part of the multivalued function $\sqrt[3]{z^3}$. (Right) The locus of points in the $z$-plane where the function has fewer than 3 values, with the color indicating which of the different values above a point come together.}
		\label{fig:Numerical_A13minus}
	\end{figure}
\end{center}

Let us consider the multivalued function $h(z)$ defined by
\begin{equation}
\boxed{\{h(z)^3=z^3\}}
\end{equation}
\noindent The multigraph of $\Re(h(z))$ is depicted in \Cref{fig:Numerical_A13minus} (left), and there are three points above a generic point in the $z$-plane. The locus of points $z\in\C$ where the sheets of $\Re(h(z))$ intersect is drawn in \Cref{fig:Numerical_A13minus} (right), where blue indicates that the bottom two sheets intersect and red indicates that the top two sheets intersect. This locus consists of a 6-valent vertex, occurring at $z=0$, and the colors of the edges at the 6-valent vertex alternate between blue and red.\\

\noindent The projection of the Lagrangian $\gr(d(\Re(h(z)))$ onto the $z$-plane is a 3-fold cover, with no branched points. Therefore this Lagrangian consists of 3 disjoint Lagrangian 2-disks, each of them embedded. As a smooth singularity, this $A_1^3$-singularity is equivalent to the origin in the union of the three coordinate planes in $\R^3$, see also \cite[Section 2.2.2]{CasalsZaslow}.

\subsubsection{Tetravalent vertex.}\label{sssec:tetravalent_local_model} 
\begin{center}
	\begin{figure}[h!]
		\centering\includegraphics[width=\textwidth]{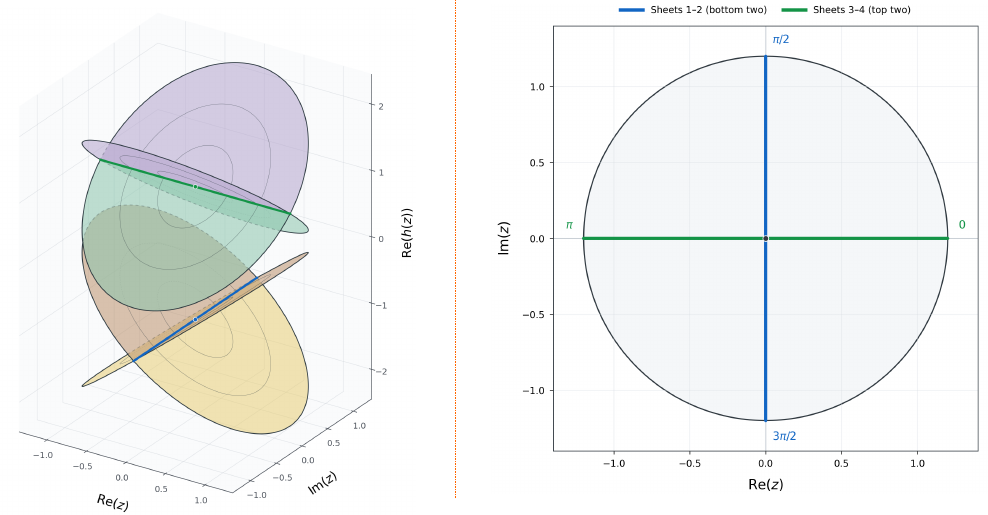}
		\caption{The local model for a tetravalent vertex of a weave, as in \Cref{sssec:tetravalent_local_model}. (Left) The graph of the multivalued function discussed in \Cref{sssec:tetravalent_local_model}. (Right) The locus of points which have fewer than four values above them.}
		\label{fig:Numerical_Tetra_Boutroux}
	\end{figure}
\end{center}

Let us consider the multivalued function $h(z)$ defined by
\begin{equation}
\boxed{\{\bigl((h(z)+1)^2-z^2\bigr)\bigl((h(z)-1)^2+z^2\bigr)=0\}}
\end{equation}
\noindent The graph of its real part near the origin is depicted in \Cref{fig:Numerical_Tetra_Boutroux} (left). Intuitively, this model consists of two copies of the local model in \Cref{sssec:weave_edge_model}, one on top of the other and the top one rotated by $\pi/2$ radians along the axis perpendicular to the $z$-plane. There are four points above a generic point in the $z$-plane. The locus of points $z\in\C$ where the sheets of $\Re(h(z))$ intersect is drawn in \Cref{fig:Numerical_Tetra_Boutroux} (right), where blue indicates that the bottom two sheets intersect and green indicates that the top two sheets intersect. This locus consists of a 4-valent vertex, occurring at $z=0$, and the colors of the edges at the 4-valent vertices alternate between blue and green.

\begin{remark}
These models above can also be parametrized directly with real functions, as was originally done in \cite[Section 2.2]{CasalsZaslow}.\sbs
\end{remark}

\begin{center}
	\begin{figure}[h!]
		\centering\includegraphics[scale=0.7]{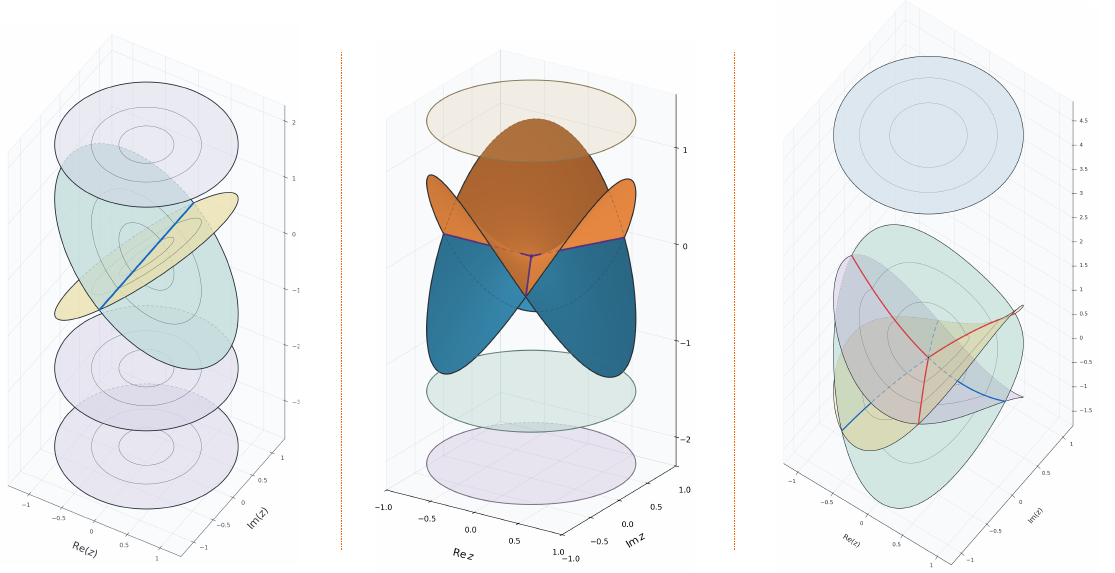}
		\caption{A way to insert the local models discussed in \Cref{sssec:trivalent_model,sssec:hexavalent_local_model,sssec:tetravalent_local_model} stacking parallel sheets above and below, according to the edge labels of a weave $\w$.}
		\label{fig:Numerical_ExampleLocalModel}
	\end{figure}
\end{center}

\noindent More generally, if we work in $\G=\SL_n$, the Lagrangian local model associated to each piece of the weave is obtained as $\gr(df(z))$ where the graph of $f$ consists of $n$ (almost) parallel sheets with only certain sheets weaving according to the models above. Let us label these sheets $1$ through $n$, bottom to top. Then, these models specifically are as follows:

\begin{enumerate}
    \item For a weave line labeled with $s_j$, we consider $\gr(df(z))$ where the graph of $f$ consists of $n$ (almost) parallel sheets except for sheets $j$ and $j+1$ intersecting according to \Cref{fig:Numerical_A12_Boutroux}. See \Cref{fig:Numerical_ExampleLocalModel} (left) for an example.\\

    \item For a trivalent weave vertex whose edges are labeled with $s_j$, we consider $\gr(df(z))$ where the graph of $f$ consists of $n$ (almost) parallel sheets except for sheets $j$ and $j+1$ intersecting according to \Cref{fig:Numerical_D4minus}. See \Cref{fig:Numerical_ExampleLocalModel} (right) for an example.\\

    \item For a hexavalent weave vertex whose edges are alternately labeled with $s_j$ and $s_{j+1}$, we consider $\gr(df(z))$ where the graph of $f$ consists of $n$ (almost) parallel sheets except for sheets $j,j+1$ and $j+2$ intersecting according to \Cref{fig:Numerical_A13minus}.\\

    \item For a tetravalent weave vertex whose edges are alternately labeled with $s_i$ and $s_{k}$, $|i-k|\geq2$, we consider $\gr(df(z))$ where the graph of $f$ consists of $n$ (almost) parallel sheets except for sheets $i,i+1$ and $k,k+1$ intersecting according to \Cref{fig:Numerical_Tetra_Boutroux}.\\
\end{enumerate}

\noindent See \cite[Figure 10]{CasalsZaslow} for another figure illustrating the multigraphs of such local functions.

%%%%%%%%%%%%%%%%%%%%%%%%%%%%%%%%%%%%%%%%%%%%%%%%%%%%%%%%
%%%%%%%%%%%%%%%%%%%%%%%%%%%%%%%%%%%%%%%%%%%%%%%%%%%%%%%%
%%%%%%%%%%%%%%%%%%%%%%%%%%%%%%%%%%%%%%%%%%%%%%%%%%%%%%%%
\subsection{Lagrangian fillings from a Demazure weave}\label{ssec:Lagrangian_filling_from_weave} Let $\w\sse\Sigma$ be an $n$-weave embedded in a smooth surface $\Sigma$. That is, $\w$ is a graph embedded into $\Sigma$ with only vertices of valence $3,4$ and $6$ whose edges are labeled by simple transpositions $\{s_1,\ldots,s_{n-1}\}$, $s_i\in W(\SL_n)$, cf.~\cite[Section 2.3]{CasalsZaslow}. If $\D^2\sse\Sigma$ is a 2-disk such that $\w\cap\D^2$ is a neighborhood of a weave edge or a weave vertex in $\w$, we constructed in \Cref{ssec:Lagrangian_filling_from_weave_local_models} a corresponding local exact Lagrangian surface $L\sse T^*\D^2$ as the multigraph of a differential of a function. This construction is globalized as follows:

\begin{fdef}[Lagrangian filling of a Demazure weave]\label{def:Lagrangian_from_weave} Let $\w\sse\Sigma$ be a $n$-weave embedded in a smooth surface $\Sigma$. By definition, the exact Lagrangian surface $L_\w\sse (T^*\Sigma,\la_\std)$ is the Lagrangian surface constructed as follows:

\begin{enumerate}
    \item Cover $\Sigma$ by a finite collection of 2-disks $\D^2_i$ such that $\w\cap\D^2_i$ consists exactly of one of the four local models of a weave, i.e.~a weave line, or a vertex of valence either 3, 4 or 6, labeled with the corresponding simple transpositions.\\

    \item On each given such 2-disk $\D^2_i$, we choose a model of the corresponding local Lagrangian surface so that for different 2-disks $\D^2_i$ and $\D^2_j$ the models coincide above the intersection $\D^2_i\cap\D^2_j$.
\end{enumerate}
The result from gluing such local Lagrangian surfaces defines a Lagrangian surface $L_\w\sse (T^*\Sigma,\la_\std)$.\sbs
\end{fdef}
\vspace{0.2cm}
\noindent In \Cref{def:Lagrangian_from_weave} we have used the notation $(T^*\Sigma,\la_\std)$ to indicate the cotangent bundle of $\Sigma$, as an exact symplectic 4-manifold, endowed with the canonical Liouville form $\la_\std$. There are choices in \Cref{def:Lagrangian_from_weave}: the cover of $\Sigma$ by the 2-disks and the corresponding model for the local Lagrangian surface above each such 2-disk. The Hamiltonian isotopy type of $L_\w$ is independent of such choices, and the latter choices exist even when requiring the local models to coincide at their intersections. A few comments on \Cref{def:Lagrangian_from_weave}:

\begin{enumerate}
    \item The Lagrangian surface $L_\w\sse (T^*\Sigma,\la_\std)$ is exact and orientable by construction. Intuitively, when we glue the local models over disks in $\Sigma$ from \Cref{ssec:Lagrangian_filling_from_weave_local_models} together: each local model is given as $\gr(df_i)$ for some multivalued function $f_i:\D^2\lr\R$ and it is possible to glue those multivalued functions as well to obtain a global multivalued function $F:\Sigma\lr\R$. The result is that $L_\w$ is also of the form $\gr(dF)\sse T^*\Sigma$ for some multivalued function $F:\Sigma\lr\R$. \footnote{Rigorously, it is the Lagrangian projection of a Legendrian surface, and thus exact.}\\

    \noindent Two examples of the multigraphs of such global multivalued $F:\D^2\lr\R$ in \Cref{def:Lagrangian_from_weave} are depicted in \Cref{fig:Numerical_Trefoil_2weave} and \Cref{fig:Numerical_34_3weave}, with the corresponding weaves $\w$ depicted to their right.\\

    \item The Lagrangian surface $L_\w\sse (T^*\Sigma,\la_\std)$ is immersed, but typically not embedded. For a Demazure weave $\w:\beta\lr\delta(\beta)$, there is a way to choose and glue the local models so as to have $L_\w\sse (T^*\Sigma,\la_\std)$ be embedded, see e.g.~\cite[Section 7.1]{casals2025spectralnetworksbettilagrangians}. The Hamiltonian isotopy class of $L_\w\sse (T^*\Sigma,\la_\std)$ is independent of such choices, cf.~\cite[Section 2]{casals2022conjugate}. Thus, given a Demazure weave $\w\sse\Sigma$, there is a uniquely well-defined embedded exact Lagrangian $L_\w\sse (T^*\Sigma,\la_\std)$, up to Hamiltonian isotopy.\\
\end{enumerate}

\noindent For specificity, we will focus on Lagrangian fillings $L_\w\sse(\D^4,\la_\std)$ associated to a Demazure weave $\w\sse\D^2$. Thus we set $\Sigma=\D^2$ and $\w:\beta\lr\delta(\beta)$ Demazure in \Cref{def:Lagrangian_from_weave} onwards. In this case, $L_\w$ is a Lagrangian filling of
$$\boxed{\dd L_\w=\La_\beta}$$
with $\La_\beta$ as in \Cref{def:Legendrian_from_braid}.

\begin{remark}
In general, even if $\w$ is a Demazure weave on the disk $\Sigma=\D^2$, the embedded exact Lagrangian $L_\w\sse(\D^4,\la_\std)$ cannot be obtained from an algebraic curve via \Cref{lem:complex_curve_real_Lagrangian}. To start, the Legendrian link $\La_\beta$ at the boundary might not even be an algebraic link. Note also that, in contrast to \Cref{lem:complex_curve_real_Lagrangian} and \Cref{ex:special_fillings}, \Cref{def:Lagrangian_from_weave} is such that $L_\w$ is automatically exact, or at least it is relatively simple to prove that. Thus, helpfully, there is no need for the finer tuning in \Cref{lem:complex_curve_real_Lagrangian} and \Cref{ex:special_fillings} to have certain periods vanish and, at least in this regard, manipulating $L_\w$ is much closer to manipulating smooth surfaces in geometric topology than it is to manipulating algebraic curves in complex algebraic geometry.\sbs
\end{remark}

%%%%%%%%%%%%%%%%%%%%%%%%%%%%%%%%%%%%%%%%%%%%%%%%%%%%%%%%
%%%%%%%%%%%%%%%%%%%%%%%%%%%%%%%%%%%%%%%%%%%%%%%%%%%%%%%%
%%%%%%%%%%%%%%%%%%%%%%%%%%%%%%%%%%%%%%%%%%%%%%%%%%%%%%%%
\subsection{Smooth topology of Lagrangian fillings $L_\w$}\label{ssec:smooth_topology_Lagrangian_filling_from_weave} The smooth topology of the Lagrangian fillings constructed in \Cref{def:Lagrangian_from_weave} can be determined as in \Cref{rmk:topology_Lagrangian_filling}. The key fact is that the restriction to $L_\w$ of the projection $\D^4\lr\R^2_{x_1,x_2}$ onto the first two coordinates is a simple $n$-fold branched cover
\begin{equation}\label{eq:branched_cover_map}
\pi:L_\w\lr\R^2_{x_1,x_2}
\end{equation}
\noindent branched over the trivalent vertices of $\w$. The two sheets where the branch point occurs in $L_\w$ over a given trivalent vertex are the $i$th and $(i+1)$st sheets if the weave edges of that trivalent vertex are labeled by $s_i$, cf.~\Cref{fig:Numerical_ExampleLocalModel} (center). In a sense:
\begin{center}
  \fbox{%
    \begin{tabular}{l}
      \itshape A weave $\w$ is a choice of branch cuts for the branched cover $L_\w\lr\R^2_{x_1,x_2}$ which is particularly suited to \\
      \itshape the contact and symplectic topology of $L_\w$ and its resulting Floer-theoretic and sheaf-theoretic invariants.
    \end{tabular}%
  }
\end{center}
\vspace{0.2cm}
\noindent This viewpoint above is particularly useful when studying both additive and monoidal categorifications of Demazure weaves. Geometrically, it also allows for the construction of Lagrangian skeleta in certain 6-dimensional exact symplectic manifolds whose Fukaya categories provide such additive categorifications, cf.~\cite{CasalsChrist}. We emphasize that it is an excessive choice, in that many more cuts are considered than minimally needed to trivialize the branched cover \eqref{eq:branched_cover_map}. That said, as motivated by \Cref{ssec:Lagrangian_filling_from_AG} and the number of new results one can prove with weaves through this viewpoint, it is certainly a useful choice of branch cuts.\\

The properties of the branched cover map \eqref{eq:branched_cover_map} imply the following result:

\begin{lemma}\label{lem:smooth_topology_weave_fillings}
Let $\beta\in\brn$ be a positive braid word and $\w:\beta\lr\delta(\beta)$ a Demazure weave. Then:\\

\begin{enumerate}[label=$(\roman*)$]
    \item The number of components of the Legendrian link $\La_\beta$ is equal to the number of cycles in the cycle decomposition of the permutation underlying $\beta\delta(\beta)$.\\

    \item The Euler characteristic of $L_\w$ is $\chi(L_\w)=n-\#(\mbox{trivalent vertices of }\w)$.\\

    \item Consider the relative first homology $\operatorname{H}_1(L_\w,\tt,\Z)$, where $\tt$ is the set of basepoints of $\dd L_\w=\La_\beta$, one per strand. Then
    $$\operatorname{rk}\left(\operatorname{H}_1(L_\w,\tt,\Z)\right)=b_1(L_\w)+n-|\pi_0(L_\w)|.$$
    \sbs
\end{enumerate}
\end{lemma}

\noindent Since $L_\w$ is orientable, if it is connected, then its number of boundary components, supplied by \Cref{lem:smooth_topology_weave_fillings}.(i), and its Euler characteristic, supplied by \Cref{lem:smooth_topology_weave_fillings}.(ii), uniquely determine the abstract diffeomorphism type of $L_\w$, as a smooth surface. As above, from now on, let us endow $\La_\beta$ with one marked point per strand, located to the right of all the crossings of $\beta$, cf.~the $t_{r,i}$ basepoints in \cite[Figure 23]{CasalsNg22}, $i\in[1,n]$, or \cite[Figure 2 (left)]{casals2026vanishinghigherlegendrianhomology}. Let $\mathfrak{t}$ be such set of marked points.

%%%%%%%%%%%%%%%%%%%%%%%%%%%%%%%%%%%%%%%%%%%%%%%%%%%%%%%%
%%%%%%%%%%%%%%%%%%%%%%%%%%%%%%%%%%%%%%%%%%%%%%%%%%%%%%%%
\subsubsection{The $\SX$ and $\SA$-lattices}\label{sssec:X_lattice} Following the viewpoint of \cite{FG2006,FG09}, cluster seeds can be interpreted as bases of a certain lattice with additional data, see e.g.~\cite[Appendix A]{GHKK18}. In the context of Demazure weaves $\w$, the lattice $M$, which corresponds to the cluster $\SX$-variety, can be interpreted in terms of the Lagrangian filling $L_\w$. Specifically, given a Demazure weave $\w:\beta\to\delta(\beta)$, the $\SX$-lattices are:
$$M:=\hH_1(L_\w,\tt;\Z)\qquad\mbox{ and }\qquad M_{un}:=\hH_1(L_\w;\Z),$$
\noindent where $\hH_1(L_\w,\tt;\Z)$ denotes the first integral relative homology group of the relative pair $(L_\w,\tt)$, with $\tt\sse\dd L_\w$ denoting the set of marked points in $\dd L_\w=\La_\b$. The natural inclusion of the absolute integral homology group $\hH_1(L_\w;\Z)\hookrightarrow \hH_1(L_\w,\tt;\Z)$ provides an inclusion of the unfrozen lattice $M_{un}$, a.k.a.~the mutable sublattice, inside of $M$.

\begin{center}
	\begin{figure}[h!]
		\centering\includegraphics[width=\textwidth]{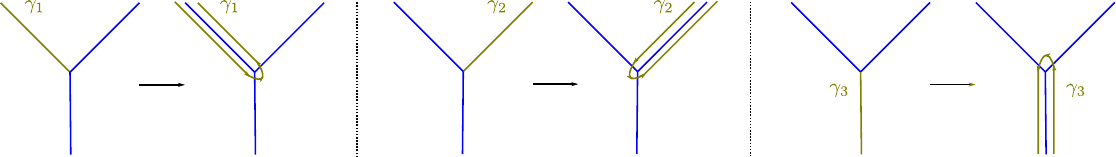}
		\caption{Perturbations of the 1-cycles associated to a weave cycle $\g_i$, $i\in[1,3]$. On the left of each of these three situations, the 1-cycle $\g_i$ projects directly to an edge of the weave $\w$. On the right, the 1-cycle has been deformed in $L_\w$ so that its projection is an embedded relative cycle. These perturbations are necessary in order to geometrically compute the intersections between weave cycles.}
		\label{fig:Cycles_Trivalent_Perturb}
	\end{figure}
\end{center}

\noindent Dually, there is also a construction for the $\SA$-lattice:
$$N:=\hH_1(L_\w\setminus \tt,\La_\beta\setminus \tt;\Z)$$
and its corresponding mutable sublattice. This $\SA$-lattice above appears naturally when studying the long exact sequence of the triple $(L_\w,\dd L_\w,\tt)$, cf.~\cite[Section 3.10]{CasalsWeng22}, as $\dd L_\w=\La_\beta$. Note that there is a non-degenerate intersection pairing between such $M$ and $N$, providing the appropriate duality.\\

In \Cref{sssec:Lusztigcycles} we introduced the Lusztig cycles over a Demazure weave $\w:\beta\to\delta(\beta)$, cf.~\Cref{def:Lusztig_cycles}. These weave cycles have key geometric meaning in the context of contact and symplectic topology. At the most basic level there is the following construction. For each Lusztig cycle $\g$ of $\w$, we can construct an element $[\gamma]\in \hH_1(L_\w,\tt;\Z)$ by drawing certain curves in the plane $\R^2_{x_1,x_2}$ where $\w$ is depicted and lifting them through the branched cover $L_\w\lr \R^2_{x_1,x_2}$. This is explained in \cite[Section 7.4]{CGGLSS25} and \cite[Section 2.4]{CasalsZaslow}. Intuitively, given a weave edge labeled by $s_i$, we can consider the two arcs above the weave edge at the sheets $i$ and $(i+1)$st. The interesting fact about the combinatorics of Lusztig cycles is that these arcs glue well at the trivalent and hexavalent vertices.

\begin{center}
	\begin{figure}[h!]
		\centering\includegraphics[width=\textwidth]{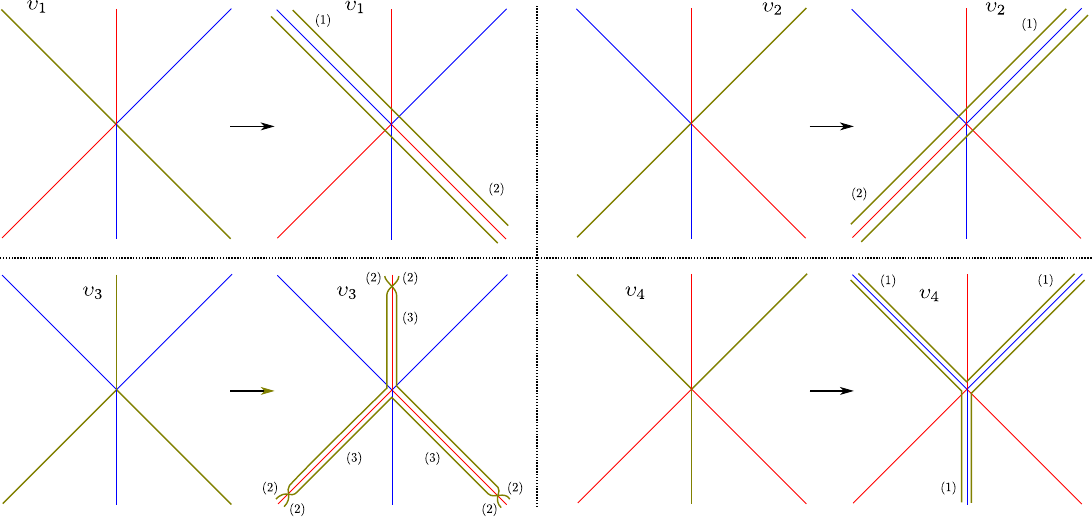}
		\caption{Perturbations of a few Lusztig cycles near the hexavalent vertex of a weave.}
		\label{fig:Cycles_Hexavalent_Perturb2}
	\end{figure}
\end{center}

\begin{center}
	\begin{figure}[h!]
		\centering\includegraphics[width=\textwidth]{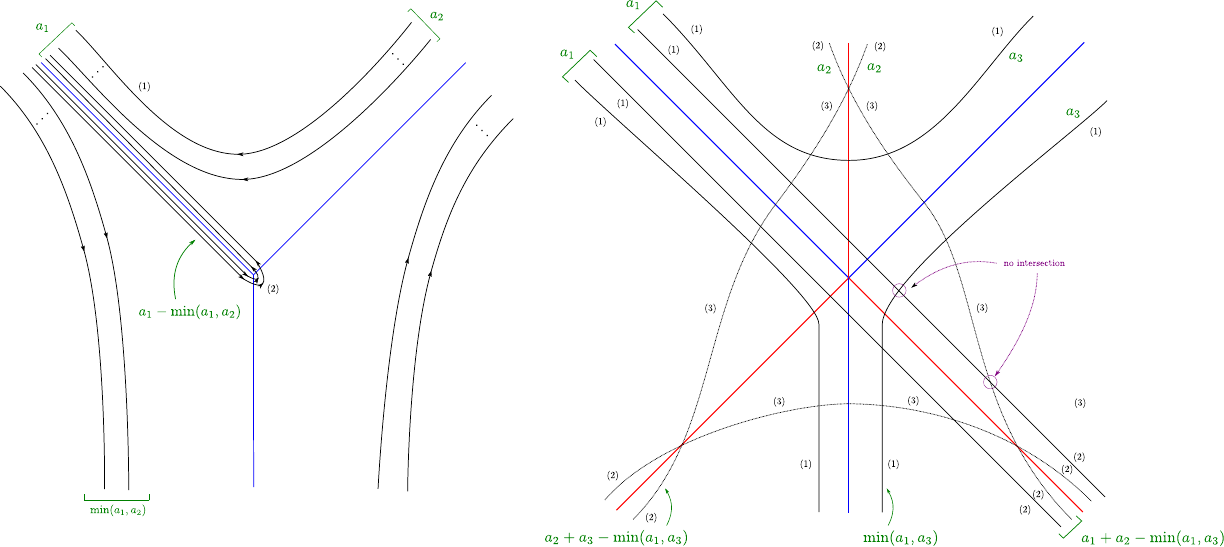}
		\caption{The geometric models explaining how the tropical propagation rules for Lusztig cycles in $\w$ arise when intersecting relative $1$-cycles in the surface $L_\w$.}
		\label{fig:Cycles_Lusztig2}
	\end{figure}
\end{center}

A way to visually understand how such gluing occurs is to perturb these arcs slightly off so they do not live on top of the weave edge under \eqref{eq:branched_cover_map}, and then directly draw their behavior via their images in $\R^2_{x_1,x_2}$, lifting those curves in $\R^2_{x_1,x_2}$ to the appropriate sheets. \Cref{fig:Cycles_Trivalent_Perturb} depicts an instance of this near a trivalent vertex with label $s_i$, where the only two sheets that are relevant are the $i$ and $(i+1)$st sheets. Similarly, \Cref{fig:Cycles_Hexavalent_Perturb2} depicts an example of this near a hexavalent vertex. The tropical Lusztig rules near trivalent and hexavalent vertices are explained in terms of these curves in \Cref{fig:Cycles_Lusztig2}.

\begin{remark}
In the case of the frozen Lusztig cycles, one needs to adjust the ends so as to lie in $\tt$: such adjustment is possible but different choices lead to interesting quasi-cluster equivalences, cf.~\cite[Section 5.5]{CGGLSS25} and \cite{casals2023demazureweavesreducedplabic}.\sbs
\end{remark}

\noindent A similar construction applies to the dual Lusztig cycles from \Cref{def:dual_Lusztig_cycles}. From a dual Lusztig cycle $\rho$ of $\w$ one can similarly construct an element $[\rho]\in\hH_1(L\setminus\tt,\La\setminus\tt;\Z)$ . Given the Lusztig cycles $\{\g_i\}$ and their duals $\{\rho_i\}$, their corresponding homology classes $\{[\g_i]\}$ and $\{[\rho_i]\}$ are Poincar\'e dual to each other, under the non-degenerate intersection pairing between $\hH_1(L_\w,\tt;\Z)$ and $\hH_1(L_\w\setminus \tt,\dd L_\w\setminus \tt;\Z)$.

\noindent Since the topology of $L_\w$ is essentially generated by the branched points of \eqref{eq:branched_cover_map}, and thus the trivalent vertices of $\w$ from a combinatorial viewpoint, it can be deduced that the Lusztig cycles of $\w$ provide a basis for the lattice. We summarize these properties in the following:

\begin{lemma}\label{lem:lustzig_cycles_basis}
Let $\w:\beta\to\delta(\beta)$ be a Demazure weave. Then the following hold:

\begin{enumerate}[label=$(\roman*)$]
    \item The homology classes associated to the mutable Lusztig cycles of $\w$ form a basis of the mutable $\SX$-lattice $\hH_1(L_\w;\Z)$. Similarly, the homology classes $\{[\g_i]\}$ associated to the Lusztig cycles $\{\g_i\}$ of $\w$ form a basis of the $\SX$-lattice $\hH_1(L_\w,\tt;\Z)$.\\

    \item The homology classes $\{[\rho_i]\}$  associated to the dual Lusztig cycles $\{\rho_i\}$ of $\w$ form a basis of the $\SA$-lattice $\hH_1(L_\w\setminus \tt,\La_\beta\setminus\tt;\Z)$.\\

    \item Let $\eta_1$ be a Lusztig cycle or a dual Lusztig cycle, and let $\eta_2$ another such cycle. Then
    $$\langle \eta_1,\eta_2\rangle_\w=\langle[\eta_1],[\eta_2]\rangle$$
    where the left hand side is the weave intersection form, as in \Cref{ssec:intersection_cycles}, and the right hand side is the algebraic intersection product in homology, induced by the cap product.\\
\end{enumerate}
\noindent In particular, the two bases $\{[\g_i]\}$  and $\{[\rho_i]\}$ are dual bases under the non-degenerate intersection pairing between $\hH_1(L_\w,\tt;\Z)$ and $\hH_1(L_\w\setminus \tt,\La_\beta\setminus\tt;\Z)$.\sbs
\end{lemma}

\noindent Note that, technically, the algebraic intersection pairing of two homology classes in $\hH_1(L_\w,\tt;\Z)$ requires a boundary push-off (or an order convention). We are implicitly making such a choice when computing in $\hH_1(L_\w,\tt;\Z)$, see \cite[Section 3.2]{CGGJ24}.

%%%%%%%%%%%%%%%%%%%%%%%%%%%%%%%%%%%%%%%%%%%%%%%%%%%%%%%%
%%%%%%%%%%%%%%%%%%%%%%%%%%%%%%%%%%%%%%%%%%%%%%%%%%%%%%%%
%%%%%%%%%%%%%%%%%%%%%%%%%%%%%%%%%%%%%%%%%%%%%%%%%%%%%%%%
\subsection{Demazure weaves and Lagrangian cobordisms}\label{ssec:Lagrangian_coboridsms} In \Cref{sec:weaves_braid_varieties}, \Cref{sec:braid_variety_cluster} and \Cref{sec:Reeb_variables} we studied the connection of a Demazure weave $\w:\beta\to\delta(\beta)$ to the braid variety $X(\beta)$ by scanning the weave $\w$ top-to-bottom, cf.~\Cref{ssec:weave_tori_braid_varieties}. Specifically, given a vertex $v\in\w$ we can consider the horizontal slices $h_\pm$ above and below $v$ and the corresponding positive braid words $\beta_{\pm}$, as in \Cref{fig:BraidVariety_DiagramMoves2} and \Cref{fig:BraidVariety_Trivalent2}:
\begin{enumerate}
    \item At hexavalent and tetravalent vertices, \Cref{lem:braid_invariance} provides isomorphisms
    \begin{equation}\label{eq:cobordism_maps_braid_varieties}R3:X(\beta_-)\stackrel{\simeq}{\lr} X(\beta_+),\quad C:X(\beta_-)\stackrel{\simeq}{\lr} X(\beta_+)
    \end{equation}

    \item At a trivalent vertex, \Cref{lem:braidvarieties_trivalent} provides an open embedding \begin{equation}\label{eq:cobordism_maps_braid_varieties2}
        \iota_v:X(\beta_-)\times\bG_m\hookrightarrow X(\beta_+).
    \end{equation}
 \end{enumerate}

\noindent The Lagrangian local models for these vertices, as presented in \Cref{ssec:Lagrangian_filling_from_weave_local_models}, actually yield Lagrangian cobordisms $L_v$ from $\La_{\beta_-}$ to $\La_{\beta_+}$, $\La_{\beta_-}$ at the concave end and $\La_{\beta_+}$ at the convex end. See for instance \cite[Section 6]{EtnyreNg22_Survey} for a survey on Lagrangian cobordisms. These Lagrangian cobordisms $L_v$ occur in the symplectization of $(\R^3,\xi_\std)$, and are exact and embedded. In particular, each such Lagrangian cobordism $L_v$ induces a dg algebra map
\begin{equation}\label{eq:cobordism_map}
\varepsilon_{L_v}:\mathbb{A}(\La_{\beta_+})\lr \mathbb{A}(\La_{\beta_-})\otimes_{\Z[C_{-*}(\Omega \La_{\beta_-})]} \Z[C_{-*}(\Omega L_v)]
\end{equation}
as in \cite[Theorem 1.6]{EHK} or \cite[Theorem 6.3]{EtnyreNg22_Survey}. See also \cite{LiCobordism1,li2024lagrangiancobordismfunctormicrolocal}.

\begin{example} The maps \Cref{eq:cobordism_map} are a refinement of \Cref{eq:augmentation_map_filling}, which corresponds to the case where $\La_{\beta_-}=\emptyset$ is empty. Indeed, a Lagrangian filling is a Lagrangian cobordism with empty concave (bottom) end: \Cref{prop:torus_from_filling} encodes precisely \eqref{eq:cobordism_map} where $\mathbb{A}(\La_{\beta_-})\cong \Z$, generated by the unit (the empty word on Reeb chords), and we identify
$$\hH_0(C_{-*}(\Omega L;\Z))^{ab}\cong\Z[\pi_1(L)]^{ab}\cong\Z[\hH_1(L;\Z)]\cong\Z[s_1^{\pm1},\ldots,s^{\pm1}_{b_1(L)}],$$
as $L$ is a punctured smooth surface, cf.~\cite{EkholmLekili2023}.\sbs    
\end{example}
\noindent In particular, $\varepsilon_{L_v}$ induces a morphism between the $\mbox{H}_0$ of their abelianizations. After applying $\Spec$ to $\hH_0$ of the abelianization of \eqref{eq:cobordism_map} and base changing to $k$, we obtain:
\begin{enumerate}
    \item At a hexavalent or tetravalent vertex $v\in\w$, the Lagrangian cobordism $L_v$ is the Lagrangian concordance given by the trace of a Legendrian Reidemeister III move (\Cref{fig:Front_Arc_Cusp_Crossing_R3_Reidemeister}) or the Legendrian isotopy given by a commutation of distant crossings, respectively. The Floer theoretic map \eqref{eq:cobordism_map} provides isomorphisms
    \begin{equation}\label{eq:Floer_maps_iso}\Spec\mbox{H}_0(\mathbb{A}^{ab}(\La_{\beta_-}))\stackrel{\simeq}{\lr} \Spec\mbox{H}_0(\mathbb{A}^{ab}(\La_{\beta_+})),\quad \Spec\mbox{H}_0(\mathbb{A}^{ab}(\La_{\beta_-}))\stackrel{\simeq}{\lr} \Spec\mbox{H}_0(\mathbb{A}^{ab}(\La_{\beta_+}))
    \end{equation}

    \item At a trivalent vertex, \eqref{eq:cobordism_map} provides an open embedding \begin{equation}\label{eq:Floer_inclusion_trivalent}
        \Spec\mbox{H}_0(\mathbb{A}^{ab}(\La_{\beta_-}))\times\bG_m\hookrightarrow \Spec\mbox{H}_0(\mathbb{A}^{ab}(\La_{\beta_+})),
    \end{equation}
    where the $\bG_m$-factor can be geometrically identified as $\bG_m\cong \Spec k[H_1(L_\w,\La_{\beta_-};\Z)]\cong \Spec k[\tau^{\pm1}]$ is generated by the relative homology class $\tau$ of the stable manifold of the unique critical point of the saddle cobordism $L_\w$, where $\w$ has a unique trivalent vertex.
 \end{enumerate}

\noindent  Here is a refinement of \Cref{prop:torus_from_filling}, in the form of a horizontal slice decomposition, in the context of Lagrangian cobordisms associated to weave vertices:

\begin{lemma}[Floer theory morphisms induce weave morphisms]\label{lem:floer_morphisms_weave_morphisms}
Let us consider the identification $X(\beta)\cong \Spec \mbox{H}_0(\mathbb{A}^{ab}(\La_{\beta}))$ of \Cref{thm:braid_varieties_from_Legendrian}.(i), see also \Cref{eq:braid_variety_dga2}. Then the weave morphisms \eqref{eq:cobordism_maps_braid_varieties} and \eqref{eq:cobordism_maps_braid_varieties2} coincide with the Floer-theoretically induced morphisms \eqref{eq:Floer_maps_iso} and \eqref{eq:Floer_inclusion_trivalent}, respectively.\sbs
\end{lemma}

\noindent \Cref{lem:floer_morphisms_weave_morphisms} follows from the description of the Floer-theoretic maps in \cite[Section 5]{CasalsNg22}. There is a result similar to \Cref{lem:floer_morphisms_weave_morphisms} for the sheaf-theoretic induced morphisms, using the other identifications in \Cref{thm:braid_varieties_from_Legendrian}. A few comments are in order:

\begin{enumerate}
    \item In the diagrammatic calculus of Demazure weave, we prominently used the Reeb coordinates, as discussed in \Cref{sec:Reeb_variables}. These variables coincide with the Floer-theoretic variables that keep track of the degree 0 Reeb chords associated to the $(-1)$-closure of $\beta\delta(\beta)$. Specifically, for such fronts, there are Reeb chords only in degrees 0 and 1, which generate the Legendrian dg algebra $\mathbb{A}(\La_{\beta})$. The Reeb chords of degree 0 correspond to the crossings of $\beta\delta(\beta)$. Intuitively, we can think of the degree 0 Reeb chords as generators and the differentials of the degree 1 Reeb chords as relations among them, cf.~\cite[Section 5]{CasalsNg22} and \cite[Corollary 1.1]{casals2026vanishinghigherlegendrianhomology}. Via the identification
    $$X(\beta)\cong \Spec \mbox{H}_0(\mathbb{A}^{ab}(\La_{\beta}))$$
    of \Cref{thm:braid_varieties_from_Legendrian}.(i), this description precisely leads to the Reeb $z$-coordinates for $X(\beta)$ and the corresponding equations from \Cref{lem:equations_braidvar}.\footnote{To be precise, from a Floer-theoretic viewpoint, there are also $z$-coordinates associated to the degree 0 Reeb chords from the crossings of $\delta(\beta)$, but the differentials of the degree 1 Reeb chords can be used to express those in terms of the $z$-coordinates of the degree 0 Reeb chords associated to the crossings of $\beta$.}\\

    \item The torus $\T$-action on $X(\beta)$ described in \Cref{sssec:braid_variety_further_symmetries} also has a symplectic interpretation. Namely, it corresponds to the $\T$-action induced by link automorphisms on the Legendrian dg algebra, cf.~\cite[Section 3.2]{CasalsNg22}. In brief, these are automorphisms induced by carrying around the basepoints of a Legendrian link around their corresponding components. In our case, the basepoints are all located to the right of all the degree 0 Reeb chords, i.e.~of all the crossings of $\beta\delta(\beta)$. Moving them past the crossings of $\beta\delta(\beta)$ and keeping track of their action on the corresponding degree 0 Reeb chords, using \cite[Section 4.1.1]{CasalsNg22}, leads to a Floer-theoretic proof of the formula \Cref{lem:torusaction_coordinates}.
\end{enumerate}
%%%%%%%%%%%%%%%%%%%%%%%%%%%%%%%%%%%%%%%%%%%%%%%%%%%%%%%%
%%%%%%%%%%%%%%%%%%%%%%%%%%%%%%%%%%%%%%%%%%%%%%%%%%%%%%%%
%%%%%%%%%%%%%%%%%%%%%%%%%%%%%%%%%%%%%%%%%%%%%%%%%%%%%%%%
\begin{remark} It should be mentioned that, in line with \Cref{ssec:Demazure_weaves_from_combinatorics}, a first connection between the combinatorics of reduced plabic graphs, and Lagrangian fillings was studied in \cite{STWZ19}. Such Lagrangian enhancements of the conjugate surfaces from \cite{GonKen13} illustrated some of the advantages of using the microlocal theory of sheaves in relation to cluster structures, especially the use of microlocal parallel transport. Note also that there exists a hybrid calculus, developed in \cite{casals2022conjugate}, combining (and transitioning between) these Lagrangian conjugate surfaces and weaves.\sbs
\end{remark}

\section{Cluster aspects of Lagrangian fillings}

The object of this last section is to bring together the cluster algebra structure on $\Z[X(\beta)]$ and the symplectic topology of $\La_\beta$ and its Lagrangian fillings $L_\w$ arising from Demazure weaves $\wb$. First, in \Cref{ssec:torus_Lagrangian_filling_from_weave} we describe the torus $T_\w\sse X(\beta)$ underlying a cluster seed in terms of the symplectic topology of $L_\w$, and interpret the associated $\SA$ and $\SX$-lattices with their intersection forms via the algebraic topology of $L_\w$. Second, the symplectic construction of the corresponding cluster $\SA$ and $\SX$-variables is presented in \Cref{ssec:X_coordinates_Lagrangian_filling_from_weave}, and \Cref{ssec:cluster_mutation_Lagrangian_filling_from_weave} explains the relation between Lagrangian disk surgeries and weave and cluster mutations.
%%%%%%%%%%%%%%%%%%%%%%%%%%%%%%%%%%%%%%%%%%%%%%%%%%%%%%%%
%%%%%%%%%%%%%%%%%%%%%%%%%%%%%%%%%%%%%%%%%%%%%%%%%%%%%%%%
%%%%%%%%%%%%%%%%%%%%%%%%%%%%%%%%%%%%%%%%%%%%%%%%%%%%%%%%
\subsection{The torus $T_\w$ associated to the Lagrangian filling $L_\w$}\label{ssec:torus_Lagrangian_filling_from_weave} Let $\beta\in\brn$ be a positive braid word and $X(\beta)$ its braid variety. Consider a Demazure weave $\w:\beta\to\delta(\beta)$. By \Cref{cor:weave_to_toric_embedding}, specifically \eqref{eq:weave_tori_map}, there exists an open embedding
$$\iota_\w:(\bG_m)^d\hookrightarrow X(\beta).$$
\noindent Independently, \Cref{prop:torus_from_filling} or \eqref{eq:cobordism_map} applied to the Lagrangian filling $L=L_\w$ provides an open inclusion
\begin{equation}\label{eq:Floer_torus}
\Spec(\hH_0(\varepsilon_{L_\w})):\Spec(\Z[\hH_1(L_\w,\tt;\Z)])\hookrightarrow \Spec\mbox{H}_0(\mathbb{A}^{ab}(\La_{\beta}))
\end{equation}
of algebraic varieties over $\Z$. \Cref{lem:floer_morphisms_weave_morphisms} and \cite[Section 5]{CasalsNg22} imply the following:
\begin{proposition}\label{prop:weave_tori_and_Lagrangian_fillings}
In the notation above, there exists a commutative diagram
\[
\begin{tikzcd}
    \Spec(\Z[\hH_1(L_\w,\tt;\Z)]) \arrow[r, hook,"\eqref{eq:Floer_torus}"] \arrow[d, "\tiny\mbox{id}"'] & \Spec(\hH_0(\mathbb{A}^{ab}(\La_\beta))) \arrow[d, "\mbox{\tiny\Cref{thm:braid_varieties_from_Legendrian}.(i)}","\cong"
'] \\
    T_\w \arrow[r, hook,"\iota_\w"',"\eqref{eq:weave_tori_map}"] & X(\beta)
\end{tikzcd}
\]
In particular, $\operatorname{im}(\Spec(\hH_0(\varepsilon_{L_\w})))=\operatorname{im}(\iota_\w)$
as open toric charts in $X(\beta)$.\sbs
\end{proposition}
\noindent \Cref{prop:weave_tori_and_Lagrangian_fillings} states that the Floer-theoretic toric chart associated to $L_\w$ coincides with the toric chart $T_\w\sse X(\beta)$ constructed via the weave $\w$ using the transversality conditions on flags given by $\w$, cf.~\Cref{ssec:weave_tori_braid_varieties}. The identification of the weave torus $T_\w$ with $\Spec(\Z[\hH_1(L_\w,\tt;\Z)])$ allows us to give a symplectic geometric description of the points in $T_\w$, as follows.\\

\noindent The $R$-points of the algebraic variety $\Spec(\Z[\hH_1(L_\w,\tt;\Z)])$, where $R$ is a commutative ring, are given by
\begin{equation}\label{eq:points_of_Tw}
\mbox{Hom}_{\Z-\mbox{alg}}(\Z[\hH_1(L_\w,\tt;\Z)],R)\cong \mbox{Hom}_{\Z-\mbox{mod}}(\hH_1(L_\w,\tt;\Z),R^\times)\cong \hH^1(L_\w,\tt;R^\times)
\end{equation}

\noindent Consequently, \color{blue}the $R$-points of a weave torus $T_\w\cong \Spec(\Z[\hH_1(L_\w,\tt;\Z)])$ parametrize isomorphism classes of rank-one $R$-local systems on the surface $L_\w$, with a common trivialization of the stalks at $\tt$\color{black}. Indeed, this follows from the isomorphism $H_1\cong \pi_1^{ab}$, and the stalks must be isomorphic to the rank-one free $R$-module $R$ because $R^\times\cong\mbox{GL}_1(R)$.\\

Recall that the cluster algebra structure on $X(\beta)$ has such weave tori $T_\w$ underlying some of its cluster seeds. \Cref{prop:weave_tori_and_Lagrangian_fillings} also explains the choice of $\SX$-lattice $M=\hH_1(L_\w,\tt;\Z)$ in \Cref{sssec:X_lattice}, as $\varepsilon_{L_\w}$ has target the group ring of the character lattice of $\operatorname{im}(\Spec(\hH_0(\varepsilon_{L_\w})))$, which is naturally isomorphic to $\hH_1(L_\w,\tt;\Z)$.\\

\begin{center}
  \fbox{%
    \begin{tabular}{l}
      \itshape The torus $T_\w\sse X(\beta)$ associated to a Demazure weave $\wb$ appears naturally\\
      \itshape  as a Hamiltonian invariant of the Lagrangian filling $L_\w$ via its first homology lattice.
    \end{tabular}%
  }
\end{center}
\vspace{0.2cm}
\noindent Specifically, following \Cref{ssec:A_and_X_schemes_weave} and \Cref{sssec:X_lattice}, there is a choice of whether we consider $X(\beta)$ to be a cluster $\SA$-variety or a cluster $\SX$-variety:

\begin{enumerate}
    \item If we are in the context of $\SX$-cluster variables, the $\SX$-torus $T_\w\sse X(\beta)$ appears as
    \begin{center}
        \fbox{$T_\w^\SX\cong \Spec\Z[\hH_1(L_\w,\tt;\Z)]$}
    \end{center}
    \noindent Equivalently, $\hH_1(L_\w,\tt;\Z)$ is the character lattice of the $\SX$-torus $T_\w^\SX$. In terms of local systems, the $R$-points of $T^\SX_\w$ parametrize $\mbox{GL}_1(R)$-local systems on $L_\w$ with a common trivialization for the stalks at $\tt$.\\
    
    \item For the $\SA$-cluster variables, the $\SA$-torus $T_\w\sse X(\beta)$ appears via
     \begin{center}
        \fbox{$T_\w^\SA\cong \Spec\Z[\hH_1(L_\w\setminus \tt,\La_\beta\setminus \tt;\Z)]$}
    \end{center}
    \noindent Equivalently, $\hH_1(L_\w\setminus \tt,\La_\beta\setminus \tt;\Z)$ is the character lattice of the $\SA$-torus $T_\w^\SA$. In terms of local systems, the $R$-points of $T^\SA_\w$ parametrize $\mbox{GL}_1(R)$-local systems on $L_\w$ with a common trivialization at the intervals $\La_\beta\setminus \tt$.\\
\end{enumerate}

\noindent As in \Cref{lem:lustzig_cycles_basis}, we have the lattice duality
\begin{equation}\label{eq:duality_lattices}
\hH_1(L_\w,\tt;\Z)^\vee\cong\hH_1(L_\w\setminus \tt,\La_\beta\setminus \tt;\Z)
\end{equation}
induced by the intersection pairing. Thus, having fixed an $\SA$-torus $T_\w\sse X(\beta)$, its cocharacter lattice is naturally identified with $\hH_1(L_\w,\tt;\Z)$.\\

In order to abstractly name $T_\w^\SX$ or $T_\w^\SA$ from $L_\w$, we only need the abstract homotopy type of $L_\w$, not its specific Lagrangian embedding into $L_w\hookrightarrow\D^4$ provided by the defining construction in \Cref{ssec:Lagrangian_filling_from_weave}. In contrast, both of the following facts require the information of the Lagrangian embedding $L_w\hookrightarrow\D^4$:

\begin{enumerate}
    \item The specific embedding $T_\w\hookrightarrow X(\beta)$ is induced by the specific embedding $L_w\hookrightarrow\D^4$. It must be Lagrangian in order to induce any morphisms for the Legendrian dg algebra of $\La_\beta$ and, from this viewpoint, embeddedness guarantees Floer-theoretic unobstructedness of such morphisms.\\

    \item The bases for the two lattices in \eqref{eq:duality_lattices} that provide cluster structures are given by the Lusztig cycles and their duals, cf.~\Cref{lem:lustzig_cycles_basis}. From a symplectic topological viewpoint, naming such Lusztig cycles intrinsically requires the specific embedding $L_w\hookrightarrow\D^4$ as well. Indeed, intuitively, these 1-cycles in $L_\w$, when mutable, are selected by the property of bounding an embedded Lagrangian disk in the complement $\D^4\setminus L_\w$, and such that the boundaries of (the closures of) such 2-disks cleanly intersect $L_\w$.
\end{enumerate}

%%%%%%%%%%%%%%%%%%%%%%%%%%%%%%%%%%%%%%%%%%%%%%%%%%%%%%%%
%%%%%%%%%%%%%%%%%%%%%%%%%%%%%%%%%%%%%%%%%%%%%%%%%%%%%%%%
%%%%%%%%%%%%%%%%%%%%%%%%%%%%%%%%%%%%%%%%%%%%%%%%%%%%%%%%
\subsection{The $\SX$ and $\SA$-coordinates associated to the Lagrangian filling $L_\w$}\label{ssec:X_coordinates_Lagrangian_filling_from_weave} At this stage, a weave $\wb$ provides an open torus chart $T_\w\sse X(\beta)$, via \Cref{cor:weave_to_toric_embedding}, and \Cref{prop:weave_tori_and_Lagrangian_fillings} above shows that this torus can be geometrically understood in terms of integral homology groups associated to the surface $L_\w$. In particular, we have geometrized in terms of the Lagrangian fillings $L_\w$ of $\La_\beta$ the following data:
\begin{enumerate}
    \item The braid variety $X(\beta)$, via the isomorphism of \Cref{thm:braid_varieties_from_Legendrian}.(i).
    $$X(\beta)\cong \Spec \mbox{H}_0(\mathbb{A}^{ab}(\La_{\beta}))$$

    \item Each weave torus $T_\w\sse X(\beta)$, via the isomorphisms in \Cref{ssec:torus_Lagrangian_filling_from_weave}, cf.~\Cref{prop:weave_tori_and_Lagrangian_fillings}:
    $$T_\w^\SX\cong \Spec\Z[\hH_1(L_\w,\tt;\Z)],\qquad T_\w^\SA\cong \Spec\Z[\hH_1(L_\w\setminus \tt,\La_\beta\setminus \tt;\Z)].$$
    \item The Lusztig cycles $\{\g_i^\w\}$ and their duals $\{\rho_i^\w\}$, geometrically realized as 1-cycles in $L_\w$ via \Cref{sssec:X_lattice}, cf.~\Cref{lem:lustzig_cycles_basis}. In particular, the weave quiver $Q_\w$ is geometrically realized via the intersection form in $L_\w$.\\
\end{enumerate}

\noindent Thus, from the viewpoint of the theory of cluster algebras, the key ingredients that remain to be geometrized are the $\SX$ and $\SA$-cluster variables. For that we recall \eqref{eq:points_of_Tw}, which allows us to interpret an $R$-point of a weave torus as an isomorphism class of a rank-one $R$-local system on $L_\w$ with certain trivializations. For an $R$-point $p\in T_\w$, we denote by $\LL_p$ the corresponding $\GL_1(R)$-local system. The $\SX$-variables can be geometrically interpreted as follows:

\begin{fdef}[Cluster $\SX$-variables via $L_\w$]\label{def:Xvariables_via_fillings}
Let $\beta$ be a positive braid word and $\wb$ a Demazure weave. Consider the vertex $i\in Q_\w$ of the weave quiver, associated to the Lusztig cycle $\g_i^\w$, and its homology class $[\g_i^\w]\in \hH_1(L_\w,\tt;\Z)$. By definition, the $\SX$-cluster variable $\SX_i^\w\in \Z(X(\beta))$ associated to $i$ is the rational function
\begin{equation}\label{eq:Xvariable_via_fillings}
\SX_i^\w:X(\beta)\lr\PP^1_\Z,\quad\mbox{s.t. } \SX_i^\w(p)=\hol_{[\g_i^\w]}(\LL_p)\quad\forall p\in T^\SX_\w,
\end{equation}
where $\hol_{[\g_i^\w]}(\LL_p)\in\GL_1(R)$ denotes the monodromy of the local system $\LL_p$ along the class of $\g_i^\w$.\sbs
\end{fdef}
\noindent In \Cref{def:Xvariables_via_fillings}, the rational function $\SX_i^\w$ is only explicitly defined in $T_\w\sse X(\beta)$. Since $X(\beta)$ is irreducible and connected, such a rational function uniquely extends to a rational function on $X(\beta)$. Note that the expression {\it monodromy of a local system} can be taken literally if $\dd\g_i^\w=\emptyset$, i.e.~it is the (conjugacy class of the) automorphism obtained by parallel transporting via $\LL_p$ along $\g_i^\w$ until we reach back a chosen basepoint.\\

Suppose otherwise that $\g_i^\w$ is represented by a relative 1-cycle with its boundary $\dd\g_i^\w\neq\emptyset$ contained in $\tt$, as it is the case for the frozen Lusztig cycles, and assume that it is a connected interval. Then the expression $\hol_{[\g_i^\w]}(\LL_p)$ refers to the automorphism of the stalk at one of its endpoints obtained by parallel transporting along the interval and then using the common trivialization for the stalk at both endpoints, which lie in $\tt$.\footnote{If there were several connected components, we take the product of these elements in $\GL_1(R)$.} By construction, these $\SX$-cluster functions $\SX_i^\w$ are non-vanishing and regular on $T_\w$. In general, they do vanish and acquire poles in the ambient cluster $\SX$-scheme $X(\beta)$.\\

\begin{fdef}[Cluster $\SA$-variables via $L_\w$]\label{def:Avariables_via_fillings}
Let $\beta$ be a positive braid word $\wb$ a Demazure weave. Consider the vertex $i\in Q_\w$ associated to the dual Lusztig cycle $\rho_i^\w$ and its homology class $[\rho_i^\w]\in \hH_1(L_\w\setminus\tt,\La_\beta\setminus\tt;\Z)$. By definition, the $\SA$-cluster variable $\SA_i^\w\in \Z(X(\beta))$ associated to $i$ is the rational function
\begin{equation}\label{eq:Avariable_via_fillings}
\SA_i^\w:X(\beta)\lr\PP^1_\Z,\quad\mbox{s.t. } \SA_i^\w(p)=\hol_{[\rho_i^\w]}(\LL_p)\quad\forall p\in T^\SA_\w,
\end{equation}
where $\hol_{[\rho_i^\w]}(\LL_p)$ denotes the monodromy of the local system $\LL_p$ along the relative class of $\rho_i^\w$.\sbs
\end{fdef}
\vspace{0.2cm}
In \Cref{def:Avariables_via_fillings}, the expression {\it monodromy of a local system} is understood as above. Namely, suppose that $[\rho_i^\w]$ is represented by a relative 1-cycle with its boundary contained in $\La_\beta\setminus\tt$, as in \Cref{sssec:X_lattice}, and assume that it is a connected interval. Then $\hol_{[\r_i^\w]}(\LL_p)$ refers to the automorphism of the stalk at one of its endpoints obtained by parallel transporting along the interval and then using the common trivialization for the local system at $\La_\beta\setminus\tt$, where both endpoints lie. As above, for several connected intervals we consider the product of the resulting functions. As before, this description makes it apparent that the rational functions $\SA_i^\w$ are regular non-vanishing functions on $T_\w^\SA$. It is a remarkable fact that, unlike the $\SX$-variables from \Cref{def:Xvariables_via_fillings}, the rational functions in \Cref{def:Avariables_via_fillings} are globally regular, which follows from the following algebraic description. Specifically:

\begin{enumerate}
    \item \Cref{def:cluster_Avariables_weave} and \Cref{def:cluster_Xvariables_weave} respectively introduced the cluster $\SA$-variables $A_i^\w$ and cluster $\SX$-variables $X_i^\w$ associated to a Demazure weave. These were defined algebraically, purely in terms of the weave $\w$ and the braid variety $X(\beta)$. Intuitively, they measure the relevant quantitative transversality among the flags given by a point of the braid variety, cf.~\Cref{sec:braid_variety_cluster}.\\

    \item \Cref{def:Avariables_via_fillings} and \Cref{def:Xvariables_via_fillings} introduced the cluster $\SA$-variables $\SA_i^\w$ and cluster $\SX$-variables $\SX_i^\w$ directly in terms of the symplectic topology of Lagrangian fillings, by using either Floer theory or microlocal sheaf techniques. By construction, these are invariants under compactly supported Hamiltonian isotopies.\\
\end{enumerate}

\noindent By \cite[Section 4]{CasalsWeng22}, and see also \cite[Section 7.2.1]{CasalsZaslow}, these two constructions coincide:

\begin{theorem}[\cite{CasalsWeng22}]\label{thm:comparison}
Let $\beta$ be a positive braid word, $X(\beta)$ its associated braid variety and $\wb$ a Demazure weave. Then
$$A_i^\w=\SA_i^\w\quad\&\quad X_i^\w=\SX_i^\w.$$
\sbs
\end{theorem}

\noindent\Cref{thm:comparison} implies Theorem \ref{thm:symplectic_intro1} and explains a few key properties from the weave side:

\begin{enumerate}
    \item It explains the invariance of $A_i^\w$ and $X_i^\w$ under weave equivalences from a geometric viewpoint. Indeed, all the Floer theoretic (or sheaf theoretic) invariants we used only depend on the compactly supported Hamiltonian isotopy class of the filling $L_\w$ of $\La_\beta$. Since a weave equivalence $\w\simeq\w'$ induces a Hamiltonian isotopy $L_\w\simeq L_{\w'}$, $\SA_i^\w$ and $\SX_i^\w$ are automatically invariant under weave equivalences.\\

    \item The global regularity of the cluster $\SA$-variables $\SA_i^\w$ can also be established purely in geometric terms, without introducing any coordinates, via \cite[Theorem 1.1]{CasalsLi24}. \Cref{thm:comparison} transports this regularity to the cluster $\SA$-variables $A_i^\w$ defined purely via weaves.\\

    \item \Cref{thm:comparison} historically explains the introduction of the Reeb $z$-coordinates to describe $A_i^\w$ and $X_i^\w$ as functions on the braid variety, as these Reeb $z$-coordinates naturally generate the Legendrian dg algebra of $\La_\beta$, from which $X(\beta)$ and subsequently $\SA_i^\w$ and $\SX_i^\w$ can be obtained by counts of pseudoholomorphic strips or microlocal holonomies.\\

    \item The cyclic invariance of $X(\beta)$, the $\T$-action and some of the additional symmetries of $X(\beta)$ and its cluster structure, cf.~\Cref{sssec:braid_variety_further_symmetries}, can all be realized through the symplectic topology related to $\La_\beta$ and its Lagrangian fillings. For instance, the cyclic symmetry is explained geometrically by a Legendrian isotopy of the Legendrian link $\La_\beta$. The invariance property of the Legendrian dg algebra then provides the necessary isomorphism describing the cyclic symmetry of $X(\beta)$ and explains why it is a quasi-cluster transformation.\\
\end{enumerate}

\begin{remark}
Geometrically, the dual Lusztig cycles $\rho_i^\w$ typically give representatives of their relative classes that are not connected. Therefore the associated cluster variable $\SA_i^\w$ is obtained as a product of holonomies, one per each interval, as explained above. We emphasize that \Cref{thm:comparison} implies that this product is an irreducible regular function, if the cluster $\SA$-variable is mutable. (Else it is a unit.) Irreducibility here does not create a contradiction with being expressed as a product, as each of the factors in the product, coming from the holonomy along each connected component of $\rho_i^\w$, is only a rational function, and typically {\it not} a regular function.\sbs
\end{remark}

%%%%%%%%%%%%%%%%%%%%%%%%%%%%%%%%%%%%%%%%%%%%%%%%%%%%%%%%
%%%%%%%%%%%%%%%%%%%%%%%%%%%%%%%%%%%%%%%%%%%%%%%%%%%%%%%%
%%%%%%%%%%%%%%%%%%%%%%%%%%%%%%%%%%%%%%%%%%%%%%%%%%%%%%%%
\subsection{Cluster mutation in the context of Lagrangian fillings}\label{ssec:cluster_mutation_Lagrangian_filling_from_weave} Let $\wb$ be a Demazure weave and suppose that $\mu_i(\w)$ is obtained from $\w$ by a weave mutation at $\g^\w_i$. The goal of this section is to describe a geometric operation on embedded Lagrangian fillings which explains how $L_\w$ and $L_{\mu_i(\w)}$ can be geometrically obtained from each other. The summarized description is:

\begin{center}
  \fbox{%
    \begin{tabular}{l}
      \itshape The Lagrangian filling $L_{\mu_i(\w)}$ is obtained from $L_\w$ by Lagrangian disk surgery at a Lagrangian disk\\
      \itshape  $\D\sse\D^4\setminus L_{\w}$ whose boundary is the Lusztig cycle $\partial\D=\g_i^\w$.
    \end{tabular}%
  }
\end{center}
\vspace{0.2cm}

In \Cref{sssec:local_model_surgery} we now describe the operation of Lagrangian disk surgery, entirely within the context of 4-dimensional symplectic topology. Its connection to weave mutation is crystallized in \Cref{sssec:disk_surgery_via_weaves}.

\subsubsection{Lagrangian disk surgery}\label{sssec:local_model_surgery} Let us proceed with a more detailed description of Lagrangian disk surgery. Let $L\sse (\D^4,\la_\std)$ be an embedded exact Lagrangian filling of a Legendrian link $\La=\dd L \sse \dd(\D^4,\la_\std)$. To construct new Lagrangian fillings that are {\it not} Hamiltonian isotopic to $L$, we could try to perform a local operation on $L$. That said, a local operation near a point $p\in L$ can be challenging to find: e.g.~\cite[Theorem 1.1.A]{EliashbergPolterovich96} implies that local knotting is {\it not} possible for Lagrangian surfaces.

\begin{center}
	\begin{figure}[h!]
		\centering
		\includegraphics[width=\textwidth]{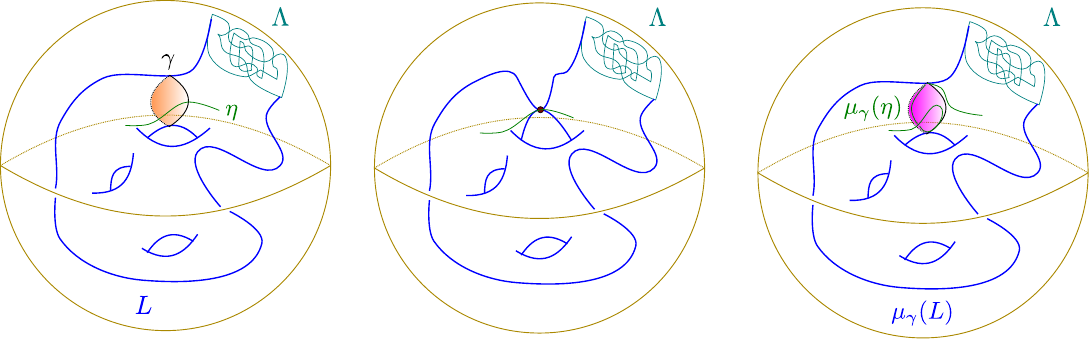}
		\caption{A Lagrangian filling $L\sse(\D^4,\la_\std)$ of a Legendrian link $(\La\sse \dd \D^4,\xi_\std)$ undergoing Lagrangian disk surgery. (Left) The Lagrangian disk $\D\sse \D^4\setminus L$ being surgered is depicted in orange, with boundary $\dd \overline{\D}=\gamma\sse L$. (Center) The critical moment of the surgery, where the filling $L$ degenerates to an immersed exact Lagrangian. (Right) The Lagrangian filling $\mu_\D(L)$, also denoted $\mu_\gamma(L)$, which is the result of the surgery. Here the new Lagrangian disk $\mu_\D(\D)\sse \D^4\setminus \mu_\g(L)$ after the surgery is depicted in purple, with boundary $\gamma\sse \mu_\g(L)$ oriented in reverse.}
		\label{fig:Talks_CartoonSurgery}
	\end{figure}
\end{center}

\begin{center}
	\begin{figure}[h!]
		\centering
		\includegraphics[scale=0.5]{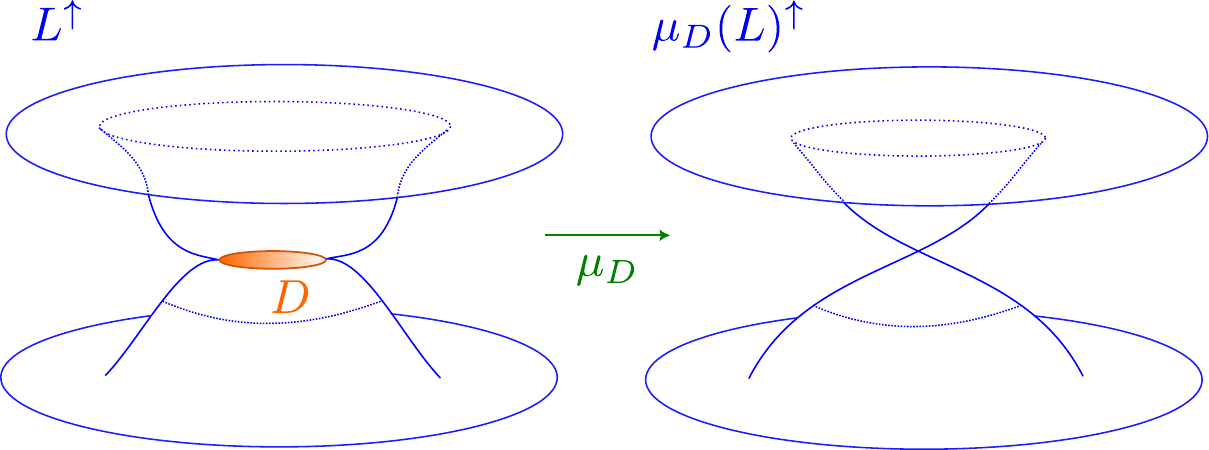}
		\caption{The effect of Lagrangian disk surgery seen in fronts for the Legendrian lifts $L^\uparrow$ and $\mL^\uparrow$ of $L$ and $\mL$. These fronts can be understood in arbitrary dimensions, as spherical spins of the cusp, for the left front, and the cone point, for the right front. The disk $\D$, in orange, is the same disk in orange as in \cref{fig:Talks_CartoonSurgery} (Left) and a neighborhood of the cone point corresponds to a neighborhood of the purple disk in \cref{fig:Talks_CartoonSurgery} (Right).}
		\label{fig:Surgery_InFront}
	\end{figure}
\end{center}

\noindent L.~Polterovich introduced a beautiful construction in \cite[Sections 3\&4]{Polterovich_Surgery} that produces a (potentially) new Lagrangian filling $\mu_D(L)$ if one has an embedded Lagrangian disk $\D\sse \D^4\setminus L$ whose boundary $\dd\D\sse L$ belongs to $L$ as a clean intersection, i.e.~$\D\cap L=\dd\D$ and as transversely as possible. See also \cite[Section 3.3]{LalondeSikorav91_LagrangianSurgery} for the case of Lagrangian surfaces. This operation is local in a neighborhood of the disk $\D$, so the Lagrangian embeddings of $L$ and $\mu_D(L)$ coincide outside of a neighborhood of $\dd \D$. Intuitively, this operation can be depicted as in \cref{fig:Talks_CartoonSurgery}: the Lagrangian disk $\D$ is contracted to a point and expanded in a {\it different} way. It is proven in \cite[Theorem 1.1]{thomson2026legendrianhopflinkexactly} that there are only two ways to locally expand such contracted point, and thus this other expansion is locally unique. For a visualization locally around $\D$, the Legendrian lifts of $L$ and the surgered $\mu_D(L)$ can be represented via fronts as in \cref{fig:Surgery_InFront}. This is justified by the computations that now follow, cf.~also \cite[Section 6.2]{CasalsMurphyPresas}.\\

\begin{figure}
     \centering
     \begin{subfigure}[b]{0.4\textwidth}
         \centering
         \includegraphics[width=\textwidth]{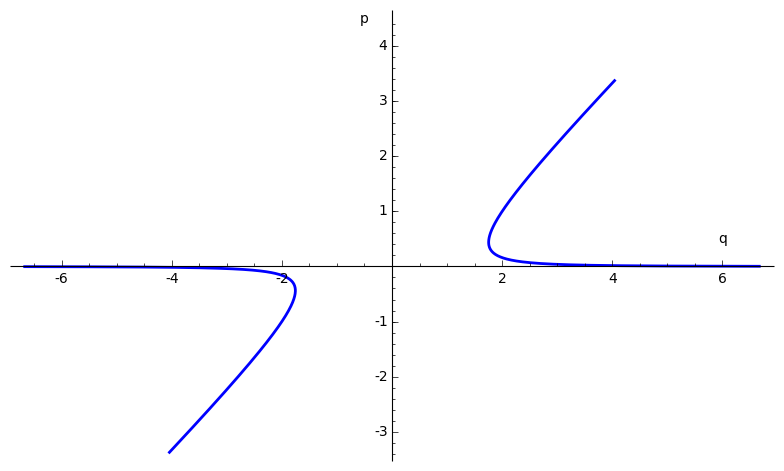}
         \caption{Lagrangian handle $\Gamma^+$.}
         \label{fig:y equals x}
     \end{subfigure}
     %\hfill
     \begin{subfigure}[b]{0.4\textwidth}
         \centering
         \includegraphics[width=\textwidth]{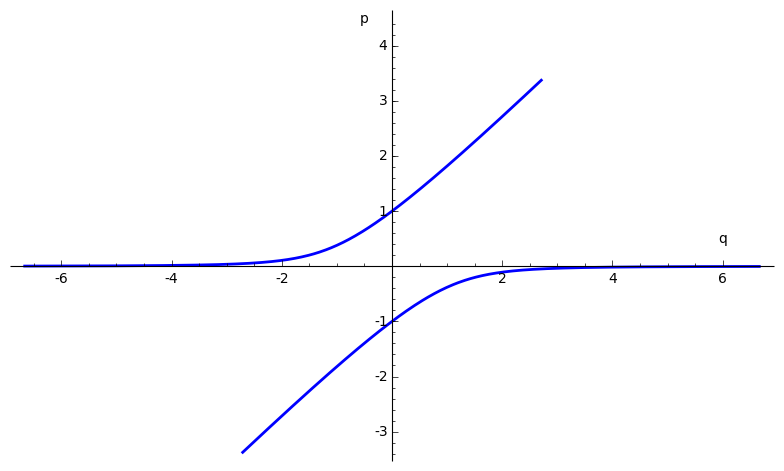}
         \caption{Lagrangian handle $\Gamma^-$.}
         \label{fig:three sin x}
     \end{subfigure}
     \caption{The Lagrangian handles $\Gamma^{\pm}\subseteq\R^{2n}(q,p)$.}
     \label{fig:Lagr_handles}
\end{figure}

The local model for the Lagrangian surgery in \cite{Polterovich_Surgery} can be described as follows. Consider local coordinates $(q_1,q_2,p_1,p_2)\in\R^{4}$ and model the two sheets of an immersed point in a Lagrangian, as the one in \cref{fig:Talks_CartoonSurgery} (Center), and the Liouville form via
$$S_1:=\{p_1=0,p_2=0\},\quad S_2:=\{q_1=p_1,q_2=p_2\},\quad \displaystyle\lambda_\std=p_1dq_1+p_2dq_2.$$
The two Lagrangian handles $\Gamma^\pm$ presented in \cite{Polterovich_Surgery} provide two ways to resolve such an immersed point, resulting in an embedded exact Lagrangian. These Lagrangian handles are depicted in Figure \ref{fig:Lagr_handles}, and in order to explicitly parametrize them we use coordinates $t=(t_1,t_2)\in\R^{2}$. First, the positive Lagrangian handle $\Gamma^+$ can be described via the parametrization $\Gamma^+:\R^{2}\setminus\{0\}\longrightarrow\R^{4}$
defined as
$$\Gamma^+(t_1,t_2)=\left((\mu+\mu^{-1})t_1,(\mu+\mu^{-1})t_{2},\mu t_1,\mu t_{2}\right)\mbox{ where }\displaystyle\mu=t_1^2+t_2^2.$$

\noindent Note that $\Gamma^\pm$ are diffeomorphic to $S^{1}\times\R$, we have the two asymptotics
$$\displaystyle\lim_{\mu\to\infty}\Gamma^+\subseteq S_2,\quad \displaystyle\lim_{\mu\to0}\Gamma^+\subseteq S_1,$$
and in this model $L$ and $\D$ are described as follows:

\begin{enumerate}
    \item The given Lagrangian $L$ is obtained by gluing the above positive Lagrangian handle $\Gamma^+$ to the Lagrangian sheet $S_1$ at the limit $\mu=0$, and to the Lagrangian sheet $S_2$ at the limit $\mu=\infty$. Thus, intuitively, locally near the boundary $\dd\D\sse L$, the Lagrangian is modeled by $\Gamma^+$.\\

    \item The Lagrangian disk $\D$ is modeled within the Lagrangian graph $\Gamma_\ell$ of the linear function
    $$\ell:\R^2_{q_1,q_2}\lr\R^2_{p_1,p_2},\quad p_i=\frac{q_i}{2},\quad i\in[1,2].$$
    Specifically, $\D$ is the image of the restriction of $\Gamma_\ell$ to $\{\|q\|\leq2\}$.
    The boundary $\dd\D\sse L$ of $\D$ in this model can be taken to be the intersection $\Gamma_\ell\cap\Gamma^+$, which can be parametrized by $\Gamma^+(t_1,t_2)$ where $(t_1,t_2)\in S^1$ belong to the 1-sphere.\footnote{Note that, generically, two Lagrangian submanifolds intersect in points, not in codimension 1 subsets of each other. This model is {\it not} a generic situation and $\Gamma_\ell\cap\Gamma^+\cong S^{1}$ being diffeomorphic to a sphere does not persist under a generic Hamiltonian isotopy applied to either Lagrangian.}
\end{enumerate}

\noindent The surgered Lagrangian $\mu_\D(L)$ is obtained by using the negative Lagrangian handle $\Gamma^-:\R^2\setminus\{0\}\longrightarrow\R^{4}$ parametrized by
$$\Gamma^-(t_1,t_{2})=\left((\mu-\mu^{-1})t_1,(\mu-\mu^{-1})t_{2},\mu t_1,\mu t_{2}\right).$$
This parametrization satisfies the asymptotics $\displaystyle\lim_{\mu\to\infty}\Gamma^-\subseteq S_2$ and $\displaystyle\lim_{\mu\to0}\Gamma^-\subseteq S_1$, and can be glued to $S_1$ and $S_2$ in the asymptotic limits, thus constructing the Lagrangian $\mu_\D(L)$:

\begin{fdef}[Lagrangian disk surgery]\label{def:disk_surgery} In the local model above, the Lagrangian disk surgery of $L$ along $\D$ is the embedded exact Lagrangian $\mu_\D(L)$ obtained by exchanging the positive Lagrangian $\Gamma^+\sse L$ by the negative Lagrangian handle $\Gamma^-$.
\end{fdef}
\vspace{0.2cm}
\noindent This Lagrangian disk surgery in \Cref{def:disk_surgery} has the following properties:\\

\begin{enumerate}
    \item $\Gamma^+$ and $\Gamma^-$ are smoothly isotopic and thus $L$ and $\mu_\D(L)$ are smoothly isotopic. In fact, since $\Gamma^+$ and $\Gamma^-$ are even Lagrangian isotopic, so are $L$ and $\mu_\D(L)$.\\

    \item $\Gamma^+$ and $\Gamma^-$ are not Hamiltonian isotopic via a compactly supported Hamiltonian diffeomorphism of $\R^{4}$. Equivalently, $\Gamma^+$ and $\Gamma^-$ are not exact Lagrangian isotopic relative to their boundaries. This follows, for instance, from the fact that $\Gamma^\pm$ are Lagrangian fillings of the Hopf link inducing different augmentations, cf.~\cite[Section 8.1]{EHK}. Alternatively, and historically preceding ibid., this follows from the fact that the Chekanov and Clifford tori, which differ by a unique Lagrangian disk surgery, are not Hamiltonian isotopic, as proven in \cite[Section 4]{Chekanov96_Torus}.\\

    \item In the model above, the new Lagrangian 2-disk $\mu_\D(\D)$ is given by
    $$\mu_\D(\D):=\{(q,p)\in\R^{4}: q_1=0,q_2=0,\|p\|\leq1\},$$
    whose boundary $\dd\mu_\D(\D)$ is indeed contained in $\mbox{im}(\Gamma^-)$, as can be parametrized by $\Gamma^-(t_1,t_2)$ where $(t_1,t_2)\in S^{1}$ belong to the unit sphere.\\

    \item Lagrangian disk surgery is an involutive operation, in the sense that
    $$\mu_{\mu_\D(D)}(\mu_\D(L))\simeq L$$
    are Hamiltonian isotopic. This follows from the fact that, in the local model above, Lagrangian disk surgery applied to $\Gamma^-$ along $\mu_\D(\D)$ results in a Lagrangian which is Hamiltonian isotopic to $\Gamma^+$.\\

    \item It is natural to wonder whether there is a third way to resolve an immersed point, not just with $\Gamma^\pm$. It follows from \cite[Theorem 1.1]{thomson2026legendrianhopflinkexactly} that $\Gamma^\pm$ are the only embedded exact Lagrangian fillings of the Hopf link, up to Hamiltonian isotopy. Thus $\Gamma^\pm$ are the only two possible Lagrangian handles we could have constructed, up to Hamiltonian isotopy.\\
\end{enumerate}

%%%%%%%%%%%%%%%%%%%%%%%%%%%%%%%%%%%%%%%%%%%%%%%%%%%%%%%%
%%%%%%%%%%%%%%%%%%%%%%%%%%%%%%%%%%%%%%%%%%%%%%%%%%%%%%%%
%%%%%%%%%%%%%%%%%%%%%%%%%%%%%%%%%%%%%%%%%%%%%%%%%%%%%%%%
\subsubsection{Lagrangian disk surgery via weaves}\label{sssec:disk_surgery_via_weaves} 
\begin{center}
	\begin{figure}[h!]
		\centering
		\includegraphics[width=\textwidth]{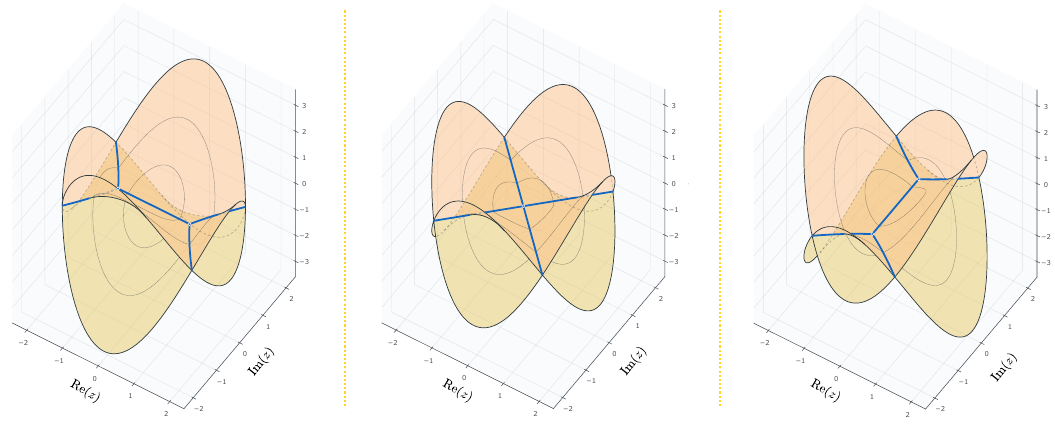}
		\caption{(Left) The weave, inside the corresponding 3-dimensional graph, associated to the exact Lagrangian cylinder $L_+$ from \eqref{eq:Lagrangian_surgery_models}. (Center) The weave corresponding to the immersed exact Lagrangian $\{w^2-z^2=0\}$. (Right) The weave associated to the exact Lagrangian cylinder $L_-$ from \eqref{eq:Lagrangian_surgery_models}.}
		\label{fig:Graphs_weave_mutation}
	\end{figure}
\end{center}

The local model for Lagrangian disk surgery, as presented in \Cref{sssec:local_model_surgery}, can be realized as a weave mutation between Demazure weaves, as follows. Consider the singularity $P(z,w)=w^2-z^2$ and, following \Cref{ssec:Lagrangian_filling_from_AG}, we can consider the two real Morsifications:
\begin{equation}\label{eq:Lagrangian_surgery_models}
P_\pm(z,w):=w^2-z^2\pm1,\qquad L_\pm:=\{(z,w)\in\C^2:w^2-z^2\pm1=0\}
\end{equation}
\noindent As in \Cref{lem:complex_curve_real_Lagrangian}, each $L_\pm$ provides a real Lagrangian cylinder in $(\R^4,\la_\std)$ and each is exact, since the residue formula implies that the unique period of $zdw-wdz$ on each of $L_\pm$ is an integer multiple of $2\pi i$. The graphs of the fronts of their Legendrian lifts of $L_\pm$ are depicted in \Cref{fig:Graphs_weave_mutation}, along with that of the immersed Lagrangian $\{P(z,w)=0\}$.

\noindent An $\L$-compressing disk $\D_\pm\sse\R^4\setminus L_\pm$ for each of $L_\pm$ can be described as follows:\\

\begin{enumerate}
    \item Consider the embedded 2-disk
    \begin{equation}\label{eq:compressing_disk_MorsificationA1_plus}
    \begin{split}
    \D_+&:=\left\{(z,w)\in\C^2: \Re(z)\in[-1,1],\Im(z)=0, w\in\left[-\sqrt{z^2-1},\sqrt{z^2-1}\right]\right\}\cong\\
        &=\{(x,iy):x,y\in\R,x^2+y^2\leq 1\}
        \end{split}
    \end{equation}
    which is Lagrangian and, by construction, $\dd\D_+\sse L_+$ provides a representative for the generator of $\hH_1(L_+,\Z)$. The associated weave $\w_+$ is depicted in \Cref{fig:Graphs_weave_mutation} (left), and it has trivalent vertices at $z=\pm1$ with the only mutable Lusztig cycle being supported at the segment $\Re(z)\in[-1,1]$ with $\Im(z)=0$. Note that the projection of $\D_+$ to $\R^2_{x_1,x_2}\cong\C_z$ is indeed given by this Lusztig cycle.\\
    
    \noindent Compressing the Lagrangian 2-disk $D_+$ as in \Cref{sssec:local_model_surgery} leads to the immersed Lagrangian $\{P(z,w)=0\}\cong \C_{\{z=w\}}\cup \C_{\{z=-w\}}\sse\C^2$, which has a unique (immersed) singular point at the origin. The front for its Legendrian lift is drawn in \Cref{fig:Graphs_weave_mutation} (center). This degeneration, given by compressing the disk, can be performed through exact Lagrangians by introducing a real parameter $\alpha\in(0,1]$ in \eqref{eq:Lagrangian_surgery_models} and considering the exact Lagrangian surfaces $\{w^2-z^2+\alpha=0\}$ from $\alpha=1$ to $\alpha=0$.\\

     \item Similarly, we can consider the embedded 2-disk
    \begin{equation}\label{eq:compressing_disk_MorsificationA1_plus2}
        \begin{split}
        \D_-&=:=\left\{(z,w)\in\C^2: \Re(z)=0,\Im(z)\in[-1,1], w\in\left[-\sqrt{z^2+1},\sqrt{z^2+1}\right]\right\}\\
                &=\{(ix,y):x,y\in\R,x^2+y^2\leq 1\}
        \end{split}
    \end{equation}
    which is an exact Lagrangian disk with $\dd\overline{\D}_-\sse L_-$. As above, the projection of $\D_-$ to $\R^2_{x_1,x_2}\cong\C_z$ is given by the Lusztig cycle \{$\Re(z)=0,\Im(z)\in[-1,1]$\} associated to the weave $\w_-$ for $L_-$. The two trivalent vertices of such a weave are located at $z=\pm i$, as illustrated in \Cref{fig:Graphs_weave_mutation} (right). Compressing $D_-$ as in \Cref{sssec:local_model_surgery} can be realized similarly via the exact Lagrangian surfaces $\{w^2-z^2-\alpha=0\}$ where $\alpha\in[0,1]$ is a real parameter.\\
\end{enumerate}

In conclusion, the Lagrangian disk surgery, as described in \Cref{sssec:local_model_surgery}, can be locally modeled by the transition from one of the two sides of \Cref{fig:Graphs_weave_mutation} to the other, through the center model. Thus we obtain

\begin{equation}\label{eq:mutation_via_fillings}
       \mu_{\D_+}(L_+)\simeq L_-\qquad\&\qquad\mu_{\D_-}(L_-)\simeq L_+
\end{equation}

\noindent In combination with \Cref{fig:Graphs_weave_mutation}, \eqref{eq:mutation_via_fillings} geometrically justifies that the diagrammatic weave mutation from \Cref{fig:weave_mutation} (left) is diagrammatically encoding a Lagrangian disk surgery. Namely, if $\w_+$ and $\w_-$ are the only two weaves $\w_\pm:\s_1^3\lr\s_1$, \Cref{fig:weave_mutation} (left), which differ by a weave mutation, then:

        \begin{equation}\label{eq:mutation_via_fillings}
        \fbox{$\mu_{\D_+}(L_{\w_+})\simeq L_{\w_-}\quad\&\quad \mu_{\D_-}(L_{\w_-})\simeq L_{\w_+}$}
        \end{equation}

\noindent This geometric realization of weave mutation, with $L_\pm$ not being Hamiltonian isotopic, was a strong reason to define weave mutation as we did in \Cref{def:weave_mutation}.

\begin{remark}
It follows from the local model above that $L_\pm$ are actually compactly supported smoothly isotopic to each other. In particular, their relative homology classes coincide. This explains why in diagrammatic Soergel calculus, cf.~\cite{EliasWilliamson2016Calculus,EMTW2020}, the move \Cref{fig:weave_mutation} (left) is instead considered an equivalence. Indeed, in Soergel calculus only the homology classes of these surfaces are relevant as a Hochschild homology is being taken, cf.~\cite{Webster_2008}.\sbs
\end{remark}

%%%%%%%%%%%%%%%%%%%%%%%%%%%%%%%%%%%%%%%%%%%%%%%%%%%%%%%%%%%%%%%%%
%%%%%%%%%%%%%%%%%%%%%%%%%%%%%%%%%%%%%%%%%%%%%%%%%%%%%%%%%%%%%%%%%

\subsubsection{$\L$-compressing systems}\label{sssec:L_compressing_systems} Let $\La\sse(T^\infty\R^2,\xi_\std)$ be a Legendrian link and $L\sse (T^*\R^2,\la_\std)$ a Lagrangian filling. We can try to apply Lagrangian disk surgeries to $L$, as discussed in \cref{sssec:local_model_surgery}, to produce Lagrangian fillings that are potentially new, i.e.~not Hamiltonian isotopic to $L$. The input of a Lagrangian disk surgery is $L$ and an embedded Lagrangian disk $\D\sse\D^4\setminus L$ with $\dd\overline{D}\sse L$. The next question thus becomes: \textcolor{blue}{given a Lagrangian $L$, how do we find such embedded Lagrangian 2-disks $\D$}?\\

\noindent The short answer is that this is not known in general. We will nevertheless present many interesting cases, some connected to singularity theory and some to algebraic combinatorics, in which these disks appear naturally. In order to have a precise language to talk about these situations, we introduce the following notion:

\begin{fdef}[$\L$-compressing systems]\label{def:Lcompressing_systems} Let $L\sse (\D^4,\la_\std)$ be a Lagrangian filling of $\dd L\sse(S^3,\xi_\std)$. By definition, an $\L$-compressing system $\D(L):=\{\D_1,\ldots,\D_s\}$ of rank $s$ is a collection of disks $\D_i\sse \D^4\setminus L$, $i\in[1,s]$, such that\\
\begin{enumerate}
    \item $\D_i$ is an embedded Lagrangian disk with $\dd\overline{\D_i}\sse L$ a clean intersection, for all $i\in[1,s]$,\\
    \item The collection of boundary classes $\{\dd\overline{\D_1},\ldots,\dd\overline{\D_s}\}$ is linearly independent in $H_1(L)$.\\
\end{enumerate}
By definition, an $\L$-compressing system $\D(L):=\{\D_1,\ldots,\D_s\}$ is said to be maximal if the boundaries of the disks $\D_i\in\D(L)$ span $H_1(L,\Z)$.\sbs
\end{fdef}
\vspace{0.2cm}
\noindent A given Lagrangian $L\sse (\D^4,\la_\std)$ might not admit a maximal $\L$-compressing system. To wit, the Legendrian representative of $m(5_2)$ in \cref{fig:Front_Arc_Cusp_Crossing_R3_Example1} (left) has two known distinct embedded Lagrangian fillings, both once-punctured tori: both of them admit an $\L$-compressing system of rank 1, but neither of them admits an $\L$-compressing system of rank 2. There are relaxations of \cref{def:Lcompressing_systems}, allowing for {\it immersed} Lagrangian disks, cf.~\cite[Section 3]{CasalsWeng22}, but we focus on maximal $\L$-compressing systems in this manuscript.

\begin{remark} (1) \cref{def:Lcompressing_systems} naturally generalizes to embedded exact Lagrangians $L\sse (\D^4,\la_\std)$ in a Liouville manifold $(W,\la_\std)$ by using the concept of relative Lagrangian skeleton. For instance, a maximal $\L$-compressing system for a given $L\sse(\D^4,\la_\std)$ is a collection $\D(L)$ of embedded Lagrangian disks $\D_i\sse \D^4\setminus L$, $i\in[1,s]$, with $\dd\overline{D_i}\sse L$ such that
$$L\cup(\D_1\cup\ldots\D_s)$$
is a relative Lagrangian skeleton of $(\D^4,\la)$, relative to the Weinstein hypersurfaces given by the Weinstein ribbon of $\dd L\sse\dd\D^4$. In the case of $(W,\la)\cong (\D^4,\la_\std)$, this generalizes the condition of $\{\dd\D_1,\ldots,\dd\D_s\}$ being a basis of $H_1(L,\Z)$.\\

\noindent (2) It is oftentimes useful to relax the condition that the boundary homology classes of the disks $\D_i$ are linearly independent. This case is particularly relevant in the study of Landau-Ginzburg models which appear on the B-side of the homological mirror symmetry correspondence for Fano varieties, see e.g.~\cite{PascaleffTonkonog20}.\sbs
\end{remark}

\begin{example}\label{ex:Lcompressing_systems} (1) The Legendrian links $\La(\bG)$ associated to plabic graphs $\bG$, as the lifts of their alternating strand diagram, come endowed with both a Lagrangian filling {\it and} an $\L$-compressing system. Intuitively, the $\L$-compressing system is given by the bounded faces of $\bG$, understood as Lagrangian disks in the zero section. See e.g.~\cite[Section 4.3]{CasalsGao24}, and see \cite[Section 3]{casals2023demazureweavesreducedplabic} for their relationship to Demazure weaves $\w$ and their Lagrangian fillings $L_\w$. These $\L$-compressing systems are maximal in many interesting cases, e.g.~for reduced plabic graphs.\\

\noindent (2) Given an isolated plane curve singularity $f:\C^2\lr\C$ with a real Morsification, the Lagrangian fillings from \cref{ex:Lagrangian_fillings}.(3) come endowed with a maximal $\L$-compressing system. The Lagrangian disks are provided by the Lefschetz thimbles, and their boundary curves in the Lagrangian Milnor fiber are the vanishing cycles, cf.~\cite[Section 2]{CasalsLagSkel}.\sbs
\end{example}

%%%%%%%%%%%%%%%%%%%%%%%%%%%%%%%%%%%%%%%%%%%%%%%%%%%%%%%%%%%%%%%%%
%%%%%%%%%%%%%%%%%%%%%%%%%%%%%%%%%%%%%%%%%%%%%%%%%%%%%%%%%%%%%%%%%

\subsubsection{$\L$-compressing systems from Demazure weaves}\label{sssec:L_compressing_systems_via_weaves} Let $\wb$ be a Demazure weave and $\g_i^\w$ a mutable Lusztig cycle. If we want to geometrically realize the weave mutation $\w\stackrel{\mu}{\sim}\w'$ via a Lagrangian disk surgery we must find a Lagrangian 2-disk $\D_i\sse\D^4\setminus L_\w$ in the complement of $L_\w$ whose boundary represents the homology class $\g_i^\w$. Recall that $L_\w\sse\D^4$ is embedded in a way that the restriction of the projection
$$\pi:\D^4\lr\R^2$$
to $L_\w$ induces a simple branched cover $\pi|_{L_\w}:L_\w\lr\R^2$, for which the weave $\w$ acts as a branch cut. In the same manner that we associated a curve in $L_\w$ to a Lusztig cycle $\g_i^\w$ by using this branched cover, cf.~\Cref{sssec:X_lattice}, we can in many cases find such 2-disk $\D_i$ in a way that it maps to the support of the Lusztig cycle $\g_i^w$. Specifically, let us consider Lusztig cycles $\g_i^\w$ whose weights are either 0 or 1 and whose support in $\w$ is a tree, seen as a graph. Such weave cycles are referred to as $\mathsf{Y}$-cycles, cf.~\cite{CasalsWeng22,CasalsZaslow}. For these Lusztig cycles one obtains the following first result:

\begin{proposition}[$\L$-compressing 2-disk above a Lusztig $\mathsf{Y}$-cycle]\label{prop:compressing_disk_Lusztig_cycle}
Let $\wb$ be a Demazure weave, $L_\w\sse(\R^4,\la_\std)$ its associated Lagrangian filling, and $\g_i^\w$ a mutable Lusztig $\mathsf{Y}$-cycle. Then there exists a properly embedded Lagrangian 2-disk $\D_i\sse\R^4\setminus L_\w$ such that:\\
\begin{enumerate}
    \item $[\dd\D_i]=[\g_i^\w]$ and, after applying a sequence of weave equivalences to $\w$, the image $\pi(\D_i)$ coincides with the support of the corresponding Lusztig cycle $\g_i^\w$.\\

    \item The Lagrangian disk surgery of $\L_w$ along $D_i$ is Hamiltonian isotopic to $L_{\mu_i(\w)}$, i.e.~
        \begin{equation}\label{eq:mutation_via_fillings2}
        \fbox{$\mu_{\D_i}(L_\w)\simeq L_{\mu_i(\w)}$}
        \end{equation}
    where $\mu_i(\w)$ denotes the weave mutation of $\w$ at $\g_i^\w$.\sbs
\end{enumerate}
\end{proposition}

\noindent The simplest instance of \Cref{prop:compressing_disk_Lusztig_cycle}, and one of the ingredients for its proof, is the local model for Lagrangian disk surgery, presented in \Cref{sssec:local_model_surgery}, as realized via weaves in \Cref{sssec:disk_surgery_via_weaves}. Indeed, the two Lagrangian 2-disks $\D_\pm$ in \eqref{eq:compressing_disk_MorsificationA1_plus} and \eqref{eq:compressing_disk_MorsificationA1_plus2} project exactly over the unique Lusztig cycle for the corresponding Demazure weaves $\w:\s_1^3\lr\s_1$.\sbs \\

In general, it is an interesting question whether one can iterate such Lagrangian disk surgeries. For instance, given a mutation sequence (e.g.~by naming a sequence of vertices in the quiver $Q_\w$), the question is whether it is possible to geometrically realize it by iterated Lagrangian disk surgeries. At its core, this problem reduces to a question of algebraic intersection numbers, which are tracked by the arrows of the quivers, and geometric intersection numbers. Intuitively, if these algebraic and geometric intersection numbers coincide, then iterating is possible.\\

A study of this phenomenon can be found in \cite{CasalsGao24}, where it is shown that any cluster seed in $\Z[X(\beta)]$ can be obtained as the seed associated to a Lagrangian filling, as stated in Theorem \ref{thm:symplectic_intro2}. The proof shows that it is possible to geometrically iterate, using Lagrangian disk surgeries, such algebraic sequences of quiver mutations. The key idea is to introduce a potential $W_\w$ for $Q_\w$ that keeps track of the difference between the algebraic and geometric intersection numbers and show that such a quiver with potential is non-degenerate. In fact, it is proven in ibid.~that the quiver with potential $(Q_\w,W_\w)$ is rigid. From a symplectic topological perspective, the quiver with potential $(Q_\w,W_\w)$ encodes the $A_\infty$-structure for the union of all the exact Lagrangian curves in the symplectic push-off of the Lagrangian filling $L_\w$. In a sense, such a union can be seen as a collection of vanishing cycles and, being in a symplectic surface, all the relevant pseudo-holomorphic disk contributions can be combinatorialized in terms of polygons, as presented in \cite[Section 2]{CasalsGao24}.

%%%%%%%%%%%%%%%%%%%%%%%%%%%%%%%%%%%%%%%%%%%%%%%%%%%%%%%%%%%%%%%%%
%%%%%%%%%%%%%%%%%%%%%%%%%%%%%%%%%%%%%%%%%%%%%%%%%%%%%%%%%%%%%%%%%
%{\bf EDITED TILL HERE\\}

%%%%%%%%%%%%%%%%%%%%%%%%%%%%%%%%%%%%%%%%%%%%%%%%%%%%%%%%
%%%%%%%%%%%%%%%%%%%%%%%%%%%%%%%%%%%%%%%%%%%%%%%%%%%%%%%%
%%%%%%%%%%%%%%%%%%%%%%%%%%%%%%%%%%%%%%%%%%%%%%%%%%%%%%%%
%\subsection{Summary of section}

\bibliography{biblio} 
\bibliographystyle{plain}

\textsc{Roger Casals: University of California Davis, Dept. of Mathematics, USA.\\}
\textit{Email address:} \texttt{casals@ucdavis.edu}
\vspace{0.2cm}

\end{document}